\documentclass[12pt]{article}
\usepackage{remark,oldgerm,amsmath,amssymb,theorem}
\usepackage{diagrams}
\theorembodyfont{\normalfont\slshape}
\font\lcirc=lcircle10 %scaled 800; from Paul Vojta
\def\cxdot{{\hskip2.5pt\raise3.8pt\hbox{\lcirc\char113}}}
\def\cx{{\textstyle{\cdot}}}
\def\mapright#1{\smash{\mathop{\longrightarrow}\limits^{#1}}}
\def\mapdown#1{\Big\downarrow\rlap
    {$\vcenter{\hbox{$\scriptstyle#1$}}$}}

\long\def\ignore#1{}

\newcommand{\openbox}{\leavevmode
  \hbox to.77778em{%
  \hfil\vrule
  \vbox to.675em{\hrule width.6em\vfil\hrule}%
  \vrule\hfil}}
\newcommand{\qedsymbol}{\openbox}
\DeclareRobustCommand{\qed}{%
  \ifmmode % if math mode, assume display: omit penalty etc.
  \else \leavevmode\unskip\penalty9999 \hbox{}\nobreak\hfill
  \fi
  \quad\hbox{\qedsymbol}}
\newcommand{\op}{{op}}
\newcommand{\ie}{\textit{i.e.}}

\newcommand{\proofname}{{ Proof:}}
\newenvironment{proof}[1][\proofname]{\par
 \normalfont
\trivlist
  \item[\hskip\labelsep\itshape
    #1
]\ignorespaces
}{%
 \qed\endtrivlist}
\def\oh#1{{\cal O}_{#1}}

\newcommand{\cHom}{{\mathcal Hom}}
\def\uar#1{\mathop{#1}\limits_{\raise3pt\hbox{$\to$}}}

\def\boxit#1{\vbox{\hrule\hbox{\vrule\kern1pt
    \vbox{\kern1pt#1\kern1pt}\kern1pt\vrule}\hrule}}

\def\varprojlim{\mathop{\vtop{\ialign{$##$\cr
 \hfil{\fam0lim}\hfil\cr\noalign{\nointerlineskip}%
 {\leftarrow}\mkern-6mu\cleaders\hbox{$\mkern-2mu{-}\mkern-2mu$}\hfill
 \mkern-6mu{-}\cr\noalign{\nointerlineskip\kern-.2326ex}\cr}}}}
\def\varinjlim{\mathop{\vtop{\ialign{$##$\cr
 \hfil{\fam0lim}\hfil\cr\noalign{\nointerlineskip}%
 {-}\mkern-6mu\cleaders\hbox{$\mkern-2mu{-}\mkern-2mu$}\hfill
 \mkern-6mu{\to}\cr\noalign{\nointerlineskip\kern-.2326ex}\cr}}}\nolimits}
\def\invlim{\varprojlim}
\def\dirlim{\varinjlim}

\def\uar#1{\vec {#1}}

 \newcommand{\bc}{{\bf C}}
 
 \newcommand{\bn}{{\bf N}}
 \newcommand{\bz}{{\bf Z}}

\newcommand{\et}{{\acute et}}

 \newcommand{\fp}{{\bf F}_p}

\newcommand{\ext}{\mbox{Ext}}
\newcommand{\Ext}{\ext}

\newcommand{\spec}{\mathop{\rm Spec}\nolimits}

\newcommand{\ep}{\epsilon}

\newcommand{\Ker}{\mathop{\rm Ker}\nolimits}
\newcommand{\cok}{\mathop{\rm Cok}\nolimits}

\def\angles#1#2{\langle{#1},{#2}\rangle}
\newcommand{\ov}{\overline}

\newcommand{\gr}{\mathop{\rm Gr}\nolimits}

\newcommand{\ca}{{\cal A}}

\newcommand{\Hom}{\mathop{\rm Hom}\nolimits}
\newcommand{\id}{{\rm id}}
\newcommand{\lotimes}{\stackrel{\bf L}\otimes}

\newcommand{\End}{\mathop{\rm End}\nolimits}

\newcommand{\dlog}{\mathop{\rm dlog}}

\newcommand{\bA}{{\bf A}}

\newcommand{\Pic}{{\rm Pic}}
\newcommand{\bT}{{\bf T}}

\newcommand{\ot}{\otimes}

\newcommand{\nat}{\natural}
\def\odiagram#1{
  \def\normalbaselines{\baselineskip20pt\lineskip3pt \lineskiplimit3pt }
   \matrix{#1}}
 \def\ldiagram#1{$$\odiagram#1 % display math mode hack
      \refstepcounter{eqncounter}\
      \edef\@currentlabel{\p@equation
 \thetheorem.\theeqncounter}}
 \def\endldiagram{\eqno  \@eqnnum$$\global\@ignoretrue}
\def\ldiagram#1{$$\odiagram#1 % display math mode hack
      \refstepcounter{eqncounter}\
      \edef\@currentlabel{\p@equation
 \thetheorem.\theeqncounter}}
 \def\endldiagram{\eqno  \@eqnnum$$\global\@ignoretrue}

\newcommand{\J}{{\cal J}}
\newcommand{\str}[1]{\stackrel{#1}{\longrightarrow}}
\newcommand{\Cee}{{\cal I}}
\newcommand{\ccx}{{[\cx]}}

\newcommand{\micn}{{MICN}}
\newcommand{\cE}{\mathcal{E}}

\newcommand{\cExt}{\cE xt}
\newcommand{\cG}{G}
\newcommand{\ccG}{\mathcal{G}}
\newcommand{\cB}{\mathcal{B}}

\newcommand{\cX}{\mathcal{X}}

\newcommand{\tX}{{\tilde X}}
\newcommand{\tpi}{{\tilde \pi}}
\newcommand{\omegap}{\omega'}
\newcommand{\tY}{{\tilde Y}}
\newcommand{\tT}{{\tilde T}}
\newcommand{\tS}{{\tilde S}}
\newcommand{\cL}{\mathcal{L}}
\newcommand{\cD}{\mathcal{D}}
\newcommand{\cA}{\mathcal{A}}
\newcommand{\cT}{\mathcal{T}}

\newcommand{\cH}{\mathcal{H}}

\newcommand{\tF}{{\tilde F}}

\newcommand{\bH}{{\bf H}}
\renewcommand{\Pic}{{\bf P}{\rm ic}}
\newcommand{\bK}{{\bf K}}

\newcommand{\tU}{{\tilde U}}
\newcommand{\bM}{{\bf M}}
\newcommand{\bJ}{{\bf J}}

\newcommand{\Mod}{Mod}
\newcommand{\cEnd}{{\cal E}nd}
\newcommand{\cZ}{{\cal Z}}

\newcommand{\cI}{{\cal I}}
\newcommand{\cC}{{\cal C}}
\newcommand{\tf}{\tilde f}
\newcommand{\cY}{{\cal Y}}

\newcommand{\hig}{HIG}
\newcommand{\hign}{HIGN}
\newcommand{\mic}{MIC}
\newcommand{\gm}{{\mathbb{G}_m}}
\newcommand{\cnv}{\circledast}

\newcommand{\fhig}{\mbox{$F$-$HIG$}}

\newcommand{\ghig}{\mbox{$\cG$-$HIG$}}

\newcommand{\gchig}{\mbox{$\cG_\cx$-$HIG$}}

\newcommand{\bQ}{{\bf Q}}
\newcommand{\ind}{{\rm ind}}
\newcommand{\cS}{\mathcal{S}}
\newcommand{\tspl}{\mathcal{TSP}}
\newcommand{\spl}{\mathcal{SP}}

\newcommand{\cF}{\mathcal{F}}
\newcommand{\bL}{{\bf L}}
\newcommand{\bP}{{\bf P}}

\newcommand{\bv}{{\bf V}}
\def\lift#1#2{\mathcal{ #1} /\mathcal {#2}}
\newcommand{\ta}{{\tilde a}}
\newcommand{\tb}{{\tilde b}}
\newcommand{\tg}{{\tilde g}}

\newtheorem{theorem}{Theorem}[section]
\newtheorem{definition}[theorem]{Definition}

\newtheorem{corollary}[theorem]{Corollary}
\newtheorem{proposition}[theorem]{Proposition}
\newtheorem{lemma}[theorem]{Lemma}
\newtheorem{claim}[theorem]{Claim}
\newtheorem{notation}[theorem]{Notation}

\newremark{example}[theorem]{Example}
\newremark{remark}[theorem]{Remark}

\numberwithin{equation}{theorem}
\title{Nonabelian Hodge Theory in Characteristic $p$. }
\author{A. Ogus and V. Vologodsky}
\begin{document}
\maketitle
\begin{abstract}
Given a scheme in characteristic $p$ together with a lifting
modulo $p^2$, we construct a functor from a category
of suitably nilpotent modules with connection to the category
of Higgs modules.  We  use  this functor to
generalize the decomposition theorem of Deligne-Illusie to  the
case of de Rham cohomology with coefficients.
\end{abstract}
\tableofcontents

\section*{Introduction}
Let $X/\bc$ be a smooth projective scheme over the complex numbers
and let $X^{an}$ be the associated analytic space. Classical Hodge
theory provides a canonical isomorphism:
\begin{equation}\label{sh.e}
H^n(X,\Omega^\cx_{X/\bc}) \cong H^n(X^{an},\bc) \cong
\bigoplus_{i+j=n} H^i(X,\Omega^j_{X/\bc}).
\end{equation}
Carlos Simpson's ``nonabelian Hodge theory''\cite{si.hbls} provides
a generalization of this decomposition to the case of cohomology
with coefficients in a representation of the fundamental group of
$X^{an}$.  By the classical Riemann-Hilbert correspondence, such a
representation can be viewed as a locally free sheaf $E$ with
integrable connection $(E,\nabla)$ on $X$.  If $(E,\nabla)$
satisfies  suitable  conditions, Simpson associates to it a
\emph{Higgs bundle} $(E',\theta)$, \ie, a locally free sheaf $E'$
together with an $\oh X$-linear map
 $\theta \colon E'\to E'\ot \Omega^1_{X/\bc}$
such that $\theta\wedge \theta \colon E'\to E'\ot \Omega^2_{X/\bc}$
vanishes. This integrability implies that the iterates of $\theta$
are zero, so that $\theta$ fits into a complex (the {\em Higgs
complex})
$$E'\ot\Omega^\cx_{X/\bc} := E' \to E'\ot \Omega^1_{X/\bc} \to E'\ot \Omega^2_{X/\bc} \cdots .$$
As a substitute for the Hodge decomposition (\ref{sh.e}), Simpson
constructs a natural isomorphism:
\begin{equation}\label{cohc.e}
H^n(X, E\otimes \Omega^\cx_{X/,\bc},d)\cong  H^n(X^{an},V) \cong
H^n(X,E'\ot \Omega^\cx_{X/\bc},\theta),
\end{equation}
In general, there is no simple relation between $E$ and $E'$, and in
fact the correspondence $E \mapsto E'$ is not holomorphic.

Our goal in this work is to suggest and investigate an analog of
Simpson's theory for integrable connections in  positive
characteristics, as well as as an extension of the
paper~\cite{di.rdcdr} of Deligne and Illusie to the case of de Rham
cohomology  with coefficients in a $D$-module. Let $X$ be a smooth
scheme over a the spectrum $S$ of a perfect field $k$, and let $F
\colon X \to X'$ be the relative Frobenius map. Assume as in
\cite{di.rdcdr} that there is a lifting $\tX$ of $X'$ to $W_2(k)$.
Our main result is the construction of a functor $C_\tX$ (the
\emph{Cartier transform}) from the category $\mic(X/S)$ of modules
with integrable connection on $X$ to the category $\hig(X'/S)$ of
Higgs modules on $X'/S$, each subject to  suitable nilpotence
conditions.

The relative Frobenius morphism $F$ and   the $p$-curvature $$\psi
\colon E \to E\ot F^*\Omega^1_{X'}$$
 of a module with  integrable  connection $(E,\nabla)$
play a crucial role in the study of  connections in characteristic
$p$. A connection $\nabla$ on a sheaf of $\oh X$-modules $E$ can be
viewed as an action of the sheaf of PD-differential operators
\cite[(4.4)]{bo.ncc} \footnote{The name ``differential operators''
is misleading: although $D_X$ acts on $\oh X$, the map $D_X \to
\End(\oh X)$ is not injective.} $D_X$ on $X$.  This sheaf of rings
has a large center ${\cal Z}_{X}$: in fact, $F_*{\cal Z}_X$ is
canonically isomorphic to the sheaf of functions on the cotangent
bundle $\bT^*_{X^{\prime}}$:
\begin{equation} \label{pcurvd.e}
    c: S^{\cdot}T_{X^{\prime}}  \cong  F_{*}{\cal Z}_{X},
 \end{equation}
and $F_*D_X$ is an Azumaya algebra over $S^\cx T_{X'}$~
\cite{bmr.lmsla}.
 The map $c$ takes a vector field $\theta$ (\ie, a derivation of $\oh X$)
to  $ \theta^p - \theta^{(p)} \in D_X$, where $\theta^{(p)} \in
Der(\oh X)$ is the $p$th iterate of $\theta$ and $\theta^p$ is the
$p$th power of $\theta$ in $D_X$. If $\nabla$ is an integrable
connection on $E$, then by definition $\psi_\theta$ is the $\oh
X$-linear endomorphism  of $E$ given by the action of
$\nabla_{c(\theta)}$.

Let $\tilde X$ be a lifting of $X$. Our construction of the Cartier
transform $C_\tX$ is based on a study of the sheaf of liftings of
the relative Frobenius morphism $F \colon X \to X'$. For each open
subset $U \subseteq X$,
 the set  ${\cal L}_{\tilde X}(U)$ of all Frobenius liftings
$\tilde F: \tilde U \to \tilde U^{\prime}$ is naturally a torsor
under the group $F^*T_{X'}$.  Key to our construction is the fact
that the  $F^*T_{X'}$-torsor $q \colon \cL_\tX \to X$ has a
canonical connection
$$\nabla: {\cal L}_{\tilde X}\to F^*T_{X^{\prime}}\otimes  \Omega^1_{X}, $$
compatible with the Frobenius descent connection on the vector
bundle $F^*T_{X'}$. If $\tilde F$ is a local section of ${\cal
L}_{\tilde X}$ , $\nabla (\tilde F)  \in
\cHom(F^*\Omega^1_{X'},\Omega^1_X)$ is given by
$$\zeta_{\tilde F}: F^* \Omega^1_{X'} \to \Omega^1_X,$$
where $\zeta_\tF := p^{-1}d\tilde F$ is the lifting of the inverse
Cartier operator defined by $\tF$. Thus the sheaf of functions
${\cal A}_{\tilde X}: = q_* {\cal O}_{{\cal L}_{\tilde X}}$ acquires
a connection, as does its $\oh X$-linear dual  ${\cal B}_{\tilde
X}$. The torsor structure on   ${\cal L}_{\tilde X}$ induces an
action of the    completed PD symmetric algebra $\hat \Gamma_{\cdot}
F^*T_{X'}$ on $\cA_\tX$ and ${\cal B}_{\tilde X}$. We show that
 the induced action of $S^{\cdot} T_{X'}$ coincides with the action of the center
 $S^{\cdot} T_{X'}\subset D_X$ defined by the $p$-curvature of the connection $\nabla$.
 Thus  ${\cal B}_{\tilde X} $ becomes a module over the algebra
 $D_{X}^{\gamma} := D_X \otimes _{S^{\cdot} T_{X'}} \hat \Gamma_{\cdot} T_{X'}$.

We define the {\em Cartier transform}
 $C_{\tilde X}$ from the category of  $D_{X}^{\gamma}$-modules to the category
of $\hat \Gamma_{\cdot} T_{X'}$-modules by the formula:
$$C_{\tilde X}(E)= \iota_*\cHom_{D^\gamma_{X}}({\cal B}_{\tX}, E),$$
where $\iota$ is the involution of $T_{X'}$ sending $\theta'$ to
$-\theta'$. In fact, $  {\cal B}_{\tilde X} $ is  {\em a splitting
module} for the Azumaya algebra $D_{X}^{\gamma}$, and from this
point of view, the Cartier transform is just the equivalence of
categories between the category of modules over a split Azumaya
algebra and the category of ${\cal O}$-modules on the underlying
space defined by the choice of a splitting module. \footnote{The
role of the involution $\iota$ is to insure that  our constructions
are compatible with the standard Cartier operator and with the
decomposition of the de Rham complex constructed by Deligne and
Illusie~\cite{di.rdcdr}.}
  In particular,  the Cartier transform gives an equivalence of categories between the category
  $\mic_{p-1}(X)$ of nilpotent $D$-modules of level less then or equal to $p-1$
  and  the category $\hig_{p-1}(X')$ of Higgs modules supported on the $(p-1)^{st}$ infinitesimal
  neighborhood of the zero section $X' \hookrightarrow  \bT^*_{X'}$.
The larger categories of $D_X^\gamma$-modules and $\hat\Gamma_\cx
T_{X'}$-modules have the advantage of being  tensor categories, and
the Cartier transform is in fact compatible with the tensor
structures.

We also obtain an analog of Simpson's isomorphism (\ref{cohc.e}):
 if $(E',\theta')$ is the Cartier transform of a module with connection $(E,\nabla)$
whose level is less than the $p$ minus the dimension of $X$, then we
construct an  isomorphism in the derived category between the de
Rham complex of $(E,\nabla)$ and the Higgs complex of
$(E',\theta')$.  This result generalizes the decomposition theorem
of Deligne-Illusie~\cite{di.rdcdr}.

% The category $\hig_{p-1}(X'/S)$ has a canonical $\gm$-action, which, by
% means of the Cartier transform, carries over to $\mic_{p-1}(X/S)$.
% The derived category of $\gm$-equivariant Higgs modules has some nice
% homological properties; in particular we can define a Cartier transform
% on a {\em full} subcategory of the derived category of $\gm$-equivariant Higgs
% bundles.  This fact allows us to formulate and prove a strong result
% about the compatibility of the Cartier transform with derived direct images.
% As a consequence, we show that  the Cartier transforms of Gauss-Manin connections are canonically $\gm$-equivariant.
% The same is true of  the connections underlying $p$-torsion Fontaine-modules
% as defined by Faltings~\cite{fa.ccpgr}, and in fact our theory can be
% used to give a simple
% description~(\ref{fontmod.d}) of the category of such objects.

Let us describe the structure and content of the paper in more
detail. We work  with a smooth morphism $X/S$ of schemes in
characteristic $p$. We shall see that the Cartier transform depends
on a lifting $\tilde X'/\tilde S$ of $X'/S$ modulo $p^2$ rather than
a lifting of $X/S$, and we write $\lift X S$ for the pair
$(X/S,\tilde X'/\tilde S)$. In Theorem~\ref{fliftc.t} of section
\ref{lfrob.ss} we construct the torsor $\cL_{\lift X S}$ of liftings
of Frobenius and compute its connection  in
Proposition~\ref{nabexp.p} and $p$-curvature in
Proposition~\ref{psicomp.p}, using the geometric language of the
crystalline site and in particular Mochizuki's geometric description
of the $p$-curvature, which we recall in Proposition~\ref{moch.p}.

We also discuss in section~\ref{frm.ss} the relationship between
$\cA_{\lift X S}$ and some more familiar constructions in the
literature.

Section \ref{chf.s} is devoted to the construction of the Cartier
transform. We begin by reviewing in Theorem \ref{azum.t} the Azumaya
property of the algebra of differential operators and the canonical
fppf splitting module described in \cite{bmr.lmsla}. Then we discuss
the global Cartier transform $C_{\lift X S}$ as well as a local
version which depends on a lifting $\tF$ of the relative Frobenius
morphism $F_{X/S}$. Theorem~\ref{dzetas.t} constructs from such a
lifting $\tF$, or just the corresponding splitting $\zeta$ of the
inverse Cartier operator,  a surjective \'etale endomorphism
$\alpha_\zeta$ of $\bT^*_{X'}$ and a splitting module $\cB_\zeta$ of
$\alpha_\zeta^*D_{X/S}$. The restriction $\hat \cB_\zeta$  of
$\cB_\zeta$ to the formal completion of $\bT^*_{X'}$ along  its zero
section splits the ring $\hat D_{X/S}:= D_{X/S}\ot_{S^\cx{T_{X'/S}}}
\hat S^\cx T_{X'/S}$ of HPD differential operators, and this
splitting module defines an equivalence, which we call in
Theorem~\ref{loccart.t}
 the \emph{local Cartier transform},   between
the category of modules over $\hat D_{X/S}$ and the category of
modules over the ring $\hat S^\cx T_{X'/S}$.
 In fact,
$\hat \cB_\zeta$ is naturally isomorphic to the dual of the divided
power envelope of $\cA_{\tX}$ along the ideal of the section of
$\cL_{\lift X S}$ defined by $\tF$. This gives the  compatibility
between the local and global  Cartier transforms.

In Theorem \ref{rh.t} we explain how the Cartier transform can be
viewed as an analog of the  Riemann-Hilbert correspondence, with the
sheaf of $\oh X$-algebras  $\cA_{\lift X S}$ playing the role of
$\oh {X^{an}}$.  We also discuss a filtered version of the
construction, in which we study filtered  $D_{X/S}^\gamma$-modules
$(E,N_\cx)$, where
\begin{equation}\label{nsat.e}
(\Gamma_j T_{X'/S}) N_kE \subseteq N_{k-j}E
\end{equation}
for all $k$ and $j$. The algebra $\cA_{\lift X S}$ has a canonical
 filtration with this property,
and we show that the filtered object $C_{\lift X S}(E,N_\cx)$, can
be computed from the tensor product filtration on ${\cA_{\lift X
S}\ot E}$, which  again satisfies (\ref{nsat.e}).
  This construction will  become important in our  analog  Theorem~\ref{DR.t}
of the cohomological theorems of Simpson and Deligne-Illusie and in
particular to our study of  the ``conjugate filtration'' in
cohomology.

% In particular
%we show  in Theorem~(\ref{rh.t}) that the  Cartier transform
% $C_{\lift X S}$ can  be written as follows:  if $E$ is a module over $D_{X/S,\gamma}$,
%$$C_{\lift X S}(E) = (E \otimes \cA_{\lift X S})^{\nabla,\gamma} := \cHom_{D_{X/S,\gamma}}(\oh X, E\otimes \cA_{\lift X S}),$$
%with quasi-inverse
%$$C^{-1}_{\lift X S}(E') = (E'\otimes \cA_{\lift X S})^{\theta} :=
%\cHom_{\hat \Gamma_\cx T_{X'/S}}(\oh {X'},E'\otimes \cA_{\lift X S}),$$
%(after a change in the sign of $\psi$.)
%The proof uses classical Cartier descent and convolution in the category of F-Higgs sheaves
%instead of the theory of Azumaya algebras, and is modeled after the techniques introduced
%by Kato and Nakayama in logarithmic geometry~\cite{kkcn.lblelsc, o.lrhc}.

Section   \ref{fct.s} investigates the compatibility of the Cartier
transform with direct and inverse images with respect to a morphism
of smooth $S$-schemes $h \colon X \to Y$.
  We begin with a review
of the construction of the Gauss Manin connection on the relative de
Rham cohomology $R^qh_*(E\ot \Omega^\cx_{X/Y})$ when $h$ is smooth
and discuss its analog for Higgs fields.  Our review culminates with
Theorem~\ref{declev.t}, which shows that $R^qh_*$ increases the
level of nilpotence of a connection by at most the relative
dimension $d$ of $h$, strengthening the  result \cite[5.10]{ka.ncmt}
of Katz. In particular, we show that if $N_\cx$ is a filtration of
$E$ such that $\gr ^N E$ has zero $p$-curvature, then the filtration
of $R^qh_*(E \ot \Omega^\cx_{X/Y})$ induced by Deligne's
``filtration d\'ecal\'ee'' $N^{dec}$ of $E\ot \Omega^\cx_{X/Y}$ has
the same property. Theorem~\ref{push.t} shows that the Cartier
transform is compatible with  direct image by constructing, given a
lifting $\tilde h'$ of of $h' \colon X' \to Y'$,   an isomorphism
in $\hig(Y'/S)$
\begin{equation}\label{intropush.e}
  R^q h'^{\hig}_*C_{\lift X S} E\cong C_{\lift Y S} R^qh^{DR}_* E
\end{equation}
if the level of $E$ is less than $p-d$; we also show that this
construction is compatible with the filtrations $N^{dec}$.  This
result can be regarded as a relative version of the cohomology
comparison Theorem~\ref{DR.t}.

The remainder of section  \ref{fct.s} is devoted to  derived
versions of these results in a certain filtered derived category of
Higgs modules. The first important ingredient of this approach is a
new construction, described in Proposition~\ref{newder.p}, of the
functors $Lh_{DR}^*$ and $Rh_{DR*}$ in characteristic $p$, due to
Bezrukavnikov and Braverman  \cite{bb.glpc}, based on the Azumaya
property of the algebra $F_*D_{X/S}$. This construction allows us to
work locally over the cotangent bundle. Another ingredient is the
{\it conjugate filtration}
\begin{equation}\label{conjfiltration.e}
  \cdots \subset \Cee _X ^i \subset \cdots  \subset
\Cee _X ^1 \subset   F_{X/S *} D_{X/S},
\end{equation}
$$ \mbox{where} \quad  \Cee _X ^i = S^i T_{X'/S}(F_{X/S *} D_{X/S}) $$
and the concept of the {\it $\Cee$-filtered derived category}
 $DF(F_{X/S *} D_{X/S}, \Cee _X)$ of modules over the {\it filtered } algebra $F_{X/S *}D_{X/S}$.
Objects of this category  are filtered complexes $(E^{\cx},
N^{\cx})$ of $F_{X/S *}D_{X/S}$-modules such that for every integer
$i$
$$\Cee_X N^iE^{\cx}\subset N^{i+1} E^{\cx},$$
or equivalently, such that the associated graded module has
vanishing $p$-curvature.
 We lift the functors $Rh^{DR}_*$ and $Lh^*_{DR}$ to functors between the $\Cee$-filtered derived categories
 and prove in Proposition~\ref{i!} that for a smooth morphism $h: X\to Y$ of
 relative dimension $d$ the functor $Rh_*^{DR}$ increases the
range of the $\Cee$-filtration at most by $d$:
$$Rh_*^{DR}(DF_{[k,l]}(F_{X/S *} D_{X/S}, \Cee _X))\subset DF_{[k-d,l]}(F_{Y/S *} D_{Y/S}, \Cee _Y).$$
We also explain in Remark~\ref{kasp.r} how our formalism combined
with the  splitting property of the Azumaya algebra $Gr _{\Cee}
F_{X/S *}D_{X/S}$ leads to a  generalization of the Katz's
formula~\cite[Theorem~3.2]{ka.asde} relating the p-curvature to the
Kodaira-Spencer mapping.
 In section~\ref{dct.ss} we explain how the Cartier transform lifts to an equivalence of triangulated categories
$$C_{{\cal X}/S }: DF_{[k,l]}(F_{X/S *} D_{X/S}, \Cee _X)\cong DF_{[k,l]}(S^{\cdot}T_{X'/S}, \J_{X'}) , \quad \mbox{  when   } l-k< p.$$
 between the category $DF_{[k,l]}(F_{X/S *} D_{X/S}, \Cee _X)$ and the  {\it $\J$-filtered derived category}
  $ DF_{[k,l]}(S^{\cdot}T_{X'/S}, \J_{X'})$ of Higgs modules, where  $\J_{X'}\subset S^{\cdot} T_{X'/S}$ is
the ideal generated by $T_{X'/S}$. We then show in Theorem (3.20)
that, for a smooth morphism $h: X \to Y$,  a lifting $\tilde h'
:\tilde X' \to \tilde Y'$ induces a quasi-isomorphism
$$ C_{{\cal X}/S } \circ Rh_*^{DR} \cong  R{h'_*}^{HIG} \circ  C_{{\cal X}/S }, \quad\mbox{for}\quad l-k+d< p. $$
The exposition of sections \ref{ddii.ss}--\ref{dct.ss} does not
depend on  sections \ref{gm.ss}--\ref{ctdr.ss}, which obtain many of
the same results on the level of cohomology by more explicit
methods.

Section 4 is devoted to applications and examples. First we give a
characterization of the local \'etale essential image of the
$p$-curvature functor from the category $\mic(X/S)$ to the category
of $F$-Higgs sheaves. We show  in Theorem \ref{locpsi.t} that if $E$
is coherent and $\psi \colon E \to E\ot F_{X/S}^*\Omega^1_{X'/S}$ is
an $F$-Higgs field, then, \'etale locally on $X$, $(E,\psi)$ comes
from a connection if and only if, \'etale locally, $(E,\psi)$
descends to $X'$.  This can be regarded as a nonabelian analog of
the exact sequence~\cite[4.14]{mi.ec}
$$ 0  \rTo \oh {X'}^* \rTo^{F_{X/S}^*}  F_{X/S*}\oh X^*\rTo^{dlog} F_{X/S*} Z^1_{X/S}
\rTo^{\pi_{X/S}^* - C_{X/S}} \Omega^1_{X'/S} \rTo 0,$$ where
$C_{X/S}$ is the Cartier operator and $\pi \colon X' \to X$ is
$\id_X\times F_S$.
 Next in Theorem~\ref{gerbel.t} and  Proposition~\ref{brauer.p}
come a comparison  of   the gerbes of liftings  of $X'$ and of
splittings of $F_{X/S*}D^\gamma_{X/S}$ and
 a cohomological formula  for the class of  $F_{X/S*} D_{X/S}$ in
the Brauer group.   We prove in Theorem~\ref{abel.t} that if $X$ is
an abelian variety, then $F_{X/S*}D_{X/S}$ always splits over the
formal completion of the zero section of its cotangent bundle, and
in section \ref{counter.ss}  we construct an example of a liftable
surface for which $F_{X/S*}D_{X/S}$ does {\em not} have this
property. Section \ref{fm.ss} contains a discussion of $p$-torsion
Fontaine modules, especially as developed in \cite{fa.ccpgr} and
\cite{o.fgthd}, from the point of view of the Cartier transform.  As
an application, we give a reduction modulo $p$ proof of the
semistability of the Higgs bundles arising from Kodaira-Spencer
mappings. Finally,  in  section \ref{bk.ss}, we  show how our
nonabelian Hodge theory can be used to give a ``reduction modulo
$p$'' proof of a celebrated recent theorem of Barannikov and
Kontsevich, answering a question of Sabbah~\cite{sa.tdc}.

We conclude with an appendix devoted to generalities about Higgs
fields, and in particular to the study of the tensor product
structure on the category of Higgs modules.  This structure can be
viewed as a convolution with respect to the additive group law on
the cotangent space and  makes sense when restricted to the formal
and divided power completions of the zero section. The tensor
category of Higgs modules has an internal Hom, and an object $F$ of
$\hig(X)$ defines what we call a  ``Higgs transform'' $E \mapsto
\cHom_{HIG}(F,E)$ from the category of Higgs modules to itself.
 Our key technical result
is Proposition~\ref{charac.t}, which shows that the Higgs transform
with respect to a character sheaf on the cotangent space defines
(after  a change of sign) an involutive autoequivalence of tensor
categories.
%This construction seems to
%describe many formulas that appear in the literature in the context
%of differential equations, in addition to the current paper.
In the last part of the appendix we introduce, using $D_{X/S}$ as a
model, the notion of a tensor structure on an Azumaya algebra $\cA$
over a group scheme.  Such a structure makes the category of
$\cA$-modules a tensor category.

Both authors would like to express their gratitude to Roman
Bezrukavnikov. The second author would like to say that he learned
the main idea of this work from him: in particular, he explained
that the ring of differential operators in characteristic p is an
Azumaya algebra over the cotangent bundle and suggested that it
might split over a suitable infinitesimal neighborhood of the zero
section. The first author was blocked from realizing his vision
(based on \cite{o.hcpc}) of a nonabelian Hodge theory in positive
characteristics until he learned of this insight. Numerous
conversations with  Roman also helped us to overcome many of the
technical and conceptual difficulties we encountered in the course
of the work. The authors also benefited greatly from Pierre
Berthelot, who in particular explained to the first author years ago
how a lifting of Frobenius makes $\hat D_{X/S}$ into a matrix
algebra.  Special thanks go to the referee who pointed out a mistake
in an early draft as well as a simplification in our argument which
offered a way around it. This led us to the realization that we
could greatly strengthen one of our main results and allowed us to
develop the filtered Cartier transform in the context of cohomology
and derived categories. We are also extremely grateful to the
referee for pointing out an enormous number of misprints and
ambiguities in an early draft. We would also like to thank Alexander
Beilinson,  Alexander Braverman, Luc Illusie, and Ofer Gabber for
the interest they showed and  the advice they provided. Finally, we
would like to alert the reader to a forthcoming work by Daniel
Schepler which extends this theory to log geometry. \footnote{Both
authors would like to acknowledge the support this collaboration
received from the Committee on Research at the University of
California at Berkeley. The second author was partially supported by
NSF grant DMS-0401164, but support for the team effort was denied by
the National Science Foundation and the Miller Institute for Basic
Research.  }

\section{The torsor of Frobenius liftings}\label{ffe.s}
\subsection{Liftings of Frobenius}\label{lfrob.ss}
If $X$ is a scheme in characteristic $p$, let $F_X$ denote its
absolute Frobenius endomorphism, \ie, the map which is the identity
on the underlying topological space and which takes each section
of $\oh X$ to its $p$th power.
  For any morphism  $f \colon X \to S$ of schemes in characteristic $p$, $F_S\circ f
= f\circ F_X$, and one has the  relative Frobenius diagram:
\begin{diagram}
X & \rTo^{F_{X/S}}&  X^{(S)}  &\rTo^{\pi_{X/S}} & X \cr
  &\rdTo^f & \dTo_{f^{(S)}} & & \dTo_f \cr
   && S &\rTo^{F_S} &S.
\end{diagram}
The square in this diagram is  Cartesian, and the map $F_{X/S}$ is
the unique morphism over $S$ such that $\pi_{X/S}\circ F_{X/S} =
F_X$. If no confusion seems likely to result, we may simplify the
notation, writing $X'$ for $X^{(S)}$, $F$ for $F_{X/S}$, etc.  We
also often write $X/S$ for the morphism $f \colon X \to S$, viewed
as an $S$-scheme.

If $f \colon X \to S$ is any morphism of schemes in characteristic
$p >0$ and  $n$ is a positive integer, by a \emph{lifting of $f$
modulo $p^n$} we shall mean a morphism ${\tilde f \colon \tilde X
\to \tilde S}$ of flat $\bz/p^n\bz$-schemes, together with a
Cartesian diagram
\begin{diagram}
X & \rTo & \tilde X \cr \dTo^f && \dTo \tilde f \cr S &\rTo &\tilde
S,
\end{diagram}
where $S \to \tilde S$ is the closed subscheme defined by $p$.
Note that if $\tilde X/\tilde S$ is such a lifting and $X/S$ is flat
(resp. smooth), then so is $\tilde X/\tilde S$.  We shall be
primarily interested in the case $n = 2$, and if $n$ is not
specified,  this is what we shall mean.  If the absolute Frobenius
endomorphism $F_S$ lifts to $\tS$, then $\tilde f' \colon \tilde
X\times_{F_{\tilde S}} \tilde S \to \tilde S$ lifts $X'/S$.  For
example, if $S$ is the spectrum of a perfect field $k$ and $\tS$ the
spectrum of its truncated Witt ring, then there is a unique such
$F_\tS$, but in general there is no reason for a lifting of $F_\tS$
or of $X'$ to exist, even locally on $S$, unless $S$ is smooth over
a perfect field.

Throughout the rest of this section, let us  fix a smooth $X/S$ as
above.  We assume that a lifting $\tilde X'/\tilde S$ of $X'/S$
modulo $p^2$ exists, and we denote the pair $(X/S,\tilde X'/\tilde
S)$ by $\lift X S$. Note that, given a lifting $\tX$ of $X$, it is
very rare for there to exist a  global lifting of $F_{X/S} \colon
\tilde X \to \tilde X'$. (For example, no such lift can exist if $X$
is a smooth proper curve of genus at least two over a perfect field,
as is well known.) However it follows from the smoothness of $X'/S$
that such lifts do exist locally, and we shall see that the sheaf of
such liftings is crystalline in nature.

Let us fix a divided power structure on the ideal $p\oh {\tilde S}$
and consider the crystalline site $Cris(X/\tilde S)$.
 If $(U,\tT)$ is an object of $Cris(X/\tS)$, let $T$
be the reduction of $\tT$ modulo $p$.  The ideal $J_T$ of  the
inclusion $i \colon U \to T$ is a divided power ideal, and so $a^p =
0$ for every local section $a $ of $J_T$.  Then the relative
Frobenius map $F_{T/S}$ factors through $U'$, and there is a unique
and canonical morphism $f_{T/S} \colon T \to X'$ such the following
diagram commutes.
\begin{equation}\label{fst.d}
\begin{diagram}
& & T &\rTo^{F_{T/S}}& T' \cr &\ldTo^{f_{T/S}}& \dTo_{f'_{T/S}} &
\ruTo_{i'} \cr X' &\lTo^{inc'} & U'
\end{diagram}
\end{equation}

Let us note for future reference that the differential of $f_{T/S}$
vanishes:
\begin{equation}\label{diff.e}
0 = df_{T/S} \colon\Omega^1_{X'/S} \to f_{T/S*}\Omega^1_{T/S}
\end{equation}
Indeed,  $df'_{T/S}\circ di' = dF_{T/S} = 0$,
and since $di' $ %\colon {i'}^*\Omega^1_{T'/S} \to \Omega^1_{U'/S}$,
is an epimorphism, $df'_{T/S} = 0$.

%%%%%%%%   Cut Here%%%%%%

If  $g \colon T_1 \to T_2$ is a morphism in $Cris(X/S)$,  then $
 f_{T_2/S}\circ g=
 f_{T_1/S}$.  Hence if $E'$ is  a sheaf of $\oh {X'}$-modules,
 there is a natural isomorphism
 $$\theta_g \colon g^* f_{T_2/S}^* E' \cong f_{T_1/S}^* E',$$
 and the  collection
 $\{f_{T/S}^*E', \theta_g \}$ defines a crystal
 of $\oh {X/S}$-modules.  The corresponding object of $\mic(X/S)$
 is $F_{X/S}^*E'$ with its Frobenius descent connection.
(This is the unique connection $\nabla$ on $F_{X/S}^*E'$
 which annihilates the sections of $F_{X/S}^{-1}E' \subseteq
F_{X/S}^*E'$.)

An extension of crystals
 \begin{equation}
 \label{exttors.e}
 0 \to E \to H \stackrel{h}{\to} \oh X\to 0
 \end{equation}
 gives rise to  a sheaf $h^{-1}(1)\subset H$ of $E$-torsors on
 $Cris(X/S)$;   this construction defines an equivalence
 between the category of $E$-torsors
 and the category of extensions (\ref{exttors.e}).
Recall that giving a crystal $E$ amounts to giving a  quasi-coherent
sheaf of $\oh X$-modules with an integrable connection
 $\nabla_E: E\to E \otimes \Omega^1_{X/S}$.
 Similarly, giving an $E$-torsor ${\cal L}$ on Cris(X/S) is equivalent to
giving an $E$-torsor ${\cal L}$ on Zar(X) together with a map
 $$\nabla _{\cal L}: {\cal L}\to E \otimes \Omega ^1_{X/S}$$ such that
 $\nabla_{\cal L} (l + e) = \nabla _{\cal L}(l) + \nabla _E (e)$ and
 such that the composition
 $${\cal L}\rTo^{\nabla _{\cal L}} E \otimes \Omega ^1_{X/S} \rTo{\nabla_E} E \otimes \Omega ^2 _{X/S}$$
 is equal to zero.

If $E$ is a locally free crystal of $\oh {X/S}$-modules, we shall
denote by ${\bf
 E}$ the corresponding {\it crystal of affine group schemes} over
 $Cris(X/S)$. That is, for each $T \in Cris(X/S)$,
$${\bf E}_T:= \spec_T S^\cx \Omega' _T,$$
where $\Omega'$ is the crystal of $\oh {X/S}$-modules dual to $E$.
In particular, a vector bundle $E'$ over $X'$ defines a crystal of
affine schemes $F^*_{X/S}{\bf E'}$.
 More generally, for an $E$-torsor ${\cal
 L}$  on $Cris(X/S)$, we denote by ${\bf
 L}$ the corresponding crystal of affine schemes, which
has a natural action
 ${\bf E}\times {\bf L}\to {\bf L}$.

Now let us fix a pair $\lift X S := (X/S,\tX'/\tS)$ as above. By a
\emph{lifting of $f_{T/S}$ to $\tT$} we shall mean a morphism $\tF
\colon \tT \to \tX'$ lifting $f_{T/S}$. The sets  of such liftings
on open subsets of $\tT$ form a sheaf $\cL_{\lift X S,\tT}$ on the
Zariski topology of $\tT$ (which coincides with the Zariski topology
of $T$). Since $\tilde X'/\tilde S$ is smooth, such liftings exist
locally,
 and by standard deformation theory,
the sheaf $\cL_{\lift X S,\tT}$  of such liftings forms a torsor
under the abelian sheaf $\cHom(f_{T/S}^*\Omega^1_{X'/S},p\oh{\tilde
T}) \cong f_{T/S}^*(T_{X'/S})$.
% (Warning:  in order for ``standard deformation theory''  to agree with
% ``crystalline deformation theory,'' it is necessary that the divided power
% structure on the ideal $(p/p^2)$ be the trivial one.  If $p$ is odd, this
% is true for the standard divided power structure, but this is not so if $p = 2$.)

\begin{theorem}\label{fliftc.t}
Let $\lift X S: = (X/S,\tilde X'/\tS) $ be as above.  Then there is
a unique crystal of $F_{X/S}^* T_{X'/S}$-torsors $\cL_{\lift X S}$
on $X/S$ with the following properties.
\begin{enumerate}
\item{For each
object $T$ of $X/S$ admitting a flat lifting $\tT \in Cris(X/\tS)$,
$\cL_{\lift X S,T}$ is the sheaf of liftings of $f_{T/S}$ to $\tT$.}
\item{For each morphism $\tilde g \colon \tT_1 \to \tT_2$ of flat objects
in $Cris(X/\tS)$ and  each  lifting $\tF \colon \tT_2 \to \tX'$ of
$f_{T_2/S}$, the transition map $\theta_g \colon g^*\cL_{\lift X
S,T_2} \to \cL_{\lift X S,T_1}$ satisfies
$$\theta_g(\tF) = \tF\circ g \colon \tT_1 \to \tX'.$$}
\end{enumerate}
We denote by ${\bf L_{\lift X S}} $ the crystal of affine schemes
$\spec\cA_{\lift X S}$ corresponding to the
$F_{X/S}^*T_{X'/S}$-torsor $\cL_{\lift X S}$; thus  $\cA_{\lift X
S}$ is a crystal of quasi-coherent $\oh {X/S}$-algebras.
%For any $\tT \in Cris_f(X/\tS)$ with reduction $T$, $\cA_{\lift X S,T}$
%can be identified with the set of morphisms
%$\alpha \colon {\cL_{\lift X S}}_{|_\tT} \to {\oh{X/S}}_{|_\tT}$.
\end{theorem}
\begin{remark}\label{mapeq.r}
We should point out that if $\tilde T_1$ and $\tilde T_2$ are two
flat liftings of an object $T$ of $Cris(X/\tilde S)$, then the set
of liftings of $f_{T/S}$ to $\tilde T_1$ and to $\tilde T_2$ can be
canonically identified.
More precisely, %\begin{corollary}\label{mapeq.c}
let $\tT_1$ and $\tT_2$ be flat objects of $Cris(X/\tS)$, and let
$\tilde g$ and $\tilde g'$ be  two morphisms $\tT_1 \to \tT_2$ with
the same reduction modulo $p$.  Then $\cL_{\lift X S }(\tilde g) =
\cL_{\lift X S}(\tilde g')$ as maps $\cL_{\lift X S}(\tilde T_2) \to
\cL_{\lift X S}(\tilde T_1) $. This follows from the proof of the
theorem, but it can also be deduced from the following elementary
argument. Let $g\colon T_1 \to T_2$ be the common reduction modulo
$p$ of $g_1$ and $g_2$. Then there is a  map $h \colon
\Omega^1_{T_2/S} \to g_*\oh {T_1}$ such that ${\tilde g}^{\prime
*}(\tilde a) = \tilde g^*(\tilde a) + [p] h(da)$ for every section
$\tilde a$ of $\oh{\tilde T_2}$ lifting a section $a$ of $\oh
{T_2}$.  Then if $\tF \in \cL_{\lift X S}(\tilde T_2)$  is any lift
of $f_{T_2/S}$ and $\tilde b$ is a section of $\oh{\tilde X'}$ with
image $b$ in $\oh {X'}$,
$$(\tF \circ \tilde g')^*(\tilde b) =
 (\tF \circ \tilde g )^*(\tilde b) + [p] h (df_{T_2/S}(db)).$$
%where
%$$d\tF \colon f_{T_2/S}^*\Omega^1_{X'/S} \to \Omega^1_{T_2/S} $$
%is  the differential of $f_{T_2/S}$.
But we saw in (\ref{diff.e}) that  $df_{T_2/S} = 0$, hence $\tF\circ
\tilde g = \tF\circ \tilde g'$.
\end{remark}

\begin{proof}[Proof of Theorem~\ref{fliftc.t}]
We will need the following easy technical result.
\begin{lemma}\label{crisf.l}
Let $Cris_f(X/\tS)$ denote the full subsite of $Cris(X/\tS)$
consisting of those objects which are flat over $\tS$. Then the
morphism of sites $a \colon Cris_f(X/\tS) \to Cris(X/\tS)$ induces
an equivalence between the respective categories of crystals of $\oh
{X/\tS}$-modules.
\end{lemma}
\begin{proof}
 Indeed, the question is local on $X$,
so we may assume the existence of a lifting $\tX/\tS$.  Then both
categories can be identified with the category of pairs
$(E,\epsilon)$, where $E$ is a quasi-coherent $\oh \tX$-module and
$\epsilon$ is an isomorphism between the two pullbacks of $E$ to the
divided power completion of $\tX\times \tX$ along the diagonal,
satisfying the cocycle condition~\cite[\S6]{bo.ncc}.
%  This is because a crystal on
% $Cris(X/\tS)$ is determined by specifying its value on a
% covering $\tT$  of the final object and descent data
% on $\tT \times \tT$, each of which in fact belongs to $Cris_f(X/\tS)$.
\end{proof}
Thus we can identify the category of  crystals of $\oh
{X/S}$-modules on $Cris(X/S)$ and the category of $p$-torsion
crystals of  $\oh {X/\tS}$-modules on $Cris_f(X/\tS)$. The same is
true for torsors over  crystals of $\oh {X/S}$-modules.

It is clear that  the family $\{ \cL_{\lift X S,\tT} : \tT \in
Cris_f(X/\tS)\}$, together with the family of transition maps
$\theta_g$ described in the theorem, forms a sheaf of sets on
$Cris_f(X/\tS)$.  Furthermore, as we saw above, this family
naturally forms a sheaf of $F_{X/S}^*T_{X'/S}$-torsors. This proves
the theorem.

% The crystal $F_{X/S}^*\bT_{X'/S}$ is, as a sheaf of sets, represented
% by the $\oh {X/S}$-algebra $F^*_{X'/S}S^\cx \Omega^1_{X'/S}$, together
% with its Frobenius descent connection.  Since $\cL_{\lift X S}$
% is locally isomorphic to  $F_{X/S}^*\bT_{X'/S}$, it too is affine,
% represented by a quasi-coherent crystal
%  of $\oh {X/S}$-algebras $\cA_{\lift X S}$.
\end{proof}

Let us record some basic facts about vector groups which we will
need later. Let $\pi_T \colon \bT \to X$ be a vector group over $X$
and let $T$ be its sheaf of sections. Thus $T$ is a locally free
sheaf of $\oh X$-modules of finite rank
 and $\bT = \spec_X S^\cx \Omega$,
where $\Omega$ is the dual of $T$. The pairing $T \times \Omega \to
\oh X$ extends to a pairing $T \times S^\cx \Omega \to \oh X$, where
sections of $T$ act as derivations of $S^\cx \Omega$.  This action
defines a map:
$$\xi  \mapsto D_\xi \colon T  \to  \pi_{T*}T_{\bT/X},$$
which identifies $T$ with the sheaf of translation invariant vector
fields of $\bT$ relative to $X$.  It also induces an isomorphism
$\pi_T^*T \to T_{\bT/X}$. Moreover, there is a canonical  pairing of
$\oh X$-modules:
$$\Gamma_n T \otimes S^{n+m} \Omega \to S^m\Omega$$
which is perfect when $m= 0$; see (section \ref{exfor.ss})  and
\cite[A10]{bo.ncc}. If we endow  $\Gamma_\cx T$ with the topology
defined by the PD-filtration of $\Gamma_\cx T$ and $S^\cx \Omega$
with the discrete topology, this action is continuous. Thus it
extends to a continuous  action
 of the completion $\hat\Gamma_\cx T$
and  identifies $\cHom_{\oh X}(S^\cx \Omega, \oh X)$ with the
completed divided power algebra $\hat \Gamma_\cx T$  of $T$.
\footnote{Thus the Cartier dual of $\bT$ is the formal scheme
$\hat\bT^*_\gamma $ associated to the $PD$-algebra $\hat \Gamma_\cx
T$ with the topology defined by the divided power filtration $\{
\prod_{j\ge n} \Gamma_j T : n \in \bn \}$.  } This action identifies
the sheaf of divided power algebras $\Gamma_\cx T$~\cite{roby.lplf}
with subring of translation invariant elements
 in the {\em full} ring  of differential operators
\cite[2.1]{bo.ncc} of $\bT$ relative to $X$.

A section $\xi $ of $T$ can be thought of as a section of the map
$\pi_T \colon \bT \to X$; let $t_\xi \colon \bT \to \bT$ be
translation by $\xi$. Then  the derivation $D_\xi$ belongs to the
divided power ideal of $\Gamma_\cx T$,  $exp({D_\xi})$ makes sense
as a differential operator of infinite order, and one has the
formula (Taylor's theorem):
\begin{equation}\label{taylor.e}
t_\xi^*(f) = (\exp { D_\xi})(f).
\end{equation}
for the action of $t_\xi^*$ on $S^\cx \Omega$. The increasing
filtration
$$N_n S^\cx\Omega := \sum_{i\le n}S^i \Omega \subseteq S^\cx \Omega$$
is invariant under   $t^*_\xi$; furthermore $t^*_\xi$ acts trivially
on the successive quotients.

Now let $\pi_\cL\colon \cL \to X$ be a $\bT$-torsor over $X$. It
follows from the translation invariance of $D_\xi$ that  the action
of $\hat \Gamma_\cx T$ on $S^\cx \Omega$ carries over to an action
on $\pi_{\cL*}\oh \cL$.  Similarly, there is a canonical  filtration
$N_\cx$ on $\pi_{\cL*}\oh \cL$ and  a canonical isomorphism
\begin{equation}\label{grn.e}
\gr^N_\cx (\pi_{\cL*}\oh \cL) \cong  S^\cx\Omega
\end{equation}
Note that $N_i \pi_{\cL*}\oh \cL$ can also be characterized at the
annihilator of $\prod_{j>i} \Gamma_j T$. The bottom level   of the
filtration $N$ corresponds to the translation invariant sections, so
there is a canonical exact sequence
$$ 0 \to \oh X \to \cE \to \Omega \to 0,$$
where $\cE := N_1 \pi_{\cL*}\oh \cL$ is the set of affine functions
on $\cL$.
%The inclusion $\oh X \to \cE$
%defines maps
%$S^{n-1}\cE \to S^n\cE$ for all $n$, and the inclusion
%$\cE \to \pi_{\cL*} \oh \cL$ defines maps
%$S^n \cE \to  q_*\oh \cL$.
%One verifies immediately that this map is injective with image $N_n \cA_{\lift X S}$
%and that one gets in the limit an isomorphism $\dirlim S^n \cE \cong
%\pi{\cL_*} \oh \cL$.

A section $\ell$  of $\cL$ determines an isomorphism $s_\ell \colon
\cL \to \bT$: $s_\ell(\ell'):= \ell' - \ell \in \bT$ for all
sections $\ell'$ over all $X$-schemes $T$. This isomorphism
determines an isomorphism
$$\sigma_\ell:= s_\ell^* \colon S^\cx \Omega \to \pi_{\cL*} \oh \cL.$$
This is the unique isomorphism of filtered $\oh X$-algebras with the
property that $\sigma_\ell(\omega)(\ell') =
\angles{\omega}{\ell'-\ell}$ for all  local sections $\ell'$ of
$\cL$ over $X$ and $\omega$ of $\Omega$. (The uniqueness comes from
the fact that any polynomial $\alpha \in A[t_1,\ldots t_d]$ of
degree less than or equal to $1$ is determined  by its values on all
$A$-valued points.)  Note in particular that, as a $\hat \Gamma_\cx
T$-modules, $\pi_{\cL*} \oh \cL$ is \emph{locally coinvertible},
\ie, its $\oh X$-linear dual is, locally on $X$, free of rank one
over $\hat \Gamma_\cx T$.

Finally, let us remark that if $T \to T'$ is an $\oh X$-linear map
of locally free sheaves, and $\cL'$ is the $\bT'$-torsor deduced
from $\cL$ by pushout, then the morphism $\cL \to \cL'$ induces an
isomorphism
\begin{equation}%\label{torspush.e}
\pi_{\cL'*}\oh {\cL'} \cong \cHom_{\hat \Gamma_\cx T}(\hat\Gamma_\cx
T',\pi_{\cL*} \oh \cL)
\end{equation}

Let us summarize these remarks for our crystal of torsors
$\cL_{\lift X S}$.

\begin{proposition}\label{torsor.p}
Let $\lift X S := (X,\tilde X')$ and $\cL_{\lift X S}$ be as above,
and let $\cA_{\lift X S}$ denote the corresponding crystal of $\oh
X$-algebras.
\begin{enumerate}
\item{There is a natural horizontal action of $\hat\Gamma_\cx F^*_{X/S}T_{X'/S}$
on $\cA_{\lift X S}$, compatible with the action  of
$F^*_{X/S}T_{X'/S}$ by translation, as described in formula
(\ref{taylor.e}) above. As a sheaf of $\hat\Gamma
F^*_{X/S}(T_{X'/S})$-modules on $X$, $\cA_{\lift X S}$ is locally
coinvertible.}
\item{There is a natural horizontal filtration $N_\cx$ on $\cA_{\lift X S}$,
invariant under the action  of $F^*_{X/S}T_{X'/S}$.   In fact
$N_i\cA_{\lift X S}$ is the annihilator of $\prod_{j>i}\Gamma_j
F^*_{X/S}T_{X'/S}$, and there is a canonical isomorphism:
$$\gr^N_\cx\cA_{\lift X S} \cong F_{X/S}^*S^\cx \Omega^1_{X'/S}.$$}
\item{Let $\tT$ be a flat object of $Cris(X/\tS)$ and let $\tF \colon \tT \to \tX'$
be a lift of $f_{T/S}$.  Then there is a unique isomorphism  of
(filtered) $\oh T$-algebras
$$\sigma_\tF \colon f_{T/S}^* S^\cx \Omega^1_{X'/S} \rTo^\cong \cA_{\lift X S,T}.$$
with the  following property.  For every section $\tilde a'$ of $\oh
{\tX'}$ lifting a section $a'$ of $\oh {X'}$,
$\sigma_\tF(f_{T/S}^*da') \in N_1\cA_{\lift X S,T}$ is the  $\oh
T$-valued  function on ${\cL_{\lift X S}}(\cT)$ such that for each
$\tF'$,
$$[p]\sigma_\tF(f_{X/S}^*da')(\tF') =  \tF^{\prime *}(\tilde a') - \tF^* (\tilde a'). $$
Furthermore $\gr_N \sigma_\tF$ is the identity.}
\end{enumerate}
\end{proposition} \qed

In particular we have a fundamental exact sequence: \footnote{The
first explicit construction of this sequence was given in
\cite{srv.ddrc}.}
\begin{equation}\label{nexc.e}
0 \to \oh X \to \cE_{\lift X S} \to F_{X/S}^*\Omega^1_{X'/S} \to 0,
\end{equation}
where
$$ \cE_{\lift X S}:= N_1\cA_{\lift X S}.$$
A section $\tF$ of $\cL_{\lift X S}$ determines as above a
homomorphism $\sigma_\tF$ which induces a splitting (not compatible
with the connections) of this sequence, and in fact the set of
splittings is bijective with the set of sections.

Since $\pi_\cL \colon \cL_{\lift X S}$ is an
$F_{X/S}^*T_{X'/S}$-torsor over $X$, there is a natural
identification $\Omega^1_{\cL/X} \cong
\pi_{\cL}^*F^*_{X/S}\Omega^1_{X'/S}$.  The following result is the
key to our theory; it shows that the $p$-curvature of the connection
on $\cA_{\lift X S}$ is very rich.

\begin{proposition}\label{psicomp.p}
The action described in part (1) of Proposition~\ref{torsor.p} of
$F_{X/S}^*T_{X'/S}\subseteq \hat\Gamma F_{X/S}^*T_{X'/S}$ on
$\cA_{\lift X S}$ is the same as the action given by the
$p$-curvature $\psi$ of the connection $\nabla$ on $\cA_{\lift X
S}$.  That is, the  diagram
\begin{diagram}
\cA_{\lift X S} & \rTo^{\psi} &  \cA_{\lift X S}\ot
F_{X/S}^*\Omega^1_{X'/S} \cr
 \dTo^\cong && \dTo_\cong \cr
\pi_{\cL*}\oh {\cL }  & \rTo^d &\pi_{\cL*} \Omega^1_{\cL/X},
\end{diagram}
where $d$ is the usual exterior derivative and $\psi$ is the
$p$-curvature of the connection on $\cA_{\lift X S}$, is
commutative.
\end{proposition}

This formula can be proved by explicit calculation (see
Remark~\ref{ez.r} below). We prefer to give here a conceptual proof
based on a geometric construction of the $p$-curvature due to
Mochizuki and communicated to us by Brian Osserman; see
\cite{oss.clex}. This construction begins with  the following
crystalline interpretation of $F_{X/S}^*\Omega^1_{X'/S}$.

\begin{proposition}\label{xip.p}
Let $X/S$ be a smooth morphism of schemes in characteristic $p$, let
$X(1) := X\times_S X$, and let $(D(1), \overline I,\gamma)$ denote
the divided power envelope of
 the ideal $I$ of the diagonal immersion $X \to X(1)$.
Then there is a unique and functorial isomorphism
$$\xi_p \colon F_{X/S}^*\Omega^1_{X'/S} \to \ov I/ (\ov I^{[p+1]}+ I\oh {D(1)})$$
such that, for every local section $a$ of $\oh X$,
$$d_p(a) :=\xi_p(d\pi^*(a)) = ((1\otimes a) - (a\otimes 1))^{[p]} \quad \pmod{\ov I^{[p+1]}+ I\oh {D(1)}}.$$
\end{proposition}
\begin{proof}
For each section $a$ of $\oh X$, let $\xi(a) := 1\otimes a - a
\otimes 1 \in I\oh {D(1)} \subseteq \ov I$. Note that $\xi(a)$
annihilates $\ov I/(\ov I^{[p+1]}+ I\oh {D(1)})$, and hence that the
actions of $(a\ot 1)$ and of $(1\ot a)$  on $\ov I/ (\ov I^{[p+1]}+
I\oh {D(1)})$ are the same. Thus this quotient can be viewed as a
sheaf of $\oh X$-modules. If $b$ is another section of $\oh X$, then
\begin{eqnarray*}
\xi(a+b)^{[p]} & = &(\xi(a)+ \xi(b))^{[p]}  = \xi(a)^{[p]} + \sum_{i=1}^{p-1}\xi(a)^{[i]}\xi(b)^{[p-i]} + \xi(b)^{[p]} \\
    & = &  \xi(a)^{[p]}+ \xi(b)^{[p]} \pmod {I\oh {D(1)}}.
\end{eqnarray*}
Furthermore, $\xi(ab) = (1\ot a) \xi(b) + (b\ot 1)\xi(a)$, so a
similar calculation show that
%\begin{eqnarray*}
%(1\ot a^p) \xi(b)^{[p]}  + \sum_{i=1}^{p-1}(a^i\ot 1)\xi(b)^{[i]}(b^{p-i}\ot 1)\xi(a)^{[p-i]} + (b^p\ot 1)\xi(a)^{[p]} \\
%    & = & a^p \xi(b)^{[p]} + b^p \xi(a)^{[p]} \pmod {I\oh {D(1)}}
%\end{eqnarray*}
$$\xi(ab)^{[p]} =  a^p \xi(b)^{[p]} + b^p \xi(a)^{[p]} \pmod {I\oh {D(1)}}.$$
Finally, if $a$ is a local section of $f^{-1}(\oh S)$, $\xi(a) = 0$.
These properties imply that $d_p$ is a derivation $\oh X \to
F_{X*}\left (\ov I/ (\ov I^{[p+1]}+ I\oh {D(1)})\right)$, and hence
that $d_p$ factors through an $\oh X$-linear $\xi_p$ as claimed.  To
see that $\xi_p$ is an isomorphism, we may work with the aid of a
system of local coordinates $t_1,\ldots t_m$ for $X/S$. Let $\xi_i
:= \xi(t_i)$, so that, using multi-index notation, $\{ \xi^{[I]} : I
\in \bn^m \}$ forms a basis for $h_{1*}\oh {D(1)}$.
  Note that $\overline I \subseteq \overline I^{[p]}  + I\oh {D(1)}$ and
that $\xi^{[I]} \in I\oh {D(1)}$ if any $I_j < p$.  It follows that
$ \ov I/ (\ov I^{[p+1]}+ I\oh {D(1)})$ is freely generated by
$\xi_1^{[p]},\ldots \xi_m^{[p]}$, and hence that $\xi_p$ is an
isomorphism.
\end{proof}

\begin{proposition}[Mochizuki]\label{moch.p}
Let  $E$ be a crystal of $\oh X$-modules on $X/S$, Let $h_1$ and
$h_2$ be the canonical maps $D(1) \to X$, and let $\epsilon \colon
h_2^*E \to h_1^*E$ be the canonical isomorphism.   Then the
$p$-curvature $\psi$ of $E$ identifies, via the isomorphism $\xi_p$
of Proposition~\ref{xip.p}, with the map sending each local section
$e$ of $E_X$ to the class of $\epsilon(h_2^*(e)) - h_1^*(e)$ in
$\overline I/(\overline I^{[p+1]} + I\oh{D(1)})\ot E$.
\end{proposition}
\begin{proof}
%
%As explained in the introduction, if $(E,\nabla)$
%is an integrable connection on $X/S$,
%then for each derivation $D \in T_{X/S}$,
%its $p$th iterate $D^{(p)}$ is again a derivation,
%and one verifies easily that $c(D) :=(\nabla_D)^p  - \nabla_{D^{(p)}}$
%is an $\oh X$-linear endomorphism of $E$.  Furthermore,
% there is a unique $\oh X$-linear map
%$\psi\colon E \to E\otimes F^*_{X/S}\Omega^1_{X'/S}$,
%the \emph{$p$-curvature} of $(E,\nabla)$,
%such that  for each section  $D$ of $T_{X/S}$, the
%following diagram commutes~\cite{ka.ncmt}:
%\begin{diagram}
%E &\rTo^\psi& E\otimes F^*_{X/S}\Omega^1_{X'/S} \cr
% &\rdTo_{\nabla_{c(D)}} & \dTo_{\id\ot \angles {\pi^*D}\ }\cr
%   & & E
%\end{diagram}
%
We verify this formula with the aid of a system of local coordinates
$(t_1,\ldots t_m)$, using the notation above. Then if $D_i :=
\partial/\partial t_i$,
$$\ep (h_2^*(e))   = \sum_I\xi^{[I]} \nabla_{D}^I h_1^*(e); $$
note that $D_i^{(p)} = 0$. Thus, modulo $\ov I^{[p+1]} + I \oh
{D(1)}$, $\ep (h_2^*(e)) - h_1^*(e))$ reduces to
$$\sum_i \xi_i^{[p]} \nabla_{D_i}^p h_1^*(e) = \sum_i \xi_p (d \pi^*(t_i)) \nabla^p_{D_i}h_1^*(e) = (\id\ot \xi_p)(\psi(e)).$$
\end{proof}

%  Recall from (\ref{splift.p})
%that the functor $\cL_{\lift X S}$ of splittings of the fundamental
% extension (\ref{exseq.p})
%is represented by $\spec \cA_{\lift X S}$, where $\cA_{\lift X S} = \dirlim
%S^n\cE_{\lift X S}$.   This scheme
% is naturally a torsor under the pullback
%$F^* {\bf T}_{X'/S}$ of the tangent bundle $\spec_{X'/S} S^\cx \Omega^1_{X'/S}$
%via $F_{X/S}$.
%Its sheaf of algebras
%  $\cA_{\lift X S}$ admits a connection $\nabla_{\cA}$
%relative to $X/S$ and a horizontal increasing
%filtration $N_n$, with $N_n\cA_{\lift X S} \cong S^n\cE_{\lift X S}$.  Let
%$\psi_\cA$ denote the $p$-curvature of the connection $\nabla_{\cA}$.

\begin{remark}\label{pcurvpull.r}
Let $h \colon X \to Y$ be a morphism of smooth $S$-schemes, let
$E\in \mic(Y/S)$ be a module with an integrable connection, and let
$\psi_Y : E\to E\otimes_{{\cal O}_Y} F^*_{Y/S}\Omega^1_{Y'/S}$ be
its $p$-curvature. Then the $p$-curvature of $h^*E$ is the
composition:
$$h^*E\mapright{h^*\psi_Y} h^*E\otimes _{{\cal O}_X}F^*_{X/S} {h'}^*
\Omega^1_{Y'/S} \mapright{id\otimes h^*} h^*E\otimes _{{\cal O}_X}
F^*_{X/S}\Omega^1_{X'/S}.$$ This follows immediately from
Proposition~\ref{moch.p}; it was first proved years ago by O.
Gabber, using an indirect method.
\end{remark}

\begin{proof}[Proof of Proposition~\ref{psicomp.p}]
Let $\tilde F \colon \tilde X \to \tilde X'$ be a local lift of
$F_{X/S}$. Let $(\tilde D,\ov J,\gamma)$ denote the PD-envelope of
the diagonal ideal $J$ of $\tilde X(1)$, let $(D,\ov I, \gamma)$
denote its reduction modulo $p$, and denote again by $\tF$ the
induced maps $\tilde X(1) \to \tilde X'(1)$ and $\tilde D(1) \to
\tilde D'(1)$. Since $\overline J$ is flat over $\tilde S$,
multiplication by $p$ induces an injective map
$$[p] \colon \overline I/I\oh{D(1)} \to \overline J/pJ\oh{\tilde D(1)}.$$
Since $\overline J/\overline J^{[p+1]}$ is flat over $\tilde S$,
$(pJ\oh{\tilde D(1)} + \overline J^{[p+1]})\cap p\oh{\tilde D(1)} =
p(J\oh{\tilde D(1)} + \ov J^{[p+1]}),$ so multiplication by $p$
induces an injective map
$$[p]\colon  \ov I/(I\oh{D(1)} + \ov I^{[p+1]}) \to \ov J/(pJ\oh{\tilde D(1)} + \ov J^{[p+1]})$$

If $\tilde a $ is a local section of $\oh {\tilde X}$, we let
$\xi(\tilde a) := 1 \otimes \tilde a - \tilde a \otimes 1$.
\begin{claim}
Let $\tilde F \colon \tilde X \to \tilde X'$ be a local lift of
Frobenius, let $a$ be a local section of $\oh X$, let $a' :=
\pi^*(a)$, and let $\tilde a'$ be a local lift of $a'$ to $\oh
{\tilde X'}$.  Then
\begin{eqnarray*}
\xi(\tilde F^*(\tilde a')) &=& -[p]\xi(a)^{[p]}\pmod {pJ \oh{ \tilde D(1)}}\\
    & = & -[p]d_p(a) \pmod{pJ\oh{\tilde D(1)} + \ov J^{[p+1]}}
\end{eqnarray*}
\end{claim}
\begin{proof}
We may prove this claim with the aid of a local lifting $\tilde a$
of $a$.  Then $\tilde F^*(\tilde a') = \tilde a^p + p\tilde b$ for
some section $\tilde b$ of $\oh {\tilde X}$.  Since $p\xi(\tilde b)
\in pJ \oh {\tilde D(1)}$,
$$\xi(\tilde F^*(\tilde a')) = 1\otimes \tilde a^p - \tilde a^p \otimes 1
\pmod{ pJ\oh{\tilde D(1)}}.$$

Now $1\otimes \tilde a = \tilde a\otimes 1 +  \xi(\tilde a)$, so
$$1\otimes \tilde a^p = \tilde a^p\otimes 1 + p \zeta +  (\xi(\tilde a))^p,$$
where $p\zeta  = \sum_{i=1}^{p-1}{p\choose i}\tilde a^i \xi(\tilde
a)^{p-i}\in pJ\oh{\tilde D(1)}$.  Since $(\xi(\tilde a))^p = p!
(\xi(\tilde a))^{[p]}$ and $(p-1)! \equiv -1 \pmod p$, this proves
the claim.
\end{proof}

Let $\sigma_\tF$ be the splitting  associated with $\tF$ described
in  Proposition~\ref{torsor.p},
 and let $\alpha := \sigma_\tF(d\pi^*(a)) \in \cA_{\lift X S}$.  Then
$\psi(\alpha) \in \cA_{\lift X S}\ot F_{X/S}^*\Omega^1_{X'/S}$, and
by the $p$-curvature formula of Proposition~\ref{moch.p}, $(\id\ot
\xi_p) \psi(\alpha)$ is the class of $h^*_2(\alpha) - h^*_1(\alpha)$
in $ \cA_{\lift X S} \ot (\ov I/I\oh{D(1)} + \ov I^{[p+1]})$. If
$\tF' \colon \tilde D(1) \to \tX'$ is any section of $\cL_{\lift X
S} (\tilde D(1))$ and $\tF_i := \tF \circ \tilde h_i$,
\begin{eqnarray*}
[p] (\id\ot \xi_p) \psi(\alpha)(\tF') & = &  [p] h_2^*(\alpha)(\tF') - [p]h_1^*(\alpha(\tF')) \\
 & = &  [p]\sigma_{\tF_2}(F_{X/S}^*da')(\tF') - [p]\sigma_{\tF_1}(F_{X/S}^*da')(\tF') \\
 & = & \left (\tF^{\prime *}(\tilde a') - \tF_2^*(\tilde a')\right) - \left ( \tF^{\prime *}(\tilde a') - \tF_1^*(\tilde a')\right )\\
  & = & \tilde h_1^*\tF^*(\tilde a') - \tilde h_2^*\tF^*(\tilde a') \\
   & = & -\xi (\tilde F^*(\tilde a')) \\
    & = & [p](d_pa) \\
    & = & [p]\xi_p(d\pi^*a)
\end{eqnarray*}
Since $[p]$ is injective, it follows that $\psi(\alpha) = d\pi(a)$.
This proves the formula for elements of the form $\alpha =
\sigma_\tF(da')$. The general case follows from the fact that
\marginpar{not good} both $\psi$ and the action described in
Proposition \ref{torsor.p}.1 annihilate $\oh X \subseteq \cA_{\lift
X S}$ and both are compatible with the algebra structure.
\end{proof}

It is helpful to have at our disposal an explicit formula for the
connection on $\cA_{\lift X S}$. Recall from  \cite{ka.ncmt} that
the \emph{inverse Cartier isomorphism} $C^{-1}_{X/S}$ is a canonical
$\oh {X'}$-linear map:
\begin{equation}\label{cart.e}
C^{-1}_{X/S}\colon \Omega^i_{X'/S} \to
\cH^i(F_{X/S*}\Omega^\cx_{X/S});
\end{equation}
if $i =1$ and $a$ is a local section of $\oh X$, then
$C_{X/S}^{-1}(d\pi_{X/S}^*(a))$ is the cohomology class of
$a^{p-1}da$. Let $Z^i_{X/S}$ denote the sheaf of closed $i$-forms on
$X/S$.  Then the \emph{Cartier operator} is the composite
$$C_{X/S} \colon F_{X/S*} Z^i_{X/S}  \to \cH^i(F_{X/S*}\Omega^\cx_{X/S}) \to  \Omega^i_{X'/S}$$
where the first map is the natural projection and the second is the
inverse of $C^{-1}_{X/S}$. Since $\Omega^1_{X'/S}$ is locally free,
locally on $X$ there exists a section  of $C_{X/S}$ (in degree one),
giving rise to a commutative diagram:
\begin{equation}\label{zeta.d}
\begin{diagram}
  &  & F_{X/S*}\Omega^1_{X/S} \cr
   &\ruTo^{\zeta} & \uTo & \luTo^\zeta\cr
\Omega^1_{X'/S} & \rTo & F_{X/S*}Z^1_{X/S} &\rTo^{C_{X/S}}
&\Omega^1_{X'/S}\cr
 & \rdTo_{C^{-1}_{X/S}} & \dTo& \ldTo_{C^{-1}_{X/S}} \cr
&& \cH^1(F_{X/S*}\Omega^\cx_{X/S})
\end{diagram}
\end{equation}

Mazur's formula~\cite{m.fhf}  shows that a lifting $\tF$ of the
relative Frobenius morphism  $F_{X/S}\colon X \to X' $ determines
such a splitting $\zeta$.
% and hence by (\ref{zetas.t}), a splitting of the
%Azumaya algebra $\cD_{X/S}$ over the
%formal completion of the cotangent space $\bT^*_{X'/S}$ of $X'/S$ along its zero section.
% In general, such liftings and splittings exist only locally
%on $X$.   We shall see that a lifting of $X'/S$
%modulo $p^2$ determines a splitting of $\cD_{X/S}$ over the nilpotent divided power
%envelope of the zero section of $\bT^*_{X'/S}$.  This global result
%depends on the crystalline nature of the sheaf of liftings of $F_{X/S}$,
%which we now  explain.
  Suppose that
$\tilde F \colon \tilde X \to \tilde X'$ is a lifting of $F_{X/S}$
modulo $p^2$. Since $$ d\tilde F \colon \Omega^1_{\tilde X'/S} \to
F_{X/S*}\Omega^1_{\tilde X/S}$$ is divisible by $p$, there is a
unique map $\zeta_\tF$ making the  following diagram commute:
\begin{equation}\label{mazcart.d}
\begin{diagram}
\Omega^1_{\tilde X'/\tilde S} & \rTo^{d\tilde F} & \tilde
F\Omega^1_{\tilde X/\tS} \cr \dTo && \uTo^{[p]} \cr \Omega^1_{ X'/S}
& \rTo^{\zeta_\tF}& F_{X/S*}\Omega^1_{X/S} \cr \cr
\end{diagram}
\end{equation}
Then $\zeta_\tF$ is a splitting of the inverse Cartier operator in
the sense of  diagram (\ref{zeta.d}).
 Let us recall the proof. Let $\tilde a$ be lift of a section $a$ of $\oh X$
and let $\tilde a'$ be a lift of $\pi^*a$.  Then $\tilde F^*(\tilde
a') = \tilde a^p + p\tilde b$ for some $\tilde b \in \oh{\tilde X}$.
Hence
$$[p]\zeta_\tF(d\pi^*a) = d\tilde F^*(\tilde a') = [p]a^{p-1} da + [p]db,$$
where $b $ is the image of $\tilde b $ in $\oh X$.  Then
$\zeta_\tF(d\pi^*a) = a^{p-1}da + db$ is closed, and its image in
$\cH^1(F_{X/S*}\Omega^\cx_{X/S})$ is the class of $a^{p-1}da$, as
required.

%
%There is a more direct link between $\cE_{\lift X S}$ and $\cL_{\lift X S}$:
%the set of liftings of Frobenius can be identified with the set of
%splittings of (\ref{nexc.e}).  Furthermore, this set of splitting is represented
%by an affine scheme over $X$, as the following result shows.

\begin{proposition}\label{nabexp.p}
Let $\lift X S$ be as above and let $\tU$ be a lift of some open
subset of $X$, let $\tF$ be an element of $\cL_{\lift X S}(\tU)$,
and let $\sigma_\tF$ be the corresponding splitting of $\cE_{\lift X
S,U}$ described in Proposition~\ref{torsor.p}. Then for any local
section $\omega'$ of $\Omega^1_{X'/S}$,
$$\nabla(\sigma_{\tilde F}(1\ot\omega')) =  - \zeta_{\tilde F}(\omega')$$
where $\nabla$ is the connection (\ref{nexc.e}) on $\cE_{\lift X
S}$.
\end{proposition}
\begin{proof}
Since both sides are linear over $\oh {X'}$, it suffices to prove
the formula if $\omega' = da'$, where $a'$ is a section of $\oh
{X'}$. Let $\tT$ be the first infinitesimal neighborhood of $\tU$ in
$\tU\times_\tS \tU$ with its two natural projection $\tilde h_i
\colon \tilde T \to \tilde U$, and let $T$ be the reduction of $\tT$
modulo $p$, so that $\oh T \cong \oh U \oplus \Omega^1_{U/S}$. The
crystal structure on $\cE_{\lift X S}$ gives us isomorphisms
$$ h_2^*\cE_{\lift X S,U} \rTo^\cong \cE_{\lift X S, T}\lTo^\cong h_1^*\cE_{\lift X S,U}$$
reducing to the identity modulo the ideal $\Omega^1_{U/S}$ of $\oh
T$. Using the resulting identifications,
$$\nabla(\sigma_\tF(da')) := h_2^*(\sigma_\tF(F_{X/S}^*da')) - h_1^*(\sigma_\tF(F_{X/S}^*da'))
 \in \cA_{\lift X S,U}\ot \Omega^1_{X/S} \subseteq \cA_{\lift X S,T}.$$
Let us evaluate this section on an arbitrary section $ \tF' \colon
\tT \to \tX'$ of $\cL_{\lift X S}(\tT)$. Let $\tF_i := \tF\circ
\tilde h_i \in \cL_{\lift X S}(\tT)$. If  $\tilde a'$  is a lift of
$a'$, then by Proposition~\ref{torsor.p},
\begin{eqnarray*}
[p]\nabla(\sigma_\tF(F_{X/S}^*da'))(\tF')  & = &  [p]h_2^*(\sigma_\tF(F_{X/S}^*da'))(\tF') - [p]h_1^*(\sigma_\tF(F_{X/S}^*da'))(\tF')\\
 & = &  [p]\sigma_{\tF_2}(F_{X/S}^*da')(\tF') - [p]\sigma_{\tF_1}(F_{X/S}^*da')(\tF') \\
 & = & \left (\tF^{\prime *}(\tilde a') - \tF_2^*(\tilde a')\right) - \left ( \tF^{\prime *}(\tilde a') - \tF_1^*(\tilde a')\right )\\
  & = & \tilde h_1^*\tF^*(\tilde a') - \tilde h_2^*\tF^*(\tilde a') \\
   & = & -d\tF^*(\tilde a') \\
    & = & -[p]\zeta_\tF(da')
\end{eqnarray*}
\end{proof}

\begin{remark}\label{ez.r}
 Somewhat more generally, let
$\zeta$ be a section of $C^{-1}_{X/S}$ as in (\ref{zeta.d}), and let
$$(\cE_\zeta,\nabla) := \oh X \oplus F_{X/S}^*\Omega^1_{X'/S},$$
where $\nabla \colon \cE_\zeta\to \cE_\zeta\otimes \Omega^1_{X/S}$
is the map
$$(f,g\otimes\omega') \mapsto (df-g\otimes\zeta(\omega'), \omega'\otimes dg). $$
Then $\nabla$ is an integrable connection on $\cE_\zeta$, and one
can simply compute
 that its $p$-curvature is the map
$$\psi \colon \cE_\zeta \to \cE_\zeta\otimes F_{X/S}^*\Omega^1_{X'/S}
\quad\quad (f,g\ot \omega') \mapsto (g,0)\ot \omega' .$$ (See for
example \cite[2.10]{o.fgt}.) If $\tF$ is a lift of Frobenius, then
$\sigma_\tF$ provides a splitting of the fundamental exact sequence
(\ref{nexc.e}) and hence an isomorphism $\cE_{\lift X S} \cong
\cE_\zeta$ inducing the identity maps on $\oh X$ and
$F_{X/S}^*\Omega^1_{X'/S}$. The formula of
Proposition~\ref{nabexp.p} shows that this morphism is horizontal,
and hence provides another proof of Proposition~\ref{psicomp.p}.
\end{remark}

\subsection{Functoriality}\label{funt.ss}

The geometric construction of $\cL_{\lift X S}$ we have given makes
it quite straightforward to check its functoriality.  Note first
that a morphism $h \colon X \to Y$ of smooth $S$-schemes induces a
morphism of schemes $h' \colon X' \to Y'$,  a morphism of $\oh
{X'}$-modules $T_{X'/S} \to h^*T_{Y'/S}$, and hence a morphism of
crystals of vector bundles:
$$\bT_{h'} \colon F_{X/S}^*\bT_{X'/S} \to h^*F^*_{Y/S}\bT_{Y'/S}.$$

\begin{proposition}\label{exfunct.p}
Let $h \colon X \to Y$ be a morphism of smooth $S$-schemes and let
 ${\tilde h}'$ be a lift of $h'$. Then the pair
${\tilde h} := (h,{\tilde h}')$  induces a morphism of crystals of
torsors:
$$\cL_{\lift X S} \rTo^{{\cL_{\tilde h}}}h^{*}\cL_{\lift Y S}, $$
compatible with the actions of $F_{X/S}^*\bT_{X'/S}$ and
$h^*F_{Y/S}^*\bT_{Y'/S}$ via the morphism $\bT_{h'}$.  This induces
an isomorphism of crystals of $h^*F_{Y/S}^*\bT_{Y'/S}$-torsors,
$$  h^* F_{Y/S}^*\bT_{Y'/S}\times_{F_{X/S}^*\bT_{X'/S}}\cL_{\lift X S} \rTo^{{\cL_{\tilde h}}}h^{*}\cL_{\lift Y S}, $$
%where $h^* F_{Y/S}^*\bT_{Y/S}\times_{F_{X/S}^*\bT_{X/S}}\cL_{\lift X S}$ is the pushout
%$of $\cL_{\lift X S}$ via the morphism $\bT_{h'}$
a horizontal  morphism of filtered $\oh X$-algebras:
\begin{eqnarray*}
\theta_{\tilde h} \colon  (h^*\cA_{\lift Y S},N_\cx) &\to& (\cA_{\lift X S},N_\cx) \\
\end{eqnarray*}
and a horizontal isomorphism of $h^*F_{Y/S}^*\Gamma_\cx
T_{Y'/S}$-algebras
$$ h^*\cA_{\lift Y S} \rTo^\cong \cHom_{F_{X/S}^*\hat\Gamma_\cx T_{X'/S} }(h^*F_{Y/S}^*\hat\Gamma_\cx T_{Y'/S},\cA_{\lift X S})  $$
%Furthermore,
%for any  flat $\tT_1 \in Cris(X/\tS)$
%and $\tT_2 \in Cris(Y/\tS)$,  any PD-morphism $\tg \colon \tT_1 \to \tT_2$,
%and any $\tF \in \cL_{\lift X S}(\tT_1)$, the following diagram commutes:
%\begin{diagram}
%\tT_1 & \rTo^\tF  & \tX' \cr
%\dTo^\tg && \dTo_{\tilde h'} \cr
%\tT_2 & \rTo^{\cL_{\tilde h}(\tF)} &\tilde Y'.
%\end{diagram}

%If $\tilde h'_2$ and $\tilde h'_1$ are liftings of $h'$, and if
%$D \colon {h'}^*\Omega^1_{Y'/S}\to \oh {X'}$ is the map which expresses their
%difference,
%then the corresponding maps $\theta_2$ and $\theta_1$  differ by the map
%$$F_{Y/S}^*(D) \colon F_{X/S}^*{h'}^*\Omega^1_{Y'/S} \cong h^*F_{Y/S}\Omega^1_{Y'/S}
%\to \oh {X}$$
\end{proposition}
\begin{proof}
Recall from \cite[6.5]{bo.ncc} and \cite[5.11]{bo.ncc} that if $E$
is a crystal on $Y/S$, then   $h^*E$ is the unique sheaf such that
for each morphism $g \colon T_1 \to T_2$ from an object in
$Cris(X/S)$ to an object in $Cris(Y/S)$, $(h^*E)_{T_1} =
g^*(E_{T_2})$.  Now if $\tilde T$ is an object of $Cris_f(X/\tS)$,
$h'\circ f_{T/S}$ is a map $T \to Y'$, and the set $\cL_{h,T}$ of
its liftings $\tilde T \to \tilde Y'$ is a torsor under
$f_{T/S}^*h'^*T_{Y'/S}$.  We claim first of all that $T \mapsto
f_{T/S}^*h'^*T_{Y'/S}$ can be identified with
$h^*F_{Y/S}^*\bT_{Y'/S}$ and that $T \mapsto \cL_{h,T}$ can be
identified with $h^*\cL_{\lift  Y S}$.  Indeed, if $g \colon T_1 \to
T_2$ is as above, then
$$(h^*F_{Y/S}^*T_{Y'/S})_{T_1} = g^*((F_{Y/S}^*T_{Y'/S})_{T_2}) = g^*f_{T_2/S}^*T_{Y/S} =
f_{T_1}^*h'^*T_{Y/S},$$ proving the first part of the claim. Suppose
further that $\tT_1 \in Cris_f(X/\tS)$, $\tT_2 \in Cris_f(Y/\tS)$,
and $\tg \colon \tT_1 \to \tT_2$ is a PD-morphism, compatible with
$h$ and let $g \colon T_1 \to T_2$ be its reduction modulo $p$. If
$\tF$ is a local section of $\cL_{\lift Y S,T_2}$, then $\tF\circ
\tilde g \colon \tT_1 \to \tilde Y'$ is a lift of $f_{T_2}\circ g =
h'\circ f_{T_1/S}$, and the sheaf of such lifts forms a $g^*
f_{T_2/S}^* \bT_{Y'/S}$-torsor. Thus $\tF \mapsto \tF\circ \tg$
defines an isomorphism of torsors from  $\cL_{\lift Y S,T_2}
\times_{T_2}T_1$ to the torsor of such  liftings, proving the second
part of the claim. Now if $\tF_1 \colon \tilde T_1 \to \tilde X'$ is
a local section of $\cL_{\lift X S,T_1}$, then $\tilde h' \circ
\tF_1$ is  such a lifting.  Thus composition with $\tilde h'$
defines a morphism  $\cL_{\lift X S,T_1} \to g^*\cL_{\lift Y
S,T_2}$, which is evidently compatible with the torsor actions.
\end{proof}

\begin{corollary}\label{morphex.c}
Let  $h \colon X \to Y$ is a morphism of smooth $S$-schemes.  Then a
lift  $\tilde h' \colon \tilde X' \to \tilde Y'$ of $h'$ induces an
exact sequence
$$ h^*\cE_{\lift Y S} \to \cE_{\lift X S} \to F_{X/S}^*\Omega^1_{X'/Y'} \to 0. $$
If $h$ is smooth, this sequence is short exact (and locally split).
\end{corollary}
\qed

\subsection{Further remarks}\label{frm.ss}

If $F_S$ lifts to $\tilde S$ and $\tilde X/\tilde S$ lifts $X/S$,
then $\tilde X' := \tilde X \times_{F_{\tilde S}}\tilde  S$ lifts
$X'$. In this case there is a lifting $\tilde \pi' \colon \tilde X'
\to \tilde X$ of $\pi \colon X'\to X$, and the following proposition
applies.

\begin{proposition}\label{deltapi.p}
Suppose that $\tilde \pi \colon \tilde X' \to \tilde X$ lifts $\pi
\colon X \to X'$. Then for each section $\tilde a$ of $\oh{\tilde
X}$, there is a unique section $\delta_\tpi(\tilde a)$ of
$\cE_{\lift X S} \subseteq \cA_{\lift X S}$ such that for every
lifting $\tF \colon \tU \to \tX'$ of $F_{X/S}$ over some open subset
$U$ of $X$,
$$[p]\delta_\tpi(\tilde a) (\tF) = \tF^*\tilde \pi^*(\ta) -\ta^p$$
on $\tU$.  Furthermore
\begin{enumerate}
\item{The following diagram commutes:
\begin{diagram}
\oh X & \rTo^{[p]} & \oh {\tilde X} & \rTo & \oh X \cr \dTo^{F_X^*}
&& \dTo_{\delta_{\tilde \pi}} && \dTo_{d\circ \pi^*} \cr \oh X &
\rTo &  \cE_{\lift X S} & \rTo &F_{X/S}^*\Omega^1_{X'/S},
\end{diagram}
where the bottom row is the fundamental extension (\ref{nexc.e}).}
\item{  The
set of all $\delta_{\tilde \pi}(\tilde a)$ for  $ \tilde a \in \oh
{\tilde X}$ generates $\cE_{\lift X S}$ as an $\oh X$-module.}
\item{If $\tilde F \colon \tilde X \to \tilde X'$ is a lift of $F$,
$$ \delta_{\tilde \pi}(\tilde a)  =  \delta_{\tilde \pi}(\ta)(\tF) + \sigma_{\tilde F}(F_{X}^*da),$$
where $\sigma_{\tilde F}$ is the splitting defined in
Proposition~\ref{torsor.p}.}
\item{For every local section $\tilde a $ of $\oh {\tilde X}$ lifting some $a
\in \oh X$,
\begin{eqnarray*}
\nabla\delta_{\tilde \pi}(a) &=& -1\ot a^{p-1}da \in \cA_{\lift X S}\ot \Omega^1_{X/S} \quad\mbox{and} \\
\psi\delta_{\tilde\pi}(a) &= & 1\ot F_X^*(da)  \in \cA_{\lift X
S}\ot F_{X/S}^*\Omega^1_{X'/S}.
\end{eqnarray*}}
\item{If $\ta$ and $\tb$ are sections of $\oh \tX$ reducing to sections $a$ and
$b$ of $\oh X$.
\begin{eqnarray*}
\delta_\tpi(\ta \tb) & = & a^p\delta_\tpi(\tb) + b^p\delta_\tpi(\ta). \\
\delta_\tpi(\ta + \tb) & = & \delta_\tpi(\ta) + \delta_\tpi(b)
+\sum_{0<i< p} {a^i \over i!} {b^{p-i} \over (p-i )!.}
\end{eqnarray*}
}
\end{enumerate}
\end{proposition}
%{\bf Warning}:  The map $\delta_{\tilde \pi}$ is not additive and does not
%factor through $\oh X$.
\begin{proof}
First of all, note that $\pi \circ F_{X/S} = F_X$, which takes any
section $a$ of $\oh X$ to $a^p$.  Hence if $\tF$ is a lift of
$F_{X/S}$ and $\ta$ is a lift of $a$, $\tF^*\tilde \pi^*(\ta) -
\ta^p$ is divisible by $p$.  Thus the formula defining
$\delta_{\tilde \pi}$ as a function $\cL_{\lift X S}(\tU) \to \oh X$
makes sense.  Now if $\tF'$ is another lift of $F_{X/S}$,
\begin{eqnarray*}
[p]\delta_{\tilde \pi}(\ta)(\tF') & = & \tF^{\prime*}\tpi^*(\ta) -\ta^p\\
  & = & \tF^{*}\tpi^*(\ta) -\ta^p + \tF^{\prime*}\tpi^*(\ta) - \tF^*\tpi^*(\ta)\\
   & = & [p]\delta_{\tilde \pi}(\ta)(\tF) + [p]\sigma_{\tF}(F_{X}^*(da)(\tF'),
\end{eqnarray*}
by Proposition~\ref{torsor.p}.  This proves that, as functions on
$\cL_{\lift X S}(\tU)$,
$$\delta_\tpi(\ta) = \delta_\tpi(\ta)(\tF) + \sigma_{\tF}(F_{X}^*da). $$
  This proves that $\delta_\tpi$ is well defined and satisfies (3).
 If $\tilde a = [p]b$ for some $b \in \oh X$, then $\tilde a^p = 0$, and
  $[p]\delta_{\tilde \pi}(\tilde a)(\tF) = F_{X/S}^*\tpi^*(pb) = [p]b^p$. This proves the
commutativity of the first square in the diagram, and shows that the
sub-$ \oh X$-module of $\cE_{\lift X S}$ generated by the image of
$\delta_{\tilde \pi}$ contains  $\oh X$.  We have already proved
(3), which implies the commutativity of the second square and the
fact the set of  images of all the $\delta_\tpi(\ta)$'s generates
$\cE_{\lift X S}$. To prove (4), we may assume that a lifting
$\tilde F$  of $F_{X/S}$ exists. Then by (3),
\begin{eqnarray*}
\nabla\delta_{\tilde \pi}(\tilde a) & = & \nabla [p]\delta_\tpi(\tilde a)(\tF) + \nabla\sigma_{\tilde F}F_X^*(da) \\
  & = & [p] d\delta_\tpi(\tilde a)(\tF) - [p]\zeta_{\tilde F}(d\pi^*(a)) \\
   & = & d\left( [p]\delta_\tpi(\tilde a)(\tF) - (\tilde F \circ \tilde \pi)^*(\ta) \right) \\
  & = & -d\tilde a^p \\
   & = &- [p]a^{p-1}da
\end{eqnarray*}
This proves the first equation in (\ref{deltapi.p}.4). The second
follows from the formula for the $p$-curvature $\psi$   in
Proposition~\ref{psicomp.p}; see also Remark~\ref{ez.r}.

The proofs of the  formulas of (\ref{deltapi.p}.5) are
straightforward calculations which we leave to the reader.
\end{proof}

\begin{remark}\label{deltapi.r}
We have seen that if   $F_{\tilde S} \colon S \to S$ lifts $F_S$ and
$\tilde X' = \tilde X\times_{F_{\tilde S}} \tilde S$, then the
projection $\tilde X' \to \tilde X$ is a natural global choice of a
lifting $\tilde \pi$ as above.  If $\tilde X'$ is some other lifting
of $X'$, then such a lift $\tilde \pi$ will exist locally on $X$.
However in general there may be no lift of $F_S$ even locally on
$S$, and consequently there may be no lift $\tilde \pi$ even locally
on $X$. However, if $\tilde a \in \oh{\tilde  X}$ is a local lift of
$a \in \oh X$, then we can choose  a local lift $\tilde a' \in
\oh{X'}$ of $\pi^*(a)$.  Then the analogs of the formulas in
Proposition~\ref{deltapi.p} hold with $\tilde a'$ in place of
$\delta_{\tilde \pi}(\tilde a)$.
\end{remark}

Let us describe another construction of the fundamental exact
sequence (\ref{nexc.e}). For each $T \in Cris(X/\tS)$, let $\Gamma$
be the graph of $f_{T/S}\colon T \to X'$, and for each lifting $\tF
\colon \tT \to \tX'$ of $f_{T/S}$ let $\tilde \Gamma$ be the graph
of $\tF$. Let $J_\tT$ be the ideal of the of the immersion
\begin{equation}\label{jgamma.e}
 j \colon T \rTo^\Gamma T \times_S X' \rTo^{inc} \tilde T \times_{\tilde S} \tilde X'.
\end{equation}
A morphism $\tg \colon \tT_1 \to \tT_2$ in $Cris_f(X/\tT)$ induces a
corresponding morphism of conormal sheaves:
$g^*J_{\tT_2}/J_{\tT_2}^2  \to J_{\tT_1}/J_{\tT_1}^2$, and so the
family $\{ J_\tT/J_\tT^2 : \tT \in Cris_f(X/\tT)\}$ forms a sheaf on
$Cris_f(X/\tS)$. If  $\tF \colon \tT \to \tX'$  is a lifting of
$f_{T/S}$ and
 $\tilde c$ is a section of $J$,  $\tilde \Gamma^*(\ta)$ vanishes
on $T$, and hence is divisible by $p$.

\begin{proposition}\label{conorm.p}
For each $\tT \in Cris_f(X/\tS)$, there is a unique morphism
$$\beta \colon J_\tT/J^2_\tT \to \cE_{\lift X S,\tT} \quad \tilde c \mapsto \tilde \beta_{\tilde c}$$
such that for every local lift $\tF$ of $f_{T/S}$ and every section
$\tilde c$ of $J_\tT$,
$$[p]\beta_{\tilde c}(\tF) = \Gamma_\tF^*(\tilde c) \in \oh{\tT}.$$
In fact, $\beta$ defines an  isomorphism of crystals of
$\oh{X/S}$-modules and fits into a commutative diagram:
\begin{diagram}
0  & \rTo &  \oh T & \rTo^{[p]} & J_\tT/J_\tT^2 & \rTo & I/I^2 &
\rTo & 0 \cr &&\dTo && \dTo_\beta && \dTo \cr 0 & \rTo & \oh T &
\rTo &\cE_{\lift X S,\tilde T} &\rTo &f_{T/S}^*\Omega^1_{X'/S} &
\rTo & 0,
\end{diagram}
where $I$ is the ideal of $\Gamma \colon T \subseteq T\times X'$ and
the bottom row is the exact sequence (\ref{nexc.e}).
\end{proposition}
\begin{proof}
Suppose for example that $\ta'$ is a local section of $\oh{\tX'}$
and that $\tilde b$ is a local section of $\oh \tT$ such that
$f_{T/S}^*(a') = inc^*(\tilde b)$. Then $\tilde c := 1\ot \tilde b
- \tilde a \ot 1$ is a section of $J$, and $J$ is locally generated
by such elements.  If $\tF$ is any local lift of $f_{T/S}$,
$[p]\beta_{\tilde c}(\tF) = \tF^*(\tilde b) - \tilde a$.  If $\xi'$
is a local section of $f_{T/S}^*T_{X'/S}$ and $\tF' = \xi' + \tF'$,
then $\beta_{\tilde c} (\tF') = \beta_{\tilde c} (\tF) + \angles
{\xi'} {da'}$. This shows that $\beta_{\tilde c}$ defines a section
of $\cE_{\lift X S ,T}$.  It is clear that $\beta_{\tilde c}$
depends only on the class of $c$ mod $J^2$, and so $\tilde c\mapsto
\beta_{\tilde c}$ defines a map $J/J^2 \to \cE_{\lift X S,T}$.

Since the global definition of $\beta$ has been given, it is enough
to prove that it is an isomorphism locally, and we may assume that
a lifting $\tF$ of $f_{T/S}$ exists. By definition $I := J/(p)$ is
the ideal of $\Gamma$. Then $I/I^2 \cong f_{T/S}^*\Omega^1_{X'/S}$,
and the image of $\sigma_\tF(da)$ in $I/I^2$ is the class of
$$1\ot a' - f_{T/S}^*(a') \ot 1 = (f_{T/S}\times \id)^* (1\ot a' - a'\ot 1),$$
which corresponds to $f_{T/S}^*(da') $ in
$f_{T/S}^*\Omega^1_{X'/S}$.
%Thus, $\sigma_\tF$  splits the sequence (\ref{nexc.e}).

Next we verify the exactness of the top row of the diagram.
 This is formal, except for the injectivity
of the map $[p]$. To see this, let $\tilde Z := \tilde T
\times_{\tilde S} \tilde X'$, and let $Z$ be its reduction modulo
$p$.  Thus $J$ is the ideal of $T$ in $\tilde Z$, and $J/(p)$ is the
ideal $I$ of  $T$ in $Z$.  Multiplication by $p$ defines  a map $\oh
{\tilde Z} \to J \to J/J^2$ which factors through $\oh T$ and whose
cokernel is $I/I^2 \cong f_{T/S}^*\Omega^1_{X'/S}$.  This gives rise
to the exact sequence
$$\oh T \rTo^{[p]}  \cE_{\lift X S ,\tT} \rTo f_{T/S}^*\Omega^1_{X/S} \to 0.$$
To prove that $[p]$ is injective, we may work locally, and in
particular we may assume that there exists a lifting $\tF \colon
\tilde T \to \tilde X'$ of $f_{T/S}$.  Let $\tilde J \subseteq \oh
{\tilde Z}$ be the ideal of the  graph of $\tF$.
  Then $\oh {\tilde Z }/\tilde J \cong \oh {\tilde T}$, which is flat over ${\tilde S}$.
Therefore
$$Tor_1(\oh{ \tilde Z}/\tilde J, \oh {\tilde Z}/p\oh {\tilde Z}) = 0,$$ \ie,
 $p\oh  {\tilde Z} \cap \tilde  J = p\tilde J$.   Furthermore, $J = p\oh{ \tilde  Z} + \tilde J$
and $p^2 =0  $, so
$$  p\tilde J \subseteq p\oh {\tilde Z} \cap J^2 =p\oh {\tilde Z} \cap(p\tilde J + \tilde J^2) \subseteq
p\oh Z  \cap \tilde J^2 \subseteq p \oh {\tilde Z} \cap \tilde J = p
\tilde J = pJ $$ It follows that $p\oh { \tilde Z}\cap J^2 = p\tilde
J  = pJ$ and hence that the map  $\oh T \to J/J^2$ induced by
multiplication by $p$ is injective.

Furthermore, if $\tF$ exists,  the ideal $J$ is locally generated by
sections of the form $\tilde c = \ta \ot \tilde b$, where $\ta$ is a
local section section of $\oh \tT$ and $\tilde b$ is a local section
of $\oh {X'}$ such that $\tF^*(\tilde b) = \ta$.
\end{proof}

\begin{remark}\label{exsphint.r}
The isomorphism class of the extension of connections in
(\ref{nexc.e}) is an element of
$Ext^1_{MIC}(F_{X/S}^*\Omega^1_{X'/S},\oh X)$, and there is a
spectral sequence with
$$E_2^{i,j} \cong H^i(X,{\cExt}^j(\Omega^1_{X'/S},\oh X)) \cong H^i(X',T_{X'/S}\ot \Omega^j_{X'/S}).$$
In particular, there is an exact sequence
$$0 \to H^1(X',T_{X'/S}) \to Ext^1_{MIC}(F^*\Omega^1_{X'/S},\oh X) \to
H^0(X',T_{X'/S}\ot \Omega^1_{X'/S}).$$ The extension (\ref{nexc.e})
has the property that its image in $H^0(X',T_{X'/S}\ot
\Omega^1_{X'/S})$ is the identity, and the above exact sequence
shows that the set of extension classes with this property is a
(pseudo)-torsor under $H^1(X',T_{X'/S})$.  Note that the same is
true of the set of isomorphism classes of liftings  of $X'/S$.   We
shall investigate this further in section~\ref{sls.ss}.
\end{remark}

It is perhaps worthwhile to elucidate the relationship between the
fundamental extension (\ref{nexc.e}) and some more familiar exact
sequences. Since the relative Frobenius morphism $F_{X/S} \colon X
\to X'$ is a homeomorphism,  (\ref{nexc.e}) remains exact when
pushed forward by $F_{X/S}$.  Pulling the resulting sequence back by
means of the canonical map $\Omega^1_{X'/S} \to
F_{X/S*}F^*_{X/S}(\Omega^1_{X'/S})$, one gets an exact sequence
\begin{equation}\label{exseqp.e}
0  \to F_{X/S*}(\oh X) \to \cE'_{\lift X S} \to \Omega^1_{X'/S} \to
0
\end{equation}
of locally free sheaves on $X'$.
%   Recall that there
%is also an exact sequence
%\begin{equation}\label{exseqr.e}
%0 \to F_{X/S*}B_{X/S} \to F_{X/S*}Z^1_{X/S} \to  F_{X/S*}\cH^1_{DR}(X/S) \to 0
%\end{equation}
%where
%$$Z^1_{X/S} := \Ker d \colon \Omega^1_{X/S}\to \Omega^2_{X/S},$$
%$$B_{X/S} := {\rm Im\ } d \colon \oh X \to \Omega^1_{X/S},$$
%and $\cH^1_{DR}(X/S)$ is the cohomology sheaf of the De Rham complex
%$\Omega^\cx_{X/S}$ of $X/S$.
% Recall also that there is a natural isomorphism
%(the inverse Cartier operator)
%$$C^{-1} \colon \Omega^1_{X'/S} \to F_{X/S*}\cH^1_{DR}(X/S).$$
Each   local section $e'$ of $\cE'_{\lift X S}$ maps to a horizontal
section of $F^*_{X/S}\Omega^1_{X'/S}$, and hence  $\nabla(e')$ lies
in $\oh X \ot \Omega^1_{X'/S} \subseteq   \cE_{\lift X S} \ot
\Omega^1_{X'/S}$.  Since $\nabla$ is integrable, in fact $\nabla(e')
\in Z^1_{X'/S}$.  Thus,
 the connection $\nabla$ on $\cE_{\lift X S}$ induces an $\oh {X'}$-linear map
$\cE'_{\lift X S} \to F_{X/S*}Z^1_{X/S}$, which fits into the
commutative diagram below:

\begin{equation}\label{cube.e}
\begin{diagram}
0 & \rTo & F_{X/S*}\oh X & \rTo & \cE_{\lift X S} & \rTo &
F_{X/S*}F_{X/S^*}\Omega^1_{X'/S} &   \rTo & 0 \cr && \uTo^{=} &&
\uTo && \uTo_{inc} \cr 0 & \rTo & F_{X/S*}\oh X & \rTo & \cE'_{\lift
X S} & \rTo & \Omega^1_{X'/S} &   \rTo & 0 \cr
  &      & \dTo^d     &      & \dTo^{C^{-1}_{\lift X S}} &   & \dTo^{C^{-1}_{X /S}} \cr
0 & \rTo &F_{X/S*}B^1_{X/S}& \rTo &F_{X/S*}Z^1_{X/S}&  \rTo &
F_{X/S*}\cH^1_{DR}(X/S) & \rTo& 0
\end{diagram}
\end{equation}

Here the middle row is the pullback of the top row along $inc$ and
the familiar bottom row  is the pushout of the middle row along $d
\colon F_{X/S*}(\oh X) \to F_{X/S*}B_{X/S}$.  Recall that the bottom
row is rarely split.  Indeed, a splitting would induce an injective
map $\Omega^1_{X'/S} \to F_{X/S*}Z^1_{X/S} \to
F_{X/S*}\Omega^1_{X/S}$ and in particular a nonzero map
$F_{X/S}^*\Omega^1_{X'/S} \to \Omega^1_{X/S}$. For example, no such
map can exist on a complete curve of genus at least two over a
field.

Note that there is also an exact sequence
\begin{equation}\label{fseq.e}
 0 \to \oh {X'} \to F_{X/S*}\oh X \to F_{X/S*}B_{X/S} \to 0
\end{equation}
When pulled back to $X$ this sequence is split by the natural map
$$s\colon F_{X/S}^*F_{X/S*}\oh X \to \oh X.$$  Thus $F_{X/S}^*F_{X/S*}(\oh X)\cong \oh X \oplus
F_{X/S}^*F_*(B_{X/S})$.  Furthermore, (\ref{nexc.e}) is the pushout
by $s$ of of the pullback by $F_{X/S}^*$ of (\ref{exseqp.e}) along
$s$.  Warning: the map $s$ is not compatible with the natural
connections on the source and target. An $S$-scheme $X/S$ for which
sequence (\ref{fseq.e}) splits is called
\emph{$F$-split}~\cite{jo.fso}.

\section{Connections, Higgs fields, and the Cartier transform}\label{chf.s}
\subsection{$D_{X/S}$ as an Azumaya algebra}
Let $X/S$ be a smooth morphism of schemes in characteristic $p > 0$.
Let $\Omega^1_{X/S}$ be its sheaf of Kahler differentials, let
$T_{X/S}$ be its dual, and let $D_{X/S}$ denote the ring of
PD-differential operators of $X/S$ \cite[\S 2]{bo.ncc}.  A section
$D$ of $T_{X/S}$ can be viewed as a  derivation of $\oh X$ relative
to $S$ and hence as a PD-differential operator of order less than or
equal to $1$, and $D_{X/S}$ is generated as a sheaf of rings over
$\oh X$  by $T_{X/S}$. If $E$ is a sheaf of $\oh X$-modules, then to
give an integrable connection $\nabla \colon E \to E\otimes
\Omega^1_{X/S}$ is the same as to give an extension of the action of
$\oh X$ on $E$ to an action of $D_{X/S}$ \cite[4.8]{bo.ncc}, which
we continue to denote by $\nabla$. The $p$th iterate $D^{(p)}$ of a
derivation is again a derivation, hence a section of $T_{X/S}$ and
an operator of order less than or equal to $1$.  This is in general
not the same as the $p$th power $D^p$ of $D$, which is an operator
of order less than or equal to $p$, even though $D^{(p)}$ and $D^p$
have the same effect on sections of $\oh X$.  For each derivation
$D$, let
\begin{equation}\label{ceq.e}
c(D) := D^p - D^{(p)}.
\end{equation}

One can show either by calculating in local coordinates
\cite{bmr.lmsla} or by means of techniques from noncommutative
algebra \cite{ka.ncmt}, that $c$ is  an $F_X^*$-linear map from
$T_{X/S}$ to the center $\cZ_{X/S}$ of $D_{X/S}$.
 By adjunction, one deduces from $c$ an $\oh {X'/S}$-linear map
\begin{equation}\label{pcurvc.e}
c' \colon T_{X'/S} \to F_{X/S*}\cZ_{X/S} \quad :\quad D' \mapsto
c'(D') :=(1\ot D')^p - (1\ot D')^{(p)}
\end{equation}

Let  $\nabla$ be an integrable connection on $E$ and  $\psi\colon E
\to E\otimes F^*_{X/S}\Omega^1_{X'/S}$ be its $p$-curvature. It
follows from the definitions that for every local section $D'$ of
$T_{X'/S}$, $\psi_{D'}$ is the endomorphism of $E$ induced by the
differential operator $c'(D')$. This mapping satisfies  the
linearity and integrability conditions of a Higgs field with
$F^*_{X/S}\Omega^1_{X'/S}$ in place of $\Omega^1_{X/S}$.  We refer
to such a map as an \emph{$F$-Higgs} field on $E$, and we denote by
$$\Psi \colon \mic (X/S)  \to \fhig (X/S)$$
the functor taking  $(E,\nabla) $ to $(E,\psi)$.

Since $c'$ maps to the center of $F_{X/S*}D_{X/S}$, it extends to a
map from the symmetric algebra $S^\cx T_{X'/S}$ to $\cZ_{X/S}$, and
in particular makes $F_{X/S*}D_{X/S}$ into a sheaf of $S^\cx
T_{X'/S}$-modules. Let $\bT_{X'/S}^* := \spec_{X} S^\cx T_{X'/S}$ be
the cotangent bundle of $X'/S$.
%Since the action of $\oh {X'}$ on $F_{X/S*}D_{X/S}$ is compatible with this
%action,
Since $F_{X/S*}D_{X/S}$ is quasi-coherent as a sheaf of $\oh
{X'}$-modules, it defines a quasi-coherent sheaf $\cD_{X/S}$  on
$\bT^*_{X'/S}$.

Recall that an \emph{Azumaya algebra}  over a scheme $Y$ is a sheaf
of associative algebras $A$ such that locally for the fppf topology,
$A$ is isomorphic to $\End_{\oh Y}(\oh Y^n)$.
  More generally, if
$Y$ is a topological space,  $R$ is a sheaf of commutative rings on
$Y$, and $A$ is a sheaf of associative $R$-algebras
 which is locally free and finite rank as an $R$-module,
we say that $A$ is an Azumaya algebra over $R$  if the canonical map
$A \ot A^{op} \to \End_R(A)$ is an isomorphism.  One can show that
if $Y$ is a scheme and $R = \oh Y$, then these definitions agree.
(See  Chapter 4 of \cite{mi.ec} for a quick review.)

Our starting point in this section is the following theorem of
\cite{bmr.lmsla}, which asserts that $\cD_{X/S}$ is an Azumaya
algebra on $\bT^*_{X'/S}$.

\begin{theorem}\label{azum.t}
Let  $X/S$ be a   smooth $S$-scheme  of relative dimension $d$.
Then the map (\ref{pcurvc.e}) induces an isomorphism:
$$  S^\cx T_{X'/S} \rTo^\cong F_{X/S*}\cZ_{X/S}.$$
This morphism makes $F_{X/S*}D_{X/S}$  an Azumaya algebra over
$S^\cx T_{X'/S}$ of rank $p^{2d}$.  The corresponding  sheaf
$\cD_{X/S}$ of
 $\oh {\bT^*_{X'/S}}$-algebras on $\bT^*_{X'/S}$
is  canonically split (isomorphic to a matrix algebra) when pulled
back via the map $\pi_{T}$ in the  diagram below:
\begin{equation}\label{btp.e}
\begin{diagram}
{\bT}^{\prime *}_{X/S}&\rTo^{:=}& X\times_{X'} \bT^*_{X'/S}
&\rTo^{\pi_T} & \bT^*_{X'/S}\cr &\rdTo&\dTo && \dTo\cr && X &
\rTo^{F_{X/S}} & X'
\end{diagram}
\end{equation}
\end{theorem}
\begin{proof}
We recall here only the main idea of the proof, referring to
\cite{bmr.lmsla} for the details. Let $M_{X/S} := F_{X/S*}D_{X/S}$
which we can view as a module over
 $\oh {{\bT}^{\prime *}_{X/S}}= F_{X/S*}\cZ_{X/S} \ot_{\oh {X'}}F_{X/S*} \oh X$
via right multiplication and the inclusion $\oh X \to D_{X/S}$ as
well as a left module over itself. These left and right actions
agree on the center $\cZ_{X/S}$,  and hence they define a
homomorphism of sheaves of rings
$$F_{X/S*}D_{X/S} \ot_{S^\cx T_{X'/S}}\oh {{\bT}^{\prime *}_{X/S}} \to \cEnd_{\oh {{\bT}^{\prime *}_{X/S}}}(M_{X/S}),$$
which one can check is an isomorphism in local coordinates.
\end{proof}
%
%
%\begin{remark}\label{pslit.r}{\rm
%Suppose that $S$ is the spectrum of a perfect field,  $x$ is a closed point of $X$,
%and $x'$ is the image of $x$ in $X'$.  Then the Azumaya algebra $F_{X/S*}D_{X/S}$
%has a canonical splitting  over the fiber of $x'$ in $\bT^*_{X'/S}$:
%$F_{X/S*}(D_{X/S}\ot_{\oh X} k(x))$ is a splitting module, since $k(x') = k(x)$.
%}\end{remark}

Observe that if $\dim X/S > 0$, then
 $\cD_{X/S}$ is not split locally in the Zariski topology
of $\bT^*_{X'/S}$.  It suffices to check this  when $S$ is the
spectrum of a field and $X$ is affine.   Then $\Gamma(X, D_{X/S})$
has no zero divisors, because its  associated graded sheaf  with
respect to the filtration by order is canonically isomorphic to the
symmetric algebra $S^\cx T_{X/S}$.  Since $\bT*_{X'/S}$ is  integral
and $\cD_{X/S}$ is locally free as an $\oh {\bT^*_{X'/S}}$-module,
it also has no zero divisors and hence is not split.

\begin{remark}\label{splitdes.r}
The power of Theorem~\ref{azum.t} can be seen from its application
to Cartier descent \cite{ka.ncmt}. Consider the action of $D_{X/S}$
on $\oh X$.  Since $D^p$ and $D^{(p)}$ agree on $\oh X$, this action
kills the ideal $S^+T_{X'/S}$ of  $S^\cx T_{X'/S}$. Thus
$F_{X/S*}(\oh X)$ can be viewed as a sheaf of $i^*\cD_{X/S}$
modules, where $i \colon X' \to \bT^*_{X'/S}$ is the zero section.
Since $i^*\cD_{X/S}$ is an Azumaya algebra over $X'$ of rank
$p^{2d}$ and $F_{X/S*}(\oh X)$ has rank $p^d$, this shows that
$i^*\cD_{X/S}$ is split, and that tensoring with the splitting
module $F_{X/S*}(\oh X)$ induces an equivalence between the category
of $\oh {X'}$-modules and the category of $D_{X/S}$-modules for
which the action of $S^+T_{X'/S}$ is zero.  This is just the
category of $\oh X$-modules endowed with an integrable connection
whose $p$-curvature is zero.
\end{remark}

Let $\cD^\ell_{X/S}$ be the commutative subalgebra of $\cD_{X/S}$
generated by the left inclusion $\oh X \to \cD_{X/S}$ and its
center.  Then $F_{X/S*} \cD^\ell_{X/S}$ defines a quasi-coherent
sheaf of algebras $\cD^\ell_{X/S}$ on $\bT^*_{X'/S}$.  In fact, it
is easy to check that the natural map $F_{X/S}^* S^\cx T_{X'/S} \to
\cD^\ell_{X/S}$ is an isomorphism, so that $\spec_{\bT^*_{X'/S}
}\cD^\ell_{X/S} \cong \bT^{\prime*}_{X/S}$~(see diagram
(\ref{btp.e})). In particular,  a sheaf $M$ of $D_{X/S}$-modules
which is quasi-coherent over $X$ can be viewed as a quasi-coherent
sheaf of $\oh {\bT^{\prime *}_{X/S}}$-modules.

\begin{proposition}\label{rankone.p}
Let $f \colon Z \to \bT^*_{X'/S}$ be a morphism and suppose $L$ is a
splitting module for $f^*\cD_{X/S}$.  Then  $L$, viewed as a sheaf
of $\oh {X\times_{X'} Z}$-modules, is locally free of rank one.
\end{proposition}
\begin{proof}
First let us prove this when $f = \pi_T$ and $L = M_{X/S}$.  Our
claim is that $M_{X/S} :=\cD_{X/S}$ is locally free of rank one over
$ \pi_T^*\cD^\ell_{X/S} = \oh X \ot S^\cx T_{X'/S} \ot \oh X$, where
the first $\oh X$ acts by multiplication on the left and the second
on the right and the tensor products are taken over $\oh {X'}$. We
may assume that we have a system of local coordinates $(t_1, \dots,
t_d)$ for $X/S$, with a corresponding set of generators $D_i$ for
$D_{X/S}$. Then  the product $D_1^{p-1} \cdots  D_d^{p-1}$ generates
$M_{X/S}$ as a module over $\pi_T^*\cD^\ell_{X/S}$, as one sees from
the fact that $[D_i,t_j] = \delta_{ij}$. This generator defines a
surjective map $\pi_T^*\cD^\ell_{X/S} \to M_{X/S}$, and since the
source and target of this map are locally free $\oh{X'}$-modules of
the same rank, it is an isomorphism.

To deduce the general statement, note that it is enough to prove the
claim about $L$ after a faithfully flat cover, and in particular
after a base extension induced by $\pi_T$.  Thus we can replace $Z$
by $Z\times_{\bT^*_{X'/S}} {\bT}^{\prime *}_{X/S} \cong Z\times_{X'}
X$. The pullback of $M_{X/S}$ to this space has the desired
property, and $L$ is necessarily locally isomorphic to $M_{X/S}$.
This concludes the proof.
\end{proof}

%
%Recall that the category   $\mic(X/S)$-modules is equipped
%with an inner tensor product $\ot$.  If $(E_i,\nabla_i)$ are objects
%of $\mic(X/S)$ with $i = 1,2$, then the underlying $\oh X$-module of
%$E_1\otimes E_2$ is just the tenosr product in the category of $\oh X$-modules,
%with the connection $\nabla = \nabla_1\ot\id + \id \ot \nabla_2$.
%One can  easily check that the $p$-curvature of $E_1\ot E_2$  is given
%by the analogous formula: $\psi = \psi_1\ot \id + \id \ot \psi_2$,
%as is well-known.  The latter formula can be interpreted as a convolution
%product in the category $\fhig$, as explained in \S(\ref{conv.ss}).  In fact
%one can also define a notion of a tensor structure on an Azumaya algebra
%which provides a framework for the compatibility between the formulas
%for $\nabla$ and $\psi$ on $E_1\ot E_2$; see \S(\ref{aags.ss}).
%Similarly, $\cHom_{\oh X}(E_1,E_2)$ has a connection $\nabla$ sending
%$h$ to  $\nabla_2\circ h - h \circ \nabla_1$, and the corresponding
%$p$-curvature has the expected form.

Let us recall that the category of left $D_{X/S}$-modules is
equipped with a tensor structure.  In section~\ref{aags.ss} we will
discuss this structure from the point of view of Azumaya algebras.

\subsection{An \'etale splitting of $D_{X/S}$}

The  proof of Theorem~\ref{azum.t} gives an explicit flat covering
of $\bT^*_{X'/S}$ which splits  $\cD_{X/S}$.  It follows from the
general theory of Azumaya algebras that there exist \emph{\'etale}
coverings over which it is split. In this section we will give an
explicit construction of such a covering, which in fact is a
surjective \'etale endomorphism of the group scheme $\bT^*_{X'/S}$.

The construction of the splitting depends on a choice $\zeta$ of a
splitting of  the Cartier operator $C_{X/S}$, as  exhibited in
diagram~\ref{zeta.d}.
%to study the $p$-curvature functor $\Psi$. % \colon \hig(X/S) \to F$-$\hig(X/S)$.
In order to express the formulas we shall encounter geometrically,
we introduce the  following notation. The map $\zeta$ induces by
adjunction a map $F_{X/S}^*\Omega^1_{X'/S} \to \Omega^1_{X/S}$ whose
dual is a map $ \phi: T_{X/S} \to F_{X/S}^*T_{X'/S}$.  Pulling back
by $\pi_{X/S}$, we find an $\oh {X'}$-linear map $\phi' \colon
T_{X'/S} \to F_{X'}^* T_{X'/S}$.
%We also have a map $F_{X'}^* T_{X'/S} \to S^p T_{X'/S}$ sending
%a section $1 \otimes \tau$ to $\tau^p$.
%This map prolongs to
%a map of $\oh {X'}$-algebras $F_{X'}^* S^\cx T_{X'/S} \to S^\cx T_{X'/S}$,
%which is nothing other than the homomorphism corresponding the relative Frobenius morphism
%$F_{\bT^*/X'}$ for the $X'$-scheme $\bT^*_{X'/S}$:
%\begin{equation}
%\begin{diagram}
%\bT^*_{X'/S}& \rTo &\bT^{*(X')}_{X'/S} & \rTo  & {\bT'}^*_{X/S} & \rTo^\pi & \bT^*_{X'/S} \cr
%&  \rdTo & \dTo && \dTo && \dTo\cr
%&& X' & \rTo^{\pi_{X/S}} & X & \rTo^{F_{X/S}}    & X'
%\end{diagram}
%\end{equation}
We let $h_\zeta$ be the composite of the map of vector bundles
induced by $\phi'$ with the relative Frobenius map for the
$X'$-scheme $\bT^*_{X'/S}$, as displayed in the diagram below.

\begin{equation}\label{hf.e}
\begin{diagram}
\bT^*_\zeta :=\bT^*_{X'/S}& \rTo^{F_{\bT^*/X'}} &\bT^{*(X')}_{X'/S}
\cr
     & \rdTo_{h_\zeta} & \dTo_{\spec \phi'}\cr
                  && \bT^*_{X'/S}
\end{diagram}
\end{equation}
This morphism is a homomorphism of affine group schemes over $X'$,
but it is not compatible with the vector bundle structures. We shall
see that $\alpha_\zeta:= h_\zeta -\id$ is surjective and \'etale and
that the Azumaya algebra $\cD_{X/S}$ splits when pulled back via
$\alpha_\zeta$.

 Recall from  Remark~\ref{ez.r} that associated
to a  splitting $\zeta$ there is an object $(\cE_\zeta,\nabla)$  of
$\mic(X/S)$, where $\cE_\zeta = \oh X \oplus
F^*_{X/S}\Omega^1_{X'/S}$.
 The connection $\nabla$
on $\cE_\zeta$ induces a connection on each $S^n \cE_\zeta$,
compatibly with the inclusion maps $S^n \cE_\zeta \to
S^{n+1}\cE_\zeta$ induced by the map $\oh X \to \cE_\zeta$, and
hence also on the direct limit $\cA_\zeta := \dirlim S^n \cE_\zeta$.
The  splitting $\sigma \colon \cE_\zeta \to \oh X$ defines an
isomorphism of $\oh X$-algebras $\cA_\zeta \cong F^*_{X/S} S^\cx
\Omega^1_{X'/S}$ and the submodule $F_{X/S}^*\Omega^1_{X'/S}$
generates an ideal $I_\zeta$ of $\cA_\zeta$, which we can identify
with $S^+ F_{X/S}^*\Omega^1_{X'/S}$. By \cite[6.2]{bo.ncc}, the
completed PD-envelope $\hat \cA_\zeta^\gamma$ of this ideal has a
natural structure of a crystal of $\oh {X/S}$-modules, so  the
connection $\nabla$ on $\cA_\zeta$ extends canonically to a
connection $\nabla_\zeta$ on $\hat \cA_\zeta^\gamma$.  Furthermore,
if $a$ is a local section of $\ov I_\zeta$, then $\nabla_\zeta
a^{[n]} = a^{[n-1]} \nabla_\zeta(a)$, and $\nabla_\zeta$ maps
$\overline I_\zeta^{[n]}$ to
 $\overline I_\zeta^{[n-1]}\ot \Omega^1_{X/S}$.
The algebra $\cA_\zeta \cong F^*_{X/S}(\Omega^1_{X'/S})$ also has a
canonical F-Higgs field $\theta$: if $\xi$ is a local
 section of $T_{X'/S}$ and $\omega$
a local section of $\Omega^1_{X'/S}$, $\theta_\xi(\omega') = \angles
\xi {\omega}$, and the action of $\theta_\xi$ on the higher
symmetric powers is determined by the Leibnitz rule.  In fact, as we
saw in  Remark~\ref{ez.r}, this F-Higgs field is also the the
$p$-curvature of the connection $\cA_\zeta \cong S^\cx F_{X/S}^*
\Omega^1_{X'/S}$.
  This  field extends to the divided power envelope $\cA_\zeta^\gamma$
and its completion $\hat \cA_\zeta^\gamma$:  the pairing
\begin{equation}\label{sgammpair.e}
 S^n F_{X/S}^*T_{X'/S} \otimes \Gamma_{n+m}F_{X/S}^*\Omega^1_{X'/S} \to
\Gamma_{m}F_{X/S}^*\Omega^1_{X'/S }
  \end{equation}
comes from the multiplication  on the symmetric algebra and the
duality between the symmetric and divided power algebras explained
for example in~\cite[A10]{bo.ncc}.\marginpar{Is this clear enough?}
In particular, if $x \in T_{X'/S}$ and $\omega \in \Omega_{X'/S}^1$,
one has
\begin{equation}\label{sgamform.e}
\xi{ \omega^{[i]}}   = \angles \xi \omega
\omega^{[i-1]}\quad\mbox{and hence} \quad \xi^p \omega^{[i]} =
\angles \xi \omega^p \omega^{[i-p]}.
\end{equation}

Let
$$\cB_\zeta := \dirlim \cHom_{\oh X}(\cA_\zeta^\gamma/\overline I^{[n]},\oh X),$$
 be the topological dual of $\hat\cA_\zeta^\gamma$,
equipped with the dual connection and F-Higgs field~(\ref{dualh.e}).
Thus $B_\zeta \cong \oplus S^nF_{X/S}^*T_{X'/S}$ as an $\oh
X$-module. Because of the sign in the definition of the dual Higgs
field, a section of $\xi $ of $T_{X'/S}$ acts on $B_\zeta$ as
multiplication by $-\xi$. The $F_{X/S}^*S^\cx T_{X'/S}$-structure of
$B_\zeta$ corresponding to this field  identifies it with
$\iota_*F_{X/S}^* S^\cx T_{X'/S}$, where $\iota \colon \bT^*_{X'/S}
\to \bT^*_{X'/S}$ is the involution $t \to -t$ of the vector group
$\bT_{X'/S}$. Note that $\nabla$ is compatible with the algebra
structure of $\cA_\zeta$ and with the divided power algebra
structure of $\hat\cA_\zeta^\gamma$.  It is not, however, compatible
with the algebra structure of $\cB_\zeta$, but rather with its
coalgebra structure.

\begin{remark}\label{geombz.r}
 If  $\zeta$ comes from a lifting
$\tF$ of $F_{X/S}$ as in (\ref{cart.e}), we can give a geometric
interpretation of the construction of $\cB_\zeta$ as follows. Let
$(\cA_\tF^\gamma, \ov I) $ be the divided power envelope of the the
ideal $I$ of the section of $\cL_{\lift X S}$  corresponding to
$\tF$. Recall from Proposition~\ref{torsor.p} that $\cA_{\lift X S}$
has a connection $\nabla$ as well as an action of $F_{X/S}^* S^\cx
T_{X'/S}$, the latter via its identification with the ring of
translation invariant PD-differential operators.
% Although the  ideal $I$ is not preserved by $\nabla$ or
%by the action  of $F_{X/S}^*S^\cx T_{X'/S}$,
Both the connection $\nabla$ and the action of $F_{X/S}^*T_{X'/S}$
extend naturally to $\cA_\tF^\gamma$ and to its PD-completion $\hat
\cA_\tF^\gamma$.
  Then $\hat \cA_{\zeta}^\gamma$
can be identified with the $\hat \cA_\tF^\gamma$ and $\cB_\zeta$
with its topological dual.  It is clear from the definitions that
these identifications are compatible with the $F_{X/S}^*S^\cx
T_{X'/S}$-module structure, and Proposition~\ref{nabexp.p} shows
that they are also compatible with the connections.
\end{remark}

\begin{proposition}\label{alphaz.p}
Let $X/S$ be a smooth morphism of schemes in characteristic $p >0$
with a splitting $\zeta$ of $C^{-1}_{X/S}$, and let $h_\zeta$ and
$\cB_\zeta:= F_{X/S}^*S^\cx T_{X'/S} $  with the connection
$\nabla_\zeta$ described above.
\begin{enumerate}
\item{The map:
$$\alpha_\zeta :=\id - h_\zeta :\bT^*_\zeta :=\bT^*_{X'/S} \to \bT^*_{X'/S}$$
is a surjective \'etale morphism of affine group schemes over $X'$.
}
\item{The action of an element $\xi'$ of $S^\cx T_{X'/S}$ on $\cB_\zeta$
defined by its $p$-curvature is multiplication by
$\alpha_\zeta^*(\xi')$.}
\end{enumerate}
\end{proposition}
\begin{proof}
We have already observed that $h_\zeta$ is a morphism of group
schemes, and consequently so is $\alpha_\zeta$.  Since $h_\zeta$
factors through the relative Frobenius map,  its differential
vanishes, and it follows that $\alpha_\zeta$ is \'etale. Then  the
images under $\alpha_\zeta$ of   the geometric fibers of
$\bT_\zeta^*/X'$ are open subgroups of the fibers of
$\bT^*_{X'/S}/X'$. Hence the image of each fiber of $\bT_\zeta^*/X'$
must contain the entire corresponding  fiber of $\bT^*_{X'/S}/X'$,
and so  $\alpha_\zeta$ is surjective. Thus $\alpha_\zeta$ is an
\'etale covering  (but not necessarily an \'etale  \emph{cover},
since it need not be a finite morphism).

  We must next compute the $p$-curvature of the divided power envelope
  $\cA_\zeta^\gamma \cong \Gamma_\cx F_{X/S}^*\Omega^1_{X'/S}$ of
  $\cA_\zeta$.  Let $\omega'$ be a local section of $\Omega^1_{X'/S}$,
  so that $x := (0,1\otimes \omega')$ belongs to the divided power
  ideal of $\cA^\gamma_\zeta$. Let $D$ be a local section of $T_{X/S}$
  and let $\xi' := \pi_{X/S}^*D \in T_{X'/S}$.  Then $\phi(D) \in
  F_{X/S}^*T_{X'/S}$, and we shall need the following formula.

\begin{claim}
  $\phi(D)^p = F_{X/S}^*h^*_\zeta(\phi'(\zeta')) \in S^\cx
  F_{X/S}^*T_{X'/S}$.
\end{claim}
To check this,  let $\bT^* := \spec_{X'} S^\cx
  T_{X'/S}$ and let $\tilde \bT^*$ denote its pullback to $X$ via the
  map $F_{X/S}$, \ie, $\tilde \bT^* = \spec_X F_{X/S}^*S^\cx T_{X'/S}$.  Then\marginpar{check diagram}
  there is a commutative diagram:
  \begin{diagram}
  \tilde \bT^* & \rTo^{F_{\tilde \bT^*/X}} & \tilde {T^*}^{(X)} \cr
\dTo^{pr} && \dTo_{pr} & \rdTo^{\pi_{\tilde \bT^*/X}} \cr
\bT^*&\rTo^{F_{\bT^*/X'}} & {\bT^*}^{(X')} &\rTo^c & \tilde \bT^*,
  \end{diagram}
where the morphism $c$ is the projection
$${\bT^*}^{(X'}) := \bT^*\times_{F_{X'}} X' \cong \tilde \bT^*\times_{\pi_{X/S}}X'
\to \tilde \bT^*.$$ Let us view $\phi(D)$ as a section of
$\oh{\tilde \bT^*}$. Then $c^*\phi(D) = \phi'(\xi')$, so
\begin{eqnarray*}
 (\phi(D))^p & = &F^*_{\tilde \bT^*}(\phi(D)) \\
       & = &  F_{\tilde \bT^*/X'}^*\pi^*_{\bT^*/X}(\phi(D)) \\
           & = & F_{\tilde \bT^*/X'}^*pr^*(\phi'(\xi')) \\
       & = & pr^* F^*_{\bT^*/X'} (\phi'(\xi'))\\
       & = & pr^*h^*_\zeta(\phi'(\xi'))
\end{eqnarray*}
Since the map $pr$ in the diagram corresponds to pullback by
$F_{X/S}$, the claim is proved.

By the definition  of the connection on
 $\cE_\zeta \subseteq\cA_\zeta^\gamma$
given in Remark~\ref{ez.r} and  of the morphism $\phi$,
\begin{eqnarray*}
\nabla_D(x)  & = &  \nabla_{D}(0,1\ot\omega')
= (\angles D {-\zeta(1\ot\omega')},0)  \\
      &   = & - (\angles{\phi(D)} {1\ot\omega'},0) \\
        &  = &-\phi(D)x \in \cA_\zeta^\gamma.
\end{eqnarray*}

The formula   \cite[6.1.1]{o.hcpc} for  the $p$-curvature of
divided powers and the computation  of the $p$-curvature of
$\cE_\zeta$
 (Proposition\ref{psicomp.p}),
then say:
\begin{eqnarray*}\label{pcurv.e}
\psi_{\xi'}(x^{[i]})& = &
        x^{[i-1]} \otimes \psi_{\xi'}(x) + x^{[i-p]}(\nabla_D(x))^p\\
          & =& (\xi' x ) x^{[i-1]} - x^{[i-p]}(\phi(D) x)^p\\
     & =&\xi'  x^{[i]} -  (\phi(D))^p (x^{[i]})\\
    & = & (\xi' -h^*_\zeta(\xi'))x^{[i]}\\
       & =&  - \alpha^*_\zeta (\xi')}{x^{[i]}.
\end{eqnarray*}

Since $\cB_\zeta \subseteq \cHom(\cA^\gamma_\zeta,\oh X)$ as a
module with connection, the formula (2) for the $p$-curvature of
$\cB_\zeta$ follows from the  formula for the $p$-curvature of the
dual of a connection; see for example Lemma~\ref{pcurvet.l}.
\end{proof}

We can now show that $\cD_{X/S}$ splits when pulled back by
$\alpha_\zeta$. Since $\bT^*_\zeta = \bT^*_{X'/S}$,
$F_{X/S*}\cB_\zeta$ can also be viewed as a quasi-coherent sheaf on
$\bT^*_\zeta$.

\begin{theorem}\label{dzetas.t}
There is a unique action of $\alpha_\zeta^*(\cD_{X/S})$ on
$F_{X/S*}\cB_\zeta$ extending the actions of
$\alpha^{-1}_\zeta(\cD_{X/S})$ and of $\oh {\bT^*_\zeta}$. The
resulting module splits the Azumaya algebra
$\alpha_\zeta^*(\cD_{X/S})$.
\end{theorem}
\begin{proof}
Proposition~(\ref{alphaz.p} shows that the actions of $S^\cx
T_{X'/S}$ on $F_{X/S*}\cB_\zeta$ defined on the one hand through the
$p$-curvature homomorphism $S^\cx T_{X'/S} \to D_{X/S}$ and through
$\alpha_\zeta$ agree, and hence that the action  of $\cD_{X/S}$
extends canonically to  an  action of $\alpha_\zeta^*\cD_{X/S}$.
Since $\cB_\zeta = F_{X/S}^*\cB'_\zeta$ and $\cB'_\zeta \cong S^\cx
T_{X'/S}$, $\cB_\zeta$ is locally free of rank $p^d$ over
$\bT^*_\zeta$. Hence it is a splitting module for the Azumaya
algebra $\alpha_\zeta^*\cD_{X/S}$.
\end{proof}

\subsection{The Cartier transform}\label{ct.ss}

In this section we explain how a lifting of $F_{X/S} \colon X \to
X'$ or just of $X'/S$ modulo $p^2$ determines   splittings of
$\cD_{X/S}$ on  suitable neighborhoods of the zero section of
$\bT^*_{X'/S}$ We then  use these splittings to define
characteristic $p$ analogs of Simpson's correspondence.

Let us begin with the global construction. Suppose we are just given
a lifting $\tX'/\tS$ of $X'/S$; and as before, let $\lift X S$
denote the pair $(X/S,\tX'/\tS)$.  The sheaf $\Gamma_\cx T_{X'/S}$
has a canonical divided power structure and can be identified with
the divided power envelope $\bT^{*\gamma}_{X'/S}$ of the zero
section of the cotangent bundle $\bT^*_{X'/S}$ of $X'/S$.  Its
completion $\hat \Gamma_\cx T_{X'/S}$ with respect to the
PD-filtration can be viewed as the sheaf of functions on the formal
scheme $\hat \bT^{*\gamma}_{X'/S}$.  The topology on the structure
sheaf is defined by the PD-filtration and is
admissible~\cite[7.1.2]{EGAIb} but  not adic, and its underlying
topological space is $X'$.
  It inherits
the structure of a formal  group scheme from the group structure of
$\bT^*_{X'/S}$, and the group law is a PD-morphism. If
$\bT^{*\gamma}_n$ is the closed subscheme defined by $\ov
I^{[n+1]}$, the group law factors through maps $
\bT^{*\gamma}_{n}\times \bT^{*\gamma}_{m}\to \bT^{*\gamma}_{n+m}$
for all $n,m$.   We shall denote by $\hig^\cx_{\gamma}(X'/S)$ the
category of ${\oh {X'}}$-modules equipped with a locally
 PD-nilpotent $\hat \bT_{X'/S}^{*\gamma}$-Higgs field.
By definition, this is the category of sheaves of $ \Gamma_\cx
T_{X'/S}$-modules with the property that each local section is
annihilated by some $\ov I^{[n]}$.   As explained in (\ref{cnv.d})
and   section~\ref{aags.ss},
 the group law on $\hat \bT_{X'/S}^{*\gamma}$
defines a tensor structure (convolution) on the category
$\hig^\cx_\gamma (X'/S)$.  If $\hig_\gamma^n$ denotes the category
of $\oh {\bT^{*\gamma}_n}$-modules, the convolution factors through
functors
$$\hig_\gamma^m(X'/S) \times \hig_\gamma^n(X'/S) \to \hig_\gamma^{m+n}(X'/S).$$
 If $E_1$ and $E_2$
are objects of $\hig^\cx_\gamma(X/S)$ and $\xi$ is a local section
of $T_{X'/S}$, then the total PD-Higgs field on the tensor product
satisfies
\begin{equation}\label{thetaconv.e}
\psi_{\xi^{[n]}} =   \sum_{i+j=n}\psi_{\xi^{[i]}}\ot
\psi_{\xi^{[j]}}.
  \end{equation}
Note that $\psi_{\xi^{[p]}}$ can be nonzero even if $E_1$ and $E_2$
have level less than $p$.  Note also that this total PD-Higgs field
commutes with the Higgs fields $\id \ot \psi$ and $\psi \ot \id$. If
$E_1 \in \hig_\gamma^m(X'/S)$ and  $E_2 \in \hig_\gamma^n(X'/S)$,
then $\cHom_{\oh X}(E_1,E_2) \in\hig_\gamma^{m+n}(X/S)$, with the
unique PD-Higgs field  satisfying:
$$\psi_{\xi^{[n]}}(h) =
   \sum_{i+j=n}(-1)^j \psi_{\xi^{[i]}} \circ h \circ\psi_{\xi^{[j]}}.$$
See section~\ref{aags.ss} for a geometric explanation of this
formula. More generally, if $E_1$ is locally PD-nilpotent, then $E_1
= \dirlim N_kE_1$, where $N_kE_1$ is the subsheaf of sections
annihilated by $\ov I^{[k+1]}$, and if $E_2 \in \hig_\gamma^n(X/S)$,
then
$$\cHom_{\oh X}(E_1,E_2) \cong \invlim\cHom_{\oh X}(N_kE_1,E_2) $$
has a natural structure of a $\hat \Gamma_\cx T_{X'/S}$ module, but
it may not be locally PD-nilpotent.

Let $D^{\gamma}_{X/S}$ denote the tensor product
$$D^{\gamma}_{X/S} := D_{X/S}\otimes_{S^\cx T_{X'/S}} \hat \Gamma_\cx(T_{X'/S})$$
via the map $S^\cx T_{X'/S}\to D_{X/S}$
 induced by the $p$-curvature mapping  $c'$ (\ref{pcurvc.e}).
The category   $\mic_\gamma(X/S)$ of $D^\gamma_{X/S}$-modules on $X$
is equivalent to the category of sheaves of $\oh X$-modules $E$
equipped with a connection $\nabla$ and a   horizontal homomorphism
$$\psi \colon \hat \Gamma_\cx(T_{X'/S}) \to F_{X/S*}\cEnd_{\oh X}(E,\nabla)$$
which extends the Higgs field
$$ S^\cx T_{X'/S} \to F_{X/S*}\cEnd_{\oh X}(E,\nabla)$$
given by the $p$-curvature of $\nabla$, such that  each local
section  of $F_{X/S*}E$ is locally annihilated by $\Gamma^iT_{X'/S}$
for $i >>0$.
  For example, $\oh X$
has an obvious structure of  a $D_{X/S}^\gamma$-module. More
generally,  if $(E,\nabla)$ is a module with integrable connection
whose $p$-curvature is nilpotent of level less than $p$,
$(E,\nabla)$ can be viewed as an object of  $\mic^\cx_\gamma(X/S)$
by letting the $p$th divided power of the ideal $\Gamma^+T_{X'/S}$
act as zero.

 The convolution product on $\hig^\cx_\gamma(X'/S)$
 allows us to make the category
$\mic^\cx_\gamma(X/S)$ into a tensor category.
 If $E_1$ and $E_2$
are objects of $\mic^\cx_\gamma(X/S)$ and $\xi$ is a local section
of $T_{X'/S}$, then the total PD-Higgs field on the tensor product
satisfies equation~\ref{thetaconv.e}. Since these endomorphisms are
horizontal and since this formula agrees with the $p$-curvature of a
tensor product when $n  = 1$, it does indeed define an object of
$\mic^\cx_\gamma(X/S)$. If $E_1 \in \mic^\cx_\gamma(X/S)$ and $E_1
\in \mic^n_\gamma(X/S)$, then $\cHom_{\oh X}(E_1, E_2) \in
\mic_\gamma(X/S)$, with the usual connection rule and the action of
$\hat\Gamma_\cx(T_{X'/S})$ defined above.

In order to keep our sign conventions consistent with other
 constructions\footnote{See for example Remark \ref{whyinv.r} below.},
we have found it convenient to introduce a twist.
  Let $\iota \colon \bT^*_{X'/S} \to \bT^*_{X'/S}$
be the inverse operation in the group law.  Then $\iota^* = \iota_*$
is an involutive autoequivalence of the tensor category
$\hig(X'/S)$.  If $(E',\psi') \in \hig(X'/S)$,
\begin{equation}\label{iotahi.e}
(E',\psi')^\iota := \iota_*(E',\psi') = \iota^*(E',\psi') =
(E',-\psi')
\end{equation}

Recall that in Theorem~\ref{fliftc.t} we constructed an algebra
$\cA_{\lift X S}$ from the torsor of Frobenius liftings  $\cL_{\lift
X S}$. We have seen in Proposition~\ref{psicomp.p} that  the
$p$-curvature of
 $(\cA_{\lift X S},\nabla_\cA)$ coincides with the action of $S^\cx T_{X'/S}$
coming from the torsor structure and hence that it extends naturally
to a continuous   divided power Higgs field $\psi_\cA$.  Thus
$\cA_{\lift X S}$ can be regarded as an element of
$\mic^\cx_\gamma(X/S)$. Let $\cB_{\lift X S}$ be its $\oh X$-linear
dual, which makes sense as an object of $\mic_\gamma(X/S)$ (although
it does not lie in $\mic^\cx_\gamma(X/S)$).

\begin{theorem}\label{cart.t}
Let $\lift X S := (X/S,\tX'/\tS)$ be a smooth morphism together with
a lift of $X'/S$ modulo $p^2$.
\begin{enumerate}
\item{The $ D_{X/S}^\gamma$-module  $\cB_{\lift X S}$ described above
is a splitting module for the Azumaya algebra $F_{X/S*}
(D_{X/S}^\gamma)$ over $\hat \Gamma_\cx(T_{X'/S})$.}
\item{The functor
$$C_{\lift X S} \colon \colon \mic_\gamma(X/S) \to \hig_\gamma(X'/S).$$
$$ E \mapsto \iota^*\cHom_{D^\gamma_{X/S}}(\cB_{\lift X S},E)$$
defines an equivalence of   categories, with quasi-inverse
$$C^{-1}_{\lift X S} \colon     \hig_\gamma(X'/S) \to \mic_\gamma(X/S)$$
$$E' \mapsto \cB_{\lift X S}\ot_{\hat \Gamma_\cx T_{X'/S}} \iota^*E' .$$
Furthermore, $C_{\lift X S}$ induces an equivalence of \emph{tensor}
categories:
$$ \mic^\cx_{\gamma}(X/S) \to \hig^\cx_{\gamma}(X'/S).$$}
\item{ Let $(E,\nabla)$ be an object of $\mic_\gamma(X /S)$, let $\psi$ be its $p$-curvature,
and let $(E',\psi') := C_{\lift X S}(E,\nabla)$. A lifting $\tF$ of
$F_{X/S}$, if it exists, induces a natural isomorphism
$$\eta_\tF \colon (E,\psi) \cong F_{X/S}^*(E',-\psi'),$$}
\end{enumerate}
\end{theorem}
\begin{proof}
To prove that $\cB_{\lift X S}$ is a splitting module for
$D^\gamma_{X/S}$, it suffices to show that it is  locally free of
rank $p^d$ over the center $\hat \Gamma_\cx T_{X'/S}$ of
$D^\gamma_{X/S}$ As we have  already observed, the action of this
center coincides with the action coming from the torsor structure as
described in Proposition~\ref{torsor.p}.
 Since $\cA_{\lift X S}$
is coinvertible by {\em op. cit.}, $\cB_{\lift X S}$ is locally free
of rank one over $F^*_{X/S}\hat \Gamma_\cx T_{X'/S}$, and hence is
locally free of rank $p^d$ over $\hat \Gamma_\cx T_{X'/S}$.  It then
follows from the general theory of matrix algebras that
$\cHom_{D^\gamma_{X/S}}(\cB_{\lift X S},\ )$ and $\cB_{\lift X S}
\ot_{\hat \Gamma_\cx T_{X'/S}}\ $ are quasi-inverse equivalences of
categories.  Since $\iota_*$ is an involutive equivalence, the
functors $C_{\lift X S}$ and $C_{\lift X S}^{-1}$ are also
quasi-inverse equivalences.

The algebra structure of $\cA_{\lift X S}$ endows $\cB_{\lift X S}$
with the structure of a coalgebra with counit.  As explained in
Proposition~\ref{spde.p}, this gives $\cB_{\lift X S}$ the structure
of a tensor splitting and makes $\cHom(\cB_{\lift X S}, \ )$ a
tensor functor;  the compatibility isomorphism
$$C_{\lift X S}(E_1) \ot C_{\lift X S}(E_2)  \rTo C_{\lift X S}(E_1)\ot C_{\lift X S}(E_2)$$
comes from the diagram:
\begin{diagram}
\cHom_{D^\gamma_{X/S}}( \cB_{\lift X S},  E_1) \ot
\cHom_{D^\gamma_{X/S}}( \cB_{\lift X S}, E_2)\cr \dTo^\ot&  \rdTo\cr
\cHom_{D^\gamma_{X/S}}( \cB_{\lift X S} \ot \cB_{\lift X S},E_1\ot
E_2) & \rTo^{\mu*}& \cHom_{D^\gamma_{X/S}}(\cB_{\lift X S},  E_1 \ot
E_2) .
\end{diagram}
Since $\iota$ is a group morphism, $\iota_*$ is also compatible with
the tensor structure.

A lifting $\tF$ of $F_{X'/S}$ defines a trivialization of the torsor
$\cL_{\lift X S}$ and hence isomorphisms of $\hat \Gamma_\cx
T_{X'/S}$-modules
$$\cA_{\lift X S} \cong F_{X/S}^*S^\cx \Omega_{X'/S}, \quad
\cB_{\lift X S}   \cong F_{X/S}^*\hat \Gamma_\cx T_{X'/S}.$$ Then
$$E \cong \iota_* E' \ot_{\hat \Gamma_\cx{ T_{X'/S}}} \cB_{\lift X S} \cong
\iota_*E' \ot_{\hat \Gamma_\cx{ T_{X'/S}}} F_{X/S}^*\hat\Gamma_\cx
T_{ X'/ S} \cong F_{X/S}^*\iota_* E',$$ as $F^*_{X/S} \hat
\Gamma^\cx T_{X'/S}$-modules.
 Statement (3) follows.
\end{proof}

\begin{corollary}\label{cart.c}
With the notation of Theorem~\ref{cart.t}, the Azumaya algebra
$F_{X/S*}D_{X/S}$ splits on the $(p-1)$st infinitesimal neighborhood
of the zero section of $\bT^*_{X'/S}$.
\end{corollary}

\begin{remark}\label{locglog.r}
Although the source and target of the  isomorphism $\eta_\tF$ in
part (3) of Theorem~\ref{cart.t},   are independent of $\tF$,
$\eta_\tF$ itself is not.  Indeed, let $\tF_2$ and $\tF_1$ be two
liftings of $F_{X/S}$, differing by a  section $\xi$ of
$F_{X/S}^*T_{X'/S}$.  Then one can form $e^\xi$ in the completed
divided power envelope $ F_{X/S}^*\hat \Gamma T_{X'/S}$.  Since $E'
\in \hig^\cx_\gamma(X'/S)$, $e^\xi$ acts naturally on $F_{X/S}^*E'$,
and we have the formula
$$\eta_{\tF_2} = e^\xi\circ\eta_{\tF_1}.$$
This follows from the fact that the isomorphism of
Theorem~\ref{cart.t} is induced by the section of  $\cL_{\lift X S}$
defined by $\tF$ and the formula (\ref{taylor.e}) for the action  by
translation of $F_{X/S}^*T_{X'/S}$ on $\cA_{\lift X S}$.
\end{remark}

A lifting $\tF$ of $F_{X/S}$, if it exists, allows us to extend the
equivalence of Theorem~\ref{cart.t} to the category $\mic^\cx (X/S)$
of all locally nilpotent connections. As explained in   \cite[4.4,
4.12]{bo.ncc}, objects of this category give rise to modules over
the ring $\hat D_{X/S}$ of hyper-PD-differential operators. This
ring   can be identified with the tensor product of $D_{X/S}$ with
the completion  $\hat S^\cx T_{X'/S}$ of $S^\cx T_{X'/S}$ along the
ideal of the zero section, and $F_{X/S*}\hat D_{X/S}$ can be viewed
as an Azumaya algebra over the sheaf of rings $\hat S^\cx T_{X'/S}$,
or equivalently, over the formal completion $\hat \bT_{X'/S}$ of the
cotangent space of $X'/S$ along its zero section. Let
$\mic_\infty(X/S)$ denote the category of sheaves of $\hat
D_{X/S}$-modules on $\oh X$, and let $\hig_\infty(X'/S)$ denote the
category of sheaves of $\hat S^\cx T_{X'/S}$-modules on $\oh {X'}$.
The subcategories $\mic^\cx(X/S)$ and $\hig^\cx(X'/S)$ are tensor
categories. The natural map $\hat S^\cx T_{X'/S} \to \hat \Gamma_\cx
T_{X'/S}$ induces a pair of adjoint functors
\begin{eqnarray*}
\gamma_* \colon \hig_\gamma (X'/S) &\to& \hig_\infty(X'/S) \\
                     %\quad \mbox{and} \quad
   \gamma^* \colon  \hig_\infty(X'/S) &\to&\hig_\gamma (X'/S),
  \end{eqnarray*}
and similarly for $\mic(X/S)$.

Let $\cA_\tF$ be the divided power envelope of  the augmentation
ideal of $\cA_{\lift X S}$ defined by the section of $\cL_{\lift X
S}$ given by $\tF$, and let $\hat \cB_\tF$ be its $\oh X$-linear
dual.  Recall from Remark~\ref{geombz.r} that it has a natural $\hat
D_{X/S}$-module structure. There are natural maps
\begin{equation}\label{locglob.e}
\cA_{\lift X S} \to \cA_\tF ; \quad
 D_{X/S}^\gamma\ot_{\hat \cD_{X/S}}\hat\cB_{\tF} \cong
%\hat \Gamma_\cx(T_{X'/S}) \ot_{\hat S^\cx T_{X'/S}} \hat \cB_{\tF} \cong
\cB_{\lift X S}
\end{equation}

\begin{theorem}\label{loccart.t}
Let  $X/S$ be a smooth morphism of schemes endowed with a lift $\tF
\colon \tX \to \tX'$ of the relative Frobenius morphism $F_{X/S}$.
\begin{enumerate}
\item{The $\hat \cD_{X/S}$-module  $\hat\cB_\tF$ described above
is a splitting module for the Azumaya algebra $F_{X/S*}\hat D_{X/S}$
over its center  $F_{X/S*} \hat \cZ_{X/S} \cong \hat S^\cx
T_{X'/S}$. }
\item{The functor
$$C_\tF \colon  \mic_\infty(X/S) \to \hig_\infty{ \,}(X'/S)$$
$$  E \mapsto \iota_*\cHom_{\hat D_{X/S}}(\hat \cB_\tF, E) \colon$$
defines an equivalence of  categories, with quasi-inverse
$$C^{-1}_\tF \colon \hig_\infty{\,}(X'/S) \to \mic_\infty (X/S)$$
$$E' \mapsto \hat \cB_\tF \ot_{\hat \cZ_{X'/S}} \iota_*E' $$
Furthermore,  $C_\tF$  induces an equivalence of  \emph{tensor}
categories
$$\mic^\cx(X/S) \to \hig^\cx(X'/S).$$}
\item{The map $\hat\cB_{\tF} \to \cB_{\lift X S}$ (\ref{locglob.e}) induces
isomorphisms of functors
$$C_\tF \circ \gamma_* \cong \gamma_*\circ C_{\lift X S} \quad\mbox{and}\quad
C_{\lift X S}\circ \gamma^* \cong \gamma^*\circ C_{\tF}$$ }
\end{enumerate}
\end{theorem}
\begin{proof}
Let $\zeta \colon \Omega^1_{X'/S} \to F_{X/S*}\Omega^1_{X/S}$ be the
splitting of Cartier associated to $\tF$ (\ref{mazcart.d}). Recall
that we constructed in Proposition~\ref{alphaz.p} a module with
connection $\cB_\zeta$ together with a horizontal action of
$F_{X/S}^*S^\cx T_{X'/S}$; as a module over this sheaf of rings,
$\cB_\zeta$ is free of rank one. As we have already noted in
Remark~\ref{geombz.r}, we can identify $\hat \cB_\tF$ with the
formal completion $\hat \cB_\zeta$ of $\cB_\zeta$; this
identification is compatible with the connections and the actions of
$ F_{X/S}^*\hat S^\cx T_{X'/S}$.  In particular, $\hat \cB_\zeta$ is
an invertible (even free) sheaf of $ F_{X/S}^*\hat S^\cx
T_{X'/S}$-modules, and hence is locally free of rank $p^d$ over
$\hat S^\cx T_{X'/S}$. Recall from Proposition~\ref{alphaz.p} that
there is a surjective \'etale group morphism $\alpha_\zeta = \id -
h_\zeta \colon \bT^*_{X'/S} \to \bT^*_{X'/S}$, and note  that its
restriction $\hat \alpha_\zeta$ to $\hat \bT^*_{X'/S}$ is an
isomorphism, with inverse
$$\hat\alpha_\zeta^{-1} =  \id + h_\zeta + h _\zeta^2 + \cdots.$$
According to Proposition~\ref{alphaz.p}, the $p$-curvature action of
$S^\cx T_{X'/S}$ on $\hat \cB_\zeta$ is given by
$\hat\alpha_\zeta^*$ followed by the standard action.  Since
$\hat\alpha_\zeta$ is an isomorphism, $\hat \alpha_{\zeta*}\hat
\cB_\zeta$ is locally free of rank $p^d$. Thus  $\hat \cB_\zeta$ is
an $F_{X/S*}\hat D_{X/S}$-module which is locally free of rank $p^d$
over the center $\hat S^\cx T_{X'/S}$, and hence is a splitting
module.  This proves (1), and (2) follows as before.  The
compatibilities stated in (3) follow immediately from the
constructions and the morphisms (\ref{locglob.e}).
\end{proof}

Let us give  a more explicit  description of the local Cartier
transform $C_\tF$. Given a splitting $\zeta$ and a Higgs module
$(E',\psi')$ we define a module with integrable connection
\begin{eqnarray}\label{bpsid.e}
\Psi_\zeta^{-1}(E',\psi') & :=& (F_{X/S}^*E',\nabla) \\
\nabla &: =  & \nabla_0  +  (\id_{E'}\otimes \zeta ) \circ
F_{X/S}^*(\psi'),
\end{eqnarray}
where $ \nabla_0$ is the Frobenius descent  connection and $
(\id_{E'}\otimes \zeta ) \circ F_{X/S}^*(\psi')$ is the ${\cal
O}_X$-linear map
$$F_{X/S}^*E' \mapright{F_{X/S}^*(\psi')} F^*_{X/S}E'\otimes  F^*_{X/S}\Omega^1_{X/S}
\mapright{ \id_{E'}\otimes \zeta } F^*_{X/S}E'\otimes
\Omega^1_{X/S}.$$ Let  $\cB'_{X/S}:= \iota_*S^\cx T_{X'/S}$, viewed
as an
 object of $\hig(X'/S)$.
\begin{lemma}  The isomorphism $\cB_\zeta  \cong\iota^*S^\cx F^*_{X/S}T_{X'/S}$
induces an isomorphism
$$ \Psi_\zeta^{-1}(\cB'_{X/S})\cong \cB_{\zeta}$$
compatible with the connections.
\end{lemma}

\begin{proof}
 For each $n$, the ideal
$\cB^{\prime > n}_{X/S}  := \oplus_{j>  n} S^j T_{X'/S}$ also
defines an object of $\hig(X'/S)$, as does the quotient $\cB'_n$ of
$\cB'_{X/S}$ by  $\cB^{\prime >n}_{X/S}$.  Let $\cA^{\prime
\gamma}_n$ denote the dual of $\cB'_n$ in $\hig(X'/S)$ and let
$\cA'_{X/S} := \dirlim \cA^{\prime \gamma}_n$. For example,
$$\cA^{\prime \gamma}_1 = \cE'_{X/S}:= \oh {X'}\oplus \Omega^1_{X'/S},$$
and if $\xi \in T_{X'/S}$, $a' \in \oh{X'}$, and $\omega' \in
\Omega^1_{X'/S}$,
$$\xi(a',\omega') = (\angles {\xi} {\omega'},0).$$
Furthermore, $\cA^{\prime \gamma}_{X/S} \cong \Gamma_\cx
\Omega^1_{X'/S}$, and if $\omega'_j \in \Omega^1_{X'/S}$ for $j =
1,\ldots r$, then each $x_j := (0,\omega'_j)$  belongs to the
divided power ideal of $\cA^{\prime \gamma}_{X/S}$, and
$$\xi( x_1^{[ i_1]} x_2^{[ i_2]} \cdots x_r^{[i_r]}) =
\sum_j \angles {\xi }{\omega'_j} x_1^{[i_1]} x_2^{[i_2]} \cdots
x_j^{[i_j-1]} x_r^{[i_r]}.$$

It follows from the definitions that $(\cE_{\zeta},\nabla)
=\Psi_{\zeta}^{-1}(\cE'_{X/S})$.
%where $(\cE_{\zeta},\nabla)$ is the connection defined in Remark~\ref{ez.r},
Then by the formula above for the action of $T_{X'/S}$ on divided
powers and the similar formula for the action of a connection on
divided powers, it follows that   $(\cA^{\gamma}_\zeta,\nabla) \cong
\Psi_{\zeta}^{-1} (\cA^{\prime \gamma}_{X/S})$. Hence by the
compatibility of $\Psi^{-1}_\zeta$ with duality,
 $(\cB_\zeta,\nabla) \cong \Psi_\zeta^{-1}(\cB'_{X/S})$.
\end{proof}

Let $(E',\psi')$ be an object of $\hig(X'/S)$.   Then the
isomorphism in the previous lemma induces
\begin{equation}\label{zeta.e}
 \Psi_\zeta^{-1}( \iota_* E')\cong  \Psi_\zeta^{-1}( E' \otimes _{S^{\cdot } T_{X'/S}}    \cB'_{X/S}) \cong
E' \otimes _{S^{\cdot } T_{X'/S}}   \Psi_\zeta^{-1}(
\cB'_{X/S})\cong  E' \otimes _{S^{\cdot } T_{X'/S}}  \cB_{\zeta}
\end{equation}
Recall  from Theorem~\ref{dzetas.t} that $F_{X/S*}\cB_\zeta$  splits
the Azumaya algebra $\alpha_\zeta^*\cD_{X/S}$ over
$$\bT^*_\zeta := \bT^*_{X'/S}\rTo^{\alpha_\zeta} \bT^*_{X'/S}. $$
This, together with  (\ref{zeta.e}),  imply the following result.

\begin{theorem}\label{zetas.t}
Let $\zeta$ be a lift of $C^{-1}_{X/S}$, let
$$\alpha^*_\zeta \colon F_{X/S*}\cZ_{X/S}\cong S^\cx T_{X'/S}\to S^\cx T_{X'/S} := \cZ_\zeta$$
be the map described in Proposition~\ref{alphaz.p}, and let
$$D_\zeta := S^\cx T_{X'/S}\ot_{\alpha_\zeta^*}F_{X/S*} D_{X/S}.$$
Let $\mic_\zeta(X/S)$ denote the category of sheaves of
$D_\zeta$-modules on $X$. For each $\xi \in T_{X'/S}$, the
$p$-curvature $\psi_\xi $    on $\Psi^{-1}_\zeta( \iota_* E' )$  is
induced  by  the action of $\alpha_\zeta^*(\xi)$ on $E'$, \ie,
$\psi_\xi = F_{X/S}^*(-\alpha_\zeta^*(\xi))$. This makes
$\Psi^{-1}_\zeta( \iota_* E' )$  a $D_\zeta$-module. Furthermore,
the functors
\begin{eqnarray*}
\hig(X'/S) & \to &\mic_\zeta(X/S)\\
(E',\psi') & \mapsto & E' \otimes_{   \cZ_\zeta   }  \cB_{\zeta}\\
(E',\psi') & \mapsto &   \Psi_\zeta^{-1}( \iota_* E')
\end{eqnarray*}
are isomorphic equivalences of categories, with quasi-inverse given
by $$E \mapsto \cHom_{\cD_\zeta}(\cB_\zeta,E).$$
\end{theorem}
\begin{corollary}\label{zettof.c}
Let $(E,\nabla)$ be an object of $\mic_\infty(X/S)$, let $\psi$ be
its $p$-curvature, and let $(E',\psi') := C_\tF(E,\nabla)$, and let
$\zeta$ be the splitting of Cartier determined by $\tF$.
 Then there is canonical isomorphism:
\begin{eqnarray*}
(E,\nabla) &\cong & \Psi^{-1}_\zeta \alpha^{-1}_{\zeta*}(E',\psi').
 \end{eqnarray*}
\end{corollary}
\begin{remark}\label{whyinv.r}
The appearance of the  involution $\iota$ in Definition \ref{cart.t}
insures  the compatibility  of the Cartier transform with the usual
Cartier operator. Let us explain this in the context of extensions.
The  group $\Ext^1_{MIC}(\oh X, \oh X)$ of isomorphism classes of
the category $EXT^1_{MIC}(\oh X,\oh X)$ of extensions of $\oh X$ by
$\oh X$ in the category $\mic(X/S)$ is canonically isomorphic to
the de Rham cohomology group $H^1_{dR}(X/S)$.   Similarly, the group
$\Ext^1_{HIG}(\oh {X'}, \oh {X'})$ of isomorphism classes of the
category $EXT^1_{HIG}(\oh {X'}, \oh {X'})$ of extensions of $\oh
{X'}$ by $\oh {X'}$ in $HIG(X'/S)$ is canonically isomorphic to
$$H^1_{Hdg}(X'/S) \cong H^1(X',\oh {X'} \oplus \Omega^1_{X'/S}).$$
The inverse  Cartier transform defines an equivalence of categories
$$ C_{\lift X S}^{-1} \colon  EXT^1_{HIG}(\oh {X'},\oh {X'})
\to EXT^1_{MIC}(\oh X,\oh X) ,$$ and hence an isomorphism of groups
$$ c_{\lift X  S }^{-1} \colon  H^1_{Hdg}(X'/S) \to  H^1_{dR}(X/S). $$
Let us consider the following diagram.
\begin{diagram}
 H^1_{Hdg}(X'/S) &\rTo& H^0(X',\Omega^1_{X'/S}) \cr
\dTo_{c^{-1}_{\lift X S}} && \dTo_{C^{-1}_{X/S}} \cr
 H^1_{dR}(X/S) & \rTo& H^0(X, \cH^1_{dR}(X/S)).
\end{diagram}
Thanks to our definition, the diagram is commutative. It suffices to
verify this when $F_{X/S}$ lifts and  for extensions $\oh {X'} \to E
\to \oh {X'}$ which split in the category of $\oh {X'}$-modules.
Then $E'$ has a basis $(e'_0,e'_1)$ such that $\psi (e'_0) = 0$ and
$\psi (e'_1) =e_0 \ot  \omega' $, where $\omega' \in
\Omega^1_{X'/S}$. Then one can check that $E := C^{-1}_{\lift X
S}(E')$ has a basis $(e_0,e_1)$ such that $\nabla(e_0) = 0$ and
$\nabla (e_1) = e_0\ot \zeta(\omega')$, where $\zeta$ is the
splitting of $C^{-1}_{\lift X S}$ defined by the lifting of
$F_{X/S}$.  This implies that the diagram commutes.
\end{remark}

\subsection{The Cartier transform as Riemann-Hilbert}
In the previous section we defined a pair of inverse
quasi-equivalences of categories:
\begin{eqnarray*}
  C_{\lift X S} \colon \mic^\cx_\gamma(X/S) \to \hig^\cx_\gamma(X'/S) \quad
&:& \quad E \mapsto \iota^*\Hom_{\cD^\gamma_{X/S}}(\cB_{\lift X S},E)\\
  C^{-1}_{\lift X S} \colon \hig^\cx_\gamma(X'/S) \to \mic^\cx_\gamma(X/S) \quad
&:& \quad E' \mapsto \cB_{\lift X S}\ot_{\hat \Gamma_\cx T_{X'/S}}
E'
\end{eqnarray*}
Our goal here is to show how the ring structure on the dual
$\cA_{\lift X S}$ of $\cB_{\lift X S}$ can be used to give an
alternative and more symmetric description of these functors.  This
viewpoint sharpens the analogy between the Cartier transform and the
Riemann Hilbert and Higgs correspondences, with the sheaf of $\oh
X$-algebras $\cA_{\lift X S}$ playing the role of the sheaf of
analytic or $C^\infty$ functions. This construction of the Cartier
transform relies on the ``Higgs transforms'' described in
(\ref{higtr.d}) and ordinary Frobenius descent instead of the theory
of Azumaya algebras.

Roughly speaking, the idea is the following.  The algebra
$\cA_{\lift X S}$ is endowed with a connection $\nabla_\cA$ and a
PD-Higgs field $\theta_\cA$. If $(E,\nabla)$ is an object of
$\mic^\cx_\gamma(X/S)$, the tensor product connection on $E
\ot\cA_{\lift X S}(E)$ commutes with the PD-Higgs field $\id
\ot\theta_\cA$.   Hence $\id \ot \theta_\cA$ induces a PD-Higgs
field on the sheaf of horizontal sections of $E \ot \cA_{\lift X
S}$, and it turns out that the corresponding object of
$\hig^\cx_\gamma(X'/S)$ is $C_{\lift X S}(E)$. Similarly, if
$(E',\theta')$ is an object of $\hig^\cx_\gamma(X'/S)$, then the
total PD-Higgs field $\theta'_{tot}$ of $E' \ot \cA_{\lift X S}$
commutes with the connection induced by $\nabla_\cA$.  Hence the
subsheaf of sections annihilated by $\theta'_{tot}$ inherits a
connection, and the corresponding object of $\mic^\cx_{\gamma}(X/S)$
is  $C^{-1}_{\lift X S}(E')$.

To make this precise, we begin with some notation and a slightly
more general setting. Let $\Omega$ be a locally free sheaf of $\oh
X$-modules, let $T$ be its dual,  and let $\bT$ be the vector group
$\spec_X S^\cx \Omega$.  Let  $\ccG$ be the group scheme $\spec_X
\Gamma_\cx T$, and let us  write ${\oh \ccG}$ for the sheaf
$\Gamma_\cx T$, $\cI$ for the divided power ideal $\Gamma_+ T$ of
$\Gamma_\cx T$,    and $\oh {\ccG_n} := \oh \ccG/\cI^{[n+1]}$ if  $n
\in \bn$.
  Recall from the discussion
preceding Proposition~\ref{torsor.p} that if $\pi_\cL \colon \cL \to
X$ is any $\bT$-torsor, then there is natural action of $\oh \ccG$
on the filtered  algebra $(\cA_\cL, N_\cx) :=(\pi_{\cL *} \oh \cL,
N_\cx)$, and that  $N_n\cA_\cL$ is the annihilator of the ideal
$\cI^{[n+1]}$.  Thus there is a natural map
\begin{equation}\label{camult.e}
\oh {\ccG_n } \times N_n \cA_\cL \to N_n \cA_\cL.
\end{equation}

We shall find it both useful and convenient to study filtered  $\oh
\ccG$-modules. We denote by $\cI^\cx$ the divided power filtration
on $\oh \ccG$ (although we should perhaps really write $\cI^\ccx$).

\begin{definition}\label{filtsat.d}
Let $E$ be an $\oh \ccG$-module. An increasing (resp. decreasing)
filtration $N_\cx$ (resp $N^\cx$) of $E$ by sub $\oh\ccG$-modules is
said to be \emph{$\cI^\cx$-saturated}, or just an
\emph{$\cI^\cx$-filtration}, if for all $j$ and $k$,
$$I^{[j]}N_kE \subseteq   N_{k-j}E, \quad\mbox{(resp.} \quad
I^{[j]}N^kE \subseteq N^{k+j}E\mbox{)}.$$
\end{definition}
For example, the filtrations $N_\cx$ of $\cA_\cL$ and $\cI^\cx$ of
$\oh \ccG$ are $\cI^\cx$-filtrations.  If $E$ is any $\oh
\ccG$-module, then the \emph{canonical filtration}
\begin{equation}\label{canfilt.e}
N_kE := \{ e \in E : \cI^{[k]}e = 0 \}
\end{equation}
is $\cI^\cx$-saturated, and $E$ is locally nilpotent if and only if
this filtration is exhaustive. If $(E_1,N_\cx)$ and $(E_2,N_\cx)$
are $\oh \ccG$-modules with $\cI^\cx$-saturated filtrations,
 then the tensor product filtration
 \begin{equation}\label{tensfil.e}
N_c(E_1 \ot E_2) := \sum_{a+b=c} Im(N_aE_1\ot N_bE_2 \to E_1\ot E_2)
    \end{equation}
is again $\cI^\cx$-saturated, because the group law induces maps
$$\mu^* \colon \cI^{[j]} \to \sum_{a+b=j} Im(\cI^{[a]}\ot \cI^{[b]} \to
\oh \ccG\ot \oh \ccG).$$

If $E$ is any $\oh {X}$-module, let $\theta_0$ denote the $\oh
\ccG$-module structure on $E$ for which the ideal $ \cI$ acts as
zero.  That is, $(E,\theta_0) = i_*E$, where $i^* \colon \oh \ccG
\to \oh X$ is restriction along the zero section. If  $E$  any $\oh
X$-module and $\theta$ is an $\oh \ccG$-module structure on $E$, let
$$E^\theta := \cHom_{\oh \ccG}(i_*\oh X,E),$$
\ie, $E^\theta$ is  the sub $\oh X$-module of $E$ consisting of all
the elements annihilated by the ideal $\cI$.

Now let $E$ be an $\oh X$-module equipped with an $\oh \ccG$-module
structure $\theta$ and an $\cI^\cx$-saturated filtration $N_\cx$.
The $\oh X$-module
$$\cA_\cL(E):= E \ot_{\oh X} \cA_\cL$$
has three natural $\oh \ccG$-module structures: the action by
transport of structure via $E$, the action by transport of structure
via $\cA_\cL$, and the convolution structure defined in
(\ref{thetaconv.e}). We shall denote these by
\begin{eqnarray*}
  \theta_E &:=& \theta_E \ot \id_\cA = \theta_E\cnv \theta_0\\
   \theta_\cA &: = & \id_E \ot \theta_\cA = \theta_0\cnv \theta_\cA\\
    \theta_{tot} & := & \theta_E\cnv \theta_\cA
\end{eqnarray*}
We endow it with the total (tensor product) filtration
$N_\cx$~(\ref{tensfil.e}). It follows from formula
(\ref{thetaconv.e}) that $\theta_\cA$ and $\theta_{tot}$ commute.
  Define
$$\cT_\cL(E) := \left(\cA_\cL(E)\right )^{\theta_{tot}},$$
with the $\oh \ccG$-structure $\theta_\cT$ induced by $\theta_\cA$
and the filtration induced by $N_\cx$  We have natural maps,
compatible with the $\oh \ccG$-structures shown:
\begin{equation}\label{ije.e}
(E,\theta_0,\theta_E)  \rTo^{i_E}
(\cA_\cL(E),\theta_\cA,\theta_{tot}) \lTo^{j_E}
(\cT_\cL(E),\theta_\cT,\theta_0),
\end{equation}
where $i_E(e) := 1\ot e$ and $j_E$ is the inclusion.  Note that
$i_E$ factors through $(\cA_\cL)^{\theta_\cA}$ and $j_E$ factors
through $(\cA_\cL)^{\theta_{tot}}$. Endow $\ca_\cL(\cT_\cL(E))$ with
the tensor product filtration, and let
$$ h \colon \cA_\cL(\cT_\cL(E)) \to \cA_\cL(E)$$
be the map defined by the commutative diagram:
\begin{equation}\label{hmap.e}
\begin{diagram}
\cA_\cL(\cT_\cL(E))& =&    \cT_\cL(E) \ot \cA_\cL & \rTo^{j_E\ot
\id} & E \ot \cA_\cL \ot \cA_\cL \cr &&     \dTo^h &  \ldTo_{\id_E
\ot m} \cr
  &&   E \ot \cA_\cL.
\end{diagram}
\end{equation}

\begin{proposition}\label{higgsh.p}
Let $E$ be an $\oh X$-module with a locally nilpotent $\oh
\ccG$-module structure $\theta$ and an $\cI^\cx$-filtration $N_\cx$
which is bounded below.
\begin{enumerate}
\item The map $i_E$ of (\ref{ije.e}) is injective
and strictly compatible with the filtrations,  and its image is
$\left( \cA_\cL(E) \right)^{\theta_\cA}$.
\item The map $j_E$ of (\ref{ije.e}) is injective
and strictly compatible with the filtrations, and its image is
$\left (\cA_\cL(E) \right)^{\theta_{tot}}$.
\item  The map $h$ of (\ref{hmap.e}) fits in  a commutative diagram
  \begin{diagram}
    \cT_\cL(\cT_\cL(E)) & \rTo^{j_{\cT_\cL(E)}} & \cA_\cL(\cT_\cL(E))\cr
\dTo^k && \dTo_h \cr E &\rTo^{i_E} & \cA_\cL(E).
  \end{diagram}
Furthermore, $h$ and $k$ are strict filtered isomorphisms,
compatible
 with $\oh \ccG$-module structures as shown:
 \begin{eqnarray*}
h\colon \left(\cA_\cL(\cT_\cL(E)),\theta_\cA,\theta_{tot} \right)
 &\to&\left (\cA_\cL(E),\theta_{tot},\theta_{\cA}\right) \\
   k\colon \left (\cT_\cL(\cT_\cL(E)),\theta_\cT\right)
 &  \to& \left (E,\theta_E\right).
 \end{eqnarray*}
\item If $s$ is a  section of $\cL$,
then $s^*\circ j_E$ induces a strict isomorphism
$$\eta_s \colon (\cT_\cL(E),\theta_\cT) \to (E,\iota_*\theta_E),$$
where  $\iota$ is the inversion mapping of the group scheme $\ccG$.
\end{enumerate}
\end{proposition}
\begin{proof}

This result can be interpreted and proved in many ways. For example,
it is a special case of the theory of Higgs transforms on affine
group schemes as explained in section~\ref{higgs.s}. Here we give a
simpler version.  Indeed, all of the statements of the proposition
can be verified locally on $X$, and so we may and shall assume
without loss of generality that $\cL$ has a section $s$ defining an
isomorphism $\cL \cong \bT$ and hence $\cA_\cL \cong S^\cx \Omega$.

Note that the $\oh\ccG$-module structure  (\ref{camult.e}) on
$\cA_\cL$  and the map
$$s^* \colon \cA_\cL \to \oh X$$
defined by the section $s$ of $\cL$ together define a perfect
pairing
\begin{equation}\label{capair.e}
\oh {\ccG_n} \times N_n\cA_\cL \to \oh  X.
  \end{equation}
% In fact, independent of any section $\tF$,
% one gets a perfect pairing
% $$ \oh \ccG \ot \gr \cA_\cL \to \oh X,$$
% since $\oh \ccG \cong \gr \oh \ccG$.
If $E$ is any $\oh X$-module, let
$$\cH_\cx(\oh \ccG, E) :=  \dirlim \cHom_{\oh X}(\oh {\ccG_n},E)
\subseteq \cHom_{\oh X}( \oh\ccG, E).$$ Then the pairing
(\ref{capair.e})  defines an isomorphism:
\begin{equation}
  \label{eq:3}
\cA_\cL(E) := \cA_\cL \ot E \cong \cH_\cx(\oh \ccG, E),
\end{equation}
Let us denote by $\theta_\cA$ and $\theta_{tot}$ the $\oh
\ccG$-module structures on $\cH_\cx(\oh \ccG,E)$ deduced from the
corresponding structures on $\cA_\cL(E)$. These can be described
explicitly as follows. If $E_1$ and $E_2$ are two $\oh
\ccG$-modules, then $\cHom_{\oh X}(E_1,E_2)$  can be give an $\oh
\ccG \otimes \oh \ccG$-module structure by the rule
$$(b_1,b_2)(\phi)(e_1) := b_2\phi(b_1 e_1).$$
Then $\theta_{tot}$ corresponds to the $\oh \ccG$-structure induced
by $\mu_*$ and $\theta_{\cA}$ to the structure induced by $pr_{1*}$.
\footnote{In the systematic treatment in the appendix, we use
$\mu'_*$ instead of $\mu_*$.} The  total filtration $N_\cx$ of
$\cA_\cL(E)$ corresponds to the filtration $N_\cx$ of $\cH_\cx(\oh
\ccG,E)$ defined by
$$N_k\cH_\cx(\oh\ccG ,E) := \{ \phi : \phi (\cI^{[j]}) \subseteq N_{k-j}E\}.$$

Now if $E$ is a locally nilpotent $\oh \ccG$-module, let us consider
the following maps:
\begin{eqnarray*}
  i_E \colon E \to \cH_\cx(\oh \ccG,E) & \quad & i_E(e)(b) := i_\ccG^*(b) e \\
   \sigma \colon \cH_\cx(\oh \ccG, E) \to E & \quad & \phi  \mapsto \phi(1)\\
\end{eqnarray*}
The map $i_\ccG$ is the identity section of $\ccG$ and the map $i_E$
corresponds to the map $i_E$ defined in (\ref{ije.e}). Similarly the
map $\sigma$ corresponds to the map $\eta_s$ of statement (4) of
Proposition~\ref{higgsh.p}.  Note that $\sigma$ is compatible with
the filtrations and that
 $ \sigma \circ i_E = \id_E$.   This shows that
$i_E$ is injective and strictly compatible with the filtrations.
  The  image of $i_E$ is just the set of
homomorphisms which factor through $\oh {\ccG_0} = i_{\ccG*}(\oh
X)$, which corresponds to $(\cA(E))^{\theta_\cA}$.  This proves (1)
of Proposition~\ref{higgsh.p}, and (2) is a tautological consequence
of the definitions.

Define
$$  \tau \colon E \to \cH_\cx(\oh \ccG,E) \quad
\mbox{by} \quad \tau(e)(b) := \iota^*(b)e,$$ where $\iota \colon
\ccG \to \ccG$ is the inverse mapping in the group $\ccG$. If $e \in
E$, then \emph{a priori} $\tau(e)$ is  just an element of
$\cHom_{\oh X}(\oh \ccG, E)$, but if $e$ is annihilated by
$I^{[n+1]}$ then $\tau(e) \in \cHom(\oh {\ccG_n}, E)$.  Thus $\tau$
is well-defined if $E$ is locally nilpotent. Note that $\sigma\circ
\tau   = \id_E$, so $\tau$ is also injective. If $e \in N_kE$ and $b
\in \cI^{[j]}$, then $\iota^*(b)e \in N_{k-j}E$, so $\tau(e) \in
N_k\cH_\cx(\oh \ccG,E)$. Thus $\tau$ is compatible with the
filtrations, and in fact is strictly compatible because $\sigma$ is
also compatible.

 It is clear that
that the image of $\tau$ consists precisely of the elements of
$\cH_\cx(\oh \ccG,E)$ which are $\iota$-linear.   We claim that
these are the elements which correspond to elements of  $\cT_\cL(E)
\subseteq \cA_\cL(E)$. Indeed, if $\phi \colon \oh \ccG \to E$ is
$\iota$-linear then it follows from the commutativity of the diagram
\begin{diagram}
 \ccG & \rTo^{pr} & X \cr
\dTo^{(\id,\iota)} && \dTo_{i_\ccG} \cr
  \ccG \times \ccG & \rTo^\mu & \ccG \cr
\end{diagram}
that $\theta_{tot}(\phi) = 0$, and the converse follows from the
fact that the diagram is Cartesian.

Thus we can write $\tau = j_E \circ \ov \tau$, where $\ov \tau
\colon E \to \cT_\cL(E)$ is an isomorphism of $\oh X$-modules,
inverse to the mapping $\sigma \circ j$. It is clear from the
definitions of $\tau$ and $\theta_\cA$ that $\tau$ takes $\theta_E$
to $\iota_*\theta_\cA$, and this proves (4) of
Proposition~\ref{higgsh.p}.

It remains for us the prove statement (3). First let us check that
$h$ is compatible with the Higgs fields as described there. As we
have observed in  equation~(\ref{ije.e}), $j_E$ takes $\theta_0$ to
$\theta_{tot}$.  More precisely, but perhaps somewhat cryptically:
$\theta_{tot} \circ j_E = j_E\circ \theta_0$, where for example we
are writing $\theta_{tot}$ for the endomorphism of $E \ot \cA_\cL$
induced by some element of $\oh \ccG$ corresponding to the $\oh
\ccG$-module structure given by $\theta_{tot}$.
\begin{eqnarray*}
  (\theta_{tot}\cnv \theta_\cA) \circ (j_E \ot \id_\cA) & = &
    (j_E \ot \id_\cA) \circ (\theta_0\cnv \theta_\cA) \\
(\id_E \ot m)\circ (\theta_{tot}\cnv \theta_\cA) \circ (j_E \ot
\id_\cA) & = &
   (\id_E \ot  m)\circ (j_E \ot \id_\cA) \circ (\theta_0\cnv \theta_\cA) \\
\theta_{tot}\circ (\id_E \ot m)  \circ (j_E \ot \id_\cA) & = &
   (\id_E\ot  m)\circ (j_E \ot \id_\cA) \circ \theta_\cA \\
\theta_{tot} \circ h & = &
    h \circ \theta_\cA.
  \end{eqnarray*}
Similarly, $(\theta_0\cnv \theta_\cA) \circ j_E = j_E \circ
\theta_\cT$, so
\begin{eqnarray*}
(\theta_0\cnv \theta_\cA \cnv \theta_\cA) \circ (j_E \ot \id_\cA)
      &=& (j_E \ot \id_\cA) \circ (\theta_\cT \cnv \theta_\cA)\\
(\id_E \ot m) \circ (\theta_0\cnv \theta_\cA \cnv \theta_\cA) \circ
(j_E \ot \id_\cA)
      &=&(\id_E \ot m)\circ (j_E \ot \id_\cA) \circ (\theta_\cT \cnv \theta_\cA)\\
 (\theta_0\cnv \theta_\cA) \circ (\id_E \ot m) \circ (j_E \ot \id_\cA)
      &=&h \circ (\theta_\cT \cnv \theta_\cA)\\
\theta_\cA \circ h &= & h \circ \theta_{tot}.
\end{eqnarray*}
Thus $h$ takes $\theta_\cA$ to $\theta_{tot}$ and $\theta_{tot}$ to
$\theta_\cA$ as claimed. Since $h$ takes $\theta_{tot}$ to
$\theta_\cA$, it induces the  map $k$:
$$\cT_\cL(\cT_\cL(E)) :=\left(\cA_\cL(\cT_\cL)(E)\right) ^{\theta_{tot}}
 \to \cA_\cL(E)^{\theta_\cA} = E. $$

  Let us check that $h$ is compatible with the filtrations.
By   definition,
$$N_c^{tot} \cA_\cL(\cT_\cL(E)) =
 \sum_a \left (Im(N_a\cT_\cL(E) \ot
N_{c-a} \cA_\cL \to \cT_\cL(E) \ot \cA_\cL(E)\right).$$
 The definition of $N_a\cT_\cL(E)$ shows
that its image under $j_E$ is contained in the sum of the images of
$N_bE\ot N_{a-b}\cA_\cL$.  Hence $h$ maps
$N_c^{tot}\cA_\cL(\cT_\cL(E))$ into the sum of the images of
$$N_b E \ot N_{a-b}\cA_\cL \ot N_{c-a} \cA_\cL \rTo^m E \ot \cA_\cL(E),$$
which is contained in $N^{tot}_c \cA_\cL(E)$.

Note that if $h$  is a strict isomorphism, then it induces a strict
isomorphism from the annihilator of $\theta_{tot}$ to the
annihilator of $\theta_\cA$, \ie, from $\cT_\cL(\cT_\cL(E))$ to $E$.
Thus  $k$ is also a strict isomorphism.

Thus it remains only to show that $h$ is a strict isomorphism.
Suppose first that $\gr_i E = 0$ for all $i \neq k$. Then the $\oh
\ccG$-structure on $E$ factors through $i_\ccG$, so $\theta_\cA =
\theta_{tot}$ and $E = \cT_\cL(E) \subseteq \cA_\cL(E)$.  Then the
map $h$  is:
$$e \ot a \mapsto e\ot 1 \ot a \mapsto e \ot a,$$
\ie, the identity map.  Now we can proceed by d\'evissage.
Statement (4) shows that the the functor $E \mapsto \cT_\cL(E)$
preserves strict exact sequences, and since $\gr \cA_\cL$ is locally
free, the same is true of the functors $E \mapsto  \cA_\cL(E)$ and
$E \mapsto \cA_\cL(\cT_\cL(E))$. Now suppose that $N_aE = 0 $ and
that $h$ induces an isomorphism for  $E' := N_{b-1}E$.  Then we have
a strict exact sequence
$$ 0 \to E' \to N_{b}E \to E'' \to 0,$$
where $E'' := N_b E/N_{b-1}E$.  We have seen above that the theorem
is true for $E''$, and it holds for $E'$ by the induction
assumption. Then it also holds for $N_bE$ by the strict exactness of
the functors $\cA_\cL(\ )$ and $\cA_\cL(\cT_\cL(\ ))$. It follows by
induction that $h$ is a strict isomorphism whenever the filtration
on $E$ is bounded, and, by taking direct limits, whenever the
filtration is bounded below and exhaustive. This completes the
proof.
% Since the filtrations of $E$ and $\cA_\cL$ are bounded
% below, the same is true of the total filtration on $\cA_\cL(E)$,
% and so it will suffice to prove that $\gr(h)$ is an isomorphism.
% To prove this, let us note first that, since
% $\gr \cA_\cL$ consists of locally free $\oh X$-modules,
% the natural map
% \begin{equation}
%   \label{grea.e}
% \gr (E) \ot \gr (\cA_\cL) \to \gr \cA_\cL(E)
% \end{equation}
% is a graded isomorphism.
% \marginpar{  The maps
% $\gr \tau$ and  $\gr i_E$:
% $ \gr E \to \gr \cA_\cL(E)$ if
% we use the grading coming from $\cA$
% or the grading from $E$ (not the
% total grading}
%    Since $\gr \cA_\cL$ is locally
% free,  the functors $\gr \cA_\cL(E)$ and $\gr _\cA(\cT(E))$
% also preserve strict exact sequences.
% Now it is now easy to deduce the theorem by d\'evissage.
% Then it follows by induction that $\gr h$ is an
% isomorphism whenever the filtration on $E$ is finite,
% and then in general by taking direct limits.
\end{proof}

\begin{remark}
  It is easy to see that the filtration of
$\cT_\cL(E)$ induced by the total filtration $N_{tot}$ on
$\cA_\cL(E)$ is the same as the filtration induced by the filtration
$\cA_\cL \ot N_\cx E$.  The total filtration has the advantage of
being again $\cI^\cx$-saturated, a fact we will exploit in our
cohomology computations in the next section.
\end{remark}

\begin{remark}\label{fhiggs.r}
  A similar result holds for standard Higgs fields
if one works with the divided power completion of $\cA_\cL$ along
the ideal of a section.  More abstractly, suppose that $T$ and
$\Omega$ be as above, let $\theta \colon E \to E \ot \Omega$ be a
locally nilpotent $T$-Higgs field on $E$.  Let $I$ be the ideal of
the symmetric algebra $S^\cx T$ generated by $ T$. Then an
$I$-saturated filtration on $E$ is just a filtration $N$ such that
$IN_kE \subseteq N_{k-1}E$. Let $\cA_\Omega$ be the divided power
algebra $\Gamma_\cx \Omega$, and define $\cA_\Omega(E) := E \ot
\cA_\Omega$ and $\cT_\Omega(E) := \left (
\cA_\Omega(E)\right)^{\theta_{tot}}$. Then the evident analog of
Proposition~\ref{higgsh.p} holds.
\end{remark}

We will sometimes want to consider \emph{graded} Higgs fields and
PD-Higgs modules, \ie, graded  modules over the graded ring $\oh
\ccG$, where  $\oh \ccG = \Gamma_\cx T$ or $S^\cx T$. There is an
evident functor $\gr$ from the category of $\cI^\cx$-filtered (resp.
$\cI$-filtered) modules to the category of graded $\oh
\ccG$-modules, compatible with the convolution tensor product.  In
particular, if $\cL$ is a $\bT$-torsor, then $\gr \cA_\cL \cong
S^\cx \Omega = S^\cx_\bT$, as a graded $\Gamma_\cx T$-modules (note
that the multiplication sends $\Gamma_a T \ot S^b\Omega$ to
$S_{b-a}\Omega$); furthermore its divided power envelope $\Gamma_\cx
\Omega$ is in a natural way a graded $S^\cx \Omega$-module.  If $E$
is an $\cI^\cx$-filtered $\oh \cG$-module,  the natural map
$$ \gr E \ot \cA_\bT \cong \gr E \ot \gr \cA_\cL \to \gr \cA_\cL(E)  $$
is an isomorphism, since $\gr \cA_\cL$ is locally free over $\oh X$,
and it is compatible both with $\theta_\cA$ and $\theta_{tot}$.   In
particular, if $E$ satisfies the hypothesis of
Proposition~\ref{higgsh.p}, the map $\cT_\cL(E) \to \cA_\cL(E)$
induces a map
$$\gr \cT_\cL(E) \to \gr(\cA_\cL(E)) \cong \gr(E) \ot  \cA_\bT(E)$$
whose image is annihilated by $\theta_{tot}$ and hence induces a map
$$\gr \cT_\cL(E) \to \cT_\bT(\gr(E))$$

\begin{corollary}\label{grhigg.c}
Let $(E,\theta,N_\cx)$ be an $\oh X$-module $E$ equipped with an
$\oh \ccG$-module structure $\theta$ and an $\cI^\cx$-filtration
$N_\cx$, as in Proposition~\ref{higgsh.p}.
\begin{enumerate}
\item The map
$\gr \cT_\cL(E) \to \cT_\bT(\gr(E))$ above is an isomorphism. In
fact there is a commutative diagram of isomorphisms:
\begin{diagram}
\cT_\bT(\gr E)\ot\cA_\bT &\lTo& \gr \cT_\cL(E) \ot  \cA_\bT & \rTo &
\gr \cA_\cL(\cT_\cL(E)) \cr
   \dTo^h && \dTo && \dTo_{\gr(h)}\cr
\gr(E)\ot \cA_\bT & \lTo^{\id} &  \gr E  \ot  \cA_\bT & \rTo& \gr
\cA_\cL(E),
\end{diagram}
compatible with the $\oh \ccG$-module structures as in \emph{op.
cit.}.
\item There is a natural isomorphism of graded $\oh \ccG$-modules:
$$\gr \cT_\cL(E) \cong \iota_* \gr E.$$
\end{enumerate}
\end{corollary}
\begin{proof}
  The existence and the commutativity of the diagram is clear, as is the fact
that the arrows are compatible with the $\oh \ccG$-module
structures.  Furthermore, it follows from Proposition~\ref{higgsh.p}
that $h$ and $\gr (h)$ are isomorphisms. It follows that the middle
vertical arrow is an isomorphism, and that the image of $\gr
\cT_\cL(E)$ in $\gr E \ot \cA_\bT$ is exactly the annihilator of
$\theta_{tot}$, \ie, $\cT_\bT(\gr(E))$. This proves (1). Then (2)
follows by applying (4) of Proposition~\ref{higgsh.p} with $E$
replaces by $\gr E$ and $\cL$ replaced by $\bT$.
\end{proof}

There is a useful cohomological complement to  the construction of
Remark~\ref{fhiggs.r}. Recall that associated to a $T$-Higgs module
$(E,\theta)$ is its Higgs (Koszul) complex
$$ E \to E \ot \Omega \to E \ot \Omega^2 \to \cdots,$$
where $\Omega^i := \Lambda^i \Omega$.

\begin{proposition}\label{cohiggs.p}
  Let $E$ be a graded $\oh X$-module with a graded
$T$-Higgs field $\theta$:
$$\theta \colon E \to E\ot \Omega,$$
where $\Omega$ is in degree $1$.   Using the notation of
Remark~\ref{fhiggs.r}, let
$$ \cA_\Omega^{i,j}(E) := \cA_\Omega (E) \ot \Omega^i\ot \Omega^j$$
and  let
\begin{eqnarray*}
d_\cA \colon:\cA_\Omega^{i,j}(E)  &\to& \cA_\Omega^{i+1,j}(E)  \\
d_{tot} \colon  \cA_\Omega^{i,j}(E) &\to& \cA_\Omega^{i,j+1}(E)
\end{eqnarray*}
be the boundary maps associated to the fields $\theta_\cA$ and
$\theta_{tot}$ respectively, tensored with the identity. Then these
maps fit into a graded double complex $\cA_\Omega^{\cx\cx}(E)$, and
the maps $i$ and $j$ of  Remark~\ref{fhiggs.r} define augmentations
of the double complex
\begin{eqnarray*}
  E \ot \Omega^\cx &\to &\cA_\Omega^{\cx\cx}(E)     \cr
  \cT_\Omega(E) \ot \Omega^\cx &\to &\cA_\Omega^{\cx\cx}(E)
\end{eqnarray*}
For each $i$, $\cA_\Omega^{i,\cx}(E)$ is a resolution of $E\ot
\Omega^i$ and for each $j$ $\cA_\Omega^{\cx,j}(E)$ is a resolution
of $\cT_\Omega(E)\ot \Omega^j$
\end{proposition}
\begin{proof}
  It is immediate to verify that the boundary maps commute
and hence define a double complex.  The fact that
$\cA_\Omega^{i,\cx}(E)$ is a  graded resolution of $E\ot \Omega^i$
follows from the filtered Poincar\'e lemma~\cite[6.13]{bo.ncc}
 for the divided power algebra
$\Gamma_\cx(\Omega)$.  Since $h$ is an isomorphism  transforming
 $\theta_{tot}$ into $\theta_{\cA}$, the second statement follows.
\end{proof}

Let us now return to our discussion of the  Cartier transform.
Recall that the center of $D^\gamma_{X/S}$ can be identified with
the divided power algebra $\Gamma_\cx T_{X'/S}$.  Let $\cI^{\cx}_X$
denote the divided power filtration of the divided power ideal
$\cI_X$ of $\Gamma_\cx T_{X'/S}$. Let $\micn_\gamma(X/S)$ denote the
category of $D^\gamma_{X/S}$-modules $E$ equipped with an
exhaustive, horizontal, and  bounded below filtration
$\cI^\cx$-filtration $N_\cx$
 (see Definition~\ref{filtsat.d}).
Similarly, let $\hign_\gamma(X'/S)$ denote the category of
$\Gamma^\cx T_{X'/S}$-modules $E'$  equipped
 with an exhaustive and bounded below
$\cI^\cx$-saturated filtration $N'_\cx$.

If $(E,\nabla,N_\cx)$ is an object of $\micn_\gamma(X/S)$, let
$$E^{\nabla,\gamma}:= \cHom_{D^\gamma_{X/S}}(\oh X, E), \mbox{ and }$$
$$E^\nabla := \Ker(E \rTo^\nabla E \ot \Omega^1_{X/S}).$$
The action of the center $\Gamma_\cx T_{X'/S}$ of $D^\gamma_{X/S}$
defines a PD-Higgs field on $F_{X/S*}E$ and hence an F-PD Higgs
field $\psi$ on $E$; note that $E^\psi$
 is invariant under
the connection $\nabla \colon E \to E\ot \Omega^1_{X/S}$.
Furthermore, $E^{\nabla,\gamma} = {E^\psi}^\nabla$, since
$D_{X/S}^\gamma$ is generated as a topological ring by  $T_{X/S}$
and   $ \Gamma_\cx T_{X'/S} $.

We endow $\cA_{\cL}(E) := E \ot \cA_\cL$ with the tensor product
$D_{X/S}^\gamma$-module structure $\nabla$ coming from the given
structures on $E$ and on $\cA_{\lift X S}$ and  with the tensor
product filtration coming from the filtrations $N_\cx$ of $E$ and
$\cA_{\lift X S}$. We also endow it with the F-PD-Higgs field
$\theta_\cA := \id \ot \theta_\cA$, where $\theta_\cA$ is the
F-PD-Higgs field of $\cA_{\lift X S}$.

\begin{lemma}\label{thetacomm.l}
The action $\theta_\cA$ of $\Gamma_\cx(T_{X'/S})$   on $\cA_{\lift X
S}(E)$ commutes with the action of $D_{X/S}^\gamma$ corresponding to
the tensor product $D_{X/S}^\gamma$-module structure.
\end{lemma}
\begin{proof}
 As we have
already observed, it follows from the formula (\ref{thetaconv.e})
that the $p$-curvature PD-Higgs field of  $\cA_{\lift X S}(E)$
commutes with $\id\ot \theta_\cA$. That is, the action of
$\Gamma_\cx T_{X'/S} \subseteq D_{X/S}^\gamma$ commutes with $\id\ot
\theta_\ca$.
 Furthermore, if $D \in T_{X/S}$,
$\xi' \in\Gamma_\cx T_{X'/S}$, $e \in E$, and $\alpha \in \cA_{\lift
X S}$,
\begin{eqnarray*}
(\id \ot {\theta_{\xi'}}) \nabla_D(e\ot a)
          & = & \nabla_D(e) \ot \theta_{\xi'}(a) + e\ot\theta_{\xi'}\nabla_D(a) \\
            & = &\nabla_D(e) \ot \theta_{\xi'}(a) + e\ot\nabla_D\theta_{\xi'}(a) \\
           & = &\nabla_D (\id \ot_{\theta_{\xi'}}) (e\ot a)
\end{eqnarray*}
Since $D_{X/S}^\gamma$ is generated by  $T_{X/S}$ and $\Gamma_\cx
T_{X'/S}$, it follows that $\nabla_D$ commutes with $\theta_{\xi'}$
for every $D \in D_{X/S}^\gamma$.
\end{proof}

Now recall  that, by definition, $\iota_*C_{\lift X S}(E) =
\cHom_{D_{X/S}^\gamma}(\cB_{\lift X S},E)$, with the $\Gamma_\cx
T_{X'/S}$-module structure coming from $E$, where $\cB_{\lift X S}
:= \cHom_{\oh X}(\cA_{\lift X S}, \oh X)$ in the category of
$D_{X/S}^\gamma$-modules. Thus when $E$ is locally nilpotent,
\begin{equation*}
\iota_*C_{\lift X S}(E) :=\cHom_{D^\gamma_{X/S}}(\cB_{\lift X S}, E)
\cong \left (E \ot \cA_{\lift X S} \right)^{\nabla, \gamma} =
 \left((\cA_{\lift XS}(E))^{\psi}\right)^\nabla.
\end{equation*}
Of course, the total PD-Higgs field on $C_{\lift X S}(E)$ is zero,
but because of the commutation of $D_{X/S}^\gamma$ and $\id\ot
\theta_{\cA}$, $\iota_*C_{\lift X S}(E)$ is stable under the
PD-Higgs field $\id \ot \theta_\cA $ of $\cA_{\lift X S}$. In fact,
the induced PD-Higgs field induced by $\theta_\cA$ on is $\iota_*$
of the PD-Higgs field induced by $\theta_E$. A geometric explanation
of this fact is given in the appendix after
Definition~\ref{higtr.d}; it can also be checked by direct
computation. Thus it follows that
\begin{equation}\label{ctxs.e}
C_{\lift X S}(E) := \iota_*\cHom_{D_{X/S,\gamma}}(\cB_{\lift X S},E)
\cong \left (\cT_{\lift X S}(E)\right)^\nabla
\end{equation}
where $\cT_{\lift X S}(E) :=  (\cA_{\lift X S}(E))^{\theta_{tot}}$
as in Proposition~\ref{higgsh.p}.

It is clear from the construction that there are natural maps,
compatible with the  connections and F-PD-Higgs fields shown:
\begin{equation}\label{ija.e}
  (E,\theta_0,\nabla)  \rTo^i (\cA_{\lift X S}(E),\theta_\cA,\nabla)
 \lTo^j (F_{X/S}^*C_{\lift X S}(E),\theta,\nabla_0)
  \end{equation}
Here $\nabla_0$ is the Frobenius descent connection on
$F_{X/S}^*C_{\lift X S}(E)$. Since  $N_\cx$ is an
$\cI^\cx$-filtration on $E$, the filtration on $F_{X/S}^*
C_{X/S}(E)$ induced  by the total filtration of $\cA_{\lift X S}(E)$
is horizontal and is also an $\cI^\cx$-filtration with respect to
the action of $\theta_\cA$.  It follows that it descends to an
$\cI^\cx$-filtration on $C_{\lift X S}(E)$. Thus we obtain a
filtered version of the Cartier transform:
\begin{equation}\label{filcart.e}
  C_{\lift X S} \colon \micn^\cx_\gamma(X/S) \to \hign^\cx_\gamma(X'/S).
  \end{equation}

On the other hand, if $(E',\theta', N_\cx')$ is an object of
$\hign^\cx_\gamma(X'/S)$,  we can endow
\begin{equation}\label{calp.d}
\cA'_{\lift X S}(E') := F_{X/S}^*E' \ot_{\oh X} \cA_{\lift X S}
\end{equation}
with the tensor product F-PD-Higgs field $\theta'_{tot}$.  It
follows as in Lemma~\ref{thetacomm.l} that $\theta'_{tot}$  commutes
with the tensor product connection on $F_{X/S}^*E' \ot \cA_{\lift X
S}$, where $F_{X/S}^*E'$ is given the Frobenius descent connection
$\nabla_0$. Thus
\begin{equation}
C'_{\lift X S}(E'):= (\cA'_{\lift X S}(E'))^{\theta'_{tot}}
\end{equation}
inherits a  nilpotent $D_{X/S}^\gamma$-module structure from
$\cA_{\lift X S}$, which we denote by $\nabla_{\cA'}$. We have
natural maps
\begin{equation}\label{ijap.e}
  F_{X/S}^*(E',\theta',\nabla_0)  \rTo^{i'} (\cA'_{\lift X S}(E'),\theta'_{tot},\nabla_{\ca'})
 \lTo^{j'} (C'_{\lift X S}(E'),\theta'_{tot},\nabla_\cA)
  \end{equation}
  where  $i'$  takes the PD-Higgs field of $E'$
%to the total Higgs field of $\cA_{\lift X S}(E')$ and
%$j'$ is horizontal.
 As before, the $\cI^\cx$-filtration  $N'_\cx$ on
$E'$ induces an $\cI^\cx$-filtration on $C'_{\lift X S}(E')$, and we
get a functor:
\begin{equation}\label{filcartp.e}
  C'_{\lift X S} \colon \hign^\cx_\gamma(X'/S) \to \micn^\cx_\gamma(X/S) .
  \end{equation}

The commutative diagram
  \begin{diagram}
 F_{X/S}^*C_{\lift  X S}(E) \ot \cA_{\lift X S}  & \rTo^{j\ot \id}
        & E \ot \cA_{\lift X S} \ot \cA_{\lift X S} \cr
\dTo^h& \ldTo_{\id \ot m} \cr
 E\ot \cA_{\lift X S}
   \end{diagram}
defines a horizontal map
\begin{equation}\label{hdef.e}
  h \colon (\cA'_{\lift X S}(C_{\lift X S}(E)),\theta_\cA,\theta_{tot},N_\cx)
 \to (\cA_{\lift X S}(E),\theta_{tot},\theta_\cA,N_\cx).
  \end{equation}
A similar construction defines a horizontal map
\begin{equation}\label{hpdef.e}
  h' \colon (\cA_{\lift X S}(C'_{\lift X S}(E')),\theta'_{\cA'},\theta'_{tot},N'_\cx) \to
 (\cA'_{\lift X S}(E'),\theta'_{tot},\theta'_\cA,N'_\cx)
  \end{equation}

\begin{theorem}\label{rh.t}
Let $\lift X S := (X/S,\tX'/\tS)$ be a smooth morphism with a
lifting of $X'$ mod $p^2$ as described above.
\begin{enumerate}
\item{  Let  $(E,\nabla,N_\cx)$ be an object of $\micn^\cx_\gamma(X/S)$ and let
$(E',\theta',N'_\cx): =  C_{\lift X S}(E,\nabla,N_\cx) $. Then the
map $h$ (\ref{hdef.e}) is a filtered isomorphism, and fits into a
commutative diagram:
\begin{diagram}
C'_{\lift X S}(E') & \rTo^{j'} &  \cA'_{\lift X S}(E') & \lTo^{i'}&
E'\cr \dTo^\cong &&\dTo^h && \dTo_{\id}\cr E& \rTo^i &   \cA_{\lift
X S}(E) & \lTo^j & C_{\lift X S}(E)
\end{diagram}}
\item{Let $(E',\theta',N'_\cx)$ be an object of $\hign^\cx_\gamma(X'/S)$, and let
$(E,\nabla,N_\cx)  := C^{-1}_{\lift X S}(E',\theta', N'_\cx) $. Then
the map $h'$ (\ref{hpdef.e}) is a filtered isomorphism and fits into
a commutative diagram:
\begin{diagram}
   C_{\lift X S}(E)& \rTo^{j} &    \cA_{\lift X S}(E)  & \lTo^i & E \cr
\dTo^\cong &&\dTo^{h'}  && \dTo_{\id} \cr E' &\rTo^{i'}& \cA'_{\lift
X S}(E') & \lTo^{j'} & C'_{\lift X S}(E')
\end{diagram}}
\end{enumerate}
Consequently, $C'_{\lift X S}$ is quasi-inverse to the Cartier
transform $C_{\lift X S}$ and is therefore isomorphic to the functor
$C^{-1}_{\lift X S}$ of Theorem~\ref{cart.t} (ignoring the
filtrations).
% Moreover, an object  $E$ of
% $\mic_\gamma(X/S)$ is nilpotent of level at most $\ell$ if and
% only if  the natural map
% $$$\left (E \ot N_\ell\cA_{\lift X S} \right)^{\nabla,\theta}
% $$\to \left(E \ot \cA_{\lift X S }\right)^{\nabla,\theta}$$
% is an isomorphism.
\end{theorem}
\begin{proof}
This theorem is an immediate consequence of
Proposition~\ref{higgsh.p} and Cartier descent. The $p$-curvature of
the connection on $\cA_{\lift X S}(E)$ is the total Higgs field
$\theta_{tot}$. Hence
$$E' := C_{\lift X S}(E) = \left ((\cA_{\lift X S}(E))^{\theta_{tot}}\right)^\nabla =
\cT_{\lift X S}(E)^\nabla,$$ in the notation of \emph{op. cit.}.
Since  the $p$-curvature of the connection
 $\cT_{\lift X S}(E)$ vanishes,
standard Cartier descent implies that the natural map
$$F_{X/S}^*E'  \to \cT_{\lift X S}(E)$$ is
a filtered isomorphism.  Thus we have a commutative diagram
\begin{diagram}
F_{X/S}^*E'\ot \cA_{\lift X S} & \rTo^\cong
     &\cT_{\lift X S}(E) \ot_{\oh X}\cA_{\lift X S}  \cr
&\rdTo_h & \dTo_{\tilde h} \cr
 & & E\ot \cA_{\lift X S}
\end{diagram}
Proposition~\ref{higgsh.p} implies  that $\tilde h$ is a filtered
isomorphism and hence so is $h$, and it is also horizontal. The
vertical left arrow in the diagram of (1) corresponds to the map $k$
of Proposition~\ref{higgsh.p} and is also a horizontal filtered
isomorphism, compatible with the PD-Higgs fields, \ie, an
isomorphism in the category $\micn^\cx_\gamma(X/S)$. A similar
argument works if we start with an object $(E',\theta', N'_\cx)$ of
$\hign^\cx_\gamma(X'/S)$. This shows that $C_{\lift X S}$  and
$C'_{\lift X S}$ are quasi-inverse equivalences.
\end{proof}

\begin{corollary}\label{grcart.c}
Let $(E,\nabla,N_\cx)$ be an object of $\micn^\cx_\gamma(X/S)$
 and let
$$(E',\theta',N'_\cx):= C_{\lift X S}(E,\nabla,N_\cx).$$
Then there is a natural isomorphism in the category of graded $
\Gamma^\cx T_{X'/S}$-modules:
$$ (\gr (E',\theta',N'_\cx)) \cong
\iota_* \left( \gr (E, \psi,N_\cx)\right )^\nabla,$$ where $\psi$ is
the action of $\Gamma^\cx_{T_{X'/S}} \subseteq  D^\gamma_{X/S}$ and
$\iota$ is the inversion involution of $\Gamma^\cx_{T_{X'/S}}$.
\end{corollary}
\begin{proof}
Using Corollary \ref{grhigg.c},  we have
 \begin{eqnarray*}
   \gr E'& \cong &
                \gr\left (\left ( \cT_{\lift X S}(E) \right) ^\nabla  \right) \\
   & \cong & \left (\gr  \cT_{\lift X S}(E)   \right)^\nabla \\
    & \cong & \left (\iota_* \gr (E) \right)^\nabla\\
\end{eqnarray*}
\end{proof}

%The statement about the level follows from
%Propostion~\ref{lev1.p}.

% As explained  in (\ref{higtra.t}), the multiplication map of the algebra $\cA_{\lift X S}$
% induces compatibility isomorphisms for the tensor structures.
%   For any two objects $E_1$ and $E_2$ of $\mic^\cx_\gamma(X/S)$, there
% is a commutative diagram
% \begin{diagram}
% C_{\lift X S} (E_1) \ot C_{\lift X S}(E_2) &\rTo^a& C_{\lift X S}(E_1\ot E_2) \cr
% \dTo && \dTo \cr
% \cA_{\lift X S} \ot E_1\ot \cA_{\lift X S} \ot E_2 &\rTo^b& \cA_{\lift X S}\ot E_1\ot E_2
% \end{diagram}
% Here $b$ is induced by the multiplication map  $\cA_{\lift X S}\ot  \cA_{\lift X S} \to
% \cA_{\lift X S}$, and the vertical arrows are injective.  Then one can check the existence
% of the map $a$ locally, as well as the fact that it is an isomorphism.  This is then immediate.
% %In the case of finite level,% we recall that $\mic_n(X/S)$ is stable under tensor product

\begin{remark}\label{higgh.r}
A similar formalism works when there is a lifting $\tF$ of
$F_{X/S}$. Let $\micn(X/S)$ denote the category of modules with
connection $(E,\nabla)$ endowed with a horizontal filtration $N_\cx$
such that $\gr^N(E)$ is constant.  We assume also that $N_\cx$ is
exhaustive and bounded below. As before,  let $\cA_\tF$ be the
nilpotent divided power
 completion  of  $\cA_{\lift X S}$ along the
ideal of the corresponding augmentation $\cA_{\lift X S} \to \oh X$.
Then if $(E,\nabla,N_\cx)$ is an object of $\micn(X/S)$, its
$p$-curvature $\psi$ gives $(F_{X/S*}E,N_\cx)$ an $I$-saturated
Higgs field as discussed in Remark~\ref{fhiggs.r}.  Then we define:
$$\cA_\tF(E) := E \ot_{\oh X} \cA_\tF \quad\mbox{and} \quad
\cT_\tF(E) := \left ( \cA_\tF(E)\right)^{\theta_{tot}}$$ where
$\cT_\tF(E)$ has the Higgs field $\theta_\cT$ induced by
$\theta_\cA$.  Then
$$ C_{\tF}(E) := (\cA_{\tF}(E))^{ \nabla} = (\cT_\tF(E))^\nabla$$
with it inherits a Higgs field  and filtration.
  Thus we obtain a functor
$$C_{\tF} \colon \micn(X/S) \to \hign(X'/S). $$ %\quad (E,\nabla,\theta) \mapsto (E',\theta').,$$
  On the other hand, if
$(E',\theta',N'_\cx)$ is an object of $\hign(X'/S)$,  let
$$\cA_{\tF}(E') := E'\ot_{\oh X'}\cA_{\tF}.$$  Then the total Higgs field $\theta'$
on $\cA_{\tF}(E')$ commutes with the connection $\id \ot
\nabla_{\cA}$.  Let
$$C^{-1}_{\tF} := (\cA_{\tF}(E'))^{\theta'},$$
 which  inherits a connection from the action of $\id\ot\nabla_{\cA}$
and a filtration $N_\cx$  from the total filtration $N'_{tot}$. Thus
$C^{-1}_{\tF}$ is a functor
$$C^{-1}_{\tF} \colon \hign(X'/S) \to  \micn(X/S).$$
These functors are quasi-inverse equivalences, compatible with the
tensor structures and with the global functors $C_{\lift X S}$
considered above.
\end{remark}

\subsection{De Rham and Higgs cohomology}\label{cohomology.ss}
Let us continue to denote by   $\lift X S$ a smooth morphism  $X/S$
of schemes in characteristic $p$, together with a lifting $\tilde
X'/\tilde S$ of $X'/S$.
 Let $(E,\nabla)$ be a module with integrable
connection on $X/S$, nilpotent of level $\ell$.
 Our goal in this
section is to compare the de Rham cohomology of $(E,\nabla)$ with
the Higgs cohomology of  its Cartier transform $(E',\theta')$.
 We shall do this by constructing
a canonical filtered double complex  $(\cA^{\cx\cx}_{\lift X S}(E),
N_\cx)$ of $\oh {X'}$-modules and quasi-isomorphisms
$$F_{X/S*}(E\ot \Omega^\cx_{X/S},d)  \rTo
 N_n\cA^\cx_{\lift X S}(E) \lTo( E'\otimes \Omega^\cx_{X'/S},\theta'),$$
whenever $\ell +d \le n < p$, where $\cA^\cx_{\lift X S}$ is the
total complex associated to the double complex $\cA^{\cx\cx}_{\lift
X S}$.

In fact,
 \begin{eqnarray*}
\cA_{\lift X S}^{ij}(E)
 &:= & F_{X/S*}\left (E \ot \cA_{\lift X S} \ot F_{X/S}^*\Omega^i_{X'/S}\ot\Omega^j_{X/S} \right)\\
 & \cong & F_{X/S*}\left (E \ot \cA_{\lift X S} \ot\Omega^j_{X/S} \right)\ot \Omega^i_{X'/S}.
    \end{eqnarray*}
with boundary maps constructed from the de Rham differentials of
$(E,\nabla)$ and the $p$-curvature of $\cA_{\lift X S}$. In the case
$(E,\nabla)= ({\cal O}_X,d)$ we obtain an isomorphism in the derived
category
$$F_{X/S*}(E\ot \Omega^\cx_{X/S},d) \sim (\Omega^\cx_{X'/S},0)$$
between the de Rham complex of $X/S$ and the Hodge complex of
$X'/S$, when $d < p$.  This is the
 result of  Deligne and Illusie \cite{di.rdcdr}
(with a loss of one dimension). For general $E$ it can be regarded
as an analog of Simpson's ``formality'' theorem~\cite{si.hbls}.

We shall find it convenient to work with filtered connections and
their de Rham complexes.
    Let $(E,\nabla)$ be a module with integrable
connection endowed with a horizontal filtration $N_\cx$ such that
$(\gr_N(E),\nabla)$ is constant, \ie, has zero $p$-curvature. We
assume that $N_{-1}E = 0$ and $N_{p-1}E = E$, so that
$(E,\nabla,N_\cx)$ defines an object of $\micn^\cx_\gamma(X/S)$. Let
$N^{tot}_\cx$ be the tensor product filtration on $E\ot \cA_{\lift X
S}$ induced by $N_\cx$ and the filtration $N_\cx$ of $\ca_{\lift X
S}$. Let $(E',N'_\cx)$ be the Cartier transform of $(E,N_\cx)$ with
the filtration induced by $N_\cx^{tot}$, as explained in
Theorem~\ref{rh.t}. For fixed $i$, the de Rham complex of the module
with connection ${\cA_{\lift X S}(E) \ot F_{X/S}^*\Omega^i_{X'/S}}$
is the complex:
\begin{equation}\label{cald.e}
 \cA_{\lift X S}^{i,0} (E) \rTo^{d^{i,0}} \cA_{\lift X S}^{i,1}(E)  \rTo^{d^{i,1}}  \cdots.
\end{equation}
Similarly, for fixed $j$, the Higgs complex of $(\cA_{\lift X
S},\theta_\cA)$ tensored with  $E\otimes\Omega^j_{X/S}$, is the
complex
\begin{equation}\label{caldp.e}
\cA_{\lift X S}^{\cx,j}(E) := \cA_{\lift X S}^{0,j}(E)
\rTo^{d'^{0,j}} \cA_{\lift X S}^{1,j}(E) \rTo^{d'^{1,j}}
 \cdots.
\end{equation}
It follows from Lemma~\ref{thetacomm.l} that the differentials $d$
and $d'$ commute.
% Let us check this in degree zero,
% on elements of the form $\alpha \ot e$,
% with $\alpha$ a section of $\cA_{\lift X S}$
% and $e$ a section of $E$.  \marginpar{add dots and remove primes}
% \begin{eqnarray*}
%   d' d(\alpha \ot e) & = & (\id\ot\theta)\nabla(\alpha\ot e) \\
%     & = & (\id \ot\theta)(\nabla_\cA(\alpha) \ot e + \alpha \ot \nabla_E(e)) \\
%    & = & (\id \ot\theta)(\nabla_\cA(\alpha) \ot e)+ \theta(\alpha) \ot \nabla_E(e).
% \end{eqnarray*}
% Since the $p$-curvature $\theta$ of the connection $\nabla_\cA$ on $\cA$ commutes
% with $\nabla_\cA$, this is
% \begin{eqnarray*}
%   \id \ot\nabla_\ca (\theta(\alpha)) \ot e + \theta(\alpha) \ot \nabla_E(e)  & = &
%     (\nabla\ot \id)(\theta(\alpha) \ot e) \\
%   & = & d d'(\alpha \ot e)
% \end{eqnarray*}
Thus we can  form the double complex $\cA_{\lift X S}^{\cx\cx}(E)$
and the associated simple complex $\cA_{\lift X S}^\cx (E)$.

For each $i$ there is a natural map from $E'  \otimes
\Omega^i_{X'/S} $ to  $\Ker(d^{i,0})$, which can be regarded as a
morphism of filtered complexes,
\begin{equation}
  \label{aughig.e}
(E'  \otimes \Omega^i_{X'/S},N'_\cx) \to (\cA_{\lift X S}^{i,\cx}
(E),N_\cx^{tot}),
\end{equation}
compatible with the Higgs boundary maps:
\begin{diagram}
  (E'  \otimes \Omega^i_{X'/S},N_\cx') &\rTo& (\cA_{\lift X S}^{i,\cx} (E),N_\cx^{tot})\\
   \dTo && \dTo\\
  (E'  \otimes \Omega^{i+1}_{X'/S},N_\cx') &\rTo& (\cA_{\lift X S}^{i+1,\cx} (E),N_\cx^{tot})\\
\end{diagram}
In the same way we find for each $j$ a morphism
\begin{equation}
  \label{augdr.e}
F_{X/S*} (E \ot \Omega^j_{X/S}, N_\cx)\to (\cA_{\lift X
S}^{\cx,j}(E),N_\cx^{tot})
\end{equation}
compatible with the de Rham boundary maps
\begin{diagram}
  F_{X/S*} (E \ot \Omega^j_{X/S}, N_\cx) & \rTo& (\cA_{\lift X S}^{\cx,j}(E),N_\cx^{tot}) \\
\dTo && \dTo \\
F_{X/S*} (E \ot \Omega^{j+1}_{X/S}, N_\cx) & \rTo &(\cA_{\lift X
S}^{\cx,j+1}(E),N_\cx^{tot})
\end{diagram}
These assemble into morphisms of filtered complexes:
\begin{equation}
  \label{drmaps.e}
    (E'\otimes \Omega^\cx_{X'/S},N'_\cx)  \rTo^{a_{\lift X S}} (\cA^\cx_{\lift X S},N^{tot}_\cx)
 \lTo^{b_{\lift X S}}  (E\otimes \Omega^\cx_{X/S}, N_\cx)
\end{equation}

If there is a lifting $\tF$ of $F_{X/S}$, we can make the analogous
construction with $\cA_\tF$ in place of $\cA_{\lift X S}$, and we
use the analogous notation.  Then there is a natural morphism of
double complexes $\cA^{\cx\cx}_{\lift X S}(E) \to
\cA^{\cx\cx}_{\cF}(E)$. Taking associated simple complexes, we find
a commutative diagram:
\begin{equation}\label{cxsf.d}
\begin{diagram}
  (C_{\lift X S}(E)\otimes \Omega^\cx_{X'/S},N'_\cx)  &\rTo^{a_{\lift X S}}
       &(\cA_{\lift X S},N^{tot}_\cx)
  & \lTo^{b_{\lift X S}} & (E\otimes \Omega^\cx_{X/S}, N_\cx)\cr
\dTo&&\dTo&&\dTo \cr
  (C_\tF(E)\otimes \Omega^\cx_{X'/S},N'_\cx) & \rTo^{a_\cF} &(\cA_{\tF},N^{tot}_\cx)
  & \lTo^{b_\cF} &  (E\otimes \Omega^\cx_{X/S}, N_\cx).
\end{diagram}
\end{equation}

Before stating the main theorem, let us recall that if $C^\cx$ is a
 complex with an
increasing filtration $N_\cx$, then as explained in \cite{de.thII},
the \emph{filtration d\'ecal\'ee} $N_\cx^{dec}$ on $C^\cx$ is
defined by
\begin{equation}
  \label{decal.e}
  N^{dec}_kC^q := N_{k-q}C^q + d (N_{k-q-1}C^{q-1})
\end{equation}

\begin{theorem}\label{DR.t}
Let $X/S$ be a smooth morphism in characteristic $p$. Let $E :=
(E,\nabla, N)$ be an object of $\micn(X/S)$ with $N_{-1}E = 0$.
\begin{enumerate}
\item   If $\lift X S$ is a lifting of $X/S$, then the
maps $a_{\lift X S}$ and $b_{\lift X S}$  (\ref{drmaps.e}) induce
filtered quasi-isomorphisms:
    \begin{eqnarray*}
F_{X/S*}(N^{dec}_{p-1}(E\ot \Omega^\cx_{X/S}),N^{dec}_\cx) &\to&
(N_{p-1}^{dec}\cA_{\lift X S}^{\cx}(E),N_\cx^{dec}) \\
(N^{\prime dec}_{p-1}(C_{\lift X
S}(E)\ot\Omega^\cx_{X'/S}),N_\cx^{\prime dec}) & \to&
(N_{p-1}^{dec}\cA_{\lift X S}^{\cx}(E),N_\cx^{dec}).
    \end{eqnarray*}
Consequently they assemble into an isomorphism in the filtered
derived category of $\oh {X'}$-modules:
$$F_{X/S*}(N^{dec}_{p-1}(E\ot \Omega^\cx_{X/S}),N^{dec}_\cx,d) \cong
(N^{\prime dec}_{p-1}(C_{\lift X
S}(E)\ot\Omega^\cx_{X'/S}),N_\cx^{\prime dec},\theta'). $$
\item  If $\tF$ is a lifting of $F_{X/S}$, then the  maps
$a_\tF$ and $b_\tF$ (\ref{cxsf.d}) induce filtered
quasi-isomorphisms:
\begin{eqnarray*}
F_{X/S*}(E\ot \Omega^\cx_{X/S},N_\cx^{dec}) &\to& (\cA_\tF^{\cx}(E), N^{dec}_\cx)\\
(C_\tF(E)\ot\Omega^\cx_{X'/S},N_\cx^{\prime dec})&\to& (\cA_\tF^{\cx}(E), N_\cx^{dec})\\
\end{eqnarray*}
These assemble into an isomorphism in the filtered derived category
$$(F_{X/S*}(E\ot \Omega^\cx_{X/S}),N^{dec}_\cx) \cong
(C_\tF(E)\ot\Omega^\cx_{X'/S},N_\cx^{\prime dec}).$$
\end{enumerate}
\end{theorem}

\begin{corollary}\label{truncate.c}
Let $(E,\nabla)$ be an object of $\mic(X/S)$ which is nilpotent of
level $\ell$.  Then a lifting $\lift X S$ induces
 isomorphisms in the derived category:
$$F_{X/S^*}\left (\tau_{< p-\ell} (E\otimes \Omega^\cx_{X/S}) \right)\cong
\tau_{< p-\ell} (C_{\lift X S}(E)\otimes \Omega^\cx_{X'/S}),$$ and
if $\ell + \dim (X/S) < p$,
$$F_{X/S^*}(E\otimes \Omega^\cx_{X/S}) \cong
 (C_{\lift X S}(E)\otimes \Omega^\cx_{X'/S}).$$
\end{corollary}

Applying (2) of Theorem~\ref{DR.t} to the canonical
filtration~(\ref{canfilt.e}) of a locally nilpotent connection, we
obtain the following result.

\begin{corollary}\label{nilp.c}
  Let $(E,\nabla)$ be an object of $\mic(X/S)$.
Assume that the connection $\nabla$ is locally nilpotent
(quasi-nilpotent in the terminology of \cite{bo.ncc}). Then a
lifting $\tF$ of $F_{X/S}$ induces isomorphisms in the derived
category
$$F_{X/S*} (E\ot \Omega^\cx_{X/S}) \cong C_\tF(E)\ot \Omega^\cx_{X'/S}.$$
\end{corollary}

Before beginning the proof of Theorem~\ref{DR.t}, let us remark that
it is not true that the  maps
\begin{eqnarray}\label{hgr.e}
 F_{X/S*} (E \ot \Omega^\cx_{X/S}, N_\cx) &\to& (\cA^\cx_\tF(E),N_\cx) \\
   (E'\ot \Omega^\cx_{X'/S}, N_\cx) & \to & (\cA^\cx_\tF(E),N_\cx)
\end{eqnarray}
are filtered quasi-isomorphisms. However, these maps induce maps of
spectral sequences, which on the $E_1$ level are maps of complexes
of sheaves:
\begin{equation}\label{hr1.e}
\cH(\gr a_\cF) \colon \left (F_{X/S*}\cH(\gr E \ot
\Omega^\cx_{X/S}),d_1 \right)
  \to \left (\cH(\gr \cA^\cx_\tF(E)),d_1\right )
\end{equation}
\begin{equation}\label{hr2.e}
 \cH(\gr b_\cF)\colon \left (\cH(\gr E'\otimes \Omega^\cx_{X'/S}),d_1\right )
   \to  \left ( \cH(\gr \cA^\cx_\tF(E)),d_1 \right ),
\end{equation}
where $d_1$ is the differential of the spectral sequences. We shall
prove that these maps are quasi-isomorphisms  (not isomorphisms),
and hence induce isomorphisms on the $E_2$-terms of the spectral
sequence.
\begin{lemma}\label{keygr.l}
In the situation of (2) in Theorem~\ref{DR.t},  the maps
(\ref{hr1.e}) and  (\ref{hr2.e})
  above are quasi-isomorphisms.
\end{lemma}
\begin{proof}
Since the $p$-curvature of $\gr E$ vanishes, the classical Cartier
isomorphism induces a canonical isomorphism:
$$E_1^\cx(E\ot \Omega_{X/S}^\cx,N)  = \cH^\cx(\gr E \ot \Omega_{X/S}^\cx)
\cong (\gr E)^\nabla \ot \Omega^q_{X'/S}.$$ Corollary (5.1.1) of
\cite{o.hcpc}, allows us to compute the differential $d_1$ of this
spectral sequence.  It asserts that the diagram below is
anticommutative,  thus identifying the (negative of) the
differential $d_1^q$ with the
 graded map $\gr(\psi)$ induced by the $p$-curvature of $E$:
 \begin{equation}\label{pcurvd1.e}
 \begin{diagram}
\cH^q(\gr_i E \ot \Omega^\cx_{X/S}) &\rTo^{d^q_1}
&\cH^{q+1}(\gr_{i-1} E \ot \Omega^\cx_{X/S}) \cr \dTo && \dTo \cr
(\gr_i E)^\nabla\ot \Omega^q_{X'/S}& \rTo^{\gr(\psi)}& (\gr_{i-1}
E)^\nabla \ot\Omega^{q+1}_{X'/S}.
 \end{diagram}
 \end{equation}
Thus there is an isomorphism of complexes
$$(F_{X/S*} E_1^\cx(E\ot \Omega^\cx_{X/S},N_\cx), d_1)
 \cong (\gr E'\ot \Omega^\cx_{X'/S},\gr(\psi)).$$

We apply the same method to analyze the $E_1$ term of the spectral
sequence of the filtered complex $(\cA_\tF^\cx(E),N_\cx)$. The total
differential of  the double complex $\cA_\tF^{\cx\cx} (E)$ induces a
map
$$N_k\cA_\tF^{i,j} \to N_{k-1}\cA_\tF^{i+1,j} \oplus N_k\cA_\tF^{i,j+1},$$
so the differential on $\gr \cA_{\tF}^\cx(E)$ is just the de Rham
differential of the module with connection
$$\gr \cA_\tF^{\cx,0}(E) = \oplus_i \gr \cA_\tF(E) \ot F_{X/S}^*\Omega^i_{X'/S},$$
Since this connection has vanishing $p$-curvature, the classical
Cartier isomorphism provides an isomorphism:
$$ H^\cx(\gr \cA^\cx_\tF(E) \ot F_{X/S}^*\Omega^i_{X'/S},d)
 \cong  (\gr (\cA_\tF(E))^\nabla \ot \Omega^\cx_{X'/S}\ot \Omega^i_{X'/S}.$$
The differential $d_1$ of the spectral sequence is then a sum of
maps
\begin{eqnarray*}
 (\gr (\cA_\tF(E))^\nabla \ot \Omega^j_{X'/S}\ot \Omega^i_{X'/S} &\to&
(\gr (\cA_\tF(E))^\nabla \ot \Omega^j_{X'/S}\ot \Omega^{i+1}_{X'/S} \\
(\gr (\cA_\tF(E))^\nabla \ot \Omega^j_{X'/S}\ot \Omega^i_{X'/S}
&\to& (\gr (\cA_\tF(E))^\nabla \ot \Omega^{j+1}_{X'/S}\ot
\Omega^i_{X'/S}
  \end{eqnarray*}
The first of these is the map induced by differential $d'$ of
$\cA^{\cx\cx}_\tF(E)$, which comes from the $p$-curvature of $\cA$,
and   \cite{o.hcpc} identifies the second as  the map coming from
the $p$-curvature of the connection $\nabla$ on $\cA_\tF(E)$. Thus
we have an isomorphism of complexes:
$$\left (E_1^\cx(\cA_\tF^\cx(E),N_\cx,d_1\right)
 \cong \left(\gr\cA_\tF(E)^\nabla\ot \Omega^\cx_{X'/S}\ot \Omega^\cx_{X'/S},d\right),$$
where the differential on the right is the  differential of the
simple complex associated to the double complex  whose term in
degree $i,j$ is $$(\gr \cA_\tF(E))^\nabla \ot \Omega^i_{X'/S} \ot
\Omega^j_{X'/S}$$ and whose differential is the graded map induced
by the Higgs fields $\theta_\cA$ and $\theta_{tot}$. In fact, by
Corollary~\ref{grhigg.c},
 $\gr \cA_\tF(E) \cong \gr E \ot \gr \cA_\tF(E)$, compatibly
with the connections and Higgs fields.  Furthermore,
$$ (\gr E \ot \gr \cA_\tF)^\nabla \cong (\gr E)^\nabla \ot (\gr \cA_\tF)^\nabla
\cong (\gr E)^\nabla \ot \Gamma_\cx \Omega^1_{X'/S}.$$ Let us write
$\Omega$ for $\Omega^1_{X'/S}$ and $T$ for its dual. According to
\ref{grhigg.c}, $\gr E'$ is the Higgs transform of $\gr E$ with
respect to the $T$-Higgs module $\Gamma_\cx \Omega$. Thus the maps
$\gr a_\tF$ and $\gr b_\tF$ become identified with maps of complexes
which term by term are the mappings
\begin{eqnarray*}
  (\gr E)^\nabla  \ot \Omega^j & \to & (\gr E)^\nabla \ot\Gamma_\cx \Omega\ot \Omega^\cx\ot \Omega^j\\
  \gr (E') \ot \Omega^i & \to & \gr(E') \ot \Gamma_\cx \Omega\ot \Omega^i\ot \Omega^\cx\\
\end{eqnarray*}
constructed in the same way as $a_\tF$ and $b_\tF$. This is exactly
the situation discussed in Proposition~\ref{cohiggs.p}, so the lemma
follows.
\end{proof}

\begin{proof}[Proof of Theorem \ref{DR.t}]
  To prove that the arrows in (1) of the theorem are isomorphisms is a local
question, so we may without loss of generality assume that there is
a lifting $\tF$ of Frobenius. For $i < p$, the map $N_i \cA_{\lift X
S} \to N_i \cA_\tF$ is an isomorphism. Furthermore, since $N_{-1}E =
0$,
$$N_i^{tot}\cA_{\lift X S}(E) = \sum_{j=0}^i N_jE \otimes N_{i-j}\cA_{\lift X S} =
\sum_{j=0}^i N_jE \otimes N_{i-j}\cA_{\tF}  = N_i^{tot}\cA_\tF(E)$$
when $i < p$.   Thus the map:
$$(\cA^{\cx\cx}_{\lift X S}(E), N^{tot}_\cx)
\to (\cA^{\cx\cx}_{\tF}(E), N^{tot}_\cx).$$
 is a filtered isomorphism when restricted to $N_{p-1}$.
Thus statement (1) will follow from statement (2).

Since the filtration  $N_\cx$ on $E$ is exhaustive
 and  formation of direct limits
in the category of sheaves on $X$ is exact, we may and shall assume
that $N_\cx$ is finite. It will suffice for us to prove that the
maps of complexes
\begin{eqnarray*}
  F_{X/S*} \gr^{N^{dec}}(E\otimes \Omega^\cx_{X/S}) & \to &
   \gr^{N^{dec}} (\cA_\tF^\cx(E)) \\
\gr^{N^{dec}} ( C_\tF(E)\ot \Omega^\cx_{X'/S}) & \to &
   \gr^{N^{dec}} (\cA_\tF^\cx(E)) \\
\end{eqnarray*}
are quasi-isomorphisms. Recall from \cite{de.thII} that there are
natural injections $H^q(\gr^{N} C^\cx) \to  \gr^{ N^{dec}} C^q$
which assemble to form a quasi-isomorphism
\begin{equation}
  \label{decquasi.e}
(E_1(C^\cx,N_\cx),d) \cong (H^\cx(\gr^{N} C^\cx),d)\to  \gr^{
N^{dec}} C^\cx \cong (E_0(C^\cx,N_\cx^{dec}),d).
\end{equation}
Thus the theorem follows from Lemma~\ref{keygr.l}.
\end{proof}

\begin{remark}\label{compat.r}
Let $(E,\nabla)$ be an object of $\mic_{\cx}(X/S)$,  suppose that
there exists a global lifting of $F_{X/S}$, and let $(E',\psi')$
denote the Cartier transform of $(E,\nabla)$.  By
Remark~\ref{locglog.r}, there is a canonical isomorphism
$F^*_{X/S}(E',\psi') \cong (E,-\psi)$, where $\psi$ is the
$p$-curvature of $\nabla$.  This induces  isomorphisms
$$F_{X/S}^*\cH^i(E',\psi') \cong  \cH^i(E,-\psi).$$
for all $i$.  Recall from \cite{o.hcpc} that the sheaves of $\oh
X$-modules $\cH^i(E,-\psi)$ carry a canonical integrable connection
$\nabla$ whose $p$-curvature is zero, induced by the given
connection on $E$ and the Frobenius descent connection on
$F_{X/S}^*\Omega^q_{X'/S}$.  It follows easily that the above
isomorphisms are horizontal and hence descend to isomorphisms of
$\oh {X'}$-modules
$$\cH^i(E',\psi') \cong \cH^i(E,-\psi)^\nabla.$$
On the other hand, (\ref{DR.t}) gives us isomorphisms
$\cH^i(E',\psi') \cong \cH^i(E,\nabla)$.  Combining these, we find
the ``generalized Cartier isomorphism''
$$\cH^i(E,\nabla) \cong \cH^i(E,-\psi)^\nabla.$$
Another construction of such an isomorphisms was given in
\cite{o.hcpc}, independent of any lifting of $X$ or $F_{X/S}$ or
nilpotence condition on $\nabla$.  One can easily see that these two
isomorphisms are the same, because they agree when $i = 0$ and
because both sides are effaceable cohomological delta functors in
the category $\mic_\cx(X/S)$.
\end{remark}

Suppose that $X$ is noetherian and $E$ is coherent. A consequence of
the isomorphisms discussed in Remark~\ref{compat.r}
 is the fact that the de Rham complex of $(E,\nabla)$
with an integrable connection $\nabla$ is determined, as an object
in the derived category, by its formal completion along a closed
subset determined by its $p$-curvature $\psi$. Recall that
$(E,\psi)$ gives rise to a coherent sheaf  $\tilde E$ on
$\bT^*_{X'/S}$. Define the \emph{essential support} of $(E,\nabla)$
to be the set-theoretic intersection of the support of $\tilde E$
with the zero section of $\bT^*_{X'/S}$. We should perhaps recall
that $F_{X/S} \colon X \to X'$ is a homeomorphism and   from
\cite[2.3.1]{o.hcpc} that  the essential support of $(E,\nabla)$
corresponds via  $F_{X/S}$ to the support in $X$ of the  Higgs
cohomology sheaves of the $p$-curvature of $(E,\nabla)$.  (In fact,
the $d$th cohomology   sheaf suffices.)

\begin{proposition}\label{forcomp.p}
Let $X/S$ be a smooth morphism of schemes in characteristic $p > 0$
of relative dimension $d$. Let $(E,\nabla)$ be a coherent sheaf with
integrable connection on $X/S$, and let $Z\subseteq X$ be a closed
subscheme containing  the essential support of $(E,\nabla)$.  Let
$i_Z \colon X_{/Z} \to X$ denote the natural map from the formal
completion of $X$ along $Z$ to $X$.  Then the natural map of de Rham
complexes:
$$a \colon E\otimes \Omega^\cx_{X/S} \to i_{Z*}E_{/Z}\otimes \Omega^\cx_{X/S}$$
is a quasi-isomorphism.
\end{proposition}
\begin{proof}
It suffices to prove that the map above induces an isomorphism on
cohomology sheaves. The generalized Cartier isomorphism
\cite{o.hcpc} is an isomorphism of sheaves of $\oh {X'}$-modules
$$\cH^q(F_{X/S*}E\otimes \Omega^\cx_{X/S}) \cong F_{X/S*}\cH^q(E\otimes F^*_{X/S}\Omega^\cx_{X'/S})^\nabla$$
where the complex on the right is the Higgs complex of the F-Higgs
field given by the $p$-curvature of $\nabla$.  Now one has a
commutative diagram
\begin{diagram}
\cH^q(F_{X/S*}E\otimes \Omega^\cx_{X/S}) &\rTo&   \cH^q(
i_{Z*}E_{/Z}\otimes \Omega^\cx_{X/S}) \cr \dTo^\cong && \dTo_\cong
\cr \cH^q(E\otimes F^*_{X/S}\Omega^\cx_{X'/S})^\nabla
      & \rTo &  \cH^q(i_{Z*}E_{/Z}\otimes F^*_{X/S}\Omega^\cx_{X'/S})^\nabla
\end{diagram}
Thus it suffices to prove that the natural map
$$\cH^q(E\otimes F^*_{X/S}\Omega^\cx_{X'/S})
    \rTo  \cH^q(i_{Z*}E_{/Z}\otimes F^*_{X/S}\Omega^\cx_{X'/S})  $$
is an isomorphism of $\oh X$-modules.  Since the completion functor
is exact, and since the cohomology sheaves $\cH^q(E\otimes
F^*_{X/S}\Omega^\cx_{X'/S})$ have support in $Z$, this is clear.
\end{proof}

Let us also remark that in the situation of
Proposition~\ref{forcomp.p}, we can define a \emph{formal Cartier
transform} as follows.  Let $I \subseteq \oh {X'}$  be an ideal of
definition of the essential support $Z$ of $\tilde E$.  For each
$n$, let $E_n := E/F_{X/S}^*I^n E$, which inherits an integrable
connection from the connection on $E$. Then the $p$-curvature of
$(E_n,\nabla)$ is nilpotent and hence, given a lifting $\tF$ of
$F_{X/S}$, it has has Cartier transform $C_\tF(E_n)$.  These Cartier
transforms are compatible with change in $n$, and they fit  together
to define a coherent sheaf on the formal scheme $X'_{/Z}$, which we
(slightly abusively) still denote by $C_\tF(E)$.  The double complex
constructions used in the proof of Theorem~\ref{DR.t} also fit
together into a formal double complex.
  The following statement is a consequence of this and the previous
  proposition.

\begin{proposition}\label{coherent.p}
Suppose that $X$ is noetherian and that $(E,\nabla)$ is a coherent
sheaf on $X$ with integrable connection. Let $\tF$ be a lifting of
$F_{X/S}$ and let $C_\tF(E)$ denote the formal Cartier transform of
$E$ described above.   Then the maps of Proposition~\ref{forcomp.p}
and statement (2) of Theorem~\ref{DR.t} fit together to define an
isomorphism in the derived category of $\oh {X'}$-modules
$$F_{X/S*} (E\ot \Omega^\cx_{X/S}) \cong C_\tF(E) \ot \Omega^\cx_{X'/S}$$
\end{proposition}
\section{Functoriality of the Cartier transform}\label{fct.s}
\subsection{Gauss-Manin connections and fields}\label{gm.ss}
In this section we review the definitions of higher direct
%  and inverse
images of modules with connections and Higgs fields.
 We show that their formation with respect to
a smooth morphism of relative dimension $d$, increases the level of
nilpotence of a connection (resp. of a Higgs field)  by at most $d$.
This  result strengthens the nilpotence theorem of
Katz~\cite[5.10]{ka.ncmt} and will be used in our discussion of the
compatibility of the Cartier transform with higher direct images.

  Recall that if
$h \colon X \to Y$ is a smooth morphism of smooth $S$-schemes and if
$(E,\nabla)$ is a module with integrable connection on $X/S$, then
the sheaves
$$R^nh^{DR}_*(E,\nabla) := R^nh_*(E\ot \Omega^\cx_{X/Y},d)$$
are endowed with a canonical connection, called the
\emph{Gauss-Manin connection}.  By the same token, if $(E,\theta)$
is a module with a Higgs field $\theta$, then the sheaves
$$R^nh_{HIG*}(E,\theta) := R^nh_*(E\ot \Omega^\cx_{X/Y},\theta)$$
are endowed with a canonical Higgs field, which we shall call the
\emph{Gauss-Manin field}. Each of these can be constructed in many
ways.  For the reader's convenience we explain one of these here; a
variant of the  ``explicit'' construction explained in
\cite[3.4]{ka.ncmt}.   We write out the details in the de Rham case
only; the Higgs case is analogous but easier.

Let $(E,\nabla)$ be a module with integrable  connection on $X/S$
and let $\xi$ be a local section of $T_{X/S}$. Then interior
multiplication by $\xi$ defines a map of graded sheaves
$$i_\xi \colon E\ot \Omega^\cx_{X/S} \to E\ot \Omega^\cx_{X/S},$$
of degree $-1$.  The \emph{Lie derivative} with respect to $\xi$ is
by definition the map
$$L_\xi := d i_\xi + i_\xi d ,$$
which has degree zero. By construction $L_\xi$ is a morphism of
complexes, homotopic to zero. Now recall that a smooth morphism $h$
induces exact sequences
\begin{eqnarray}
  \label{omegs.e}
0 \to h^*\Omega^1_{Y/S} \to &\Omega^1_{X/S}& \to \Omega^1_{X/Y} \to 0\\
  0 \to T_{X/Y} \to &T_{X/S}& \to h^*T_{Y/S} \to 0
\end{eqnarray}

Pull the second of these sequences
 back via the map $h^{-1}T_{Y/S} \to h^*T_{Y/S}$ to obtain an
exact sequence of sheaves of $f^{-1}(\oh Y)$-modules:
\begin{equation}
  \label{tseq.e}
  0 \to T_{X/Y} \to T^Y_{X/S} \to h^{-1}T_{Y/S} \to 0.
\end{equation}
Let us note that $T^Y_{X/S} \subseteq T_{X/S}$ is closed under the
bracket operation and that the inclusion $T_{X/Y} \to T^Y_{X/S}$ is
compatible with the bracket operations.   Moreover, if $g$ is a
local section of $h^{-1}(\oh Y)$ and $\xi$ is a local section of
$T^Y_{X/S}$, then $\xi(g)$ also belongs to $h^{-1}(\oh Y)$, and if
$\eta$ is a local section of $T_{X/Y}$, then
$$[\eta,\xi](g) = \eta(\xi(g)) - \xi (\eta(g)) = 0. $$
It follows that $[\eta,\xi] \in T_{X/Y}$, so that $T_{X/Y}$ is an
ideal in  the Lie algebra $T^Y_{X/S}$ and the map $T^Y_{X/S} \to
h^{-1}T_{Y/S}$ is a Lie algebra homomorphism.

\begin{lemma}\label{liekoz.l}
If $\xi$ is a local section of $T_{X/S}^Y$, then $L_\xi $ preserves
the Koszul filtration $K^\cx$ of $E\ot \Omega^\cx_{X/S}$ induced by
the exact sequence (\ref{omegs.e}). In particular, $L_\xi$ induces a
morphism of complexes
$$L_\xi \colon E\ot \Omega^\cx_{X/Y} \to E\ot \Omega^\cx_{X/Y}.$$
Furthermore, if $\xi$ and $\xi'$ are local sections of $T^Y_{X/S}$,
then $[L_{\xi},L_{\xi'}] = L_{[\xi,\xi']}$ in $\End(E\ot
\Omega^\cx_{X/Y})$.
  \end{lemma}
  \begin{proof}
 By definition,
$$K^i (E\ot \Omega^q_{X/S}) =
Im \left (h^*\Omega^i_{Y/S} \ot E\ot\Omega^{q-i}_{X/S} \right)\to
E\ot \Omega^q_{X/S}.$$ Let $\xi$ be a local  section of $T^Y_{X/S}$.
Since $L_\xi$ acts as a derivation with respect to multiplication by
$\Omega^\cx_{X/S}$, it suffices to check that if $\omega$ is a local
section of $h^*\Omega^1_{Y/S}$, then $L_\xi(\omega)$ also belongs to
$h^*\Omega^1_{Y/S}$. Again using the fact that $L_\xi$ is a
derivation, we see that it suffices to check this when $\omega$ lies
in $h^{-1}\Omega^1_{Y/S}$. But if $\omega \in h^{-1}\Omega^1_{Y/S}$
and if the image of $\xi$ in $h^*T_{Y/S}$ lies in $h^{-1} T_{Y/S}$,
%$i_\xi(\omega) \in h^{-1} \oh Y$
$L_\xi(\omega) = di_\xi(\omega) + i_\xi d\omega \in
h^{-1}\Omega^1_{Y/S}$.

The fact that  the action  of $T_{X/S}$ on $E$ by Lie derivative is
compatible with the bracket follows from the integrability of
$\nabla$, and it is well-known that the same is true for its action
on $\Omega^\cx_{X/S}$.
  Since $L_\xi$, $L_{\xi'}$ and $L_{[\xi,\xi']}$ act
as derivations with respect to multiplication by forms, it follows
that $[L_\xi,L_{\xi'}] = L_{[\xi,\xi']}$ on  $E \ot
\Omega^\cx_{X/S}$ and hence also on $E\ot \Omega^\cx_{X/Y}$
  \end{proof}

Now let
$$T^\cx_{X\to Y} := T_{X/Y} \rTo^{d} T^Y_{X/S},$$
regarded as a complex in degrees $-1$ and $0$,
 where the boundary map is the inclusion.
We can give $T^\cx_{X \to Y}$ the structure of a differential graded
Lie algebra by defining $[\eta,\eta'] := 0$ if $\eta, \eta' \in
T_{X/Y}$, $[\eta,\xi] := [d\eta,\xi] \in T_{X/Y}$
 if $\eta \in T_{X/Y}$
and $\xi \in T^Y_{X/S}$, and $[\xi,\xi']$ the usual bracket if $\xi,
\xi' \in T^Y_{X/Y}$. The exact sequence (\ref{tseq.e}) defines an
isomorphism in the derived category of $f^{-1}\oh Y$-modules:
\begin{equation}
  \label{txytot.e}
T^\cx_{X \to Y} \to h^{-1} T_{Y/S}.
\end{equation}
which is compatible with the bracket structure on $h^{-1} T_{X/Y}$.

If $\eta$ is a local section of $T_{X/Y}$, then $i_\eta$ defines a
section of degree $-1$ of the complex $\cEnd(E\ot
\Omega^\cx_{X/Y})$, which we denote by $\nabla^{-1}(\eta)$.   If
$\xi$ is a local section of $T^Y_{X/S}$, then Lemma~\ref{liekoz.l}
tells us that $L_\xi$ defines a section $\nabla^0(\xi)$ of degree
$0$ of $\cEnd(E\ot \Omega^\cx_{X/Y})$. Let us observe that
$\nabla^{-1}$ and $\nabla^0$ assemble into a morphism of complexes:
$$\nabla^\cx \colon T^\cx_{X\to Y} \to \cEnd(E\ot \Omega^\cx_{X/Y}).$$
Indeed, if $\xi \in T^Y_{X/S}$, then $\nabla^0(\xi) $ is a morphism
of complexes, so it is annihilated by the total differential of
$\cEnd(E\ot \Omega^\cx_{X/Y})$. If $\eta \in T_{X/Y}$, then
$\nabla^{-1}(\eta)$ has degree $-1$, so
$$d \nabla^{-1}(\eta)  =  d\circ \nabla^{-1}(\eta) + \nabla^{-1}(\eta) \circ d
 =  d\circ i_\eta + i_\eta \circ d = L_{\eta}     =  \nabla^0(d\eta).$$
Let us also check that $\nabla^\cx$ is a morphism of differential
graded Lie algebras.  If $\xi,\xi' \in T_{X/S}^Y$, then  we saw in
Lemma~\ref{liekoz.l} that
$$[\nabla^0(\xi), \nabla^0(\xi')] := [L_\xi,L_{\xi'}] = L_{[\xi,\xi']}
= \nabla^0( [\xi,\xi']).$$ We must also check that if $\xi \in
T^Y_{X/S}$ and $\eta \in T_{X/Y}$, then
$$[\nabla^{-1}(\eta),\nabla^0(\xi)] = \nabla^{-1}([\eta,\xi]),\ie, $$
that $[i_\eta, L_\xi] = i_{[\eta,\xi]}$.  Observe first that both
sides are derivations of $E \ot \Omega^\cx_{X/Y}$ of degree $-1$
with respect to multiplication by forms, and in particular are $\oh
X$-linear.  Thus it suffices to check the formula for closed
$1$-forms.  In fact, if  $\omega \in \Omega^1_{X/S}$ is closed, then
\begin{eqnarray*}
[\nabla^{-1}(\eta),\nabla^0(\xi)](\omega) & = & i_\eta L_\xi(\omega) - L_\xi i_\eta(\omega) \\
             & = & i_\eta(d \angles \xi \omega) - L_\xi \angles \eta \omega \\
             & = & \eta \angles \xi \omega - \xi \angles \eta \omega \\
            & = & \angles {[\eta,\xi]} \omega \\
       & = &  \nabla^{-1}([\eta,\xi]),
\end{eqnarray*}
as required.  Finally, let us observe that $\nabla^\cx$ is a
derivation with respect to multiplication by sections of $h^{-1}(\oh
Y)$.

 \begin{definition}\label{gmanin.d}
   Let $h \colon X \to Y$ be a smooth morphism
of smooth $S$-schemes and let $(E,\nabla)$ (resp.$(E,\theta))$ be a
module with integrable connection  (resp. Higgs field) on $X/S$.
Then the \emph{Gauss-Manin connection} (resp. \emph{Higgs field}) on
$R^nh_*(E\otimes \Omega^\cx_{X/Y})$ is the map
$$T_{Y/S}  \to \End R^nh_*(E\ot \Omega^\cx_{X/Y})$$
obtained by composing the adjunction map
$$ T_{Y/S} \to   h_* h^{-1}T_{Y/S}  = R^0h_* h^{-1} T_{Y/S}$$
with the inverse of the   isomorphism $ R^0h_* T^\cx_{X \to Y} \to
R^0h_* h^{-1}T_{Y/S}$ defined by  (\ref{txytot.e}) and  the maps
$$R^0h_*(\nabla^\cx) \colon R^0h_*(T^\cx_{X \to Y} ) \to R^0h_*\cEnd(E\ot \Omega^\cx_{X/Y})
\to \End R^nh_*(E\ot \Omega^\cx_{X/Y}).$$

 \end{definition}

 \begin{remark}\label{gmder.r}
The integrability of the Gauss-Manin connection defined here follows
from the compatibility of the maps (\ref{txytot.e}) and $\nabla^\cx$
with the bracket operations.  A similar construction defines the
Gauss-Manin Higgs field, and thus we obtain sequence of  functors
   \begin{eqnarray*}
     R^nh^{DR}_* \colon \mic(X/S) &\to& \mic(Y/S) \\
     R^nh^{HIG}_* \colon \hig(X/S) &\to& \hig(Y/S)
   \end{eqnarray*}
It is straightforward to check that these fit into  sequences of
exact effaceable $\delta$-functors and hence are derived functors.
This makes it easy to compare this construction with the many others
which appear in the literature and in particular with the derived
category constructions appearing in section~\ref{ddii.ss}.
 \end{remark}

Now suppose that $N_\cx$ is an increasing filtration on $E$ which is
stable under the connection (resp. Higgs field). Then the
filtrations $N_\cx$ and $N^{dec}$  of $E \ot \Omega^\cx_{X/Y}$ are
stable under the action of $T^\cx_{X\to Y}$, and hence the higher
direct images of the corresponding filtered pieces  and  the graded
objects inherit Gauss-Manin connections.

 \begin{theorem}\label{declev.t}
Let $h \colon X \to Y$ be a smooth morphism
 of smooth $S$-schemes.
 Let $E$ be a sheaf of $\oh X$-modules endowed with
 an integrable connection $\nabla$
 (resp. a Higgs field $\theta)$.  Suppose that $N_\cx$ is a filtration
 on $E$ such that $\gr^N \nabla$ is constant
 (resp., such that $\gr^N \theta = 0$).  Then for each $n$ and $i$,
 the action of the Gauss-Manin connection (resp. field)
 on $R^nh_*(\gr_i^{N^{dec}} (E\ot \Omega^\cx_{X/Y}))$
 is constant (resp. trivial).
 \end{theorem}
 \begin{proof}
   If $\theta$ is a Higgs field such that
 $\gr^N(\theta) = 0$ , then
$\theta$ maps  $N_iE$ to $N_{i-1}E\ot \Omega^1_{X/S}$.  It follows
that the actions of $T_{X/Y}$ and $T^Y_{X/S}$ on $E\ot
\Omega^\cx_{X/Y}$ by interior multiplication and Lie derivative map
$N^{dec}_i$ to $N^{dec}_{i-1}$. Hence $T^\cx_{X\to Y}$ acts
 trivially on $\gr_i^{N^{dec}}(E\ot\Omega^\cx_{X/Y})$.

Now suppose that $\nabla$ is a connection on $E$ and $N$ is a
horizontal filtration on $E$.  Recall that we have a natural
quasi-isomorphism  (\ref{decquasi.e}) of complexes
$$ a \colon (E^{\cx,j}_1(E \ot \Omega^\cx_{X/Y},N),d_1) \to
\gr_j^{N^{dec}}(E\ot \Omega^\cx_{X/Y},d).$$ Here
 $E^{i,j}_1(E\ot \Omega^\cx_{X/Y},N) = \cH^{j-i}(\gr_i^N E\ot \Omega^\cx_{X/Y})$.
Note that if  $\xi \in  T_{X/Y} \subseteq  T^Y_{X/S}$, then $i_\xi$
is well-defined on $E\ot \Omega^\cx_{X/Y}$, and hence $L_\xi = d
i_\xi + i_\xi d$ acts as zero on $\cH^q(E\ot \Omega^\cx_{X/Y})$.
Thus the action of $T^\cx_{X \to Y}$ factors through
$h^{-1}T_{Y/S}$; the boundary maps  $d_1$ are compatible with this
action.
 Thus  $R^nh_*(E^{\cx,j}_1,d_1)$  has  a connection also,
and we claim that $R^nh_*(a)$ is compatible with the connections. To
see this, it is convenient to recall the ``dual'' version of the
filtration d\'ecal\'ee:
$$N^*_i (E\ot \Omega^q_{X/Y}): = N_{i-q} E \ot \Omega^q_{X/Y}
\cap d^{-1} (N_{i-q-1}E\ot\Omega^q_{X/Y}).$$ Then there is also a
natural quasi-isomorphism
$$ a^*\colon   (\gr_j^{N^*} (E\ot \Omega^\cx_{X/Y}),d) \to
(E^{\cx,j}_1(E\ot \Omega^\cx_{X/Y},N), d_1) $$ Then $a^*$ and $a
a^*$ are  compatible with the actions of $T^\cx_{X \to Y}$. Although
$a$ is not compatible with the action of $T^\cx_{X\to Y}$ on the
level of complexes, it follows that  it \emph{is} compatible with
the induced action of $T_{Y/S}$ on  higher direct images.

Now suppose that $\gr E :=\gr^N E$ is constant. The theorem will
follow if we prove that  the Gauss-Manin connection on
$R^nh_*(E^\cx_1(E\ot \Omega^\cx_{X/Y},N_\cx),d_1)$ is constant.  Let
us consider the relative Frobenius diagram:
\begin{equation}\label{relrel.d}
\begin{diagram}
X & \rTo^{F_{X/Y}} & X^{(Y)} & \rTo^{\pi_{X/Y/S}} &  X' \cr
&\rdTo_h&\dTo_{h^{(Y)}} && \dTo_h' \cr &&Y &  \rTo^{F_{Y/S}} & Y'.
\end{diagram}
\end{equation}
Here $F_{X/S} = \pi_{X/Y/S}\circ F_{X/Y}$ and he square is
Cartesian, so $\Omega^q_{X^{(Y)}/Y} \cong
\pi_{X/Y/S}^*\Omega^q_{X'/Y}$.

The morphism of filtered complexes
$$ (E\ot \Omega^\cx_{X/S}, N_\cx) \to
(E\ot \Omega^\cx_{X/Y}, N_\cx)$$ induces a morphism of spectral
sequences, which on the $E_1$-level corresponds to the top row of
the following commutative diagram:
\begin{diagram}
  \cH^q(\gr E \ot \Omega^\cx_{X/S}) & \rTo & \cH^q(\gr E \ot \Omega^\cx_{X/Y})  \cr
   \uTo^{c_S} && \uTo_{c_Y} \cr
    \cH^0(\gr E \ot \Omega^\cx_{X/S})\ot \Omega^q_{X'/S} & \rTo
        & \cH^0(\gr E \ot \Omega^\cx_{X/Y}) \ot \Omega^q_{X^{(Y)}/Y} \cr
     & \rdTo_a & \uTo_b \cr
    &&     \cH^0(\gr E \ot \Omega^\cx_{X/S})\ot \Omega^q_{X'/Y'}
\end{diagram}
The vertical maps $c_S$ and $c_Y$ induced by the inverse Cartier
isomorphism are isomorphisms because $\gr E$ is constant, the map
$a$ is surjective, and the map $b$ is injective.  Thus $\cH^0(\gr E
\ot \Omega^\cx_{X/S})\ot \Omega^q_{X'/Y'}$ can be identified with
the image of the arrow at the top of the diagram.  Since the
differentials of the spectral sequence leave this image invariant,
they induce maps
$$\cH^0(\gr E \ot \Omega^\cx_{X/S})\ot \Omega^q_{X'/Y'} \to
\cH^0(\gr E \ot \Omega^\cx_{X/S})\ot \Omega^{q+1}_{X'/Y'}.$$ and
define a complex $\cH^0(\gr E\ot \Omega^\cx_{X/S})\ot
\Omega^\cx_{X'/Y'}$ of sheaves of $\oh {X'}$-modules on $X'$. Since
the natural map
$$\pi_{X/Y/S}^* \left (\cH^0(\gr E \ot \Omega^\cx_{X/S})  \ot \Omega^q_{X'/Y'} \right)  \to
 \cH^0(\gr E \ot \Omega^\cx_{X/Y}) \ot \Omega^q_{X^{(Y)}/Y}$$
is an isomorphism, we see that the complex $E_1^\cx(E,N_\cx)$
descends to a complex of $\oh {X'}$-modules on $X'$.

Note that if $\xi \in T_{X/Y} \subseteq T^Y_{X/S}$, then $i_\xi$ is
well-defined on $E\ot \Omega^\cx_{X/Y}$, and hence $L_\xi = di_\xi+
i_\xi d$ acts as zero on $\cH^q(E\ot \Omega^\cx_{X/Y})$.  Thus the
action of $T^\cx_{X \to Y}$ on $\cH^q(E\ot \Omega^\cx_{X/Y})$
factors through $h^{-1}T_{Y/S}$. For the same reason,
$T^{-1}_{Y/S}$ acts as zero on the image of $\cH^q(\gr E \ot
\Omega^\cx_{X/S})$  in $\cH^q(\gr E \ot \Omega^\cx_{X/Y})$. and it
follows that the action of $h^{-1}T_{Y/S}$ on $E_1^\cx(E,N_\cx)$ is
nothing but the Frobenius descent connection. It follows that the
Gauss-Manin connection on  $R^nh^{(Y)}_*(\cH^\cx\gr^N( E\ot
\Omega^\cx_{X/Y}),d_1)$ is the Frobenius descent connection.
 \end{proof}

The following result is an improvement of the  result
\cite[5.10]{ka.ncmt} of Katz,  which gives a multiplicative instead
of an additive estimate for the level of nilpotence of higher direct
images.

 \begin{corollary}\label{level.c}
In the situation of the previous theorem, suppose that $h \colon X
\to Y$ has relative dimension $d$ and denote by $\micn_\ell(X/S)$
the category of objects of $\micn(X/S)$ of level $\ell$, \ie, such
that there exists an integer $k$ such that $N_kE = 0$ and
$N_{k+\ell}E = E$. Then for each $q$, $R^qh^{DR}_*(E\ot
\Omega^\cx_{X/Y},N_\cx^{dec})$ lies in $\micn_{d+\ell}(Y/S)$, and
the analogous statement for Higgs modules also holds.
\end{corollary}

\begin{remark}
   In the case of connections, we can use the
diagram~\ref{pcurvd1.e}, which computes the boundary maps of the
complex $E_1^{\cx,j}(E\ot \Omega^\cx_{X/Y},N_\cx)$, to see that
$$R^nh_*\gr_i^{N^{dec}}( E\ot \Omega^\cx_{X/Y}, d) \cong
F_{Y/S}^*R^nh'_*\gr_\cx^N( E \ot \Omega^\cx_{X'/Y'},\ov \psi)$$
where $\ov \psi$ is the map induced by the $p$-curvature.
\end{remark}

\begin{example}\label{badlevel.e}
{\rm Let $k$ be a field of characteristic $p$, $S := \spec k$,  $Y
:= \spec k[t]$.  If $d$ is a positive integer,  let $m := d+2$,
assume $(p,m) = 1$,
 and consider the hypersurface  $X$
in $\bP^{n+1}$  over $S$ defined by $X_0^m + X_1^m + \cdots
X_{d+1}^m + t X_0X_1 \cdots X_{d+1}$.
 Once $Y$ is replaced by a suitable affine neighborhood
of the origin, $X/Y$ will be smooth, and the iterated
Kodaira-Spencer mapping
$$(\kappa_{\partial/\partial t})^d :
H^0(X,\Omega^d_{X/Y} ) \to H^d(X,\oh X)$$ is  an isomorphism
\cite[3.4]{o.gtcc}.  Then Katz's formula~\cite[Theorem 3.2]{ka.asde}
  implies that the  iterated $p$-curvature mapping
$$(\psi_{\partial/\partial t'})^d \colon
H^0(X,\cH^d(\Omega^\cx_{X/Y})) \to H^d(X, \cH^0(\Omega^\cx_{X/Y}))$$
is also an isomorphism. This implies that the level of the
Gauss-Manin connection on $R^dh_*(\Omega^\cx_{X/Y})$ is $d$.
Moreover, if $d >p$, the action of the center of $D_{Y/S}$ on
$R^dh_*(\Omega^\cx_{X/Y})$ does not factor through the divided power
neighborhood of the zero section. }\end{example}

\subsection{The Cartier transform and de Rham direct images}\label{ctdr.ss}

Let $h \colon X/S \to Y/S$ be a smooth morphism of smooth
$S$-schemes, endowed with liftings $\tX'/\cS$ and $\tY'/\cS$.  We
shall explain how a lifting $\tilde h' \colon \tX' \to \tY'$ of $h'$
defines a compatibility isomorphism between the Cartier transform of
the de Rham direct image of a module with connection and the Higgs
direct image of its Cartier transform.

It is convenient to work  with filtered categories as described in
Corollary~\ref{level.c}. If $\ell < p$, an   object $(E,\nabla,
N_\cx)$  of $\micn_\ell(X/S)$ can be  viewed as an object of
$\micn^\cx_\gamma(X/S)$ and we apply the filtered Cartier transform
of Theorem~\ref{rh.t} to obtain an object $(E',\theta',N'_\cx)$ of
$\hign_\ell(X'/S)$.
% If $h \colon X \to Y$ is a smooth morphism of smooth
% $S$-schemes, let
% $$R^qh^{DR}_*(E,\nabla,N):= (R^qh_*(E\ot \Omega^\cx_{X/Y}),\nabla, N_\cx)$$
% where $\nabla$ is
% the Gauss-Manin connection \ref{} and the $N_\cx$
% is the filtration  the filtration
% d\'ecal\'ee $N^{dec}$ of $E \ot \Omega^\cx_{X/Y}$.
% Theorem~\ref{} implies that this defines an object
% of $\micn_{\ell+d}(Y/S)$, where
% $d$ is the relative dimension of $X/Y$;
% and we define a filtered version of $R^qh'_{HIG*}$ analogously.
% Thus we get functors:
% \begin{eqnarray*}
%   R^qh^{DR}_* \colon \micn_{\ell}(X/S) & \rTo \micn_{\ell+d}(Y/S) \\
%   R^qh'_{HIG*} \colon \hig_{\ell}(X'/S) & \rTo \hig_{\ell+d}(Y'/S) \\
% \end{eqnarray*}

\begin{theorem}\label{push.t}
  Let $h \colon X/S \to Y/S$ be a smooth morphism
of smooth $S$-schemes, endowed with liftings $\tX'/\cS$ and
$\tY'/\cS$.  Let
 $\ell$ be an integer less than $p-d$,
where $d$ is the relative dimension of $h$. Then a lifting  $\tilde
h' \colon \tilde X'/ S \to \tilde Y'/S $ of $h' \colon X'/S \to
Y'/S$ induces an isomorphism of functors (made explicit below):
$$\Theta^q_{\tilde h'}\colon  R^qh_*^{\prime \hig}\circ C_{\lift X S} \Rightarrow
C_{\lift Y S} \circ R^qh_*^ {DR}$$
  making the diagram below 2-commutative:
\begin{diagram}
\micn_{\ell}(X/S) & \rTo^{C_{\lift X S}}&\hign_{\ell}(X'/S) \cr
\dTo^{R^qh_*^{DR}} && \dTo_{R^qh_*^{\prime\hig}} \cr
\micn_{\ell+d}(Y/S) & \rTo^{C_{\lift Y S}}&\hign_{\ell+d}(Y'/S).
\end{diagram}
\end{theorem}

We shall construct the compatibility isomorphism of
Theorem~\ref{push.t} from a canonical  filtered double complex,  a
relative version of the double complex we used in the construction
of the comparison isomorphism in Theorem~\ref{DR.t}. For any
$(E,\nabla,N) \in \micn(X/S)$, define
\begin{eqnarray*}
\cA_{\lift X  Y /\cS}^{ij}(E)
& := & F_{X/S*}\left (E \ot \cA_{\lift X S} \ot F_{X/S}^* \Omega^i_{X'/Y'}\ot\Omega^j_{X/Y} \right)\\
&\cong & F_{X/S*}\left (E \ot\cA_{\lift X S} \ot \Omega^j_{X/Y}
\right) \ot\Omega^i_{X'/Y'}.
\end{eqnarray*}

The de Rham and Higgs boundary maps then form
 a double  complex
$(\cA_{\lift X  Y /S}^{\cx\cx}(E),d',d)$, which we endow with the
total filtration $N_\cx := N_\cx^{tot}$. There is a canonical
morphism
$$(\cA^{\cx\cx}_{\lift X S}(E), N_\cx^{tot}) \to
(\cA^{\cx\cx}_{\lift X  Y/ \cS}(E), N_\cx^{tot})$$ Let us recall
from the diagram \ref{relrel.d} that we have a morphism $h^{(Y)}
\colon X^{(Y)} \to Y$ and a homeomorphism $\pi_{X/Y/S} \colon
X^{(Y)} \to X',$ which we will sometimes allow ourselves to view as
an identification to simplify the notation. The  terms of the
complex $\cA^{\cx\cx}_{\lift X Y /\cS}(E))$ are $F_{X/S*}\oh
X$-modules and the boundary maps are $\pi_{X/Y/S*}\oh
{X^{(Y)}}$-linear.

Recall from  Proposition~\ref{exfunct.p} that the lifting $\tilde
h'$ of $h$ defines a morphism of filtered algebras with connection
$$\theta_{\tilde h'}\colon  (h^*\cA_{\lift Y S},N_\cx) \to
(\cA_{\lift X S},N_\cx).$$ Then we have a morphism of filtered
relative de Rham complexes:
$$(E\ot h_{DR}^*\cA_{\lift Y S}\ot \Omega^\cx_{X/Y},N^{tot}_\cx)
\to (E\ot\cA_{\lift X  S}\ot \Omega^\cx_{X/Y} ,N_\cx^{tot}).$$ Since
$h_{DR}^*\cA_{\lift Y S}$ comes from $Y$, its $p$-curvature relative
to $Y$ vanishes, so   for each $j$, the map
\begin{equation}\label{aug1.e}
E\ot h^*\cA_{\lift Y S}\ot \Omega^j_{X/Y} \to E\ot \cA_{\lift X S}
\ot \Omega^j_{X/Y}.
\end{equation}
is annihilated by the differential:
$$ d'\colon E\ot \cA_{\lift X S} \ot \Omega^j_{X/Y} \to
 E\ot \cA_{\lift X S} \ot \Omega^j_{X/Y}\ot F_{X/S}^*\Omega^1_{X'/Y'}$$
Let $\cA_{\lift Y S}(E) : = E \ot h_{DR}^*\cA_{\lift Y S} \in
\mic(X/S)$. It follows that the maps (\ref{aug1.e}) define a
morphism of filtered complexes:
\begin{equation}\label{coqis.e}
a \colon F_{X/S*}( \cA_{\lift Y S}(E)\ot
\Omega^\cx_{X/Y},N^{tot}_\cx) \to (\cA^\cx_{\lift X Y
/\cS}(E),N^{tot}_\cx)
\end{equation}

Let $E'$ be the Cartier transform of $E$. Since formation of
$p$-curvature is compatible with de Rham pullback (see
Remark~\ref{pcurvpull.r}), the map $\theta_{\tilde h'}$ is also
compatible with the $F$-Higgs fields.  Thus we have a morphism of
filtered relative F-Higgs complexes:
$$(E' \ot h^* \cA_{\lift Y S}\ot F_{X/S}^*\Omega^\cx_{X'/Y'},N^{tot}_\cx) \to
(E\ot \cA_{\lift X S} \ot F_{X/S}^*\Omega^\cx_{X'/Y'},N^{tot}_\cx)
$$

Note that there is an isomorphism of $\oh {X^{(Y)}}$-modules
$${\pi_{X/Y/S*}h^{(Y)}}^*\cA_{\lift Y S} \cong {h'}^*F_{Y/S*}\cA_{\lift Y S}.$$
Since ${h^{(Y)}}^*\cA_{\lift Y S}$   and the Cartier transform $E'$
of $E$ are both annihilated by the relative de Rham differential
$\cA_{\lift X S}(E) \to \cA_{\lift X S}(E) \ot \Omega^1_{X/Y}$, the
same is true of the tensor product
$$\cA_{\lift Y S}(E')
 := E' \ot\pi_{X/Y/S*} h^{(Y)*}\cA_{\lift Y S}.$$
Thus we find a morphism of filtered complexes:
\begin{equation}\label{coqhis.e}
b \colon (\cA_{\lift Y S}(E')\ot \Omega^\cx_{X'/Y'},N^{tot}_\cx) \to
(\cA^\cx_{\lift X Y /\cS}(E),N^{tot}_\cx).
\end{equation}

We shall deduce Theorem~\ref{push.t} from the following result on
the level of complexes.

\begin{theorem}\label{filteredrel.t}
Suppose that $E \in \micn(X/S)$ and $N_{-1}E = 0$.
  Then the morphisms $a$ and $b$
above induce filtered quasi-isomorphisms
\begin{eqnarray*}
a \colon F_{X/S*}(N^{dec}_{p-1}(\cA_{\lift Y S}(E)\ot
\Omega^\cx_{X/Y}),N^{dec}_\cx)
&\to& (N^{dec}_{p-1}\cA^\cx_{\lift X Y /\cS}(E),N^{dec}_\cx)\\
b \colon (N^{dec}_{p-1} (\cA_{\lift Y S}(E')\ot
\Omega^\cx_{X'/Y'}),N^{dec}_\cx) &\to& (N^{dec}_{p-1}\cA^\cx_{\lift
X Y /\cS}(E),N^{dec}_\cx).
\end{eqnarray*}
The map $a$ is compatible with the Gauss-Manin connections and the
map $b$ is compatible with the Gauss-Manin Higgs fields defined in
(\ref{gmanin.d}). Moreover, the Gauss-Manin connection annihilates
the map
$$E'\ot \Omega^\cx_{X'/Y'} \to \cA^\cx_{\lift X Y /\cX} (E)$$
and the Gauss-Manin Higgs field annihilates the map
$$F_{X/S*}E\ot \Omega^\cx_{X/Y} \to \cA^\cx_{\lift X Y /\cX} (E).$$
\end{theorem}
\marginpar{for which $j$?}
\begin{proof}
The compatibilities with the Gauss-Manin connections and fields are
straightforward.  To prove that the maps in the theorem are filtered
quasi-isomorphisms, it will suffice to show that the maps of
complexes
\begin{eqnarray*}
a \colon (E^{\cx,j}_1(\cA^\cx_{\lift Y S}(E),N_\cx),d_1) &\to&
(E^{\cx,j}_1(\cA^\cx_{\lift X Y /\cS}(E),N_\cx),d_1)\\
b \colon ( E^{\cx,j}_1(\cA^\cx_{\lift Y S}(E'),N_\cx),d_1) &\to&
(E^{\cx,j}_1 (\cA^\cx_{\lift X Y /\cS}(E),N_\cx),d_1)\\
\end{eqnarray*}
are quasi-isomorphisms.  As in the proof of Theorem~\ref{DR.t}, we
find that these become maps
\begin{eqnarray*}
 \gr a \colon ( \gr E \ot S^\cx \Omega_{Y'/S}\ot \Omega^\cx_{X'/Y'}, d_1) & \to &
 \gr( E \ot S^\cx \Omega_{X'/S}\ot \Omega^\cx_{X'/Y'}, d_1) \\
 \gr b \colon ( \gr E' \ot S^\cx \Omega_{Y'/S}\ot \Omega^\cx_{X'/Y'}, d_1) & \to &
 \gr (E  \ot S^\cx \Omega_{X'/S}\ot \Omega^\cx_{X'/Y'}, d_1)
\end{eqnarray*}
Working locally on $X$, we may
%assume that the sequence
% $ 0 \to h'^*\Omega^1_{Y/S} \to \Omega^1_{X'/S} \to \Omega^1_{X/Y} \to 0$
% splits.  Then we can
identify $S^\cx \Omega^1_{X'/S}$ with the tensor product $h'^*S^\cx
\Omega^1_{Y'/S} \ot S^\cx \Omega^1_{X'/Y'} $ Then the result follows
from the filtered Poincare\'e lemma, as in
Proposition~\ref{cohiggs.p}.
\end{proof}

\begin{proof}[Proof of Theorem~\ref{push.t}]
We may assume without loss of generality that $N_{-1}E = 0$. For
each $q$, let $(E^q_{DR},N_\cx) := R^qh^{DR}_* E$  with the
filtration induced by the filtration $N_\cx^{dec}$ of $E\ot
\Omega^\cx_{X/Y}$, and let
 $(E'^q_{HIG},N_\cx) := R^q{h'}^{HIG}_*E'$ with the filtration
induced by ${N'}^{dec}_{\cx}$.

 Since the pieces of $\gr \cA_{\lift Y S}$ consists of locally
free sheaves of finite rank,
%  the natural  map
% $$ \bigoplus_{a+b = c} N_a \cA_{\lift Y S} \ot R^qh^{DR}_*(N^{dec}_b
% E\ot \Omega^\cx_{X/Y}) \to
% R^qh^{DR}_*(\cA_{\lift Y S}(E), N_c^{dec} E\ot \Omega^\cx_{X/Y})$$
% is an isomorphism, so both sides have the same image in
% $$R^qh^{DR}_* \cA_{\lift Y S}(E) \ot \cA_{\lift Y S} \ot R^qh^{DR}_*(E).$$
the projection formula gives filtered isomorphisms
$$(\cA_{\lift  Y S} \ot E_{DR}^q, N_\cx^{tot})
 \cong (R^qh^{DR}_* \cA_{\lift Y S}(E),N_\cx) $$
$$(F_{Y/S*}\cA_{\lift  Y S} \ot E_{HIG}'^q, N_\cx^{tot})
 \cong (R^qh^{HIG}_*\cA_{\lift Y S}(E'),N_\cx) $$
where the filtrations on the right are induced by the filtration
$N_\cx^{dec}$. Furthermore, these maps are compatible with the Higgs
fields and connections.    Theorem~\ref{filteredrel.t} then gives us
an isomorphism
$$\left (N_{p-1}^{tot} (\cA_{\lift Y S} \ot E_{DR}^q), N_\cx^{tot}\right)
\cong \left (N_{p-1}^{tot} (\cA_{\lift Y S} \ot E'^q_{HIG}),
N_\cx^{tot}\right )$$ compatible with the filtrations, connections,
and Higgs fields. Since $E^q_{DR}$ and has level at most $p-1$, its
Cartier transform is obtained by taking the horizontal sections of
$N_{p-1}\cA_{\lift Y S}(E^q)$, which by the above isomorphism is
$(E'^q_{HIG},N_\cx)$.

\end{proof}

\begin{remark}
 When $Y = S$, the categories $\mic_\ell(Y/S)$ and $\hig_\ell(Y'/S)$
 reduce to the category of $\oh S$-modules, and Theorem~\ref{push.t}
 above reduces to Theorem~\ref{DR.t}.
\end{remark}

\subsection{Derived direct and inverse images}\label{ddii.ss}

Let $S$ be a noetherian scheme of finite Krull dimension, $h \colon
X\to Y$  a morphism of smooth schemes over $S$. Let ${\bf
T}^*_{X'\to Y'}$  be the pullback of $\bT^{*}_{Y'/S}$ to $X'$, which
fits into the following diagram:
 $$
\def\normalbaselines{\baselineskip20pt
\lineskip3pt  \lineskiplimit3pt}
\def\mapright#1{\smash{
\mathop{\to}\limits^{#1}}}
\def\mapdown#1{\Big\downarrow\rlap
{$\vcenter{\hbox{$\scriptstyle#1$}}$}}
\begin{matrix}
 {\bf T}^*_{X'\to Y'} & \str{i'_h}   & {\bf T}^*_{X'}  \cr
 \mapdown{\overline h'}   && \cr
 {\bf T}^*_{Y'}    &  &
\end{matrix}
$$
Note that $i'_h$ is a closed embedding if and only if $h'$ is
smooth.

Let $\hig(X'\to Y')$ denote the category of sheaves of ${h'}^*
S^{\cx} T_{Y'/S}$-modules on $X'$. We define the \emph{derived
inverse image }
  $$L{h'}^*_{HIG}: D(HIG(Y'/S)) \to D(HIG(X'/S)) $$
to be the composition
$$D(HIG(Y'/S)) \str{L\overline {h'}^* }D(\hig(X'\to Y')) \str{i'_{h *}} D(HIG(X'/S)) .$$
Since $h'$ is a morphism between smooth $S$-schemes, $\overline
{h'}^*$ has bounded cohomological dimension dimension, and so takes
$D^b(\hig(Y'/S))$ to $D^b(\hig(X'/S))$.

Similarly, for a {\it smooth } morphism $h$,  the \emph{derived
direct image}
$$R{h'}^{HIG}_*: D(HIG(X'/S))\to D(HIG(Y'/S) )$$
is  the composition
$$ D(HIG(X'/S)) \str{R{i'_h}^{!}} D(\hig(X'\to Y')) \str{R\overline h'_*}D(HIG(Y'/S)) ,$$
where  $R{i'}^{!}$ sends a complex $E$  in $D(\hig(X'/S))$ to
$$R{i'}^{!}(E)= R{\cal H}om_{S^{\cx} T_{X'/S}} ( { S^{\cx} {h'}^*T_{Y'/S} }, E).$$
It is again true that this functor takes bounded complexes to
bounded complexes. Note that $R{h'}^{HIG}_*$ is right adjoint to
$L{h'}^*_{HIG}$.

Let us pass to the direct and inverse images of $D$-modules.
Proposition~\ref{newder.p} below is a reformulation, based on the
Azumaya property of the algebra of differential operators in
characteristic $p$, of the usual definition of the functors
$$L{h}^*_{DR }: D(MIC(Y/S)) \to  D(MIC(X/S)) $$
$$Rh^{DR}_* : D(MIC(X/S)) \to  D(MIC(Y/S)).  $$
Recall that $\cD_{X/S}$ is the sheaf of algebras on the cotangent
space of $X'/S$ associated to $F_{X/S*} D_{X/S}$.
 We first need the following result.

\begin{theorem}[\cite{bb.glpc}]\label{hxty.t}
Let $h \colon X  \to Y$ be a morphism of smooth $S$-schemes. Then
the $D_{X/S}\otimes h^{-1}D^{op}_{Y/S}$-module $D_{X\to Y}:=
h^*D_{Y/S}$ induces an equivalence of Azumaya algebras on ${\bf
T}^*_{X'\to Y'}$ :
$${i'}_h^*\cD_{X/S} \sim \overline {h'}^*\cD_{Y/S}.$$
\end{theorem}
\begin{proof}
To prove the theorem consider $D_{Y/S}$ as a left module over
itself. Remark~\ref{pcurvpull.r} shows that the left action of
$D_{X/S}$ on $h^*D_{Y/S}$ and the right action of $h^{-1}D_{Y/S}$
together define a left action of
$$F_{X/S*} D_{X/S} \otimes_{S^{\cx} T_{X'/S}} {h'}^*F_{Y/S*}D_{Y/S}^{op},$$
where $S^{\cx} T_{X'/S}$ acts on ${h'}^*F_{Y/S*}D^{op}_{Y/S}$ via
${i'}^*_h$ and the evident action of ${h'}^*S^{\cx} T_{Y'/S}$. This
gives us a module over the Azumaya algebra ${i'}_h^*{\cal D}_{X/S}
\otimes_{ {\cal O} _{{\bf T}^*_{X'\to Y'}}}   (\overline {h'}^*
{\cal D}_{Y/S})^{op}$. A local computation shows that this module is
locally free over ${\cal O} _{{\bf T}^*_{X'\to Y'}}$ of rank $ind
({i'}_h^*{\cal D}_{X/S} )\cdot  ind (\overline {h'}^* {\cal
D}_{Y/S})$.
\end{proof}

As a corollary, we get an equivalence of categories:
\begin{equation}\label{cxtoy.e}
C_{X'\to Y'}: \Mod({i'_h}^*{F_{X/S*} D}_{X/S})\simeq \Mod({{h'}^*{
F_{Y/S*}}D}_{Y/S}),
\end{equation}
where $\Mod(\cA)$  denotes the category of $\cA$-modules.
 Note that, since
$F_{X/S}$ is a homeomorphism, the functor
$$MIC(X/S)=\Mod(D_{X/S}) \str{F_{X/S *}} \Mod(F_{X/S*}D_{X/S})$$
is an equivalence of categories.   Thus the following result
determines $Lh^*{DR}$ and $Rh_{DR}^*$.

\begin{proposition}[\cite{bb.glpc}]\label{newder.p}
For any morphism $h: X\to Y $ there is a canonical isomorphism:
$$F_{X/S *}\circ L{h}^*_{DR } \cong {i'}_{h*} \circ C^{-1}_{X'\to Y'} \circ L{h'}^* \circ F_{Y/S *}.$$
If $h$ is smooth we also have
$$F_{Y/S *}\circ Rh^{DR}_*  \cong Rh'_*\circ C_{X'\to Y'}\circ R{i'_h}^! \circ F_{X/S *}.$$
\end{proposition}
\begin{proof}
We shall just explain the second formula. By definition, for any
$E\in  D(MIC(X/S)) $, we have
$$ C_{X'\to Y'}\circ R{i'_h}^!(F_{X/S *} E) \cong $$
$$ {\cal H}om_{{i'}_h^*F_{X/S *} D_{X/S}}(F_{X/S *}h^* D_{Y/S},
R{\cal H}om_{F_{X/S *}D_{X/S}}({i'}_h^*F_{X/S *}D_{X/S},F_{X/S *}E))
\cong $$
  $$ R{\cal H}om_{F_{X/S *}D_{X/S}}(F_{X/S *}h^*D_{Y/S}, F_{X/S *}E).$$
It follows then that
$$ Rh'_*\circ C_{X'\to Y'}\circ R{i'_h}^!(F_{X/S *} E)  =  F_{Y/S *} Rh_*(R{\cal H}om_{D_{X/S}}(h^*D_{Y/S}, E)) .$$
When $h$ is smooth this is the standard definition of $F_{Y/S
*}\circ Rh^{DR}_*$.
\end{proof}

As an application of the new construction of $Rh^{DR}_*$ let us
observe that
 if $E\in MIC(X/S)$ and the Zariski closure  of
$supp \, F_{X/S *} E \subset {\bf T}^*_{X'}$ does not intersect
${\bf T}^*_{X'\to Y'}\subset {\bf T}^*_{X'}$, then $Rh^{DR}_* E =0$.
(This follows also from Proposition~\ref{forcomp.p}).

\subsection{The conjugate filtration on $F_{X/S *} D_{X/S}$}
The algebra of differential operators in characteristic $p$, besides
the order filtration, has another natural filtration by ideals:
\begin{equation}\label{conjfiltration}
  \cdots \subset \Cee _X ^i \subset \cdots  \subset
\Cee _X ^1 \subset   F_{X/S *} D_{X/S},
\end{equation}
$$  \Cee _X ^i = S^i T_{X'}(F_{X/S *} D_{X/S}).   $$
We shall call (\ref{conjfiltration}) the \emph{conjugate filtration}
since, as we will explain in (\ref{conjf}) below, it induces the
conjugate filtration on the de Rham cohomology groups. The
associated graded algebra $\gr(F_{X/S *} D_{X/S})$ is a canonically
split tensor Azumaya algebra. In this section we shall study a
certain filtered derived category of modules over the filtered
algebra $ F_{X/S *} D_{X/S}$. We will see how the splitting property
of
 $\gr(F_{X/S *} D_{X/S})$ together with some general results in homological algebra lead to generalizations
 and simple proofs of some of the fundamental results of Katz, including the p-curvature formula for the
 Gauss-Manin connection. Our main application is the functoriality of the Cartier transform with respect to the direct images.

The following construction plays a central role in this subsection.
\begin{definition}\label{cfa.d}
Let  ${\cal A}$ be  a sheaf of algebras  over a scheme $Z$ and $\Cee
\subset {\cal A}$ be a two-sided ideal. Denote by $CF({\cal A},
\Cee)$ the category of (unbounded) filtered complexes of ${\cal
A}$-modules
$$\cdots \subset  (N^{i+1}E^{\cx}, d) \subset  (N^{i}E^{\cx}, d) \subset \cdots  \subset ( E^\cx, d), $$
satisfying the following conditions:
\begin{enumerate}
\item $\bigcup_{i \in \mathbb{Z}}  N^i E^j = E^j$ and $\bigcap_{i \in \mathbb{Z}}  N^i E^j = 0$,
\item{The filtration $N^\cx$ on each $E^j$ is an $\Cee$-filtration, that is:
\begin{equation}\label{Cee}
\cI N^iE^j \subseteq N^{i+1}E^j.
\end{equation}
(see also Definition~\ref{filtsat.d}).}
\end{enumerate}
The \emph{$\Cee$-filtered derived category $DF({\cal A}, \Cee)$} is
the Verdier quotient of the homotopy category
 $Ho( CF({\cal A}, \Cee))$  of  $CF({\cal A}, \Cee)$ by the subcategory
$Ho( CF^{ac}({\cal A}, \Cee))$ of acyclic complexes.
\end{definition}

In the context of this  definition, a filtered complex is said to be
acyclic if for every $i$ the complex $ (N^{i}E ^{\cx}, d) $ is
acyclic. Recall that by  the definition of the Verdier quotient
there is a triangulated functor
$$L:  Ho( CF({\cal A}, \Cee))\to  DF({\cal A}, \Cee),$$
such that $ L \, (Ho( CF^{ac}({\cal A}, \Cee)))=0$. The pair
$(DF({\cal A}, \Cee), L)$ has the following universal property: for
every triangulated category $T$, the composition with $L$ induces an
equivalence of categories between the full subcategory of
triangulated functors $ \Phi: Ho(CF({\cal A}, \Cee))\to T $, such
that $\Phi \, (Ho(CF^{ac}({\cal A}, \Cee)))=0$, and the category of
triangulated functors  $\Phi ' : DF({\cal A}, \Cee) \to T$.
Explicitly, $DF({\cal A}, \Cee)$ can be constructed as the category
whose objects are those of $Ho( CF({\cal A}, \Cee))$ and morphisms
$Hom _{DF({\cal A}, \Cee)}(X,Y)$ are represented by diagrams
$$ X\str {\alpha} Y' \stackrel {s}{\longleftarrow} Y,$$
where $\alpha $ and $s$  are morphisms in $Ho(CF({\cal A}, \Cee))$
and $cone \, s \in Ho( CF^{ac}({\cal A}, \Cee))$. We refer the
reader to \cite{nee.tc} for a detailed discussion. In the case when
$\Cee =0$, the filtered derived category $DF({\cal A}):= DF({\cal
A}, 0) $ was first considered by Illusie in his thesis
\cite{ill.cctd}.

Given  a filtered  complex $E^\cx$,  we denote by  $E^\cx(r)$ the
same complex  but with the shifted filtration: $N^i(E^\cx(r))=
N^{i+r}E^\cx$. Let  $CF_{\leq l}({\cal A}, \Cee)$ be the full
subcategory of $CF({\cal A}, \Cee)$ whose objects are filtered
complexes with $  N^{l+1}E^\cx=0 $, let  $CF_{\geq k}({\cal A},
\Cee)$  be the full subcategory whose objects satisfy $N^{k}E^\cx=
E^{\cdot }$, and let be $  CF_{[k,l]}({\cal A}, \Cee)  $ the
intersection of $CF_{\leq l}({\cal A}, \Cee)$
 and  $CF_{\geq k}({\cal A}, \Cee)$.
 Denote by   $DF_{\leq l}({\cal A}, \Cee)$, $DF_{\geq k}({\cal A}, \Cee)$,  and
 $DF_{[k,l]}({\cal A},\Cee)$ the quotients of the corresponding homotopy categories.

\begin{lemma}
The functor  $c_{\geq k}:  DF_{\geq k}({\cal A}, \Cee )  \to  DF({\cal A}, \Cee) $ \\
has a right adjoint functor
$$w_{\geq k}:   DF({\cal A}, \Cee )  \to  DF_{\geq k}({\cal A}, \Cee) \quad : \quad
 w_{\geq k}(E^\cx) = N^{k}E^\cx.$$
The functor   $c_{\leq l}:    DF_{\leq l}({\cal A}, \Cee)  \to
DF({\cal A}, \Cee) $ has a left adjoint functor
$$w_{\leq l}:     DF({\cal A}, \Cee)  \to  DF_{\leq l}({\cal A}, \Cee) \quad
 : \quad  w_{\leq l}(E^\cx) =  E^\cx/ N^{l+1}E^\cx . $$
Moreover, $ w_{\geq k} c_{\geq k}\simeq Id , \,  w_{\leq l} c_{\leq
l}\simeq Id$.
\end{lemma}
The proof is straightforward.

\begin{corollary} \label{strogiy}
The functors $c_{\geq k}$,  $c_{\leq l}$  and    $c_{[k,l]}: DF_{[k,
l]}({\cal A}, \Cee)  \to   DF({\cal A}, \Cee)$ are fully faithful.
The essential image of $c_{\ge k}$ consists of those objects
$(E^\cx,N^\cx E^\cx)$ such that each $N^jE^\cx \to E^\cx$ is a
quasi-isomorphism for all $j\le k$, and the essential image of
$c_{\le l}$ consists of those objects such that $N^jE^\cx$ is
acyclic for all $j > l$.
\end{corollary}
\begin{proof}
Indeed, for $E^\cx$, ${E'}^{\cx} \in  DF_{\geq k}({\cal A}, \Cee)$
we have
\begin{eqnarray*}
Hom_{DC_{\geq k}({\cal A}, \Cee)}(E^\cx, {E}^{\prime\cx}) &\simeq &
                Hom_{DF_{\geq k}({\cal A}, \Cee)}(E^\cx, w_{\geq k} c_{\geq k} {E}^{\prime\cx})\\
& \simeq & Hom_{DF({\cal A}, \Cee )}(c_{\geq k}E^\cx, c_{\geq k}
{E}^{\prime\cx}),
 \end{eqnarray*}
where the first isomorphism is induced by $w_{\geq k} c_{\geq
k}\simeq Id $ and the second one comes from the adjointness property
from the lemma. The proofs for $c_{\le l}$ and $c_{[k,l]}$
  are similar.  If $(E^\cx,N^\cx E^\cx )$ is an object of $DF(A,\Cee)$ and
each $N^jE^\cx \to E^\cx $ is a quasi-isomorphism for all $j \le k$,
then the natural map $c_{\ge k}\,w_{\ge k}(E^\cx,N^\cx E^\cx) \to
(E^\cx,N^\cx E^\cx )$ is an isomorphism in $DF(A,\Cee)$, so that
$(E^\cx,N^\cx E^\cx )$ is in the essential image of $c_{\ge k}$. The
proof for $c_{\le l}$ is similar.
\end{proof}

Let $p: {\bf V} \to Z$ be a vector bundle over a scheme $Z$, $V$
the corresponding  sheaf of sections (thus, $V$ is a locally free
sheaf of ${\cal O}_Z$-modules),  and  let ${\cal A}$ be a flat
 sheaf of algebras  over  $p_*{\cal O}_{{\bf V}} \cong S^{\cx} V^*$. Let
${\cal I}$ be the sheaf of ideals in $\cA$ generated by $V^*$.
Denote by $\gr\, {\cal A} =\oplus _{j\geq 0} {\cal I}^j/{\cal
I}^{j+1}$ the sheaf of graded algebras over $S^{\cx} V^*$. Since
${\cal A}$ is flat over $S^{\cx}V^*$ the morphism:
\begin{equation}\label{gra}
 S^{\cdot }V^* \otimes _{{\cal O}_Z} {\cal A}/{\cal I} \to \gr\, {\cal A}, \, f\otimes a \longrightarrow f a
\end{equation}
is an isomorphism and $\gr\, {\cal A} $ is a flat $S^{\cx}
V^*$-module. Denote by $ D(Mod ^{\cxdot}(\gr\, {\cal A}))$ the
derived category of graded $\gr\, {\cal A}$ -modules. We then have a
functor:
\begin{eqnarray*}
\gr: DF({\cal A}, \Cee) &\to&  D(Mod ^{\cxdot}(\gr\, {\cal A}))\\
 (E^\cx) & \mapsto &
      \bigoplus_{-\infty < j< + \infty} N^jE^\cx/N^{j+1}E^\cx.
  \end{eqnarray*}

Let $q: {\bf W} \to Z$ be another vector  bundle over $Z$ and $i:
{\bf W}\hookrightarrow {\bf V}$ a linear embedding. Set ${\cal B}=
{\cal A} \otimes _ {S^{\cx} V^*} S^{\cx} W^*$ and $\Cee'=   W ^*
{\cal B} \subset {\cal B}$. Then ${\cal B}$ is a sheaf of algebras
over $S^{\cx} W^*$ and $\Cee'\subset {\cal B}$ is a subsheaf of
ideals.

\begin{proposition}\label{i!}
 Assume that $Z$ is a noetherian scheme of finite Krull dimension.
 \begin{enumerate}
 \item  The functor
$i_* : DF({\cal B}, \Cee ')  \to  DF({\cal A}, \Cee )$ has a right
adjoint
 $$ Ri^!:   DF({\cal A}, \Cee) \to DF({\cal B}, \Cee ')$$
and the functor
 $i_* : D(Mod^{\cxdot}( \gr \, {\cal B}))  \to  D(Mod^{\cxdot}(\gr \, {\cal A}))$
has a right adjoint:
$$Ri^! : D(Mod^{\cxdot}(\gr , {\cal A}))  \to  D(Mod^{\cxdot}( \gr \, {\cal B})).$$
\item The functor $Ri^!$ takes the essential image of $DF_{\leq l}({\cal A}, \Cee)$
into  the essential image of $DF_{\leq l}({\cal A}, \Cee ')$ and the
essential image of $DF_{\geq k}({\cal A}, \Cee)$ into the essential
image of $DF_{\geq k -d}({\cal B}, \Cee ')$, where $d:= rk\, {\bf V}
- rk \,{\bf W}$.
\item For every $\Cee$-filtered ${\cal A}$ complex $E^\cx$, the morphism $\gr R i^! E^\cx \to R i^! \gr E^\cx $ defined by adjunction:
$$Id\in Hom (E^\cx, E^\cx) \to Hom (\gr i_*Ri^! E^\cx, \gr E^\cx) \to
Hom (\gr Ri^! E^\cx, Ri^! \gr E^\cx)$$ is an isomorphism.
\item $Ri^!$ commutes with the forgetful
functors
$$\Psi: DF({\cal A}, \Cee )\to D(Mod({\cal A}))\, {\textstyle and} \, \Psi' : DF({\cal B}, \Cee ' )\to
D(Mod({\cal B})),$$ i.e. the canonical morphism $\Psi ' Ri^! E^\cx
\to Ri^! \Psi E^\cx, $ defined by adjunction,  is an isomorphism.
\end{enumerate}
\end{proposition}
\begin{proof}
For (1) we use the technique from \cite{nee.gd}. The Brown
Representability Theorem ({\it loc.cit.}, Theorem 4.1.) asserts that
the existence of the adjoint functor $ Ri^!:   DF({\cal A}, \Cee)
\to DF({\cal B}, \Cee ')$ would follow if we prove that
\begin{enumerate}
\item the
 categories $DF({\cal A}, \Cee)$ and  $ DF({\cal B}, \Cee ')$ have arbitrary direct sums
\item the functor $i_*$ commutes with arbitrary direct sums
\item the category  $ DF({\cal B}, \Cee ')$ is compactly generated \footnote{Recall that an object
$X \in DF({\cal B}, \Cee ')$  is called \emph{compact} if for every
set of objects $\{Y_{\alpha} \}$ one has $\oplus Hom(X, Y_{\alpha})
\simeq Hom(X, \oplus Y_{\alpha})$. A category is said to be
\emph{compactly generated} if there exists a set $T$ of compact
objects such that for every nonzero $Y\in  DF({\cal B}, \Cee ')$
there exists $X\in T$ such that $Hom(X, Y)\ne 0$ .}.
\end{enumerate}
The first two properties are immediate.  Let us check the third.
Given an open subset $j: U\hookrightarrow Z$, denote by ${\cal B}_U
$ the filtered ${\cal B}$-module such that $N^i {\cal B}_U = j_!
{\cal B}$ for $i\leq 0$ and
 $N^i {\cal B}_U = j_! (S^i W^*{\cal B})$ for $i>0$. For any $E^\cx\in DF({\cal B}, \Cee ')$,
one has
$$Hom_{DF({\cal B}, {\cal J} )}({\cal B}_U (l), E^\cx[j])\simeq R^j\Gamma (U, N^{-l} E^\cx).$$
It follows that  $DF({\cal B}, {\cal J})$ is generated by objects of
the form ${\cal B}_U (l)$. It is known that for any noetherian space
$U$ of finite Krull dimension the functor $R\Gamma (U, \, )$
commutes with arbitrary direct sums (see, for example
\cite{voev.htss}). Thus the objects $ {\cal B}_U (l) $ are compact.

The second claim in (1) is proven by a similar argument.

For (2),  let $E^\cx \in DF_{\leq l}({\cal A}, \Cee )$.
 We want to show that $w_{\geq l+1} Ri^!(c_{\leq l}\,  E^\cx) =0$.
Indeed, for every $E^{\prime\cx}\in  DF_{\geq l+1} ({\cal B}, \Cee '
)$ we have
\begin{eqnarray*}
Hom(E^{\prime\cx}, w_{\geq l+1}\, Ri^!(c_{\leq l}\, E^\cx))
&\simeq & Hom( c_{\geq l+1}\,  E^{\prime\cx},  Ri^!(c_{\leq l}\, E^\cx)) \\
&\simeq & Hom(w_{\leq l}i_*( c_{\geq l+1}  E^{\prime\cx}),   E^\cx) \\
  &\simeq & Hom( w_{\leq l}\, c_{\geq l+1} \, i_*( E^{\prime\cx}),  E^\cx)\\
  & = & 0.
\end{eqnarray*}

To prove the second statement consider the forgetful functor
$$\Phi:  DF({\cal B}, \Cee ' ) \to DF({\cal O}_Z): = DF({\cal O}_Z, 0) $$
to the filtered derived category of ${\cal O}_Z$-modules. By
Corollary~\ref{strogiy}, we will be done if we show that
 $\Phi Ri^! (E^\cx)\in DF_{\geq k -d}({\cal O}_Z)$
for every object $E^\cx$ of $ DF_{\geq k}({\cal A}, \Cee)$. Consider
the Koszul complex
\begin{equation}\label{kosrez}
0\to   \Lambda ^d  T \otimes _{{\cal O}_Z}  {\cal A}(-d)
     \to \cdots \to
 T\otimes _{{\cal O}_Z}  {\cal A}(-1)  \to   {\cal A} \to  i_*i^* {\cal A} \to 0,
\end{equation}
where $T:= ker\, (V^*\stackrel{i^*}{\longrightarrow} W^*)$, and
where  the $\Cee $-filtration on
 $\Lambda ^m  T \otimes _{{\cal O}_Z}  {\cal A}(-m) $ is defined by
$$ N^i (\Lambda ^m  T \otimes _{{\cal O}_Z}  {\cal A}(-m))= \Lambda ^m  T \otimes _{{\cal O}_Z} \Cee ^{i-m}. $$
Then (\ref{kosrez}) is an acyclic complex in $CF({\cal A}, \Cee )$.
It yields a functorial isomorphism
\begin{equation}\label{kc}
\Phi Ri^! (E^\cx)\simeq \cHom_{{\cal A}}(\Lambda ^{\cx} T \otimes
_{{\cal O}_Z}  {\cal A}(- \cdot), E^\cx)
\end{equation}
This is the filtered complex $C^\cx$ whose term in degree $i$ is
$$C^i:= \bigoplus_{p+q=i}  \Lambda ^{p} T^* \otimes  _{{\cal O}_Z} E^q (p).$$
Since $E^q(p) \in DF_{\ge k-p}(\cA,\Cee)$ and $T^*$ has rank $d$,
this completes the proof.

For the last two statements, it will be
 enough to prove that $\gr R i^! E^\cx \to R i^! \gr E^\cx $ (resp.
$\Psi Ri^! E^\cx \to Ri^! \Psi E^\cx $ ) becomes an isomorphism
after the projection to the derived category of graded ${\cal
O}_Z$-modules (resp. the derived category of ${\cal O}_Z$-modules.)
In turn, this follows from the Koszul computation in (2).
\end{proof}

Let $h:X\to Y$ be a smooth morphism  of relative dimension $d$ of
smooth schemes over a noetherian scheme $S$ of finite Krull
dimension. We shall apply the above construction to the linear
morphism
$${\bf T}^*_{X'\to Y'}  \stackrel{i'_h}{\longrightarrow}    {\bf T}^*_{X'}, $$
and to ${\cal A}\supset \Cee $ being either $S^{\cx}T_{X'}\supset
\J_{X'}:= \bigoplus _{i>0} S^i T_{X'}$ or  the Azumaya algebra  $
F_{X/S *} D_{X/S}\supset \Cee _X : = T_{X'}(F_{X/S *} D_{X/S})$.
 We then have the filtered derived image functors
$$R{h'}^{HIG}_*= Rh'_* \circ R{i'_h}^!: DF( S^{\cx}T_{X'}, \J _{X'}  )\to DF(S^{\cx}T_{Y'}, \J _{Y'} )$$
$$Rh^{DR}_* = Rh'_*\circ C_{X'\to Y'}\circ R{i'_h}^!  :
DF(F_{X/S *} D_{X/S}, \Cee _X  ) \to  DF(F_{Y/S *} D_{Y/S}, \Cee _Y)
$$ and by the previous proposition
$$R{h'}^{HIG}_* : DF_{[k,l]}( S^{\cx}T_{X'}, \J _{X'})\to DF_{[k-d, l]}(S^{\cx}T_{Y'}, \J _{Y'}  )$$
$$Rh^{DR}_* : DF_{[k,l]} ( F_{X/S *} D_{X/S}, \Cee _X  ) \to  DF_{[k-d,l]}(   F_{Y/S *} D_{Y/S}, \Cee _Y ).  $$

In particular, this gives another proof of Corollary~\ref{level.c}.

\begin{example}\label{conjf}
 Consider the $D_{X/S}$-module ${\cal O}_X \in DF(F_{X/S *} D_{X/S}, \Cee _X  )$
endowed with the trivial filtration. Then the filtration on
$$Rh^{DR}_* {\cal O}_X \in   DF_{[-d, 0]}(F_{Y/S *} D_{Y/S}, \Cee _Y  ) $$
coincides with the ``conjugate'' filtration. Indeed, we will
construct a canonical quasi-isomorphism in the filtered derived
category $DF(i^{\prime *}_h (F_{X/S *} D_{X/S}))$:
$$(Ri^{\prime !}_h (F_{X/S *}{\cal O}_{X}), N^\cx) \simeq
 (Ri^{\prime !}_h (F_{X/S *}{\cal O}_{X}), T^\cx),$$
where for any complex $C^\cx$,
$$T^iC^q := \begin{cases}
    C^q &\text{if $q \le -i$} \\
          Im (d^q)  & \text{if $ q = -i+1$} \\
             0 & \text{if $ q > -i+1$.}
  \end{cases}$$
That is,
 $T^iC^\cx =
           \tau_{\le -i} C^\cx$, where $\tau_\le $ is the
             canonical filtration.
Note that by   (\ref{kc}),
$$\gr^{-m} Ri^{\prime !}_h (F_{X/S *}{\cal O}_{X}) \simeq \Omega ^m_{X'/Y'}\otimes F_{X/S *}{\cal O}_{X}[-m].  $$
Thus the result follows from  the following lemma, whose proof is
straightforward.

 \begin{lemma}\label{conj.l}
Let $(E^\cx,N^\cx)$
 be a filtered complex in an abelian category. Assume that
the filtration is finite and that  for every $m$
$$H^i(\gr ^{-m} \, E^\cx) =0, \quad\mbox{for every}\quad i\ne  m. $$
For each $i$, let $T_N^iE^\cx := T^iN^iE^\cx \subseteq N^i E^\cx$.
 Then the morphisms
$$(E^\cx, N^\cx)\leftarrow (E^\cx, T_N^\cx) \to (E^\cx, T^\cx)$$
are filtered quasi-isomorphisms.
 \end{lemma}
\end{example}

Observe that the graded Azumaya algebra
$$\gr\, F_{X/S *}D_{X/S} \simeq ( F_{X/S *}D_{X/S}/ \Cee _X) \otimes _{{\cal O}_{X'}} S^{\cx} T_{X'}$$
 over $S^{\cx} T_{X'}$    splits canonically:
$F_{X/S *}{\cal O}_X \otimes _{{\cal O}_{X'}} S^{\cx} T_{X'}$ is the
graded splitting module. This defines an equivalence of categories :
\begin{eqnarray*}
 C^{\cxdot -1}_{X/S}:  D(HIG^{\cxdot}( X'/S)) &\to&     D(Mod ^{\cxdot}(\gr\,F_{X/S *}  D_{X/S}))\\
C^{\cxdot -1}_{X/S}(E^{\cdot \cxdot}) & = &   E^{\cdot \cxdot}
\otimes _{S^\cx T_{X'/S}} (F_{X/S *}{\cal O}_X
 \ot S^\cx T_{X'/S}) \\
   & \cong &  E^{\cdot \cxdot} \otimes _{{\cal O}_{X'}} F_{X/S *} \oh
   X.
\end{eqnarray*}
 Observe that $C^{\cxdot -1}_{X/S}$ and its quasi-inverse
$C^\cxdot_{X/S}$ commute with $Ri_h^{\prime!}$.   By part (3) of
Proposition \ref{i!} we have a functorial isomorphism
\begin{equation}\label{lkatz}
 C^\cxdot_{X/S}\, \gr \, Ri^{\prime !}_h  (E^\cx) \simeq   Ri^{\prime !} _h  (C^\cxdot_{X/S}\, \gr \, E^\cx)
\end{equation}
and its direct image to $Y'$
\begin{equation}\label{katzformula}
 C^\cxdot_{Y/S}\,  \, Rh^{DR}_* (E^\cx) \simeq Rh^{HIG}_*   (C^\cxdot_{X/S}\,\gr \, E^\cx).
\end{equation}

Let $E^\cx$ be an object of $ DF_{[k,l]} ( F_{X/S *}D_{X/S}, \Cee_X
) $. Then the filtered complex $Rh^{DR}_* (E^\cx)$ yields a spectral
sequence:
 $$E_1^{p,q}= H^{p+q}(\gr^p  Rh^{DR}_* (E^\cx))
\Rightarrow H^*( Rh^{DR}_* (E^\cx) ).$$ We shall call it the
\emph{conjugate spectral sequence} \footnote{Let us note that, when
$E= \oh X $ the $E_r$-terms of our spectral sequence correspond to
the $E_{r+1}$ terms of the usual conjugate spectral sequence, after
a suitable renumbering.}
 (\emph{c.f.} Example \ref{conjf}).

Assume that the conjugate spectral sequence degenerates at $E_1$.
Then the quasi-isomorphism (\ref{katzformula}) induces an
isomorphism of graded Higgs modules:
\begin{equation}\label{ckf}
C^\cxdot_{Y/S}\, \gr \, R^jh^{DR}_* (E^\cx) \simeq R^jh^{HIG}_*  (
C_{X/S}^\cxdot\, \gr \, E^\cx).
\end{equation}
\begin{remark}\label{kasp.r}
Take $E^\cx = \oh X$ and assume that the Hodge  spectral sequence as
well as the conjugate spectral sequence for the de Rham direct image
degenerate at $E_1$. Then (\ref{ckf}) yields an isomorphism of Higgs
modules
\begin{equation}\label{grkatz.e}
C_{Y/S}^\cxdot\, \gr _N \, R^jh^{DR}_* ({\cal O}_X) \simeq (\gr_F
R^jh^{DR}_* (\oh {X'}), \kappa),
\end{equation}
where $\gr_F $ denotes the associated graded object with respect to
the Hodge filtration on $R^jh^{DR}_* (\oh{X'})$ and $\kappa $ is the
Kodaira-Spencer operator viewed as a Higgs field on $\gr_F
Rh_*^{DR}(\oh {X'})$. This is the Katz's \footnote{In {\it loc.cit.}
Katz considers also the log version of his formula. We shall not do
so here.}
 p-curvature formula \cite[Theorem~3.2]{ka.asde}.
See  Example~\ref{conjf} for an explication of the left side which
relates it to Katz's original formulation.
% $$\gr^p_G  Rh^{DR}_*  \oh X \cong Rh^{p+q}_*\gr^p_T\Omega^\cx_{X/S} \cong$$
% $$R^qh_*\cH^p(\Omega^\cx_{X/S}) \cong R^qh'_*( \Omega^p_{X'/S});$$
% a renumbering changes what we call the $E_r$-terms of this spectral sequence into $E_{r+1}$-terms
% appearing in the usual ``second spectral sequence'' of hypercohomology.
% The degeneration of the conjugate spectral sequence at (our) $E_1$
%  implies that the Hodge spectral sequence  also degenerates at $E_1$.\marginpar{check references
% and reorganize}
 We refer the reader to section~\ref{fm.ss} for a generalization of this remark.
\end{remark}

 \begin{remark}
 Example \ref{conjf} can be generalized as follows.
Let ${\cal A}$ be sheaf of algebras flat over $S^{\cx} V^*$ and
$i:{\bf W}\hookrightarrow {\bf V}$ a linear embedding. Consider the
functors
 $$
\def\normalbaselines{\baselineskip20pt
\lineskip3pt  \lineskiplimit3pt}
\def\mapright#1{\smash{
\mathop{\to}\limits^{#1}}}
\def\mapdown#1{\Big\downarrow\rlap
{$\vcenter{\hbox{$\scriptstyle#1$}}$}}
\begin{matrix}
DF({\cal A}, \Cee ) & \str {i^!}  & DF(i^*{\cal A}, i^*\Cee )  \cr
\mapdown{ } &  & \mapdown {} \cr
 DF({\cal A} ) & \str {i^!}  & DF(i^*{\cal A} ).
\end{matrix}
$$
This diagram is not commutative. However, we will show that for
every ${\cal A}$-module $E$ with a finite $\Cee$-filtration $(E= N^0
E \supset \cdots \supset N^n E \supset N^{n+1}E=0)$ the $\Cee$-
filtration $(R i^! E= N^{-d} R i^! E \supset \cdots \supset N^n R
i^! E \supset N^{n+1} R i^!E=0)$ is the \emph{filtration
d\'ecal\'ee} of $(R i^! E = R i^!N^0 E \supset \cdots \supset  R
i^!N^n E \supset Ri^! N^{n+1}E=0)$. To see this we, first, recall an
interpretation of the \emph{filtration d\'ecal\'ee} convenient for
our purposes.

Let $DF({\cal C})$ be the filtered derived category of an abelian
category ${\cal C}$, and let
 $DF^{\leq k}({\cal C}) \subset DF({\cal C})$  (resp. $DF^{\geq k}({\cal C}) \subset DF({\cal C})$) be the full subcategory whose objects are  filtered complexes
$(E^\cx, F^{\cx}E^\cx)$ such that, for every integer $n$,  $\gr^n
E^\cx$ has vanishing cohomology in degrees greater then $n+k$ (resp.
less then $n+k$). It is known \cite[Appendix]{beil.dcps}, that the
subcategories $DF^{\leq k}({\cal C})$ and  $DF^{\geq k}({\cal C})$
define a $t$-structure on $DF({\cal C})$ whose heart is the abelian
category of complexes $C({\cal C})$. In particular, the embedding
 $DF^{\leq k}({\cal C}) \to DF({\cal C})$ has a right adjoint functor
$$\tau_{\leq k}:  DF({\cal C}) \to DF^{\leq k}({\cal C}). $$
Explicitly,
$$ F^m  (\tau_{\leq k} (E^\cx, F^{\cx}E^\cx ))^i =  F^{m+i-k}E^i + d ( F^{m+i-k -1} E^{i-1}), $$
 if $i>k$ and
$$ F^m  (\tau_{\leq k} (E^\cx, F^{\cx}E^\cx))^i =  F^{m}E^i $$
otherwise.

The canonical filtration
$$ \cdots \subset \tau_{\leq k} (E^\cx, F^{\cx}E^\cx ) \subset \cdots \subset (E^\cx, F^{\cx}E^\cx ).$$
makes $E^\cx$ a bifiltered complex. We shall denote this bifiltered
complex by
$$(E^\cx, F^{\cx}E^\cx)^{dec} := (E^\cx, N^{\cx} F ^{\cx}E^\cx),$$
so that $(E^\cx, N^{-k} F^{ \cdot}E^\cx)= \tau_{\leq k} (E^\cx,
F^{\cx}E^\cx ) $. We then have the following generalization of
Lemma~\ref{conj.l}.
\begin{lemma} Let
$ (E^\cx, N^{\cx}F^{\cx} E^\cx)  $ be a bifiltered complex. Assume
that the filtration $N$ is finite, i.e. there exist integers $a$ and
$b$ such that $N^{a}F^{\cx} E^\cx=0 $ and $N^{b}F^{\cx} E^\cx=
N^{b-i}F^{\cx} E^\cx$ for every $i\geq 0$. Set $F^{\cx}E^\cx:=
N^{b}F^{\cx} E^\cx$. Assume also that, for every $m$,
$$\gr _N ^{-m} (E^\cx, N^{\cx} F^{\cx} E^\cx) \in
 DF^{\leq m}({\cal C})\cap DF^{\geq m}({\cal C}),$$
i.e. $H^j(\gr _F ^k \gr _N ^{-m} (E^\cx, N^{\cx} F^{\cx} E^\cx))=0$,
for every $j\ne k+m$. Then the canonical morphism
   $$ (E^\cx, N^{\cx}F^{\cx} E^\cx)\to (E^\cx, F^{\cx} E^\cx)^{dec}  $$
defined as in Example \ref{conjf} is a bifiltered quasi-isomorphism.
\end{lemma}
We omit the proof.

We apply the Lemma to the bifiltered complex  $(R i^! E,
N^{\cx}F^{\cx}R i^! E)$, where $ N^{k}F^{m}R i^! E  = N ^k R i^!  c
_{\geq m } w_{ \geq m } E $ \footnote{ Precisely, $(R i^! E,
N^{\cx}F^{\cx}R i^! E)$  is defined as $Ri^! (E,  N^{\cx}F^{\cx}E)$,
$ N^{k}F^{m}E = N^{max(k,m)}E $  in the bifiltered derived category
of ${\cal A}$-modules $(E^{\cx}, N^{\cx}F^{\cx}E )$ such that $\Cee
N^k F^m E^{\cx} \subset N^{k+1} F^m E^{\cx}$.}. By (\ref{kc}),
$$\gr _F ^k \gr _N ^{-m} R i^! E \simeq \wedge ^{k+m} T^* \otimes _{{\cal O}_Z} \gr ^k E [-k-m] .$$
Thus we get a canonical bifiltered quasi-isomorphism
$$  (R i^! E, N^{\cx}F^{\cx}R i^! E )\simeq (R i^! E,  F^{\cx} R i^! E  )^{dec}.$$
 \end{remark}

\subsection{The derived Cartier transform.}\label{dct.ss}

Let ${\cal X}/S $ be a lifting. For any $k $ and $l$ with $l-k <p$,
the Cartier transform yields equivalences of categories
$$ DF_{[k,l]} (F_{X/S *}D_{X/S}, \Cee_X )  \str {\stackrel{C_{{\cal X}/S }}{\sim}}
 DF_{[k,l]}(S^{\cx}T_{X'}, \J_{X'})    $$
$$ C_{\lift X S}(E^\cx, N^{\cx}E^\cx)= (C_{{\cal X}/S }E^\cx, C_{{\cal X}/S }N^{\cx}E^\cx).$$

\begin{theorem}\label{gmeqpush.t}
a) Let $h:X\to Y$ be a morphism of smooth schemes over $S$. Then,
for any integers $k$ and $l$ with $l-k <p$,  a lifting $\tilde h':
\tilde X' \to \tilde Y'$ of $h'$ induces an isomorphism:
$$Lh_{DR}^* \circ C^{-1}_{\lift Y S }\cong C^{-1}_{\lift X S } \circ Lh^{\prime *}_{HIG}:
DF_{[k,l]} (S^{\cx}T_{Y'}, \J_{Y'}) \to  DF_{[k,l]} (F_{X/S
*}D_{X/S}, \Cee_X ). $$ b) If in addition $h$ is smooth of relative
dimension $d$ and $l-k -d <p$, then
$$Rh_*^{DR} \circ C^{-1}_{\lift X S } \cong C^{-1}_{\lift Y S } \circ R{h'_*}^{HIG}:
DF_{[k,l]} (S^{\cx}T_{X'}, \J_{X'}) \to  DF_{[k-d,l]} (F_{Y/S
*}D_{Y/S}, \Cee_Y ).$$
 \end{theorem}
\begin{proof}
a) Define an equivalence of categories
$$ DF_{[k,l]}(i^{\prime *}_h (S^{\cx}T_{X'}), \J'_{X'})
\str {\stackrel{(C^Y_{\lift X S })^{-1}}{\sim}}  DF_{[k,l]}
(i^{\prime *}_h (F_{X/S *}D_{X/S}), \Cee '_X )  $$ to be the
composition
$$ (C^Y_{\lift X S })^{-1}:= {\cal M}_{\lift X S } \circ  \iota _*,$$
where $\iota_*:  DF_{[k,l]}(i^{\prime *}_h (S^{\cx}T_{X'}),
\J'_{X'})\to DF_{[k,l]}(i^{\prime *}_h (S^{\cx}T_{X'}), \J'_{X'})$
is the involution defined in  (\ref{iotahi.e}) and  ${\cal M}_{{\cal
X}/S }$ is the tensor product with the splitting module $F_{X/S
*}{\cal B}_{\lift X S}$:
$${\cal M}_{{\cal X}/S } (E^{\cx}, N^{\cx}E^\cx ) =
(E^{\cx}\otimes _ {\hat \Gamma T_{X'/S}} F_{X/S *}{\cal B}_{\lift X
S}, N^{\cx}E^{\cx}\otimes _ {\hat \Gamma T_{X'/S}} F_{X/S *} {\cal
B}_{\lift X S}). $$ Similarly, the splitting module $h^{\prime
*}F_{Y/S *}{\cal B}_{\lift Y S}$ yields an equivalence of categories
$$ DF_{[k,l]}(h^{\prime *}(S^{\cx}T_{Y'}), h^{\prime *} \J_{Y'})
\str {\stackrel{(C^{X}_{\lift Y S })^{-1}}{\sim}}
  DF_{[k,l]} (h^{\prime *}(F_{Y/S *}D_{Y/S}), h^{\prime *} \Cee _Y ).  $$
\begin{lemma}\label{invcomp.p}
A morphism $(h,\tilde h') \colon \lift X S\to \lift Y S$ induces an
isomorphism of functors
$$ C^Y_{\lift X S } \simeq C^X_{\lift Y S } \circ C_{X' \to Y'}$$
%  $$
% \def\normalbaselines{\baselineskip20pt
% \lineskip3pt  \lineskiplimit3pt}
% \def\mapright#1{\smash{
% \mathop{\to}\limits^{#1}}}
% \def\mapdown#1{\Big\downarrow\rlap
% {$\vcenter{\hbox{$\scriptstyle#1$}}$}}
% %
% \begin{matrix}
%  DF_{[k,l]} (i^{\prime *}_h (F_{X/S *}D_{X/S}), \Cee '_X )  & \str{C_{X'\to Y'} } &
%   DF_{[k,l]} (h^{\prime *}(F_{X/S *}D_{Y/S}), h^{\prime *} \Cee _Y )    \cr
% \mapdown{C^Y_{{\cal X}/S } } &  & \mapdown {C^X_{{\cal Y}/S }} \cr
%  DF_{[k,l]}(i^{\prime *}_h (S^{\cx}T_{X'}), \J'_{X'}) & \str {=} &
%  DF_{[k,l]}(h^{\prime *}(S^{\cx}T_{Y'}), h^{\prime *} \J_{Y'}).         .
% \end{matrix}
% $$
\begin{diagram}
 DF_{[k,l]} (i^{\prime *}_h (F_{X/S *}D_{X/S}), \Cee '_X )  & \rTo^{C_{X'\to Y'} } &
  DF_{[k,l]} (h^{\prime *}(F_{X/S *}D_{Y/S}), h^{\prime *} \Cee _Y )    \cr
\dTo^{C^Y_{\lift X S } } &  & \dTo_ {C^X_{\lift Y S }} \cr
 DF_{[k,l]}(i^{\prime *}_h (S^{\cx}T_{X'}), \J'_{X'}) & \rTo^{=} &
 DF_{[k,l]}(h^{\prime *}(S^{\cx}T_{Y'}), h^{\prime *} \J_{Y'}).
\end{diagram}
\end{lemma}
\begin{proof}

Recall from Proposition~\ref{exfunct.p} that a morphism $ (h,\tilde
h') \colon   \lift X S \to \lift Y S$ induces an isomorphism
$$ h^*\cA_{\lift Y S} \str{\cong} {\cal H}om_{F_{X/S}^*\Gamma^{\cx} T_{X'/S} }
(h^*F_{Y/S}^*\Gamma^{\cx} T_{Y'/S},\cA_{\lift X S}).  $$ Dualizing
this isomorphism, we find an isomorphism of $D_{X/S}$-modules
$$ F_{Y/S}^* h^{\prime *} \hat \Gamma T_{Y'/S}\otimes_{F^*_{X/S} \hat\Gamma T_{X'/S}} \cB_{\lift X S}
\cong h^*\cB_{\lift Y S}.  $$ With the notations of Theorem
\ref{hxty.t}, we have
$$F _{X/S *} h^*\cB_{\lift Y S} \simeq  F_{X/S *}D_{X\to Y} \otimes
_{ h^{\prime *} F_{Y/S *}D_{Y/S} } h^{\prime *} F_{Y/S *} \cB_{\lift
Y S}.$$
 Thus we get  an isomorphism of splitting modules for
$ h^{\prime *} \hat \Gamma  T_{Y'/S}\otimes _{S^{\cx}T_{X'/S}}
F_{X/S *} D_{X/S}$:
\begin{equation}\label{f} h^{\prime *} \hat \Gamma  T_{Y'/S}\otimes_{\hat\Gamma T_{X'/S}}  F_{X/S *} \cB_{\lift X S} \cong
 F_{X/S *}D_{X\to Y} \otimes
_{ h^{\prime *} F_{Y/S *}D_{Y/S} } h^{\prime *} F_{Y/S *} \cB_{\lift
Y S}.
\end{equation}
By definition, the functor  $(C^Y_{{\lift X S} })^{-1}$ is the
composition of the involution $\iota_* $ and the tensor product over
$ h^{\prime *} \hat \Gamma  T_{Y'/S}$ with the left-hand side of
(\ref{f}), and the functor $(C^X_{\lift Y S } \circ C_{X' \to
Y'})^{-1}$ is the composition of $\iota_* $ and the tensor product
with the right-hand side of (\ref{f}). Thus, (\ref{f}) induces the
desired  isomorphism $(C^Y_{{\lift X S} })^{-1}\simeq (C^X_{\lift Y
S } \circ C_{X' \to Y'})^{-1}$.
\end{proof}

Let us return to the proof of the theorem. Observe the natural
isomorphisms of functors:
\begin{eqnarray*}
Lh^{\prime *} \circ (C_{\lift Y S })^{-1} &\simeq& (C^X_{\lift Y S
})^{-1} \circ Lh^{\prime *} \quad\mbox{ and}
\end{eqnarray*}
\begin{eqnarray*}
  i'_{h *} \circ (C^Y_{{\lift X S} })^{-1} &\simeq& (C_{{\lift X S} })^{-1} \circ i'_{h *}.
\end{eqnarray*}
Hence, by Lemma \ref{invcomp.p}
$$ Lh_{DR}^* \circ
C^{-1}_{\lift Y S } \simeq i'_{h *}\circ  (C_{X'\to Y'})^{-1} \circ
Lh^{\prime *} \circ C^{-1}_{\lift Y S } \simeq  i'_{h *}\circ
(C_{X'\to Y'})^{-1} \circ (C^X_{\lift Y S })^{-1} \circ Lh^{\prime
*} $$
$$ \simeq  i'_{h *} \circ (C^Y_{{\lift X S} })^{-1} \circ Lh^{\prime *}
\simeq C^{-1}_{{\lift X S} } \circ Lh^{\prime *}_{HIG}.$$
This proves part a) of the Theorem. \\
b) By Lemma \ref{invcomp.p} it remains to construct an isomorphism
of functors:
\begin{equation}\label{comp}
R{i'_h}^! C_{{\lift X S}  } \cong C^Y_{{\lift X S} }R{i'_h}^! :
DF_{[k,l]} (F_{X/S *}D_{X/S}, \Cee_X )
 \to  DF_{[k-d,l]} (i^{\prime *}_h (S^{\cx}T_{X'}), \J'_{X'} ).
\end{equation}
Let $E^\cx\in DF_{[k,l]} (F_{X/S *}D_{X/S}, \Cee_X )$ and
 ${E}^{\prime \cx} \in   DF_{[k-d,l]} (i^{\prime *}_h (S^{\cx}T_{X'}), \J'_{X'})$. We then have functorial
 isomorphisms
$$Hom( {E}^{\prime \cx}, R{i'_h}^! C_{{\lift X S} }(E^\cx))\simeq  Hom(  C_{{\lift X S} }^{-1}  i'_{h*}
{E}^{\prime \cx},  E^\cx) \simeq $$
$$Hom( i'_{h *} (C^Y_{{\lift X S} })^{-1}   {E}^{\prime \cx},  E^\cx) \simeq
Hom(  {E}^{\prime \cx},  C^Y_{{\lift X S} }R{i'_h}^! (E^\cx)).$$ By
the Yoneda lemma this yields (\ref{comp}).
\end{proof}

\begin{remark}
 In the absence of the lifting of $ h'$ the theorem can be modified as follows.
Let ${\cal L}_{h'} $ be the ${h'}^*T_{Y'/S}$-torsor of liftings of
$h'$ and let $ exp\, {\cal L}_{h'}$ be the pushforward of ${\cal
L}_{h'} $ via the homomorphism
           $$ exp: {h'}^*T_{Y'} \to ({h'}^* \hat \Gamma T_{Y'})^* .$$
Thus $ exp\, {\cal L}_{h'}$ is a $({h'}^* \hat \Gamma
T_{Y'})^*$-torsor. We denote by ${\cal K}_{h'}$ the corresponding
 invertible module over ${h'}^* \hat \Gamma  T_{Y'}$.
Define an autoequivalence  $$\tau_{h'}: Mod({h'}^* \hat \Gamma
T_{Y'}) \to Mod({h'}^* \hat \Gamma  T_{Y'})$$
 $$\tau_{h'}( E ) = E \otimes _{{h'}^* \hat \Gamma  T_{Y'} } {\cal K}_{h'}.$$
Then, with the notations from the proof of Theorem \ref{gmeqpush.t},
one has
 $${\cal K}_{h'} \otimes_{\hat\Gamma T_{X'/S}}  F_{X/S *} \cB_{\lift X S} \cong
 F_{X/S *}D_{X\to Y} \otimes
_{ h^{\prime *} F_{Y/S *}D_{Y/S} } h^{\prime *} F_{Y/S *} \cB_{\lift
Y S},$$
\begin{eqnarray*}
 Lh_{DR}^* \circ C^{-1}_{\lift Y S }&\cong& C^{-1}_{{\lift X S} } \circ i'_{h *} \circ \tau_{h'}^{-1}\circ Lh^{\prime *} \\
DF_{[k,l]} (S^{\cx}T_{Y'}, \J_{Y'}) &\to&  DF_{[k,l]} (F_{X/S
*}D_{X/S}, \Cee_X )
\end{eqnarray*}
and if $h$ is smooth of relative dimension $d$ and $l-k+d< p$
\begin{eqnarray*}
Rh_*^{DR} \circ C^{-1}_{\lift X S } &\cong& C^{-1}_{\lift Y S } \circ R{h'_*} \circ \tau_{h'}\circ  Ri^{\prime !}_h \\
DF_{[k,l]} (S^\cx T_{X'}, \J_{X'}) &\to&  DF_{[k-d,l]} (F_{Y/S
*}D_{Y/S}, \Cee_Y ).
\end{eqnarray*}
\end{remark}

\begin{corollary}
Let $h:X\to Y$ be a smooth morphism of relative dimension $d$ and
let $E^{\cx}$  be an object of $D(HIG_{[k,l]}^{\cxdot}( X'/S))$.
Assume that $l-k -d <p$ and that there exists $\tilde h': \tilde X'
\to \tilde Y'$. Then the conjugate spectral sequence for $H^*(
Rh^{DR}_* (C^{-1}_{{\lift X S} } E^{\cx}) )$ degenerates at $E_1$.
\end{corollary}
\begin{proof}
We have
$$Rh^{HIG}_* (E^{\cx}) \simeq Rh^{HIG}_* (\gr\, E^{\cx})\simeq \gr\, Rh^{HIG}_* (E^{\cx}).$$
Here the first isomorphism comes from the grading on $E^{\cx}$ and
the second  one from  (3) of Proposition~\ref{i!}. It follows that
the spectral sequence of the filtered complex $Rh^{HIG}_* (E^{\cx})$
degenerates at $E_1$. Then by Theorem \ref{gmeqpush.t} the same is
true for $Rh^{DR}_* (C^{-1}_{{\lift X S} } E^{\cx}) $.
 \end{proof}

\section{Applications  and examples}
\subsection{Local study of the $p$-curvature}\label{lpc.ss}
Let $X/S$ be a smooth morphism of schemes in characteristic $p >0$
and let
\begin{equation}
\Psi \colon \mic(X/S) \to \fhig(X/S)
\end{equation}
denote the functor taking a module with integrable connection to the
corresponding module with $F$-Higgs field. This functor is not an
equivalence or even fully faithful. For example, the category of
pairs $(\oh,\nabla)$ with vanishing $p$-curvature is equivalent to
the category of invertible sheaves $L$ on $X'$ together with a
trivialization $F_{X/S}^*L \cong \oh X$. However, we show that if
$(E_1,\nabla_1)$ and $(E_2,\nabla_2)$ are two noetherian objects of
$\mic(X/S)$ with isomorphic images in $\fhig(X/S)$, then Zariski
locally on $X$, $(E_1,\nabla_1)$ and $(E_2,\nabla_2)$ are
isomorphic.  Moreover, we can characterize the image of the functor
$\Psi$, \'etale locally on $X$: if $\psi$ is an $F$-Higgs field on a
coherent $E$, then \'etale locally on $X/S$, $\psi$ comes from a
connection if and only if $(E,\psi)$ descends to a Higgs field on
$X'/S$. Taken together, these results can be interpreted as a
nonabelian analog of the well-known exact
sequence~\cite[4.14]{mi.ec}.
$$ 0  \rTo \oh {X'}^* \rTo^{F_{X/S}^*}  F_{X/S*}\oh X^*\rTo^{dlog}  F_{X/S*}Z^1_{X/S}
 \rTo^{\pi^* - C_{X/S}}
 \Omega^1_{X'/S} \rTo 0,$$
where $C_{X/S}$ is the Cartier operator and $\pi \colon X' \to X$
the projection. Indeed, one can recover this sequence by considering
the category of connections of the form $(\oh X,d + \omega)$, where
$\omega$ is a closed one-form, and recalling that the $p$-curvature
of such a connection is precisely $\pi^*(\omega) - C_{X/S}(\omega)$.

\begin{theorem}\label{locpsi.t}
Let $X/S$ be a smooth morphism of noetherian schemes in
characteristic $p$.
\begin{enumerate}
\item{Let $(E_i,\nabla_i)$, $i = 1,2$,  be objects of $\mic(X/S)$, with $E_i$
coherent, and let $\psi_i$ denote their $p$-curvatures. Suppose that
there exists an isomorphism $h\colon  (E_1,\psi_1) \to (E_2,\psi_2)$
in $\fhig(X/S)$.  Then Zariski locally on $X$, $(E_1,\nabla_1)$ and
$(E_2,\nabla_2)$ are  isomorphic in $\mic(X/S)$.}
\item{Let $E$ be a coherent sheaf with an F-Higgs field $\psi \colon E \to E \ot F_{X/S}^*\Omega^1_{X'/S}$.
Then \'etale locally on $X$, the following are equivalent:
\begin{enumerate}
\item{ There exists a connection on $E$ whose
$p$-curvature is $\psi$.}
\item{There exist a coherent sheaf  with a Higgs field $(E',\psi')$
on $X'/S$ and an isomorphism $(E,\psi) \cong F^*_{X/S}(E',\psi')$.}
\end{enumerate}}
\end{enumerate}
\end{theorem}
\begin{proof}
To prove (1), let $H := \Hom(E_1,E_2)$, with the internal Hom
connection and $p$-curvature. Let $H^\psi \subseteq H$ be the
subsheaf  annihilated by $\psi$, and let $F_{X/S*}(H^\nabla)$ be the
subsheaf  annihilated by $\nabla$. Then by Cartier descent, the
natural map $F_{X/S}^*F_{X/S_*}H^\nabla \to H^\psi$ is an
isomorphism of $\oh X$-modules.

 Let $x$ be a  point of $X$, and let $x'$  be its image in $X'$.
Then $k(x)$ is a finite and purely inseparable extension of $k(x')$.
The fiber   $V' := (F_{X/S*} H^\nabla)(x')$ of $F_{X/S*}H^\nabla$ at
$x'$ is  a finite dimensional $k(x')$-vector space, the fiber
 $V := H^\psi(x)$ of $H^\psi$ at $x$ is a finite
dimensional $k(x)$-vector space, and  the natural map $k(x)
\otimes_{k(x')} V' \to V$ is an isomorphism. There is also a
natural map $V \to \Hom_{k(x)}(E_1(x), E_2(x))$. Let $\bv$ be the
affine space over $k(x)$ corresponding to the $k(x)$-vector space
$V$, and let ${\bf U}$ denote the Zariski open subset of $\bv$
corresponding to those elements which define isomorphisms $E_1(x)
\to E_2(x)$. The isomorphism $h$ lies in $H^\psi$ and hence its
image $h(x)$ in $V$ corresponds to a $k(x)$-rational point of ${\bf
U}$. Let $\bv'$ be the affine space over $k(x')$ corresponding to
$V'$. Then $\bv$ is the  base change of $\bv'$ to $\spec k(x)$, and
since $k(x') \to k(x)$ is a purely inseparable extension, the
projection mapping $\bv \to \bv'$ is a homeomorphism and the image
${\bf U}'$ of ${\bf U}$ in $\bv'$ is a nonempty open subset. If
$k(x')$ is infinite, it follows  that the $k(x')$-rational points of
${\bf V'}$ are Zariski dense, so ${\bf U}'$ has a $k(x')$-rational
point.  If $k(x')$ is finite, it is perfect, and it follows that
$k(x) = k(x')$. Thus in either case there is an element $v'$ in $V'$
which induces an isomorphism $E_1(x) \to E_2(x)$. Then there exists
an element $g'$ in the stalk of the $\oh {X'}$-module
$F_{X/S*}H^\nabla$ at $x'$ whose image in $V'$  is $v'$.   Let $h'
:= F_X^*(g')$, which  defines a horizontal morphism $E_1 \to E_2$ in
some neighborhood of $x$.  The fiber of $h'$ at $x$ is an
isomorphism.  We know that $E_{1,x}$ and $E_{2,x}$ are isomorphic as
$\oh {X,x}$-modules, and in particular their reductions module any
power of the maximal ideal have the same finite length. It follows
from Nakayama's lemma that $h'$ is surjective modulo any power of
the maximal ideal, and hence is also an isomorphism modulo any such
power. Then it follows that $h'$ is an isomorphism in a neighborhood
of $x$. This proves (1).

We should remark that (1) could also have been proved from the
theory of Azumaya algebras; we preferred to explain the elementary
proof above.  We do not know of such an elementary proof of (2).
Note first that since (2) is a local statement, we may assume that
there exists a spitting $\zeta$ of $C^{-1}_{X/S}$ as in
(\ref{zeta.d}).

%
%\begin{lemma}\label{splitloc.l}
%Let $E'$ be a coherent sehaf on  $\bT^*_{X'/S}$ whose
%projection to $X'$  is again coherent.
%Then \'etale locally on $X'$, there exists a coherent sheaf $E''$ on $\bT^*_{X'/S}$
%and an isomorphism $E' \cong \alpha_{\zeta*}(E'')$.
%Furthermore, $\cB_\zeta$
%splits the Azumaya algebra $\cD_{X/S}$ on $\bT^*_{X'/S}$  when pulled back
%to the s
%\end{lemma}
%\begin{proof}
%Let $i'\colon Z' \to \bT_{X'/S}$ be the closed immersion
%defined by the annihilator of $E'$ in $\oh {\bT^*_{X'/S}}$.
% Since $E'$ is coherent as a sheaf of $\oh {X'}$-modules,  $Z'$ is finite over $X'$,
%and hence the  \'etale covering
%$\alpha_\zeta  \colon \bT_{X'/S} \to \bT_{X'/S} $ splits over $Z'$,
% \'etale locally on $X'$.  Thus, after
%replacing $X'$ by an \'etale open neighbhorhood of any point, we may assume
%that there exists a map $j' \colon Z' \to \bT^*_{X'/S}$
%such that $\alpha_\zeta\circ j' = i'$.  Let $E'' := j'_*{i'}^*E'$.
%Then $E' \cong i'_* i^{\prime*} E' \cong \alpha_{\zeta*} j'_* {i'}^*E' \cong \alpha_{\zeta *}E''$.
%\end{proof}

%We can now easily prove (2).
Suppose that $(E',\psi')$ is an object of $HIG(X'/S)$, with $E'$
coherent as an $\oh {X'}$-module.  Let $\tilde E'$ denote the
coherent sheaf on $\bT^*_{X'/S}$ corresponding to $(E',\psi')$ Let
$i'\colon Z' \to \bT_{X'/S}$ be the closed immersion defined by the
annihilator of $\tilde E'$ in $\oh {\bT^*_{X'/S}}$.
 Since $E'$ is coherent as a sheaf of $\oh {X'}$-modules,  $Z'$ is finite over $X'$,
and hence the  \'etale covering $\alpha_\zeta  \colon \bT_{X'/S} \to
\bT_{X'/S} $ splits over $Z'$,
 \'etale locally on $X'$.  Thus, after
replacing $X'$ by an \'etale localization, we may assume that there
exists a map $j' \colon Z' \to \bT^*_{X'/S}$ such that
$\alpha_\zeta\circ j' = i'$.  Let $\tilde E'' := j'_*{i'}^*\tilde
E'$, which corresponds to an object $(E'',\psi'')$ of $\hig(X'/S)$.
Then $\tilde E' \cong i'_* i^{\prime*} \tilde E' \cong
\alpha_{\zeta*} j'_* {i'}^*\tilde E' \cong \alpha_{\zeta *}\tilde
E''$. Let $(E,\nabla) := \Psi^{-1}_\zeta(E'',\psi')$ (see
Theorem~\ref{zetas.t}).  By {\em op. cit.},
 the $p$-curvature of $(E,\nabla)$ is
$F_{X/S}^*\alpha_{\zeta*}(E'',\psi'') \cong F_{X/S}^* (E',\psi')$.

Conversely, suppose that $(E,\nabla)$ is an object of $\mic(X/S)$,
with $E$ coherent as an $\oh X$-module.  Its $p$-curvature defines
an object $(E,\psi)$ of $F$-$\hig(X/S)$, and hence a
 coherent sheaf $\tilde E$ on $\bT^{*(X')}_{X'/S} := \bv F_{X/S}^* T_{X'/S}$
(see diagram (\ref{btp.e})).  The  claim  is that there exists a
coherent sheaf $\tilde E'$ on $\bT_{X'/S}$ such that
$\pi_{\bT}^*(\tilde E') \cong \tilde E$. Since $F_{X/S*}E$ is
coherent as an $\oh{X'}$-module, the scheme-theoretic support $Z'$
of $\pi_{\bT*}$  is finite over $X'$, and  there exists a section
$j'$ of $\alpha_\zeta$ over $Z'$.  If we view $\tilde E$ as a module
over $S_\zeta^\cx T_{X'/S}$ via $j^{\prime\sharp}$, then the action
of $S_\zeta^\cx T_{X'/S}$ agrees with the action of $S^\cx
T_{X'/S}$, and so the action of $D_{X/S}$ on $E$ extends to an
action of $D_\zeta$. Let $\tilde E' := \cHom_{\cD_\zeta}(\cB_\zeta ,
E)$, corresponding to an object $(E',\psi') \in \hig(X'/S)$.  Then
$(E,\nabla) \cong \Psi^{-1}_\zeta(E',\psi')$, so $(E,\psi) \cong
F_{X/S}^*(\alpha_{\zeta*}(E',\psi'))$ in $F$-$\hig(X/S)$, by
Theorem~\ref{zetas.t}.
% and there exists a section
%$j'$ of $\alpha_\zeta$ over $Z'$.  Then the Azumaya algebra
%
%  Consider the object $\cH := \cHom_{\oh X}(\cB_\zeta\,E)$ with
%it induced connection $\nabla$, and let $E' := \cH^\nabla$, regarded as a sheaf
%of $\oh {X'}$-modules.
%
% By Cartier descent, $F_{X/S}^*E' \cong \cH^\psi$, which
%can be identified with $\cH^\psi = \cH_{\cB_\zeta}(E)$ in the notation of (\ref{}).
%
%
%
%
%
%
%  Then $\psi$ defines the structure
%of a sheaf of $\bT^*_{X'/S}$-modules $E_\psi$ on $F_{X/S*}E$, which is coherent
%as a sheaf of $\oh {X'}$-modules.  Let $E''$ be the sheaf of $\oh {\bT^*_{X'/S}}$-modules
%as in Lemma (\ref{splitloc.l}), after some \'etale localization of $X$.
%  Then $E'' \cong F_{X'/S*}(E)$ as a sheaf of $\oh{ X'}$-modules
%Thus the action $\cZ_{X/S} $ on $E$ given by $\psi$
%extends to an action of $S^\cx T_{X'/S}$ through $\alpha_\zeta$, and
%the action of $D_{X/S}$ on $E$ extends to an action of $D_\zeta$.
%Then $(E,\nabla) \cong \Psi^{-1}_\zeta(E',\psi')$, where $(E',\psi') =
%\cHom_{D_\zeta}(\cB_\zeta, E)$. In particular, $(E,\psi) \cong F^*(E', \psi'')$,
%where $\psi'' = \alpha^*_\zeta(\psi')$.
\end{proof}

\subsection{Stacks of liftings and splittings}\label{sls.ss}
In this subsection we discuss  relationships between and geometric
interpretations of some of the liftings and splittings used in our
constructions. In particular, we show that there is a natural
equivalence  between the gerbe of liftings of $X'$ and the gerbe of
{\it tensor} splittings of $\cD_{X/S}$ over the completed divided
power envelope  $\hat \bT^{* \gamma}_{X'/S}$ of the zero section of
$\bT^*_{X'/S}$.

First we  shall study the gerbe of splittings of the Azumaya algebra
${\cal D}_{X/S}$ on $\bT^*_{X'/S}$.  Recall from \cite{mi.ec} and
\cite{de.tc} that the equivalence  class of this gerbe can be viewed
as the image of ${\cal D}_{X/S}$ in the cohomological Brauer group
$H^2(\bT^*_{X'/S}, {\cal O}^*_{ \bT^*_{X'/S}})$. Our first goal is
to provide a simple description of this cohomology class.

Recall from  \cite[4.14]{mi.ec} that for any smooth $Y/S$ there is
an exact sequence of \'etale sheaves on $Y'$:
\begin{equation}\label{ohy*.e}
 0  \rTo \oh {Y'}^* \rTo^{F_{Y/S}^*} F_{Y/S *} \oh Y^*\rTo^{dlog} F_{Y/S *} Z^1_{Y/S} \rTo^{\pi_{Y/S}^* - C_{Y/S}}
 \Omega^1_{Y'/S} \rTo 0.
 \end{equation}
Here $F_{Y/S *}Z^1_{Y/S} \subset F_{Y/S*}\Omega^1_{Y/S} $ is the
subsheaf of closed one-forms, $C_{Y/S}$ is the Cartier operator, and
$\pi_{Y/S}: Y'\to Y$ is the morphism induced by the Frobenius on
$S$. As we observed in  section \ref{lpc.ss}, the  morphism
$\pi_{Y/S}^*-C_{Y/S}: F_{Y/S *} Z^1_{Y/S}\to \Omega^1_{Y'/S}$ can be
viewed as the map sending the line bundle ${\cal O}_Y$ with
integrable connection $\nabla = d +\omega$ to its p-curvature. The
exact sequence \ref{ohy*.e} induces  a morphism:
     $$\phi: H^0(Y', \Omega^1_{Y'/S})\to H^1(Y', F_{Y/S *}({\cal O}^*_Y )/{\cal O}^*_{Y'})\to
H^2(Y', {\cal O}^*_{Y'})=Br(Y').$$ As we shall  recall  below, the
cotangent bundle of $X'/S$ has a canonical global one-form (the
``contact form''). We shall see in  Proposition~\ref{brauer.p} below
that the Brauer class of $\cD_{X/S}$ can be identified with the
image of this one-form under the map $\phi$. We begin with the
following  convenient geometric description of  the map $\phi$.

\begin{proposition}\label{phid.p}
Let $\omegap \in  H^0(Y', \Omega^1_{Y'/S})$ be a one-form. For each
\'etale  $U' \to  Y'$, let  $U := F_{Y/S}^{-1}(U') \to Y$ and let
${\cal P}^\nat_{\omegap}(U')$ be the groupoid of invertible sheaves
with integrable connection on $U$ whose p-curvature is equal to
$\omegap$.  Then,  ${\cal P}^\nat_{\omegap}$ forms a fibered
category which is in fact a gerbe under ${\cal O}^*_{Y'}$ on $Y'$.
The class of ${\cal P}^\nat_{\omegap}$ in  $H^2(Y', {\cal
O}^*_{Y'})$ is equal to $\phi(\omegap)$.
\end{proposition}
\begin{proof}
It is clear that ${\cal P}^\nat_{\omegap}$ forms a stack
 and that the automorphism group of each object
is $\oh {Y'}^*$.  The  local surjectivity of $\pi_{Y/S}^*- C_{Y/S}$
implies that the class of objects of ${\cal P}^\nat_{\omegap} $ is
locally not empty. If $L_1$ and $L_2$ are two objects of ${\cal
P}^\nat_{\omegap} $ over some $U'$,  then the $p$-curvature of
$\cHom(L_1,L_2)$ is zero, and hence locally has a horizontal basis.
This implies that any two objects of ${\cal P}^\nat_{\omegap} $ are
locally isomorphic, so that ${\cal P}^\nat_{\omegap} $ is indeed a
gerbe.

 The boundary map associated to the exact sequence
$$ 0 \to F_{Y/S*}(\oh Y^*)/\oh {Y'}^* \rTo^{dlog}  F_{Y/S *}Z^1_{Y /S} \rTo^{\pi_{Y/S}^* - C_{Y/S}}
 \Omega^1_{Y'/S} \rTo 0$$
takes $\omegap$ to the $(F_{Y/S*}\oh Y^*)/\oh {Y'}^*$-torsor  ${\cal
T}_{\omegap}$ of closed one-forms $\eta$ such that $\pi^*_{Y/S}\eta
- C_{Y/S}(\eta) = \omegap$. The boundary map associated to the exact
sequence
$$ 0 \to \oh {Y'}^* \to F_{Y/S *}(\oh Y ^*) \to F_{Y/S*}(\oh Y^*)/\oh {Y'}^* \to 0$$
takes ${\cal T}_{\omegap}$ to the  gerbe $\cG_{\omegap}$ of
$F_{Y/S*}(\oh Y^*)$-torsors $\cL$ equipped with an isomorphism
$\alpha \colon \overline \cL \to {\cal T}_{\omegap}$, where
$\overline \cL$ is the $(F_{Y/S*}\oh Y^*)/\oh {Y'}^*$-torsor
associated to $\cL$. Hence $\phi(\omegap) = \cG_{\omegap}$, and it
will suffice to prove that $\cG_{\omegap}$ is equivalent to ${\cal
P}^\nat_{\omegap}$. Let $\cL$ be an object of $\cG_{\omegap}$ over
$U'$, let $L$ be the associated invertible sheaf over $U$, and let
$e$ be a local section of $\cL$, \ie, a  basis for  $L$ on some open
subset $V$ of $U$.  There is a unique connection $\nabla$ on $L$
such that $\nabla(e) = e \otimes \alpha(e)$. It follows from the
fact that $\alpha$ is a morphism of torsors that $\nabla$ is
independent of the choice of $e$, and it is clear that the
$p$-curvature of $\nabla$ is $\omegap$.  This construction defines a
functor from the gerbe $\cG_{\omegap}$ to ${\cal P}^\nat_{\omegap}$,
which  is easily seen to be an equivalence.
\end{proof}

\begin{remark}{\rm
In the context of the above proposition, the form $\omegap$ gives a
morphism $i: Y' \to T^*_{Y/S'}$, and ${\cal P}^\nat_{\omegap} $ is
the gerbe of splittings of the Azumaya algebra $i^*\cD_{Y/S}$ on
$Y'$. }\end{remark}

Let us write $\bT^*$ for $\bT^*_{X/S}$, and recall that there is an
exact sequence
\begin{equation}\label{omegato.e}
 0 \to pr^*\Omega^1_{X/S} \to \Omega^1_{\bT^*_{X/S}/S} \to \Omega^1_{\bT^*_{X/S}/X} \to 0
\end{equation}
Furthermore,  $\bT^* = \spec_X S^\cx T_{X/S}$, so that there is a
canonical global section of $pr_*pr^*\Omega^1_{X/S} \cong
\Omega^1_{X/S}\ot S^\cx T_{X/S}$, corresponding to the identity
element of ${\Omega^1_{X/S} \ot T_{X/S}} \cong \End T_{X/S}$.
 The image of this section in $\Omega^1_{\bT^*/S}$ is the well-known
``contact form'' on the cotangent bundle.

\begin{proposition}[\cite{bb.glpc}]\label{brauer.p}
Let $\omegap \in \Gamma(\bT^*_{X'/S}, \Omega^1_{\bT^*_{X'/S}})$ be
the  contact form and let $\mathcal{P}^\nat_{\omegap}$ be the
corresponding
 $\gm$-gerbe on  $\bT^*_{X'/S}$ described in  Proposition~\ref{phid.p}.  Then the gerbe
 ${\cal P}^\nat_{\omegap}$ is equivalent to the gerbe
$\mathcal{S}$ of splittings of the Azumaya algebra
 ${\cal D}_{X/S}$ on $\bT^*_{X'/S}$.  In particular, the class of ${\cal D}_{X/S}$ in $Br(\bT^*_{X'/S})$
is $\phi(\omegap)$.
\end{proposition}
\begin{proof}
We have a diagram:
 \begin{diagram}
\bT^*_{X/S} & \rTo^{F_{\bT^*/X}}&{\bT}^{\prime *}_{X/S}&\rTo^{\pi_T}
& \bT^*_{X'/S}  & \rTo &\bT^*_{X/S} \cr &\rdTo&\dTo&&\dTo && \dTo\cr
 &&X & \rTo^{F_{X/S}} & X'&\rTo^\pi & X,
\end{diagram}
in which both squares are Cartesian and $F_{\bT^*/S} = \pi_\bT\circ
F_{\bT^*/X}$.
 We identify the pullback
of $\bT^*_{X/S}$  by $F_S$ with $\bT^*_{X'/S}$ and use
abbreviations:
$$ \bT^* := \bT^*_{X/S} , \quad \bT^{*\prime} := \bT^*_{X'/S}, \quad \bT^{\prime *}
 := \bT^{\prime *}_{X/S}.$$
%  Note from the fact that $Y^{(X)} = X\times_{X'} Y'$
%that $\oh {\bT^{\prime*}}$ has a natural structure of a $\cD_{X/S}$-module.
Let $U' \to  \bT^{*\prime}$ be \'etale, let $U \to \bT^*$ (resp.
$U''$) be its pullback via $F_{\bT^*/S}$, (resp. via $\pi_T$). Let
$(L,\nabla)$ be an object of ${\cal P}^\nat_{\omegap}(U')$, \ie, an
invertible sheaf with integrable connection on $U/S$ whose
$p$-curvature is $\omegap$. The connection $\nabla$ defines an
action of $D_{\bT^*/S}$ and hence of the subalgebra $D_{\bT^*/X}$ on
$L$. Since the projection of $\omegap$ to
$\Omega^1_{\bT^{*\prime}/X'}$ is equal to $0$, the $p$-curvature of
the corresponding object of $\mic(\bT^*/X)$ vanishes. Let
$$L' := \cH^0_{dR}(L\ot \Omega^\cx_{\bT^*/X}) := Ker \left(
L \rTo^\nabla L \ot \Omega^1_{\bT^*/S} \to L \ot
pr^*\Omega^1_{\bT^*/X} \right).$$ Then $L'$ has a natural structure
of a sheaf of $\oh {\bT^{\prime*}}$-modules on $U''$,  and it
follows from Cartier descent that the natural map $F_{\bT^*/S}^*L'
\to L$ is an isomorphism.  Furthermore, $\nabla$ induces a map
$\nabla' \colon L' \to L' \ot pr^*\Omega^1_{X/S}$, which defines a
$pr^{-1}D_{X/S}$-module structure on $L'$. (This is essentially the
Gauss-Manin connection for the morphism $\bT^* \to X$.) The
$p$-curvature of this module is still given by the contact form
$\omegap$, which means that the action of sections of $\oh
{\bT^{*\prime}}$ via the $p$-curvature is the same as the action via
the map $\bT^{\prime*} \to \bT^{*\prime}$ and the given $\oh
{\bT^{\prime*}}$-structure. This means that we can safely view the
$pr^{-1} D_{X/S}$-module structure and the $\oh
{\bT^{\prime*}}$-module structure as defining a $\cD_{X/S}$-module
structure on $L'$.
 Since $L'$ is an invertible sheaf on $\bT^{\prime*}$,
it has rank $p^d$ over $\bT^{*\prime}$, and thus defines a splitting
module for the Azumaya algebra $\cD_{X/S}$. Thus we have defined a
functor ${\cal P}_{\omegap} \to {\cal S}$.  It is clear that this
functor is fully faithful, since the automorphisms of objects in
either category are just give by units in $\oh {\bT^{*\prime}}$. On
the other hand, suppose that $M$ is a splitting module for
$\cD_{X/S}$. Then viewing $\oh X \to \cD_{X/S}$ via the action on
the left, we can view $M$ as a module over $\bT^{\prime*}$, and by
Proposition~\ref{rankone.p} it then becomes an invertible sheaf of
$\oh{ \bT^{\prime*}}$-modules.
  Since the $\oh{\bT^{\prime *}}$-module structure of $M$
comes from its $p$-curvature, the $p$-curvature of $M$  is just the
contact form $\omegap$. A local calculation shows that
 there is a unique
extension of the action of $\cD_{X/S}$ on $M$ to an action of
$D_{\bT^*/S}$  on $F_{\bT^*/X}^*M$ with the property that $M$ is the
annihilator of $D_{\bT^*/X}$.  This shows that the functor ${\cal
P}_{\omegap}^\nat \to \cS$ is an equivalence. The statement about
the Brauer group then follows, as explained in \cite{mi.ec}.

\end{proof}

In this following discussion we will assume that the reader is
acquainted with the notion of tensor structure on an Azumaya algebra
introduced in section~\ref{aags.ss}. In particular, we explained
there that the algebra $ {\cal D}_{X/S}$ has a canonical symmetric
tensor structure.
  Let us consider the following stacks on $X'_\et$.
\begin{enumerate}
\item{The stack $\cL$ of liftings of $X'$}.
\item{The stack $\tspl$ of tensor splittings of $\cD_{X/S}$ over the
 completed divided power envelope  $\hat \bT^{* \gamma}_{X'/S}$ of the zero section
of $\bT^*_{X'/S}$. \footnote{Note that the \'etale topologies of
$X$,  $X'$,  and  $\hat \bT^{*\gamma}_{X'/S}$ are the same.}}
\item{The stack $\spl_1$ of pairs $(M_1,\alpha)$, where $M_1$ is a splitting
of $\cD_{X/S}$ over the first infinitesimal neighborhood $\bT^*_1$
of the zero section of $\bT^*_{X'/S}$ and
 $\alpha: i^* M_1 \simeq F_{X/S *} {\cal O}_X $ is
 an isomorphism between the restriction of $M_1$ to the zero section and the canonical splitting  over $X'$. }
\item{The stack $\mathcal{EX} $ of extensions of $F_{X/S}^*\Omega^1_{X'/S}$ by $\oh X$ in $\mic(X/S)$
such that the  graded $p$-curvature mapping $\psi\colon
F_{X/S}^*\Omega^1_{X'/S} \to  \oh X \otimes
F_{X/S}^*\Omega^1_{X'/S}$ is the identity.}
\end{enumerate}

In the  discussion preparing  for Theorem~\ref{cart.t}  we
constructed a functor $\cB$ associating a  tensor splitting
$\cB_{\lift X S}$ to a lifting $\tX'$ of $X'$. Furthermore, recall
that $\tilde X'$ determines an extension (\ref{nexc.e}) as in (4),
so that we also have a functor $\cE \colon \cL \to \mathcal{EX}$.
Recall that for any tensor splitting  $M$  there is a canonical
isomorphism $\alpha: i^*M \simeq F_{X/S *} {\cal O}_X $, and hence
there is a restriction functor $i_1^*: {\tspl} \to  {\spl}_1$.  The
dual of an extension in $\mathcal{EX}$ is an object of $\spl_1$, so
there is also a functor from $\mathcal{EX} \to \spl_1$.  This
functor is easily seen to be an equivalence.  The following theorem,
shows that in fact all the above functors are equivalences.

\begin{theorem}\label{gerbel.t}
The stacks above are in fact gerbes, and the functors
$$ \cB \colon \cL \to \tspl, \quad i_1^*\colon \tspl \to \spl_1, \ \mbox{and} \quad \cE \colon \cL \to \mathcal{EX}$$
are equivalences.
\end{theorem}
\begin{proof}
It is clear that $\cL$ and $\spl_1$ are gerbes. The fact that
$i^*_1$ is an equivalence is proven in Proposition~\ref{tsp.p}, and
it follows that $\tspl$ is also a gerbe. Thus, it suffices to prove
that the composition $i^*_1\circ {\cB}:  \cL \to \spl_1$ is an
equivalence.
 Let us show that, for any lifting $\tilde X ^{\prime}$, the group of automorphisms
of  $\tilde X ^{\prime}$ reducing to the identity on $X'$
 maps isomorphically to the group of automorphisms of
$(i^*_1  {\cal B}_{\lift X S} ,\alpha) $. Indeed, the first group
can be identified with the group of vector fields on $ X ^{\prime}$
, and  the second one with  the group of invertible functions on
$X^{\prime}_1$  equal to $1$ on $X'$, and the map is the obvious
isomorphism between this two groups. The following easy and well
known result completes the proof.
\begin{lemma} Let $F: {\cal M} \to {\cal N}$ be a morphism of
gerbes on  $Y_\et$. Assume that for every \'etale morphism $U\to Y$
and every object $C\in {\cal M}(U)$ the induced map
 $$F_* : Aut (C) \to Aut(F(C))$$
is an isomorphism. Then $F$ is an equivalence of gerbes.
\end{lemma}

\end{proof}

Let $\theta \in   T_{X^{\prime}/S}(U)$ be a vector field on
$U\subset X^{\prime}$. We may view $\theta$ as a linear function on
the cotangent space $ {\bf T}^*_{U/S}$. Then the exponential
$exp(\theta)= \sum \frac{\theta ^i}{i!}$  makes sense as an
invertible  function on the  completed PD envelope
 $   {\bf T}^{* \gamma}_{U/S}    \subset  {\bf T}^{* \gamma}_{X'/S}  $. Thus we get a homomorphism of sheaves:
$$exp:  T_{X^{\prime}/S } \to     {\cal O}^*_{  {\bf T}^{* \gamma}_{X'/S}  }=(\hat \Gamma T_{X'/S})^* .         $$
This, in turn, gives a map:
$$exp: H^*_{et}(X^{\prime}; T_{ X'/S } ) \to
 H^*_{et}( X^{\prime}; {\cal O}^*_{  {\bf T}^{* \gamma}_{X'/S}  }) .   $$

In the following corollary we use $\hat \Gamma F^*_{X/S} T_{X'/S}$
-module structure on $\cB_{\lift X S}$ as introduced in
subsection~\ref{chf.s}.

\begin{corollary}
\begin{enumerate}
 \item Let $\theta \in  H^0( X', T_{X^{\prime}/S})$ be an automorphism of a lifting
$\lift X S$ reducing to the identity on $X'$. Then the induced
morphism
$$\theta_*: \cB_{\lift X S}\to \cB_{\lift X S}$$
is the multiplication by  $F^*_{X/S} (exp \, \theta) \in (\hat
\Gamma F^*_{X/S} T_{X'/S})^*$.
\item Let $(\lift X S)_1$, $(\lift X S)_2$
be liftings,  and let ${\cal L}_{Id}$ be the $T_{X'/S}$-torsor of
isomorphisms between $\tilde X_1$ and $\tilde X'_2$ reducing to the
identity on $X'$. Denote by $exp \, {\cal L}_{Id}$ the corresponding
${\cal O}^*_{ {\bf T}^{* \gamma}_{X'/S} }$-torsor and by ${\cal
K}_{Id}$ the corresponding invertible sheaf on ${\bf T}^{*
\gamma}_{X'/S} $. Then the isomorphism of
$F^*_{X/S}T_{X'/S}$-torsors ${\cal L}_{{\lift X S}_1}\otimes
_{F^*_{X/S} T_{X'/S}} F^*_{X/S}{\cal L}_{Id} \simeq {\cal L}_{{\lift
X S}_2}$ induces a tensor isomorphism of splitting modules
$$\cB_{(\lift X S)_1} \otimes _{ \hat \Gamma
F^*_{X/S} T_{X'/S}}  F^*_{X/S} {\cal K}_{Id} \simeq \cB_{(\lift X
S)_2}.$$
\item The class of the Azumaya algebra $ {\cal D}_{X/S}$ restricted
to $ {\bf T}^{* \gamma}_{X'/S} $ in the cohomological Brauer  group
$Br( {\bf T}^{* \gamma}_{X'/S} )= H^2_{et}( X^{\prime}; {\cal O}^*_{
{\bf T}^{* \gamma}_{X'/S}  })$ is equal to $exp\, \delta$, where
$\delta \in   H^2_{et}(X^{\prime}; T_{X^{\prime}/S })$ is the
obstruction to lifting of $X^{\prime}$ over $\tilde S $.
\end{enumerate}
\end{corollary}
\begin{proof}  Since $\theta _*$ and $exp \, \theta $ are tensor
automorphisms of $\cB_{\lift X S}$, by Theorem \ref{gerbel.t} it is
enough to check that $\theta _*$ and $exp \, \theta $ are equal when
restricted to $i^*_1 \cB_{\lift X S}$. In turn, this follows from
the fact that the automorphism of ${\cal L}_{\lift X S}$ induced by
the automorphism of the lifting $\lift X S$ coincides with the
translation by $F^* _{X/S}\theta \in H^0(X, F^*_{X/S} T_{X'/S})$.
This proves (1).
 The proof of the second claim is similar, and  the
 last claim follows from Proposition~\ref{azexp.p}.
 \end{proof}

\begin{remark}
The construction of  the tensor splittings in the proof of
Proposition~\ref{azexp.p}
 can be viewed in the present setting as follows.
 Let $exp \, {\cal L}_{\lift X S}$ be the pushforward of the $F_{X/S}^*T_{X'/S}$-torsor
${\cal L}_{\lift X S}$ via the homomorphism
$$ exp:  F_{X/S}^*T_{X'/S} \to (\hat \Gamma F_{X/S}^* T_{X'/S})^*.  $$
The $(\hat \Gamma F_{X/S}^* T_{X'/S})^*$-torsor
 $exp \, {\cal L}_{\lift X S}$ acquires the induced connection, as does the
associated invertible $\hat \Gamma F^* T_{X'/S}$-module $exp \,
{\cal L}_{\lift X S}\otimes _{ (\hat \Gamma F_{X/S}^* T_{X'/S})^*}
\hat \Gamma F_{X/S}^* T_{X'/S}$. We then have a horizontal
isomorphism
 $${\cal B}_{\tilde X} \simeq exp \, {\cal
L}_{\lift X S}\otimes _{ (\hat \Gamma F_{X/S}^* T_{X'/S})^*}  \hat
\Gamma F_{X/S}^* T_{X'/S}.$$
\end{remark}

Let us end by explaining the relationships between the various
liftings, splittings, and extensions we have been considering.
Consider the exact sequence of $\oh {X'}$-modules:
$$0 \rTo  F_{X/S*}B^1_{X/S} \rTo F_{X/S*}Z^1_{X/S}  \rTo  F_{X/S*}\cH^1_{DR}(X/S)  \rTo 0.$$
A splitting of this sequence amounts to lifting $\zeta$ of
$C^{-1}_{X/S}$ as in (\ref{zeta.d}). Let $\bM_{X/S}$ denote the
sheaf on $X$ which to every open set $U$ assigns the set of liftings
of $C^{-1}_{X/S}$ over $U$. If $U$ is an open subset of $X$, let
$\bL_{X/S}(U)$ denote the category whose objects are morphisms $\tF
\colon \tU \to \tU'$ lifting the relative Frobenius morphism
$F_{U/S} \colon U \to U'$ and whose morphisms $\tF_1 \to \tF_2$  are
commutative diagrams
\begin{diagram}
\tU_1 & \rTo^{\tF_1} & \tU'_1 \cr \dTo^f && \dTo^{f' } \cr \tU_2 &
\rTo^{\tF_2} & \tU'_2,
\end{diagram}
where $f$ and $f'$ reduce to the identity modulo $p$.  In
particular, $f$ and $f'$ are necessarily isomorphisms, and
$\bL_{X/S}$ defines a stack over $S$. As a variant, consider the
stack $\bJ_{X/S}$ which over each $U$ is the category whose objects
are pairs $(\tU',s)$, where $\tU'$ is a lift of $U'$ and $s$ is a
section of the  torsor $\cL_{\lift U S}$ defined by $\tU'$  as in
Theorem~\ref{fliftc.t}
 and whose morphisms
are those reducing to the identity and compatible with $s$. If we
are given a fixed lifting $\tilde X'/\tS$ of $X'/S$, then we can
also consider the fibered category $\bL_{\lift X S}$ which to every
open set $U$ in $X$ assigns the category of pairs $(\tU, \tF)$,
where $\tU$ is a lift of $U$ and $\tF \colon \tU \to \tX'$ is a lift
of $f_{U/S}$.  Morphisms in this category are diagrams as above, in
which $f'$ is the identity.  If $\tU$ is a fixed lifting of $U$,
recall that $\cL_{\lift X S}(\tU)$ is the set of all liftings of
$f_{U/S}$, so there is a natural map from $\cL_{\lift X S}$ to the
sheaf of objects of $\bL_{\lift X S}$. Finally, if $\tF_S \colon \tS
\to \tS$ is a lift of the Frobenius endomorphism of $S$ we can
define a  more rigid version of $\bL_{X/S}$. If $U$ is an open
subset of $X$, let $\bK_{X/S}(U)$ denote the subcategory of
$\bL_{X/S}(U)$ whose objects are liftings  $\tF \colon \tU \to \tU'$
of $F_{U/S}$ with $\tU'  = \tS \times_{\tF_S} \tU$ and whose
morphisms are diagrams as above with $f' = f \times_{\tF_S} \id_S$.

\begin{proposition}\label{stacks.p}
Let $\ov \bL_{X/S}$ denote the sheaf associated to the presheaf of
isomorphism classes of objects of $\bL_{X/S}$, and use the analogous
notation for $\bL_{\lift X S}$.
\begin{enumerate}
\item{The stack $\bJ_{X/S}$ is rigid, and the natural map $\bL_{ X /S} \to \bJ_{X/S}$
induces an isomorphism $\ov\bL_{X/S}\to \bJ_{X/S}$.}
\item{The map (\ref{mazcart.d}) $\tF \mapsto \zeta_\tF$ induces an isomorphism
$\ov \bL_{X/S} \to \bM_{X/S}$ and hence also $\bJ_{X/S} \cong
\bM_{X/S}$.}
\item{ The natural map
$\cL_{{\lift X S},X} \to \overline \bL_{{\lift X S}}$ is an
isomorphism.}
\item{If $\tF_\tS$ lifts $F_S$, then $\bK_{X/S}$ is rigid, and if $S$ is the spectrum
of a perfect field, then  $\tF \mapsto \zeta_\tF$ induces an
isomorphism $\bK_{X/S} \to \bM_{X/S}$.}
\end{enumerate}
\end{proposition}
\begin{proof}

 The following lemma
follows from standard deformation theory and Remark~\ref{mapeq.r};
we omit its proof.

\begin{lemma}\label{defstack.l}
Let $\tX$ and $\tX'$ be liftings of $X$ and $X'$ respectively.  Then
\begin{enumerate}
\item{The sheaf of liftings
$\tF \colon \tX \to \tX'$ of $F_{X/S}$ is a torsor under
$F_{X/S}^*T_{X'/S}$, under the standard action.}
\item{If $\tF_1$ and $\tF_2  \colon \tX \to \tX'$ lift $F_{X/S}$ and differ
by a section  $h'$ of $f_{X/S}^*T_{X'/S}$, then $\tF_1$ is
isomorphic to $\tF_2$  in $\bL_{X/S}$ if and only if $ h'$ comes
from a section of $T_{X'/S}$.}
\item{If $f$ is an automorphism of $\tX$ lifting the identity, then $\tF \circ f = \tF$;
if $f'$ is an automorphism of $\tX'$ lifting the identity such that
$f'\circ \tF = \tF$, then $f' = \id$.}
\item{The sheaf of automorphisms of an object $\tF$ of $\bL_{X/S}$ is canonically isomorphic
to $T_{X/S}$. }
\end{enumerate}
\end{lemma}\qed

Suppose that $(\tU',s)$ is a section of $\bJ_{X/S}$ over $U$.  Then,
locally on $U$, there exists a lift $\tU$ of $U$ and a lift $\tF
\colon \tU \to \tU'$ inducing $s$.  Then an automorphism $f$ of
$(\tU',s)$ corresponds to an automorphism of $\tU'$ reducing to the
identity and such that $f \circ \tF = \tF$.  By
Lemma~\ref{defstack.l}, $f$ is the identity, \ie, $\bJ$ is rigid.
It follows that the natural functor $\bL_{X/S} \to \bJ_{X/S}$
factors through $\ov \bL_{X/S}$, and the above argument makes it
clear that this morphism is surjective.  The injectivity follows
from the definitions.

It follows from the lemma that $\ov \bL_{X/S}$ is a torsor under
$$  C :=\cok ( T_{X'/S}\to F_{X/S*} F_{X/S}^* T_{X'/S}) \cong \cHom (\Omega^1_{X'/S} ,B^1_{X'/S})$$
where
$$ (F_{X/S*} \oh X)/ \oh {X'} \cong B^1_{X/S} \subseteq F_{X/S}^*\Omega^1_{X/S}$$
 is the sheaf of boundaries.
The sheaf  $\bM_{X/S}$ is also naturally a torsor under
$\cHom(\Omega^1_{X'/S},F_{X/S*}B^1_{X/S})$, and the map $\tF \mapsto
\zeta_\tF$ factors through $\ov \bL_{X/S}$:
$$ \ov \bL_{X/S} \to \bM_{X/S} : \tF \mapsto \zeta_\tF.$$
This map is a morphism of torsors, hence a bijection. Now suppose
that $\tF_\tS$ exists and suppose that $\tF$ is an object of
$\bK_{X/S}(U)$.  Then an automorphism of $\tF$ is a an automorphism
$\tf$ of $\tU$ lifting $\id_U$  such that $\tf' \tF = \tF \tf$.
where $\tf' := \tf\times_{\tF_\tS} \id$. But then it follows from
the lemma that $\tf = \id$, so $\bK_{X/S}$ is rigid and its presheaf
of isomorphism classes is a sheaf. Let $\tF_1$ and $\tF_2$ be two
objects of $\bK_{X/S}(U)$.  After shrinking $U$, $\tU_1$ and $\tU_2$
become isomorphic; let us assume they are equal. Then $\tF_1$ is
isomorphic to $\tF_2$ if and only if there exists a lifting $\tf$ of
the identity such that $\tF_2 = \tf'  \tF_1 f^{-1}$. But $\tf' \tF_1
f^{-1}  = \tf' \tF_1$, and if $\tf$ corresponds to an element $D$ of
$T_{X/k}$, $\tf' \tF_1$  differs from $\tF_1$ by the action of
$\pi^*D$. This shows that $\bK_{X/S}$ is a torsor under the cokernel
of the map
$$\pi^{_1}T_{X/S} \rTo^{\pi*} T_{X'/S} \to F_{X/S*} F_{X/S}^* T_{X'/S}.$$
When $S$ is the spectrum of a perfect field, $\pi^*$ is an
isomorphism, and it follows that $\bK_{X/S}$ is also a torsor under
$C$.

Statement (3) can be checked at the stalks. Let $x$ be a point of $U
\subseteq X$ and let $\tU_1$ be a liftings of $U$. Then the stalk of
$\cL_{\lift X S,U} = \cL_{\lift X S,\tU_1}$ at $x$ is the set germs
at $x$ of lifts of $f_{U/S}$ to $\tU_1$, and the stalk of $\overline
\bL_{\lift X S }$ at $x$ is the set of germs of isomorphism classes
of of lifts $(\tU_2,\tF)$ of $f_{U/S}$. Let $\tF \colon \tU_2\to
\tX'$ be a lift of $f_{U/S}$ in some neighborhood of $x$. Then there
is an isomorphism $\tU_1 \cong \tU_2$ near $x$, and this shows that
the map is surjective.  For the injectivity, observe that if $\tF$
and $\tF'$ are elements of $\cL_{\lift X S}(\tU_1)$ which become
equal in $\overline \bL_{\lift X S,U}$, then there is an
automorphism of $\tU_1$ which is the identity mod $p$ and which
takes $\tF$ to $\tF'$. But then by Remark~\ref{mapeq.r}, $\tF =
\tF'$.  This shows the injectivity.
\end{proof}

\subsection{Line bundles with connection}
We use the following notation. If $X$ is a scheme over a field $k$,
$E$ is a coherent sheaf of $\oh X$-modules on  $k$, and $S$ is a
$k$-scheme,
$$\bH^i(X,E)(S) := H^0(S,\oh S) \ot_k H^i(X,E).$$
If $H^i(X,E)$  is finite dimensional,  the functor $\bH^i(X,E)$ is
represented by the (vector) $k$-scheme $\spec S^\cx H^i(X,E)^\vee$.

Let $X/k$ be a smooth proper  geometrically connected scheme over a
perfect field of characteristic $p >0$, with a $k$-rational point
$x_0$. Let $\Pic^\nat_X(S)$ denote the set of isomorphism classes of
triples $(L,\nabla,\alpha)$, where $L$ is an invertible sheaf on
$X\times S$, $\nabla$ is an integrable connection on $L$ relative to
$S$, and $\alpha$ is an isomorphism $L \cong \oh X$ over $x_0\times
S$. Forgetting $\nabla$ defines a morphism $b$ from $\Pic^\nat_{X}$
to the Picard scheme $\Pic_X$ of $X$.  If $L$ is an invertible sheaf
on $X\times S$, the set of integrable connections on $L$ is either
empty or a torsor under the group   $\bH^0(X,Z^1_{X\times S/S})$ of
closed one-forms on $X\times S/S$. Note that formation of the latter
commutes with base change and that $H^0(X, Z^1_{X/k}) \cong
H^0(X',F_*Z^1_{X/k})$.  Thus $H^0(X,Z^1_{X\times S/S}) \cong
\bH^0(X',F_*Z^1_{X/k})(S)$. The Chern class map $\dlog\colon \oh X^*
\to Z^1_{X/k}$ defines a morphism $c \colon \Pic_X \to
\bH^1(X,Z^1_{X/k})$, and there is thus an exact sequence:
$$ 0 \to \bH^0(X,Z^1_{X/k}) \rTo \Pic_X^\nat \rTo^b
\Pic_X \rTo^c \bH^1(X,Z^1_{X/k}).$$ The proof of the following is
then immediate (and well-known).

\begin{proposition}\label{reppic.p}
The above sequence is exact as a sequence of sheaves in the flat
topology.  Furthermore, the functor $\Pic^\nat_X$ is representable,
and its tangent space at the origin is  canonical isomorphic to
$H^1_{dR}(X/k)$.
\end{proposition} \qed

If $(L,\nabla)$ is an object of $\Pic_X^\nat(S)$, its $p$-curvature
can be viewed as an element of $H^0(X'\times S,\Omega^1_{X'\times
S/S})$. This defines a  morphism  of group schemes $\psi \colon
\Pic_X^\nat \to \bH^0(X',\Omega^1_{X'/k})$. If $L'$ is an invertible
sheaf on $X'\times S$ trivialized along $x'_0\times S$, then
$(F_{X/S}\times \id_S)^*L$ is an invertible sheaf on $X\times S$,
and we can equip it with its canonical Frobenius descent connection
to obtain an element of $\Pic^\nat_X(S)$.  This defines a morphism
of group schemes $\phi \colon \Pic_{X'} \to \Pic^\nat_{X}$. An
element in the kernel of $b$ is given by an integrable  connection
on $\oh{ X\times S}$, relative to $S$ \ie, a closed one-form $\omega
\in \Omega^1_{X\times S/S}$, and the $p$-curvature of the
corresponding connection is $\pi_{X\times S/S}^*(\omega) -
C_{X\times S/S}(\omega)$, where $C_{X\times /S}$ is the Cartier
operator~\cite[]{ka.asde}. \marginpar{check citation} Thus there is
a commutative diagram:
\begin{equation}\label{picn.e}
\begin{diagram}
&& \Pic_{X'} \cr && \dTo_\phi & \rdTo^{F_{X/S}^*}\cr
\bH^0(X,Z^1_{X/k}) &\rTo^a &\Pic_X^\nat &\rTo^b & \Pic_X \cr
 &\rdTo_{\pi_{X/k}^*- C_{X/k}}&\dTo_\psi  \cr
&&\bH^0(X',\Omega^1_{X'/k}),
\end{diagram}
\end{equation}
where $\pi_{X/k}^*$ is the composition:
\begin{multline}\label{pixk.e}
\bH^0(X,Z^1_{X/k}) \to \bH^0(X,\Omega^1_{X/k}) \to \\
\bH^0(X',\Omega^1_{X'/k}) = \bH^0(X,\Omega^1_{X/k})\times_{F_k^*} k
= \bH^0(X,\Omega^1_{X/k})'.
\end{multline}
Here the map $\bH^0(X,\Omega^1_{X/k}) \to \bH^0(X,\Omega^1_{X/k})'$
is the relative Frobenius map of the $k$-scheme
$\bH^0(X,\Omega^1_{X/k})$. The map $C_{X/k}$ in the diagram
 is the map of group schemes induced by
the linear map of vector spaces
$$ C_{X/k} \colon H^0(X,Z^1_{X/k}) \to H^0(X',\Omega^1_{X/k}).$$

Recall that there are two spectral sequences converging to de Rham
cohomology: the Hodge spectral sequence, with $E_1^{i,j} =
H^j(X,\Omega^i_{X/k})$, and the conjugate spectral sequence, with
$E_2^{i,j} = H^i(X,\cH^j_{dR}) \cong H^i(X',\Omega^j_{X'/k})$.

\begin{lemma}\label{difpic.l}
In the corresponding diagram of tangent spaces:
\begin{diagram}
&& H^1(X',\oh X') \cr && \dTo^{d\phi} &\rdTo^{F_{X/k}^*}\cr
H^0(X,Z^1_{X/k}) &\rTo^{da} & H^1_{dR}(X/k) &\rTo^{db} & H^1(X,\oh
X) \cr
 &\rdTo_{-C_{X/k}}&\dTo_{d\psi}  \cr
&&H^0(X',\Omega^1_{X'/k})
\end{diagram}
$db$ (resp. $-d\psi$) is the edge homomorphism coming from the Hodge
(resp. conjugate) spectral sequence, and $d\psi \circ da = -
C_{X/k}$.
\end{lemma}
\begin{proof}
Since $\pi^*_{X/k}$ in diagram~\ref{picn.e} factors through the
relative Frobenius map
 in formula~\ref{pixk.e} above,
its differential is zero.  Since $C_{X/k}$ is $k$-linear, it follows
that $d\psi\circ da = - C_{X/k}$. To compute $d\psi$, let $S :=
\spec k[\ep]$, let $\eta \in  H^1_{DR}(X/k)$ and  let $(L, \nabla)$
be the corresponding line bundle with connection over $X\times S$.
Then $d\psi(\eta)$ is a section of $H^0(X',\Omega^1_{X'/k})$ and is
determined by its restriction to any nonempty open subset of $X'$.
We can choose an open subset on which $L$ is trivial, and hence
reduce to the previous calculation.  This proves the claim.
\end{proof}

As we have seen, $F_{X/k*}D_{X/k}$ defines  an Azumaya algebra
$\cD_{X/k}$ over $\bT^*_{X'/k}$; we shall study the splitting of the
pullback of this Azumaya algebra along the canonical map $q\colon
X'\times \bH^0(X',\Omega^1_{X'/k} ) \to \bT^*_{X'/k}$. The universal
$(L,\nabla)$ defines a module over the pullback of $\cD_{X/k}$ to
$X'\times \Pic^\nat_X$, and since it is locally free of rank
$p^{\dim X}$, it is a splitting module.  More generally, suppose we
are given a morphism of $k$-schemes $ f\colon Z \to
\bH^0(X'/k,\Omega^1_{X'/k})$ and a splitting  module $L$ over the
pullback of $\cD_{X/k}$ to $X'\times Z$ via the map $\id_{X'}\times
f$.  Then $L$ is a coherent sheaf on $X'\times Z$ equipped with an
action of the differential operators $F_{X\times Z/Z*}(D_{X\times
Z/Z})$, and in particular can be regarded as a coherent sheaf  with
integrable connection on $X\times Z/Z$ whose $p$-curvature is equal
to the section of $\Omega^1_{X'\times Z/Z}$ defined by $f$. By
Proposition~\ref{rankone.p}, $L$ is an invertible sheaf on $X\times
Z$. By a \emph{rigidified splitting of $\cD_{X/k}$ along $f$} we
shall mean a pair $(L,\alpha)$, where $L$ is a splitting module for
$(\id_{X'}\times f)^*q^*\cD_{X/k}$ and $\alpha$ is a trivialization
of  the restriction of $L$ to  $x_0\times Z$. Thus the universal
$(L,\nabla,\alpha)$ is a rigidified splitting of $\cD_{X/k}$ along
$\psi$.

\begin{proposition}\label{lifspl.p}
Let $ f\colon Z \to \bH^0(X',\Omega^1_{X/k})$ be a morphism and let
$(L,\nabla,\alpha)$ be the universal rigidified line bundle with
connection on $X \times \Pic^\nat_X$.
\begin{enumerate}
\item{The map $f \mapsto f^*(L,\nabla,\alpha)$ is a
bijection between the set of isomorphism classes of rigidified
splittings of $\cD_{X/k}$
 and the set of  maps $\tilde f \colon Z \to \Pic^\nat_X$ such that $\psi\tilde f = f$.}
\item{If $f$ as above is a morphism of commutative group schemes, then  under the
the bijection above, the tensor splittings of $(\id \times
f)^*q^*F_*D_{X/k}$ correspond to the group morphisms $\tilde f$ with
$\psi\tilde f = f$.}
\end{enumerate}
\end{proposition}
\begin{proof}
We have seen that a rigidified splitting of $(\id \times
f)^*q^*F_*D_{X/k}$ gives an invertible sheaf $(M,\nabla,\alpha)$
with connection on $X\times Z$ whose $p$-curvature is given by $f$
and a trivialization of $M$ on $x_0\times Z$.  Hence there is a
unique map $\tilde f\colon Z \to \Pic^\nat_X$
 such that $\tilde f^*(L,\nabla,\alpha) \cong (M,\nabla,\alpha)$,
and necessarily $\psi\tilde f = f$.  This completes the proof of
(1), and (2) follows immediately.
\end{proof}

\subsection{Abelian varieties}
\begin{theorem}\label{abel.t}
 Let $A$ be an abelian variety over a perfect field $k$ of characteristic $p$.
\begin{enumerate}
\item{The Azumaya algebra $F_*D_{A/k}$ splits (non-canonically)
over the formal completion    $\hat \bT^*_{A'/k}$
 \footnote {This result is due to Roman Bezrukavnikov.}}
\item{There exists a tensor splitting of $F_* D_{A/k}$ over  $\hat \bT^*_{A'/k}$
 if and only if $A$ is ordinary. For an ordinary $A$,
 the tensor splitting is unique. }
\end{enumerate}
\end{theorem}
\begin{proof}
It is known \cite{mm.ue}  that Hodge and   conjugate spectral
sequences for $A$ degenerate
 and that $\Pic^\nat_A$ is smooth.  Thus
 Lemma~\ref{difpic.l} implies
that the differential of $\psi \colon \Pic^\nat_A \to \bH^0(A',
\Omega^1_{A'/k})$ is surjective, and it follows that $\psi$ is
smooth.  This implies that $\psi$ has a lifting over the formal
completion of $\bH^0(A',\Omega^1_{A/k})$ at the origin, and
therefore by Proposition~\ref{lifspl.p} that $\cD_{A/k}$ splits over
$A'\times \hat  \bH^0(A',\Omega^1_{A'/k}) \cong \hat \bT^*_{A'/k}$.

It follows from  Proposition~\ref{lifspl.p}  that giving a {\em
tensor} splitting of $F_*D_{A/k}$ over  $\hat A^{\prime}$ is
equivalent to giving a group homomorphism
 $$\tilde \psi \colon \hat\bv H^0(A',\Omega^1_{A'}) \to \Pic_A^\nat$$
such that $\psi \circ \tilde \psi = \id$.  The map $b\circ \tilde
\psi$ necessarily factors through $\Pic^0_A$, and since the latter
is $p$-divisible, $b\circ \tilde \psi = 0$.  Hence $\tilde \psi$
factors through $a$  in diagram (\ref{picn.e}) and can be viewed as
a morphism $\hat \bv H^0(A,\Omega^1_{A}) \to \bH(A,Z^1_{X/k})$.
 Since  $H^0(A,Z^1_{X/k}) = H^0(A,\Omega^1_{X/k})$,
the groups $\hat \bv H^0(A,Z^1_{A/k})$ and $\hat\bv
H^0(A',\Omega^1_{A'})$ are smooth of the same dimension.  Thus the
existence of $\tilde \psi$ is equivalent to the differential of
$\psi\circ a$ at $0$  being an isomorphism. It follows from
Lemma~\ref{difpic.l} that this restriction is the negative of the
Cartier operator
$$ C_{A/k} \colon H^0(A,\Omega^1_{A/k}) \to   H^0(A',\Omega^1_{A'/k}).$$
One of the equivalent definitions of an ordinary abelian variety is
that   $C_{A/k}$ is an isomorphism. This proves that lifting $h$
exists if and only if $A$ is ordinary. Moreover, for an ordinary $A$
the morphism $\psi:\hat\bv  H^0(A,\Omega^1_{A})\to  \hat\bv
H^0(A',\Omega^1_{A'})$ is an isomorphism. Thus, in this case, the
lifting is unique. We could also remark that an ordinary abelian
variety over a perfect field of characteristic $p$ has a canonical
lifting, together with a lifting of $F$, and this gives a tensor
splitting of $\hat \bT^*_{A'/k}$ by Theorem~\ref{loccart.t}.
\end{proof}

\subsection{A counterexample:  $\cD_{X/k}$ need not split on $\hat\bT_{X'/S}$}\label{counter.ss}
In this section, we will construct an example of a smooth proper
surface  $X/k$ over a perfect field $k$ which lifts to $W(k)$ but
such that $\cD_{X/k}$ does not split over the formal completion of
$T^*_{X'/k}$ along the zero section, or even over the formal
completion of $X'\times \bH^0(X',\Omega^1_{X'/k})$ along the zero
section.

\begin{lemma}\label{picred.l}
Let $X/k$ be a smooth and proper scheme with a rational point $x_0$.
Assume that  the following properties hold:
\begin{enumerate}
\item{$\dim H^0(X,\Omega^1_{X/k} ) = \dim H^1(X,\oh X)$,}
\item{$F_X$ acts as zero on $H^1(X,\oh X)$,}
\item{The Hodge spectral sequence of $X/k$ degenerates at $E_1$,}
\item{$q^*F_*\cD_{X/k}$ splits over the formal completion of $X'\times\bH^0(X',\Omega^1_{X'/k})$
along the zero section.}
\end{enumerate}
Then $\Pic_X$ is reduced.
\end{lemma}
\begin{proof}
 By assumption (3) the Hodge and conjugate spectral sequences of $X/k$ \marginpar{ref for this}
degenerate at $E_1$ and $E_2$ respectively, and so the row and
column of the  commutative diagram of Lemma~\ref{difpic.l} are short
exact.   The map $h := F_{X/k}^*$ in the diagram below vanishes by
assumption (2). This implies that the map $d\phi$ factors  as shown
below through $da$.  By (1) and (3) the induced map $h'$ is an
isomorphism, and it follows that $C_{X/k}$ is zero, and hence that
$d\psi$ factors through an arrow $h''$ as shown.

\begin{diagram}
&&& H^1(X',\oh {X'})\cr &&\ldDotsto^{h'} &\dTo^{d\phi} & \rdTo^{h}
\cr &H^0(X,Z^1_{X/k}) & \rTo^{da} & H^1_{dR}(X/k) & \rTo^{db} &
H^1(X,\oh X) \cr &&\rdTo_{C_{X/k}} &\dTo^{d\psi} & \ldDotsto{h''}
\cr &&&H^0(X',\Omega^1_{X'/k} )\cr
\end{diagram}

Now suppose that $F_{X/k*}D_{X/k}$ splits over $X'\times \hat
\bH^0(X',\Omega^1_{X'/k})$. Choosing a rigidification of the
splitting module, we get a lifting $\tilde \psi$ of $\psi$ over
$\hat \bH^0(X',\Omega^1_{X'/k})$, so $d\psi \circ d\tilde \psi =
\id$. Then $h'' \circ db \circ  d\tilde \psi = \id$, so $db\circ
d\tilde \psi \colon H^0(X',\Omega^1_{X'/k}) \to H^1(X,\oh X)$ is
injective. By (1), the source and target have the same dimension, so
the differential of the morphism $b\circ \tilde \psi$ is an
isomorphism.
 Since $\hat \bv H^0(X',\Omega^1_{X'/k})$ is smooth,
this implies that $\Pic_X$ is smooth.

\end{proof}

Let $k$ be a perfect field of odd characteristic and let $W$ be its
Witt ring. We construct an example of a smooth projective surface
$\tilde X/W$ whose special fiber $X$ over $k$ satisfies (1)--(3),
but whose Picard scheme is not reduced, using the technique of Serre
and its generalization by Raynaud \cite[4.2.3]{ray.psp}.  Let $E$ be
an elliptic curve over $W$ with supersingular reduction and denote
by $G$ the kernel of multiplication by $p$ in $E$. By [{\em op.
cit.}], there exists a projective complete intersection $\tilde Y$,
flat of relative dimension two over $W$, with a free action of $G$
and whose quotient $\tilde X := \tilde Y/G$ is smooth over $W$.  By
the weak Lefschetz theorem, $\Pic^0_{\tilde Y} = 0$, and it follows
that $\Pic^0_{\tilde X}$ is the Cartier dual of $G$, which can be
identified with $G$ itself. Since $\Pic$ commutes with base change,
the Picard scheme of the special fiber  $X$ is the special fiber
$G_0$ of $ G$.  In particular $G_0$ is not smooth. Replacing $k$ by
a finite extension, we may assume that $X$ has a rational point.
  Thus to produce our counterexample, it will
suffice to prove that $X$ satisfies (1)--(3) of
Lemma~\ref{picred.l}.
 The degeneration of
the Hodge spectral sequence of $X/k$ follows from its liftability.
The endomorphism of  $H^1(X,\oh X)$ induced by $F_X$ corresponds via
its identification with the tangent space of $G_0$ to the Cartier
dual of the endomorphism induced by the Frobenius of $G_0$ which in
our case vanishes. Since $\Pic_X \cong G_0$, $H^1(X,\oh X)$ is
one-dimensional, and so to prove (1), it will suffice to prove that
$H^1_{dR}(X/k)$ is two-dimensional. We use Faltings' comparison
theorem \cite[5.3]{fa.ccsc} which relates the de Rham cohomology
$H^1_{dR}(X/k)$ to the \'etale cohomology $H_\et^1(\tilde X_{\ov
K},\fp)$.  In particular, this theorem implies that these have the
same dimension. Since $\tilde X_{\ov K} = \tilde Y_{\ov K}/G_{\ov
K}$, $G_{\ov K} \cong \fp\oplus \fp$, and $\tilde Y_{\ov K}$ is
simply connected, $H_\et^1(\tilde X_{\ov K},\fp) \cong \fp\oplus
\fp$. Thus $H^1_{dR}(X/k)$ is two dimensional, and (2) follows.

\subsection{Fontaine modules}\label{fm.ss}

Throughout this section we assume that $S$ is a smooth scheme over a
field of characteristic $p$.
\begin{definition}\label{fontmod.d}
Let $X/S$ be a smooth scheme and let ${\cal X}/{\cal S}$ be a
lifting.  Fix  integers $k\leq l$ with $l-k<p$. Then a
\emph{Fontaine module} on $\lift X S$ consists of a coherent sheaf
with integrable  connection $(M,\nabla)$ of $\mic(X/S)$ and a
\emph{Hodge filtration}
$$0 = F^{l+1}M \subseteq F^{l}M \subseteq \cdots \subseteq  F^k M = M$$
 satisfying Griffiths
transversality, together with an isomorphism
\begin{equation}\label{frobenius}
\phi: C^{-1}_{\lift X S}\pi_{X/S}^*(Gr_F^{\cx} M,\kappa) \cong
(M,\nabla),
 \end{equation}
where the Higgs field $\kappa $ is given by {\it the Kodaira-Spencer
operator}
$$Gr \, \nabla: Gr^i_F \, M \to Gr ^{i-1}_F \, M \ot \Omega^1_{X/S} .$$
\end{definition}
We will denote the category of Fontaine modules by ${\cal M}{\cal
F}_{[k,l]}({\cal X}/{\cal S})$. Although we shall not do so here,
one can check that if $S$ is the spectrum of a perfect field $k$,
and the lifting $\tilde X' \to \tilde S = \spec \, W_2(k)$ comes
from a smooth formal scheme $X_{W(k)}$ over $W(k)$,   the category
${\cal M}{\cal F}_{[k,l]}({\cal X}/{\cal S})$ is equivalent to the
full subcategory of $p$-torsion objects in Faltings' category ${\cal
M}{\cal F}_{k,l}^{\nabla}(X_{W(k)})$ \cite{fa.ccpgr}.

The formula $N^m (Gr_F^{\cx} M,\kappa) = \oplus _{i\leq -m} Gr^i_F M
\subset Gr_F^{\cx} M $ defines a {\it  $\J_X$-filtration} on
$(Gr_F^{\cx} M,\kappa)$. Applying the isomorphism $\phi$ we obtain
an {\it  $\Cee_X$-filtration } on $M$:
$$N^m M= C^{-1}_{\lift X S}\pi_{X/S}^*(\oplus _{i\leq -m} Gr^i_F M) \subset M.$$
\begin{theorem}Let $ (M, \nabla, F^{\cx}M,
\phi)$  and $(M', \nabla ', F^\cx, \phi')$ be Fontaine modules over
$X$. Then
\begin{enumerate}
\item  For every integer $i$, the
${\cal O}_X$-module $Gr _F^i \, M $ is locally free. In particular,
$M$ is a locally free ${\cal O}_X$-module
 \cite[Theorem 2.1]{fa.ccpgr}.
\item  Every morphism $f:M\to M'$ of
Fontaine modules is {\it strictly} compatible with the Hodge
filtration $F^{\cx}$. In particular, the category ${\cal M}{\cal
F}_{[k,l]}({\cal X}/{\cal S})$ is {\it abelian} [\emph{op. cit}].
\item  Let $h: X\to Y$ be a smooth proper
morphism of relative dimension $d$,let  $\tilde h': \tilde X' \to
\tilde Y'$ be a lifting of $h'$, and let $ (M, \nabla, F^{\cx},
\phi)\in {\cal M}{\cal F}_{[k,l]}({\cal X}/{\cal S})$ be a Fontaine
module. Assume that $l-k+d<p$. Then, the Hodge  spectral sequence
for $Rh^{DR}_* (M, F^{\cx}M)$ degenerates at $E_1$. Thus, by Theorem
\ref{gmeqpush.t} b), for every integer $i$, we have a canonical
isomorphism
\begin{equation}\label{pushf}
\phi: C^{-1}_{\lift Y S}\pi_{Y/S}^*(Gr_F^{\cx} R^i h^{DR}_*
M,\kappa) \cong (R^i h^{DR}_* M,\nabla),
\end{equation}
 which makes $(R^i h^{DR}_*M, \nabla, F^{\cx}R^i h^{DR}_* M, \phi)$ a Fontaine module
over $Y$. In particular, if $d<p$, the $D_{Y/S}$-module $R^i
h^{DR}_*{\cal O}_X$ is a Fontaine module~\cite{fa.ccpgr}.
\item The Chern classes $c_i(M)\in H_{et}^{2i}(X, \mathbb{Q}_l(i))$,
$l\ne p,\,  i>0$ are all equal to $0$.
\end{enumerate}
\end{theorem}
\begin{proof}
The key to parts a) and b) is the following general result, whose
proof can be found in \cite[8.2.3]{o.fgthd}.

\begin{lemma}\label{ogus} Let $Z$ be a smooth scheme over a field  of
characteristic $p$ and let
$$0 = F^{l+1}M^\cx \subseteq F^{l}M^\cx \subseteq \cdots \subseteq  F^k M^\cx = M^\cx $$
be a bounded filtered complex of coherent ${\cal O}_Z$-modules.
Assume that there exists a  (not necessarily filtered)
quasi-isomorphism
$$F_Z^* Gr_F M^\cx \simeq M^\cx. $$
Then the differential $M^\cx \to M^{\cx +1}$ is strictly compatible
with the filtration and, for every pair of integers $i$ and $j$, the
${\cal O}_Z$-module $H^j(Gr _F ^i M^\cx ) \simeq Gr_F ^i H^j(M^\cx)$
is locally free.
\end{lemma}
Let us return to the proof of the theorem. Since the claims in parts
(1)  and (2) are local on $X$ we may assume that there exists a
lifting $\tilde F$ of the Frobenius $F_{X/S}$. By
Theorem~\ref{loccart.t}, such a lifting induces a natural
isomorphism of ${\cal O}_{X}$-modules
\begin{equation}\label{lg}
\eta_{\tilde F}:  C^{-1}_{\lift X S} (E)\simeq F^*_{X/S} E,
\end{equation}
 for every $E\in HIG _{p-1} (X'/S)$. Composing this with
(\ref{frobenius}) we obtain an isomorphism of ${\cal O}_{X}$-modules
$$F^*_X Gr_F M \simeq M . $$
Then the statements (1)  and (2) follow from the lemma.

By Theorem \ref{gmeqpush.t} the lifting $\tilde h' $ induces an
quasi-isomorphism \begin{equation}\label{dirfon}
 C^{-1}_{\lift Y
S}\pi_{Y/S}^*Rh^{HIG}_* (Gr_F^{\cx}M, \kappa) \cong R h^{DR}_* M.
\end{equation}
Applying (\ref{lg}) we obtain locally on $Y$ an isomorphism in the
derived category of ${\cal O}_Y$-modules
$$ F^*_Y Rh^{HIG}_* (Gr_F^{\cx}M, \kappa)\cong R h^{DR}_* M.$$
We can compute $Rh^{HIG}_* (Gr_F^{\cx}M, \kappa)$ as follows. Endow
the relative de Rham complex $\Omega ^\cx _{X/Y} \otimes M $ with
the filtration
$$F^i (\Omega ^\cx _{X/Y} \otimes M) = (F^i M \to \Omega ^1 _{X/Y}\otimes F^{i-1} M
\to \cdots \to  \Omega ^d _{X/Y}\otimes F^{i-d} M), $$ and let
$(Rh^{DR}_* M, F^\cx)$ be the filtered derived direct image of
$(\Omega ^\cx _{X/Y} \otimes M, F^\cx)$. We then have an isomorphism
in the derived category of ${\cal O}_Y$-modules
$$ Gr_F Rh^{DR}_* M \simeq Rh^{HIG}_* (Gr_F^{\cx}M, \kappa).$$
Thus by Lemma \ref{ogus},  applied to the filtered complex of
coherent $\oh Y$-modules $(Rh^{DR}_*M, F^\cx)$, the Hodge spectral
sequence for $Rh^{DR}_* (M, F^\cx)$,
%\ie,  the spectral sequence of filtered complex
$(Rh^{DR}_* M, F^\cx)$, degenerates at $E_1$. Hence we get a
canonical isomorphism of ${\cal O}_Y$-modules
\begin{equation}\label{degh}
Gr_F R^i h^{DR}_* M \simeq R^i h^{HIG}_* (Gr_F^{\cx}M, \kappa)
\end{equation}
It is well known\footnote{This fact should be compared with Katz's
formula (\ref{ckf}). A conceptual proof of this result can be
obtained using an appropriate filtered derived category of
D-modules. See, for example~\cite{saito.hsd}.} that this isomorphism
is compatible with the Higgs fields. Thus passing to the cohomology
sheaves in (\ref{dirfon}) we obtain the desired isomorphism
(\ref{pushf}).
 This completes the proof of statement (3).

For  statement (4), we will first prove that for any ${\cal
O}_{X}$-coherent $N\in HIG_{p-1}(X/S)$,
$$[ C^{-1}_{{\cal X}/S}\pi_{X/S}^* N]= F^*_X [N],$$
where $[\, ]$ denotes the class of a coherent ${\cal O}_X$-module in
$K'_0(X)=K_0(X)$. Indeed, choose any filtration $N=N^0\supset N^1
\supset \cdots \supset N^m=0$ by Higgs submodules such that
$N^i/N^{i+1}\in HIG_{0}(X/S)$. Then
$$ (C^{-1}_{{\cal X}/S}\pi_{X/S}^* N^i)/ (C^{-1}_{{\cal X}/S}\pi_{X/S}^* N^{i+1})\simeq
 C^{-1}_{{\cal X}/S}\pi_{X/S}^* (N^i/N^{i+1}) \simeq$$
$$F^*_X (N^i/N^{i+1})\simeq F^*_X N^i/F^*_X N^{i+1}.$$
 This implies the claim.

In particular, for a Fontaine module $(M, \nabla, F^{\cx}M, \phi)$
it follows that $[M]= F^*_X [M]$. Thus
$$c_i([M])= c_i(F^*_X [M])= p^i c_i([M]), $$
and we are done.\\
\end{proof}

\begin{proposition}\label{sshigg.p}
Let $X$ be a smooth projective curve of genus $g$ over a field $spec
\, k = S$, $\tilde X' \to \tilde S$ a lifting, and let
$(M,\nabla,F^\cx,\phi)$ belong to the category
$\mathcal{MF}_{[0,n]}$. Assume that
$$n\, (rk\, M  -1) \, max\,  \{2g-2, 1\} <p-1.$$
 Then $(Gr_F^{\cx}M, \kappa)$ is a semistable Higgs
bundle.
\end{proposition}
\begin{proof}  We have to show that $(Gr_F^{\cx}M, \kappa)$ has no Higgs subbundles
$$(L, \theta) \hookrightarrow (Gr_F^{\cx}M, \kappa)$$
of positive degree. Replacing $(L, \theta)$ by $ \wedge ^{rk \,
L}(L, \theta)$ and $ M$  by $\wedge^{rk \, L} M$ (this is again a
Fontaine module) we reduce Proposition to the following claim:

For any Fontaine module $(M, \nabla, 0=F^{n+1}M \subset F^{n}M
\subset \cdots \subset F^0 M=M, \phi)$, with $n \, (2g-2) < p-1$ the
Higgs bundle $(Gr_F^{\cx}M, \kappa)$ does not have one-dimensional
Higgs subbundles
\begin{equation}\label{subbun}
(L, 0) \hookrightarrow (Gr_F^{\cx}M, \kappa)
\end{equation}
of positive degree.

Assume that this is not the case and consider such  an $L$ of the
largest possible degree $d>0$. Then any morphism $(L', 0) \to
(Gr_F^{\cx}M, \kappa)$, where $L'$ is a line bundle of degree $>d$,
is equal to zero. Consider the morphism
$$F_X^* L \simeq C^{-1}_{{\cal X}/S}\pi_{X/S}^* (L, 0) \hookrightarrow C^{-1}_{{\cal X}/S}\pi_{X/S}^* (Gr_F^{\cx}M, \kappa)
\stackrel{\phi}{\simeq} M$$ induced by (\ref{subbun}). We will prove
by induction on $m$ that the composition
\begin{equation}\label{prss}
F_X^* L \to M \to M/F^m M
\end{equation}
is $0$. Let us, first, check this for $m=1$. Observe that the Higgs
field $\kappa$ restricted to $M/F^1 M \hookrightarrow   Gr_F^{\cx}M$
is $0$. Thus
$$(F_X^* L, 0) \to  (M/F^1 M, 0) \hookrightarrow   (Gr_F^{\cx}M, \kappa )$$
is a morphism of Higgs bundles. Since $deg\, F^*_X L = pd>d$, this
morphism must be equal to zero.

Assume that the composition $F_X^* L \to M \to M/F^{m-1} M$ is $0$.
Then (\ref{prss}) factors through $F_X^* L \to F^{m-1}M/F^m M$. For
any $j$, $0\leq j < m$, consider the composition
$$\rho _j: F_X^* L \to F^{m-1}M/F^m M \stackrel{\kappa ^j}{\longrightarrow} F^{m-1-j}M/F^{m-j} M \otimes
( \Omega ^1_{X/S})^j,$$ and let $j_0 $ be the smallest integer less
then $m$, such that $\rho_{j_0}\ne 0$. Then $\rho_{j_0}$ induces a
nonzero map of Higgs bundles
$$ (F_X^* L \otimes ( T_{X/S})^{j_0}, 0)  \to (\bigoplus _{i\geq j_0} F^{m-1-i}M/F^{m -i}M, \kappa)
 \hookrightarrow  (Gr_F^{\cx}M, \kappa).$$
However
$$deg \,( F_X^* L \otimes ( T_{X/S})^{j_0})= pd - j_0(2g-2)\geq pd - n(2g-2)> d.$$
 This contradiction completes
the proof.
\end{proof}

\begin{remark}
 Let  $h: Y\to X$ be a smooth proper morphism of relative dimension $d$,
  and let $\tilde h': \tilde Y' \to \tilde X'$ be a lifting.
Then, for $d<p$,  $M= R^n h^{DR}_* {\cal O}_{Y}$ is a Fontaine
module on $X$. Thus, by Proposition~\ref{sshigg.p} if $n\, (rk\,M
-1) \, max\, \{2g-2, 1\} <p-1 $,  $ (Gr_F^{\cx}M, \kappa)$ is
semistable. By the standard technique this implies the following
result over the complex numbers.

\begin{theorem}\label{cxss.t}
Let $X$ be a smooth projective curve over $\bc$  and let  $h: Y \to
X $ be a smooth proper morphism. Then $ (Gr_F^{\cx}R^n h^{ DR}_*
{\cal O}_{Y} , \kappa)$ is a semistable Higgs bundle.
\end{theorem}

This result was proved by analytic methods (for any polarizable
variation of Hodge structure) by Beilinson and Deligne (unpublished)
and later, in a greater generality, by Simpson~\cite{si.hbls}  using
a similar technique.
\end{remark}

\subsection{Proof of a theorem of Barannikov and Kontsevich}\label{bk.ss}

Let us recall the following striking result of Barannikov and
Kontsevich, of which the only published proof we know is due to
Sabbah~\cite{sa.tdc}.

\begin{theorem}\label{bk.t}
Let $X/\bc$ be a quasi-projective smooth scheme over $\bc$. Suppose
that $f \in \Gamma(X,\oh X)$ defines a proper morphism to
${\bA^1/\bc}$.  Then the hypercohomologies of the complexes
\begin{diagram}
\oh X &\rTo^{d +\wedge df} & \Omega^1_{X/\bc} & \rTo^{d + \wedge df}
& \Omega^2_{X/\bc} & \cdots &\quad \mbox{and} \cr \oh X
&\rTo^{-\wedge df}    & \Omega^1_{X/\bc} & \rTo^{-\wedge df}    &
\Omega^2_{X/\bc} & \cdots& \cr
\end{diagram}
 have the same finite dimension in every degree.
\end{theorem}

We shall show how our version of nonabelian Hodge theory can be used
to give a proof of this theorem by the technique of reduction modulo
$p$. Since any pair $(X/\bc,f)$  as in Theorem~\ref{bk.t} comes from
some ``thickened'' situation, it is clear that the following result
implies Theorem~\ref{bk.t} by base change $R \to \bc$.

\begin{theorem}\label{bko.t}
Let $\cS = \spec R$ be an affine, integral, and smooth scheme over
$\bz$, let $\cX/\cS$ be a smooth quasi-projective  $S$-scheme, and
let $\tf $ be a global section of $\oh \cX$ which defines a proper
morphism: $\cX \to {\bf A}^1_\cS$.
   Then, after replacing $\cS$ by some \'etale neighborhood of its
generic point,  the following results are true.
\begin{enumerate}
\item{The hypercohomology groups
$$H^*(\cX,\Omega^\cx_{\cX/\cS},d+d\tf) \quad\mbox{and} \quad
H^*(\cX,\Omega^\cx_{\cX/\cS},-\wedge d\tf)$$
 are finitely generated  free $R$-modules whose formation commutes with  base change.}
\item{ Let $p$ be a prime, let
$X/S$ denote the reduction of $\cX/\cS$ modulo $p$, and  let $ X
\rTo^{F_{X/S}} X' \rTo^\pi X$ be the usual factorization of $F_X$.
Then for every $p$, the complexes of $\oh {X'}$-modules
\begin{diagram}
F_{X/S*}\oh {X} &\rTo^{d +\wedge df} & F_{X/S*}\Omega^1_{X/S} &
\rTo^{d + \wedge df} & F_{X/S*}\Omega^2_{X/S} & \cdots \cr \oh{ X'}
&\rTo^{-\wedge d\pi^*f} & \Omega^1_{X'/S} & \rTo^{-\wedge d\pi^*f} &
\Omega^2_{X'/S} & \cdots \cr
\end{diagram}
are quasi-isomorphic.}
\end{enumerate}
\end{theorem}

The rest of this section will be devoted to a proof of
Theorem~\ref{bko.t}.
%{\bf Note that the dimensions (even the isomorphism class) of the cohomology
%groups of the Higgs complex are unchanged if $d\pi^*f$ is replaced
%by $-d\pi^*f$.}
Along the way we shall prove some auxiliary results which may be of
independent interest, for example the finiteness criterion given in
Proposition~\ref{noncr.p} and Corollary~\ref{finzer.c}. We begin
with a ``cleaning'' lemma.

\begin{lemma}\label{prepz.l}
With the notation of  Theorem~\ref{bko.t}, let $\cZ\subseteq \cX$
be the  reduced zero locus of $df$. Then after replacing $\cS$ by
some \'etale neighborhood of its generic point, the following
conditions are satisfied.
\begin{enumerate}
\item{The morphism $\cZ \to \cS$ is proper, flat, and generically smooth, and for every $p$
the reduction modulo $p$ of $\cZ$ is reduced.}
\item{ The
restriction of $f$ to  each connected component $\cZ'$ of $\cZ$ lies
in the image of the map $\Gamma(\cS,\oh \cS) \to \Gamma (\cZ',\oh
{\cZ'})$.}
\end{enumerate}
\end{lemma}
\begin{proof}
Note that formation of $\cZ$ commutes with \'etale base change $\cS'
\to \cS$, so that our statement is not ambiguous.
  Let $\sigma$ be the generic point of $\cS$.
The statements are trivial if the generic fiber of $\cZ_\sigma$ of
$\cZ/\cS$
 is empty, so let us
assume that this is not the case. By the theorem of generic
flatness~\cite[6.9.1]{EGAIV2}, we may assume that $\cZ$ is flat over
$\cS$. Then the map from each irreducible component $\cZ_i$ of $\cZ$
to $\cS$ is dominant and the generic fiber of $\cZ_i$ is an
irreducible component of $\cZ_\sigma$. Localizing further if
necessary, we may assume that if $\cZ_i$ and $\cZ_j$ intersect, then
so do their generic fibers. There is a finite extension $k'$ of
$k(\sigma)$ such that all the connected components of $\cZ_{k'}$ are
geometrically connected and have a $k'$-rational point.  Replacing
$\cS$ by an \'etale neighborhood of $\sigma$, we may  assume that
$k' = k(\sigma)$.  Since $\cZ_\sigma$ is  reduced and $k(\sigma)$ is
a field of characteristic zero, $\cZ_\sigma/\sigma$ is generically
smooth. Since the differential of $f_{|_{\cZ_\sigma}}$ vanishes, its
restriction to the smooth locus $\cZ_\sigma^{sm}$ of $\cZ_\sigma$ is
locally constant.  Thus for each irreducible component $\cZ_i$ of
$\cZ$, there exists an element $c_i$ in $k(\sigma)$  (the value of
$f$ at a rational point) such that $f = c_i$ on
$\cZ_{i\sigma}^{sm}$.  Since $\cZ_i$ is reduced, this holds on all
of $\cZ_i$. If  $\cZ_i$ and $\cZ_j$ intersect, so do $\cZ_{i\sigma}$
and $\cZ_{j\sigma}$, and
 it follows that $c_i = c_j$.  Thus $c_i$ depends only
on the connected component of $\cZ_\sigma$ containing $\cZ_i$.
Furthermore, localizing on $\tS$, we may assume that each $c_i$
belongs to $R$. Thus (2) is proved. Now if $\cZ'$ is a connected
component of $\cZ$,  the composite $\cZ' \to \cX \to  \bA^1_\cS$
factors through the  section of $\bA^1_{\cS}/\cS$ defined by the
appropriate element of $R$, Since $X \to \bA^1_\cS$ is proper, so is
each  $\bz' \to \cS$ and hence the same is true of $\cZ \to \cS$.

We have now attained all the desired properties of $\cZ$, except for
the reducedness of its reductions modulo $p$, which is a consequence
of the following (probably standard) lemma.

\begin{lemma}\label{reduced.l}
Let $\cZ$ be a reduced scheme of finite type over $\spec \bz$.  Then
for almost all primes $p$, the reduction modulo $p$ of $\cZ$ is
reduced.
\end{lemma}
\begin{proof}
In the course of the proof, we may without loss of generality
replace $\cZ$ by the open subset defined localization by any
positive integer.  In particular, by the theorem of generic
flatness,  we may assume that $\cZ$ is flat over $\bz$. Since
$\cZ_\bQ/\bQ$ is reduced and of finite type,   it is generically
smooth over $\bz$.
 Let $ \eta \colon \cY \to \cZ$ be the normalization
mapping. Then $\cZ$ is also generically smooth over $\bz$. Thus
each irreducible component $\cY^0$ of $\cY$ contains a proper closed
subscheme $\cY^1$ such that $\cY^0 \setminus \cY^1 \to \cS$ is
smooth.  For almost all $p$, the reduction modulo $p$ of  $\cY^1$
has strictly smaller dimension than that of the reduction modulo $p$
of $\cY^0$, and we may assume this is true for all $p$.  Then the
map $\cY \to \spec \bz$ remains generically smooth modulo $p$ for
every $p$. By \cite[7.7.4]{EGAIV2}, $\eta$ is finite, and hence the
cokernel $Q$ of $\eta^\sharp \colon \oh \cZ \to \eta_* \oh \cY$ is a
coherent sheaf of $\oh \cZ$-modules.  Again by the lemma of generic
flatness, $Tor_1^\bz(Q,\fp) = 0$ for all but finitely many $p$.
Shrinking, we may assume that this is true for all $p$.  It then
follows that the reduction modulo $p$ of $\eta^\sharp$ remains
injective for all $p$.  Since $\cY$ is normal, it satisfies Serre's
condition $S_2$, and since each $p$ defines a nonzero divisor on
$\cY$ the fiber $Y$ of $\cY$ over $p$ satisfies $S_1$. Since $Y$ is
generically smooth over ${\bf F}_p$, it is generically reduced, and
since it satisfies $S_1$, it is reduced. Since $\eta^\sharp$ is
injective mod $p$, the fiber $Z$ of $\cZ$ over $p$ is also reduced.
\end{proof}

%
%Since  $f_{|_\cZ} \colon \cZ \to \cS$ is dominant, the map
%$R:= \Gamma(\cS, \oh \cS) \to \Gamma(\cZ,\oh \cZ)$ is injective.
%Let $g$ be the image of $f$ in $\Gamma(\cZ,\oh \cZ)$, let
%$R'$ be the subring of $\Gamma(\cZ,\oh \cZ)$  generated by $R$ and $g$,
%and let $\cS' := \spec R'$.  Then the natural map $h \colon \cZ \to \cS'$
%is dominant and its differential $h^*\Omega^1_{\cS'/\cS}  \to \Omega^1_{\cZ/\cS}$
%is generically injective.  Since $df$ maps to zero in $\Omega^1_{\cZ/\cS}$,
%the same is true of $dg$, and hence the image of $dg$
%in $h^*\Omega^1_{\cS'/\cS}$ also vanishes at the generic point of $\cZ$.
%Hence $dg$ is generically zero in $\Omega^1_{\cS'/\cS}$.   But $dg$
%generates $\Omega^1_{\cS'/\cS}$, so   $\Omega^1_{\cS'/\cS}$ is generically zero and
%$\cS'/\cS$  is generically unramified. Since it is generically flat, it is
%generically \'etale.  Thus we may replace $\cS$ by $\cS'$, so that (2) holds.
%Then the map $\cZ \to \cX \to  \bA^1_\cS$ factors through the  section
%of $\bA^1_{\cS}/\cS$ defined by $g$, and it follows that $\cZ/\cS$ is also proper.
\end{proof}

Let $\cE := (\oh  \cX, d +d\tf) \in \mic(\cX/\cS)$ and let $\cL :=
(\oh \cX,d\tf) \in \hig(\cX/\cS)$; we denote by just $E$ and $L$
their respective reductions modulo a prime $p$ of $\bz$. Let
$J\subseteq \oh \cX$  be the ideal of the scheme-theoretic zero
locus of $d\tf$. This is just the ideal generated locally by the set
of partial derivatives of $\tf$ in any set of local coordinates for
$\cX/\cS$. The Higgs complex  $\cL \otimes \Omega^\cx_{\cX/\cS}$ of
$\cL$ can be locally identified with the Koszul complex of this
sequence of partials, and it follows that the  cohomology sheaves of
$\cL \otimes \Omega^\cx_{\cX/\cS}$ are annihilated by
$J$~\cite[17.14]{eis.ca}.
  Since the closed subscheme of $\cX$ defined by the radical of $J$ is $\cZ$,
which is proper over $\cS$,   the hypercohomology groups $H^i(\cL
\otimes \Omega^\cx_{\cX/\cS})$ are finitely generated $R$-modules.
Since $R$ is reduced, they are free in some neighborhood of the
generic point of $\cS$, which we may assume is all of $\cS$. Since
the terms in the complex  $\cL \otimes \Omega^\cx_{\cX/\cS}$ are
flat over $\cS$, the formation of its hypercohomology will then
commute with all base change. This completes the proof of
Theorem~\ref{bko.t}.1 for the Higgs complex.

 The proof
for the de Rham complex is more difficult; in general, the de Rham
cohomology  groups of a coherent sheaf  with integrable connection
on a smooth scheme of finite type over $\bz$ are not finitely
generated. (For example, the de Rham cohomology of the trivial
connection on $\bA^1_\bz$ is not finitely generated.)  We will use
the technique of logarithmic geometry to study the irregularity of
the connection $d +df$ to obtain the finiteness we need.

Let $ Y /S$ be a smooth morphism of fine saturated and noetherian
log schemes. We just write $\Omega^\cx_{Y/S}$ for the logarithmic de
Rham complex of $Y/S$~\cite{kato.lsfi}. If  $m $ is a section of
$M_Y$,  the set $ Y_m$ of all $y \in Y$ such that $m_y \in
M_{Y,y}^*$ is open in $Y$.  In fact, since $ \alpha \colon M_Y \to
\oh Y$ is a log structure, $y \in Y_m$ if and only if $\alpha_Y(m)
\in \oh{Y,y}^*$. Let us assume that $\alpha_Y(m)$
 is a nonzero divisor of $\oh Y$, so that it defines
a Cartier divisor $D$ of $Y$ and $Y_m = Y\setminus D$.
 Suppose we are given a torsion free coherent sheaf
$E$ on $Y$ and an integrable connection $\nabla$ on $j^*E$, where $j
\colon Y_m \to Y$ is the inclusion. Then $\nabla$ induces a
connection on $j_*j^*E \cong  E(*) :=\dirlim E(nD)$. If $\nabla$
maps $E$ to $E \ot \Omega^1_{Y/S}$, then $E$ has regular singular
points along $D$; we wish to measure the extent to which this fails.
 Since $E$ is coherent, $\nabla$
maps $E$ to $E\ot \Omega^1_{Y/S}(nD)$ for some $n$; replacing $m$ by
$m^n$ we may assume that $n= 1$. Since $da \in I_D\Omega^1_{Y/S}(D)$
for all $a \in \oh Y$, the map
$$\theta_D \colon  E\ot \oh D \to E \ot\Omega^1_{Y/S}(D)_{|_D}$$
induced by $\nabla$ is  $\oh D$-linear.
%We call $\theta_D$ the \emph{$D$-irregularity of $(E,\nabla$)}.
It follows from the integrability of $\nabla$ that $\theta_D$
defines an action of  the symmetric algebra $S^\cx(I_D T_{Y/S})$ on
$E_{|_D}$, so that $E_{|_D}$ can be viewed as a module over $\bv
(I_DT_{Y/S})$.

The following result is  inspired by  of a result of
Deligne~\cite[II, 6.20]{de.edprs} which was pointed out to us by H.
Esnault.

\begin{proposition}\label{noncr.p}
Suppose that in the above situation  $\theta_D$ is noncritical, \ie,
that the support of the $\bv (I_D T_{Y/S})$-module $E_D$ defined by
$(E_D,\theta_D)$ is disjoint from the zero section. Suppose further
that  $Y /S$ is proper and that $S = \spec R$.  Then for every $i$,
$H^i(Y\setminus D,E\ot \Omega^\cx_{Y/S})$ is a finitely generated
$R$-module.
\end{proposition}
\begin{proof}
Let $\Omega_{Y/S}^q(*) := j_*j^*\Omega^q_{Y/S}$ and for each natural
number $n$, let
$$F_n (E\otimes \Omega_{Y/S}^q)(*) :=
E\ot \Omega^q_{Y/S}((n+q)D) \subseteq  j_*j^*(E\otimes
\Omega^q_{Y/S}).$$ Then $F_\cx$ defines an exhaustive filtration of
the complex $E\ot \Omega^\cx_{Y/S}(*)$ by coherent sheaves.  Since
$Y/S$ is proper, for each $n$ and $i$, $H^i(F_n E\ot
\Omega^\cx_{Y/S}(*))$ is finitely generated over $R$.  Thus it will
suffice to show that for  each $n \ge 0$, the natural map
$$F_n E\ot \Omega^\cx_{Y/S}(*) \to E\ot \Omega^\cx_{Y/S}(*)$$
is a quasi-isomorphism, and for this it will suffice to prove that
for each $n \ge 0$, the map
$$F_0 E\ot \Omega^\cx_{Y/S}(*) \to F_n E\ot \Omega^\cx_{Y/S}(*)$$
is a quasi-isomorphism. This will follow by induction if for every
$n > 0$, $\gr_n^F  E\ot \Omega^\cx_{Y/S}(*)$ is acyclic.

Multiplication by $g^n$ defines an isomorphism $F_n E(*) \to F_0
E(*) $ which induces an isomorphism
$$ E(nD)_{|_D} = \gr_n^F E(*) \to
\gr_0^FE(*) = E_{|_D}.$$
  If $e \in F_n E,$  then
$\nabla(e) \in F_nE\otimes \Omega^1_{Y/S}(D)$ and $\nabla(g^n e) \in
E\otimes \Omega^1_{Y/S}(D)$. Since $g = \alpha(m)$ and $dg = g \dlog
m$,
$$\nabla(g^n e) = ng^ne \ot  \dlog(m) + g^n \nabla( e) \in E\otimes \Omega^1_{Y/S}(D).$$
Since $g^ne\ot \dlog(m) \in E\ot \Omega^1_{Y/S}$,  $\nabla(g^n e)$
reduces to  $g^n\nabla(e)$ in $E \ot \Omega^1_{Y/S}(D)_{|_D}$.
 Thus multiplication
by $g^n$ identifies $\gr_n^F \nabla$ with  $\theta_D$ for all $n \ge
0$.  This identification extends to an isomorphism of complexes
$$ \gr_n^F ( E\otimes \Omega^\cx_{Y/S}(*)) \cong \gr_0^F ( E\otimes \Omega^\cx_{Y/S}(*))$$
  But $\gr^0_F (E\ot \Omega^\cx_{Y/S}(*))$ is just the
Higgs (Koszul) complex of $\theta_D$, whose cohomology sheaves can
be identified with $Ext^*_{\oh {\bv}}(i_*\oh X,E_D)$, where $i
\colon X \to \bv$ is the zero section of $\bv := \bv (I_D T_{Y/S})$.
These vanish since $\theta_D$ is noncritical.
\end{proof}

The following corollary then completes the proof of statement (1) of
Theorem~\ref{bko.t}:
 take  $(E,\nabla)$ to be the constant connection on $\cX/\cS$.

\begin{corollary}\label{finzer.c}
Let $X/S$ be a smooth quasi-projective scheme over $S = \spec R$,
where $R$ is a flat and finitely generated $\bz$-algebra. Let
$(E,\nabla)$ be a coherent sheaf with integrable connection on $X/S$
whose restriction to the generic fiber of $X/S$ has regular
singularities at infinity. Suppose that $f \in \oh X(X)$ is a global
function which defines a proper morphism $X \to {\bf A}^1_S$, and
let $(E',\nabla')$ be the $df$-twist of $(E, \nabla)$: $E' = E$, and
$\nabla' := \nabla + \wedge df$. Then after replacing $S$ by some
affine neighborhood of the generic point of $S$, the de Rham
cohomology $H^*(X,E'\ot \Omega^\cx_{X/S})$ is finitely generated and
free over $R$.
\end{corollary}
\begin{proof}
Let $\sigma$ be the generic point of $S$. We may find a projective
compactification $\ov X_\sigma$ of $X_\sigma$, and after blowing up
$\ov X_\sigma$ outside of $X_\sigma$ we may assume that $f$ extends
to a morphism $\ov X_\sigma \to {\bf P}^1_\sigma$, which we still
denote by $f$.  After a further blowing up outside of $X_\sigma$,
we may assume that $\ov X_\sigma$ is smooth over $\sigma$ and that
the complement of $X_\sigma$ in $\ov X_\sigma$ is a divisor with
strict normal crossings. Then the log scheme  $Y_\sigma$ obtained by
endowing $\ov X_\sigma$ with the log structure corresponding to the
inclusion $X_\sigma \to \ov X_\sigma$ is (log) smooth.  Furthermore,
$f$ extends to a morphism of log schemes $Y_\sigma \to P^1_\sigma$,
where  $P^1_\sigma$ is the
 log scheme $P^1_\sigma$ obtained by endowing ${\bf P}^1_\sigma$ with
the log structure corresponding to the
 inclusion ${\bf A}^1_\sigma \to {\bf P}^1_\sigma$,

Let $t$ be the coordinate of ${\bf A}^1_\sigma$ and let $s :=
t^{-1}$, which is a local generator of the ideal of $\infty$.  There
is a unique local section  $m$ of the sheaf of monoids
$M_{P^1_\sigma}$ over $V$ with  $s := \alpha_{P^1_\sigma}(m)$, and
$\dlog m$ is  basis for the stalk of  $\Omega^1_{P^1_\sigma}$ at
$\infty$. Let $y$ be a point of $D :=f^{-1}(\infty)$.  Then in an
\'etale neighborhood of $y$, there exists a system of coordinates
$(t_1,\cdots t_n)$ and  natural numbers $r, e_1,\ldots e_r$ such
that such that $f^*(s) = t_1^{e_1}\cdots t_r^{e_r}$. Then $f^*(dm) =
\sum_i e_i \dlog t_i$, which is nonvanishing in the fiber of
$\Omega^1_{Y_\sigma/\sigma}$ at $y$. (This implies that $f$ is log
smooth at $y$.) Since $(E,\nabla)$ has regular singularities at
infinity, there is   a coherent (even locally free) extension $\ov
E$ of $E$ to $Y_\sigma$ and a log connection
 $\nabla \colon \ov E \to \ov E \ot \Omega^1_{Y_\sigma/\sigma}$
extending $\nabla$. Now  $df = f^*(dt) = -s^{-2}ds = -
s^{-1}f^*\dlog m$.  Thus  $\nabla'$ maps $\ov E$ to
 $\ov E\ot \Omega^1_{Y_\sigma/\sigma}(D)$,
and $\theta_D$ is the map $\ov E_{|_D} \to \ov E_{|_D}\ot
\Omega^1_{Y_\sigma/\sigma}$ sending $e$ to $-e\wedge s^{-1}\dlog m$.
This is an isomorphism, so $\theta_D$ is noncritical. There exists
an affine neighborhood of the generic point of $S$ over which all
this remains true, and without loss of generality we may assume they
are true for $Y/S$. Then Proposition~\ref{noncr.p} implies that the
de Rham cohomology groups of $(E,\nabla)$ over $Y\setminus D = X$
are finitely generated over $R$; shrinking further we may assume
they are free.
\end{proof}

We now turn to the proof of statement (2) of  Theorem~\ref{bko.t}.
Assume that $\cX/\cS$ satisfies the conditions in (1) of
Theorem~\ref{bko.t} and  in Lemma~\ref{prepz.l}.  Fix a prime $p$,
let $X/S$ be the reduction of $\cX/\cS$ modulo $p$,  and let $\cS_1$
be the reduction of $\cS$ modulo $p^2$.
  Since $\cS/\bz$ is smooth and affine,
there exists a lifting $F_\cS$ of the absolute Frobenius
endomorphism of $S$ to  $\cS_1$  and hence a Cartesian square:
\begin{equation}\label{xs1.e}
\begin{diagram}
\cX'_1 &\rTo^{\tilde \pi}& \cX_1 \cr \dTo  & &\dTo \cr \cS_1 &
\rTo^{F_{\cS}} &\cS_1 .
\end{diagram}
  \end{equation}
We shall abuse notation and write $C_{\lift X S}$ for the Cartier
transform defined by the lifting $\cX'_1/\cS_1$ of $X'/S$.

Let $(E,\nabla)$ be the restriction of $(\cE,\nabla)$ to $X/S$.
According to \cite[]{ka.asde}, the $p$-curvature \marginpar{find
precise reference here and also in section 4.2}
 $\psi  \colon E \to E \otimes F_{X/S}^*\Omega^1_{X'/S}$
is multiplication by
$$F_X^*(df) - F^*_{X/S}C_{X/S}(df) = F_{X/S}^*\pi^*(df).$$
 Since this is not nilpotent, we cannot
apply our Cartier transform to it directly. Our approach will be to
approximate $E$ by nilpotent connections, and we shall see that the
Cartier transform of these approximations approximate $L$.

In general, if $(E,\nabla)$ is a connection on a smooth $X/S$ in
characteristic $p$, $F_{X/S*}(E)$ becomes an $S^\cx T_{X'/S}$ module
via the $p$-curvature $\psi$, and since $\psi$ acts horizontally,
the quotient $E_{(n)}$ of $E$ by the $n$th power of the ideal $S^+
T_{X'/S}$ of $S^\cx T_{X'/S}$ inherits a connection. In fact, this
quotient is the maximal quotient of $E$ on which the connection is
nilpotent of level $n-1$. In the situation at hand, we can be quite
explicit. Let $J \subseteq \oh X$ be the ideal of the zeroes of
$df$, \ie, the ideal generated by the partial derivatives of $f$ in
any local system of coordinates. Then $E_{(n)}$ is the quotient of
$E$ by $F_X^*(J^n)$. Our  next goal is the computation of the
Cartier transform of a suitable quotient of $E_{(n)}$.

%(Since $\cX$ is noetherian, there exists a constant
%$a$ such that $\tilde I^a \subseteq \tilde J$, but we don't really need
%this unless we want to be explicit.)

\begin{proposition}\label{carcal.p}
Suppose that $\cX/\cS$  and $\tilde f$ satisfy the conditions of
Lemma~\ref{prepz.l}, and let $Z$ be the reduction of $\cZ$ modulo
$p$. Let $n$ be a natural number and
\begin{eqnarray*}
E &:=& (\oh X, d+ df) \in \mic(X/S)\\
E_n &:=& (E/F_X^*(I^n_{Z})E, d + df) \in \mic(X/S)\\
L &:=& (\oh X, - df) \in \hig(X/S)\\
L_n &:= & L/I_Z^n L \in \hig(X/S)\\
L'_n &:=& \pi^*L_n \in \hig(X'/S).
\end{eqnarray*}
Finally, let $r$ be the maximum codimension of $\cZ$ in $\cX$. Then
if $p > rn$, the Cartier transform $C_{\cX/S}(E_n)$  of $E_n$ with
respect to $\cX'_1/\cS_1$ is isomorphic to $L'_n$.
\end{proposition}
\begin{proof}
Note that, by definition, $I_\cZ$ is the radical of the ideal $J$,
so $E_n$ is indeed a quotient of $E_{(n)}$ and $C_{\lift X S}(E_n)$
is defined. It is enough to prove the proposition after restricting
to each connected component of $\cZ$.  To simplify the notation, we
shall assume that $\cZ$ is connected. Replacing $\tf$ by $\tf-\tilde
c$, for a suitable  $\tilde c \in \Gamma(\oh \cS)$ as in (2) of
Lemma\ref{prepz.l},
 we may assume that the restriction of $\tf$
to $\cZ$ vanishes.

Recall from Proposition~\ref{deltapi.p}  that the lifting $\tilde
\pi$ of $\pi \colon X' \to X$ determines  a map  $\delta_{\tilde
\pi} \colon \oh \cX  \to \cA_{\lift X S}$.

\begin{claim}\label{alph.cl}
Let $\alpha := \delta_{\tilde \pi}(\tilde f) \in \cA_{\lift X S}$
and let
$$\beta := 1  + \alpha + {\alpha^2  \over 2!} + \cdots {\alpha^{p-1}\over (p-1)!} \in \cA_{\lift X S}.$$
  Then:
\begin{enumerate}
\item{$\alpha^{r} \in F_X^*(I_Z)\cA_{\lift X S}$, }
\item{$\psi_{\cA} (\beta) =   \bigl(\beta - {\alpha^{p-1}\over (p-1)!} \bigr) \ot F_X^*df $, and}
\item{$\nabla_\cA(\beta) = - (\beta - {\alpha^{p-1}\over (p-1)!} ) f^{p-1}\ot df$.}
\end{enumerate}
\end{claim}
\begin{proof}
By (1) of Lemma~\ref{prepz.l}, $Z$ is is reduced and in particular
satisfies Serre's condition $S_1$. Since $X$ is regular, its
absolute Frobenius endomorphism is flat, and hence the inverse image
$Z^{(p)}$ of $Z$ by $F_X$ still satisfies $S_1$. (To see this, let
$j \colon U\to Z$ be the inclusion of any dense open subset of $Z$
and observe that the map $\oh Z \to j_*j^*\oh Z$ is injective, and
remains so after pullback by $F_X$.) Since $\spec_X \cA_{\lift X S}$
is smooth over $X$, the  inverse image of $Z^{(p)}$ in $\spec_X
\cA_{\lift X S}$ also satisfies $S_1$.   Thus it suffices to check
(1)  at the generic points of $Z$, and since $Z/S$ is generically
smooth,  we may assume that it is smooth. We  may work in a
neighborhood of a point of $Z$ with the aid of a system of local
coordinates $(\tilde t_1,\ldots \tilde t_n)$ for $\cX/\cS$  such
that $I_\cZ = (\tilde t_1,\ldots \tilde t_{s})$. Let  $\tF \colon
\cX_1 \to  \cX'_1$  be the  lift of $F_{X/S}$ sending $\tilde \pi^*
\tilde t_i$ to $\tilde t_i^p$ for all $i$.  This defines a splitting
of the fundamental exact sequence (\ref{nexc.e}), and hence an
isomorphism
$$ N_1 \cA_{\lift X S} = \cE_{\lift X S} \cong \oh X \oplus F_{X/S}^*\Omega^1_{X'/S}.$$
Proposition~\ref{deltapi.p} says that, in terms of this splitting, $
\alpha = (g,F_X^*df)$, where $\tF^*\tilde \pi^*(\tilde f) = \tilde
f^p + [p] g$. Since $\tilde f \in I_\cZ$ $\tilde f^p$ belongs to
$I_\cZ^p$, and since $\tF^*\tilde \pi^*$ maps $I_\cZ$ to $I_\cZ^p$,
it follows that $[p]g \in I_\cZ^p$. It follows from the smoothness
of $\cZ$ and $\cX$ over $\cS$ that the closed subscheme of $\cX$
defined by $I_\cZ^p$ is flat over $\cS$, and hence that $g \in
I_Z^p$.  Then $g^s \in I_Z^{ps}$, and since $I_Z$ has $s$
generators, $I_Z^{ps} \subseteq F_X^*I_Z\oh X$,
 so in fact $g^s \in F_X^*I_Z\oh X$.
Since $df \in  I_Z\Omega^1_{X/S}$ by hypothesis, $F_X^*df \in
F_X^*(I) F_{X'/S}^*(\Omega^1_{X'/S})$. Thus $\alpha^s \in F_X^*(I_Z)
\cA_{\lift X S}$, and since $s \le r$, the same is true of
$\alpha^r$.  This proves (1).

Recall from Proposition~\ref{deltapi.p} that  $\psi_\cA (\alpha) =
F_X^*df$.  Hence
\begin{eqnarray*}
\psi_{\cA} (\beta) &=&   \bigl(1 + {\alpha \over 1!} +
{\alpha^2\over 2!}
         + \cdots +{\alpha^{p-2}\over(p-2)!}\bigr) \psi( \alpha)\\
         & = &  \bigl(\beta - {\alpha^{p-1}\over (p-1)!} \bigr) \ot F_X^*df. \\
\end{eqnarray*}
This proves (2). Proposition~\ref{deltapi.p} also says that
$\nabla_\cA (\alpha) = -f^{p-1}df$, so a similar calculation  proves
(3).
\end{proof}

 Recall from Theorem~\ref{rh.t} that
$$C^{-1}_{\lift X S}(L'_n) :=  \left (L'_n \ot F_{X/S*}\cA_{\lift X S}\right )^{\psi_{tot}} .$$
As an $\oh X$-module, $L'_n \ot F_{X/S*} \ca_{\lift X S} \cong
\ca_{\lift X S }/F_X^*(I_Z^n) \ca_{\lift X S}$. Since  $\alpha^{rn}
\in F_X^*(I_Z^n)$ by (1) of the claim, this module is annihilated by
$\alpha^{p-1}$ if $p > rn$.  Hence
$$\psi_\cA(\beta) = \beta \ot F_X^*df  \quad\mbox{in} \quad L'_n\ot \cA_{\lift X S} \ot F_{X'/S}^*\Omega^1_{X'/S}.$$
  Hence if we view $\beta$
as a global section of $L'_n\otimes \cA_{\lift X S}$, we find
$$\psi_{tot}(\beta) = \psi_L(1) \beta + \psi_\cA(\beta)
        = -\beta \ot \pi^*df + \beta\ot \pi^*df = 0.  $$
Thus $\beta \in  C^{-1}_{\lift X S} (L'_n) =  (L'_n\ot \cA_{\lift X
S})^{\psi_{tot}}$, and in fact $\beta$ is a basis for $C^{-1}_{\lift
X S}(L'_n)$ since it is a unit modulo $I$. Furthermore, it follows
from  (3) of the claim that
$$\nabla_{tot} (\beta) = -  \beta \ot f^{p-1}df \in E'_n \otimes \Omega^1_{X/S}. $$
 Now  consider  the Artin-Hasse exponential of $f$, which is given formally  by
$$g :=  \exp (f + f^p/p + f^{p^2}/p^2 + \cdots ),$$
and which in fact has $p$-adically integral coefficients. Then
$$g^{-1}dg = (1 + f^{p-1} + f^{p^2-1} + \cdots )df.$$
Since $f \in I_Z$ and ${p} > n$, $f^{p^2-1} = F_X^*(f^{p-1})
f^{p-1}\in   F_X^*(I_Z^n)$, so
$$dg = g(1 + f^{p-1}) df \qquad\bmod F_X^*(I_Z^n).$$
 Since $g$ is a unit,
$e := g\beta$ is also a basis for $C^{-1}_{\lift X S}(L'_n)$, and
\begin{eqnarray*}
\nabla( e) &= &g\nabla(\beta) + \beta \otimes dg \\
      & = & - g\beta \ot f^{p-1}df +  g(1+ f^{p-1}df)  \beta\otimes df\\
  & = & e\ot df
\end{eqnarray*}
In other words, $C^{-1}_{\lift X S}(L'_n)$ is isomorphic to $E_n$,
as claimed.
\end{proof}

% Proposition (\ref{carcal.p}) and Theorem (\ref{DR.t}) will allow
% us to compare  the complexes
% $F_{X'/S*}(E_n\otimes \Omega^\cx_{X/S})$
% and $L'_n\otimes \Omega^\cx_{X'/S}$.  There still remains the task
% of comparing the  complexes $F_{X'/S*}(E) \otimes \Omega^\cx_{X/S}$
% and  $L'\otimes \Omega^\cx_{X'/S}$.
%As a first step, we compare
%these to their  formal completions.
%This one can do quite generally.
% First of all, note that by (\ref{forcomp.p}), we may restrict to the formal
% completion of $X$ along $Z$.
We shall also need the following general result about morphisms in
the derived category.

\begin{proposition}\label{derar.p}
Let $X$ be a noetherian scheme or formal scheme, let $K^\cx$ be a
perfect complex of coherent sheaves of $\oh X$-modules, and let $J$
be a sheaf of ideals annihilating the cohomology sheaves of $K$.
Then there exists a natural number  $n$ such that
 for all $m > 0$, the map in the derived category
$$J^{n+m}\lotimes K^\cx \to J^m\lotimes K^\cx$$
induced from the inclusion $J^{n+m} \to J^n$
 is zero.
\end{proposition}
\begin{proof}
First we prove the statement for the induced maps on cohomology
sheaves. We may cover $X$ by a finite number of open affines on each
of which $K^\cx$ is quasi-isomorphic to a bounded complex  $\tilde
K^\cx $of locally free $\oh X$-modules. and it suffices to prove the
local statement on each of these open sets.
 Thus we may assume that $X = \spec A$ and replace $K^\cx$
by $\Gamma(X, \tilde K^\cx)$.  Then $J^m \lotimes K^\cx \cong
J^m\otimes K^\cx$ for all $m$. Let $B^q \subseteq Z^q \subseteq K^q$
be the boundaries, (resp. cycles, resp. chains) of $K^\cx$ in degree
$q$.
  By the Artin-Rees lemma,
there exists an integer $r$ such that $Z^q\cap J^{m+r}K^q \subseteq
J^{m}Z^q$ and $B^q\cap J^{m+r}K^q\subseteq J^{m}B^q$ for all $m \ge
0$. The hypothesis on $J$ implies that $JZ^q \subseteq B^q$.
 Hence if $n > r$,
$Z^q\cap J^{m+n}K^q \subseteq JZ^{q} \subseteq B^q$, so
$$ Z^q\cap J^{m+n}K^q \subseteq  B^q \cap J^{m+n} K^q \subseteq J^{m} B^q.$$
Since $K^q$ is free, $Z^q(J^{m+n}K^q) = Z^q\cap J^{m+n} K^q$ and
$J^{m}B^q = B^q(J^{m}K)$.  It follows that the map $H^q(J^{m+n}K)
\to H^q(J^{m}K)$ is zero.

The following lemma then completes the proof of the proposition.
\begin{lemma}
Let $K_0 \rTo^{f_0} K_1 \rTo^{f_1} K_2 \rTo^{f_2} \cdots K_{n+1}$ be
a sequence of morphisms in the derived category of an category.
Suppose that each $K_i$ has cohomological amplitude in $[a, a+n]$
and that the maps $H^*(K_i) \to H^*(K_{i+1})$ are all zero.  Then
the composition $K_0 \to K_{n+1}$ is zero.
\end{lemma}
\begin{proof}
The proof is by induction on $n$.  If $n = 0$, there is nothing to
prove, since $K_i \cong H^a(K_i)$ for all $i$. Let $\tau_{<}$ denote
the canonical filtration \cite{bbd.fperv}, let $f := f_1 f_2
\ldots,f_{n+1}$, and consider the following diagram:
\begin{diagram}
H^{a+n}(K_0)[-a-n] & \rTo^{H^{a+n}(f_0)} & H^{a+n}(K_1)[-a-n] \cr
\uTo & & \uTo^\alpha \cr K_0 &\rTo^{f_0}  & K_1 & \rTo^f & K_{n+1}
\cr \uTo && \uTo^\beta && \uTo \cr \tau_{<a+n}K_0 & \rTo &
\tau_{<a+n}K_1 & \rTo^{f'} &\tau_{<a+n}K_{n+1} .
\end{diagram}
Since $H^{a+n}(f_0) = 0$, $\alpha f_0 = 0$ and since $\beta$ and
$\alpha$ comprise a distinguished triangle, it follows that $f_0$
factors through $\beta$.  The induction hypothesis implies that $f'=
0$, and it follows that $f f_0= 0$. This proves the claim.
\end{proof}
\end{proof}

\begin{corollary}\label{dersp.c}
let $\cL := (\oh \cX,-df) \in \hig(\cX/\cS)$, let $\hat \cL$ denote
its formal completion along $\cZ$, and let $ b\colon \cL \to \cL_n$
denote the projection to the restriction of $\cL$ to the  $n$th
infinitesimal neighborhood of $\cZ$.  Then for sufficiently large
$n$,  there exists a map $s$ in the derived category making the
diagram below commute.
\begin{diagram}
\hat\cL \otimes \Omega^\cx_{\cX/\cS}& \rTo^{b}& \cL_n \otimes
\Omega^\cx_{\lift X S}  \cr \dTo^{\id} & \ldTo_s \cr \hat\cL \otimes
\Omega^\cx_{\lift X S}
\end{diagram}
\end{corollary}
\begin{proof}
Let us write $\hat \cL^\cx$ for  the complex $\hat \cL\ot
\Omega^\cx_{\lift X S}$, and consider for each natural number $n$
the exact sequence of complexes
$$ 0 \to I^n_\cZ   \hat \cL^\cx \to  \hat \cL^\cx \to \cL_n^\cx \to 0.$$
There is then a corresponding exact sequence of abelian groups
$$ \Ext^0(\cL_n^\cx,\hat  \cL^\cx) \to \Ext^0 (\hat \cL^\cx, \hat \cL^\cx)
 \to \Ext^0 ( I^n_\cZ \hat \cL^\cx, \hat \cL^\cx),$$
where $\Ext^0$ means hyperext, or equivalently, the group of
morphisms in the derived category.  It will thus suffice to prove
that the identity element of $\Ext^0 (\hat \cL^\cx, \hat \cL^\cx)$
maps to zero in $\Ext^0 ( I^n_\cZ \hat\cL^\cx, \hat \cL^\cx)$. But
the image of the identity element  is just the class of the
inclusion mapping, which vanishes for $n$ sufficiently large by
Proposition~\ref{derar.p}.
\end{proof}

\begin{proof}[Proof of Theorem~\ref{bko.t}]
Choose  $n$ as in Corollary~\ref{dersp.c} and localize $\cS$ so that
all primes less than the maximum of $rn$ and $n + \dim(\lift X S)$,
become invertible. Let $X/S$ and $Z/S$ denote the reductions of
$\cX/\cS$ and $\cZ/\cS$ modulo one of the remaining primes $p$.  Let
$X_{/Z}$ denote the formal completion of $X$ along $Z$, let $\hat E
:= E_{/Z}$ and $\hat L' := L'_{/Z}$, and consider the following
diagram:

\begin{diagram}
F_{X/S*}(E\otimes \Omega^\cx_{X/S}) & \rTo^a & F_{X/S*}(\hat
E\otimes \Omega^\cx_{X/S}) & \rTo^p & F_{X/S*}(E_n\otimes
\Omega^\cx_{X/S})  \cr
 && \dTo^h && \dTo^{c_n} \cr
L'\otimes \Omega^\cx_{X'/S} & \rTo^{a'} & \hat L'\otimes
\Omega^\cx_{X'/S} & \pile{\rTo^{b'}\\ \lTo_{s'}} & L'_n\otimes
\Omega^\cx_{X'/S}
\end{diagram}
Here $a,p,a',$ and $b'$ are the obvious maps, $c_n$ is the
quasi-isomorphism coming from Theorem~\ref{DR.t}, $s'$ is the
pullback via $\pi$ of the map $s$ of Corollary~\ref{dersp.c}, and
$h:= s'c_np$. Note that we do not know if  $b' h = c_n p$.
 The arrow $a$ is a quasi-isomorphism by Proposition~\ref{forcomp.p} and $a'$ is a
quasi-isomorphism by a similar (easier) argument.
 We shall show that  $h$ is a quasi-isomorphism, completing the proof of Theorem~\ref{bko.t}.

 Since our statement is local, we may restrict to an open
affine subset $U$ of $X$ and then choose a lifting $\tF$ of
$F_{X/S}$ mod $p^2$. Let $C_\tF(\hat E)$ be the formal Cartier
transform of $\hat E$ described in Proposition~\ref{coherent.p} with
respect to this lifting.

{\bf Claim:} There exists an invertible sheaf $\Lambda$ on $X_{/Z}'$
such that $C_\tF(\hat E) \cong L\otimes_{\oh {X'}} \Lambda$, where
$\Lambda$ is given the trivial Higgs field.

Indeed, the $F$-Higgs module corresponding to the $p$-curvature of
$C^{-1}_\tF(\hat L')$ is $F_{X/S}^*\hat L'$, and hence the
$p$-curvature of $\cHom(\hat E, C^{-1}_\tF(\hat L'))$ is zero.
Hence there exists an invertible sheaf $\Lambda$ on $X'$ such that
$\cHom(\hat E, C^{-1}_\tF(\hat L')) \cong F_{X/S}^*\Lambda$ with the
Frobenius descent connection. Then $C_\tF(\hat E) \cong
L'\otimes_{\oh {X'}} \Lambda$, where $\Lambda$ has the trivial Higgs
field.

By the compatibility of $C_\tF$ and $C_{\lift X S}$, the isomorphism
$\alpha_n \colon C_{\lift X S} (E_n) \cong L'_n$ of
Proposition~\ref{carcal.p} defines a trivialization of $\Lambda_n$.
Restricting to smaller affine if necessary,  we may assume that
$\Lambda$ is trivial, and choose an extension $\alpha$ of $\alpha_n$
to an isomorphism $C_\tF(\hat E) \to \hat L'$. Now consider the
commutative diagram of maps in the derived category:
\begin{diagram}
 F_{X/S*}(\hat E\otimes \Omega^\cx_{X/S}) & \rTo^p & F_{X/S*}(E_n\otimes \Omega^\cx_{X/S})  \cr
 \dTo^e && \dTo_c \cr
 \hat C_\tF(\hat E) \otimes \Omega^\cx_{X'/S} & \rTo^{b''} & C_{\lift X S}(E_n)\otimes \Omega^\cx_{X'/S} \cr
 \dTo^\alpha && \dTo_{\alpha_n} \cr
 \hat L'\otimes \Omega^\cx_{X'/S} & \rTo^{b'} & L'_n\otimes \Omega^\cx_{X'/S}
\end{diagram}
The arrows $e$ and $c$ are quasi-isomorphisms by
Proposition~\ref{coherent.p} and Theorem~\ref{DR.t}, respectively,
and $\alpha$ and $\alpha_n$ are quasi-isomorphisms by construction.
Furthermore,  $c_n = \alpha_n  c$, so
$$h = s'c_n p = s'\alpha_nc p = s'\alpha_nb'' e =
s'b'\alpha e = \alpha e$$ and hence is a quasi-isomorphism.
\end{proof}

\section{Appendix: Higgs fields and Higgs transforms}\label{higgs.s}
\subsection{Higgs fields over group schemes}
Let $X/S$ be a smooth morphism of schemes, let $\Omega_{X/S}$ be its
sheaf of Kahler  differentials and $T_{X/S}$ the dual of
$\Omega_{X/S}$. Recall that  a \emph{Higgs field} on a sheaf $E$ of
$\oh X$-modules is any of the following equivalent sets of data:
\begin{enumerate}
\item{an $\oh X$-linear map $\theta \colon E \to E\ot \Omega_{X/S}$ such that the
composition  of $\theta$ with the map $E \ot  \Omega_{X/S}\to E\ot
\Lambda^2\Omega_{X/S}$ induced by $\theta$ vanishes}
\item{a linear map $\theta \colon T_{X/S} \to \cEnd_{\oh X}(E)$ with the property that
the endomorphisms associated to any two sections of $T_{X/S}$
commute}
\item{an extension $\theta$ of the $\oh X$-module structure on $E$ to an $S^\cx T_{X/S}$-module
structure.}
\end{enumerate}
  If $E$ is quasi-coherent, then associated to the $S^\cx T_{X/S}$-module $E$
is a quasi-coherent sheaf $\tilde E$ of $\oh {{\bf T}^*}$-modules on
the cotangent bundle  ${\bf T}^*_{X/S}$ of $X/S$. Conversely, if
$\tilde E$ is such a sheaf, its direct image on $X$ is a
quasi-coherent sheaf of $\oh X$-modules equipped with a Higgs field.

These definitions make sense with any locally free sheaf $T$ in
place of $T_{X/S}$ and with the vector bundle $\bv T:= \spec_X S^\cx
T$ in place of cotangent bundle.  In fact, it will be useful for us
to work in an even  more general context, in which the vector bundle
${\bf T}^*_{X/S}$ is replaced by any commutative affine group $\cG$
scheme over $X$. Abusing notation, we shall denote by $\oh \cG$ the
sheaf of $\oh X$-bialgebras on $X$ corresponding to $\cG$.

\begin{definition}\label{gfield.d}
Let $\cG$ be a commutative flat affine group scheme over $X$ and let
$E$ be a sheaf of $\oh X$-modules on $X$. A \emph{$\cG$-field} on
$E$ is a structure $\theta$ of an $\oh \cG$-module on $E$,
compatible with the given $\oh X$-module structure via the map $\oh
X \to \oh \cG$.
\end{definition}

We denote by $\ghig$ the category whose objects are sheaves  of $\oh
X$-modules $E$ equipped with a $\cG$-field $\theta$ and whose
objects are morphisms compatible with the $\cG$-fields.  We will
often omit the $\theta$ from the notation when no confusion seems
likely to result.
% when $\cG$ is the cotangent bundle of $X/S$ we write
%instead $\hig(X/S)$.
As before, there is an evident equivalence between the category of
quasi-coherent objects in this category and the category of
quasi-coherent sheaves on $\cG$.  Since we will have to deal with
sheaves which are not quasi-coherent, we will not make use of the
topological space $\spec_X \oh \cG$.
%and or even the quasi-coherence of
%$\oh \cG$ as a sheaf of $\oh X$-algebras.
Nevertheless we will try to use geometric notation whenever
possible. Thus, if $\cA$ is a sheaf of $\oh X$-algebras, we denote
by $\Mod(\cA)$ the category of sheaves of $\cA$-modules on the
topological space $X$. If $\gamma^\sharp \colon  \cA \to \cB$ is a
homomorphism of sheaves of $\oh X$-algebras, we have functors:
\begin{eqnarray*}
\gamma^* \colon \Mod(\cA) \to \Mod(\cB) : &\quad M & \mapsto \cB \ot_\cA M\\
\gamma_* \colon \Mod(\cB) \to \Mod(\cA) : &\quad N & \mapsto N, \quad\mbox{with $a n := \gamma^\sharp(a)n$} \\
\gamma^! \colon \Mod(\cA) \to \Mod(\cB) : &\quad M & \mapsto
\cHom_\cA(\gamma_*\cB, M), \quad\mbox{with $(bh)(b') := h(bb')$},
\end{eqnarray*}
together with the standard adjunction isomorphisms:
\begin{eqnarray*}
\cHom_\cA(M, \gamma_*N) & \cong \gamma_*\cHom_\cB(\gamma^*M, N) \\%\quad\mbox{in $\Mod(\cA)$} \\
\cHom_\cA(\gamma_*N, M) & \cong \gamma_*\cHom_\cB(N, \gamma^!M) %\quad\mbox{in $\Mod(\cA)$}
\end{eqnarray*}
Note that even if $\cA$ and $\cB$ are quasi-coherent,  the functor
$\gamma^!$ does not preserve quasi-coherence, in general.

In our context we shall consider the following morphisms of
$X$-schemes and the corresponding morphisms of sheaves of $\oh
X$-algebras.  Here all fiber products are taken in the category of
$X$-schemes and all tensor products in the category of $ \oh
X$-modules.

\begin{notation}\label{morphs.e}
\begin{eqnarray*}
 p_i \colon \cG \times \cG \to \cG :  & (g_1,g_2)\mapsto g_i,  & p_i^{\sharp}\colon \oh \cG \to \oh \cG \otimes \oh \cG \notag \\
\sigma \colon \cG \times \cG \to \cG\times \cG :  & (g_1,g_2)\mapsto (g_2,g_1),  & \sigma^{\sharp}\colon \oh \cG\ot \oh\cG \to \oh \cG \otimes \oh \cG \notag \\
\iota \colon \cG \to \cG : &g \mapsto g^{-1}, & \iota^{\sharp} \colon \oh \cG \to  \oh \cG \\
\mu \colon \cG \times \cG \to \cG :& (g_1,g_2) \mapsto g_1g_2, &\mu^{\sharp} \colon \oh \cG\to  \oh \cG \otimes \oh \cG \\
\mu' \colon \cG \times \cG \to \cG: &(g_1,g_2)\mapsto g_2g_1^{-1}, &\mu^{\prime\sharp} \colon \oh \cG\to  \oh \cG \otimes \oh \cG \\
%\mu'' \colon \cG \times \cG \to \cG: &(g_1,g_2)\mapsto g_1^{-1}g_2, &\mu^{\prime \prime\sharp} \colon \oh \cG\to  \oh \cG \otimes \oh \cG \\
i\colon X \to \cG :& x \mapsto 0 , & i^{\sharp} \colon \oh \cG \to \oh X \\
p \colon \cG^n   \to X: & (g_1,\ldots g_n) \mapsto p(g_i), &  p^\sharp  \colon \oh X \to \oh {\cG^n}\\
j \colon G \to G: & g \mapsto 0 & j^\sharp \colon \oh \cG \to \oh
\cG
\end{eqnarray*}
\end{notation}

These are   the  projections $p_i$,  the inversion mapping $\iota$,
the group law $\mu$,
the  twisted group law%s $\mu'' := \mu \circ (\iota \times \id)$ and
$\mu' :=\mu\circ \sigma \circ(\iota\times \id)$, the augmentation
given by the zero section of $\cG$, the  structure map $\cG^n \to
X$, and the map $p\circ i$. Note that since $\iota^2 = \id_\cG$,
$\iota_* = \iota^*$.  If $E$ is any object of $\ghig$, we let
$E^\iota := \iota_* E = \iota^*E$.

\subsection{Convolution}\label{conv.ss}
\begin{definition}\label{cnv.d}
Let $(E_1,\theta_1)$ and $(E_2,\theta_2)$ be two objects of $\ghig$.
Then
\begin{enumerate}
\item{$E_1\boxtimes E_2 := p_1^*E_1\otimes_{\oh {\cG\times \cG}}p_2^*E_2$, as an object of $\cG\times \cG$-$HIG$.}
\item{$E_1\cnv E_2 := \mu_*(E_1\boxtimes E_2)$,  as an object of $\ghig$.}
\end{enumerate}
\end{definition}

For example, if $(E_1,\theta_1)$ and $(E_2,\theta_2)$ are objects of
$\hig(X/S)$, then ${E_1\cnv E_2}$ is the tensor product of $E_1$ and
$E_2$ in the category of $\oh X$-modules, with the  Higgs field
$\theta$ defined by
$$\theta = \theta_1\ot \id_{E_2} + \id_{E_1}\ot\theta_2.$$
Geometrically, the object $(E_1\ot E_2,\theta)$ corresponds to the
convolution of $E_1$ and $E_2$ with respect to the group structure
of the cotangent space of $X/S$.

The associative law for $\cG$ implies that the standard isomorphism
$${(E_1\ot E_2)\ot E_3} \cong E_1\ot (E_2\ot E_3)$$
induces an isomorphism
$$(E_1\cnv E_2)\cnv E_3 \cong E_1\cnv ( E_2\cnv E_3).$$
Similarly, the commutativity of $\cG$ implies that the standard
isomorphism $E_1\ot E_2 \cong E_2 \ot E_1$ induces an isomorphism
$$E_1 \cnv E_2 \cong E_2\cnv E_1.$$
Furthermore, if we let
$$U := i_* \oh X \in \ghig,$$
then the fact that $i$ is the identity section implies that the
natural isomorphism $\oh X \ot_{\oh X} E \cong E$ induces an
isomorphism in $\ghig$:
$$U\cnv E \cong E.$$
Thus $\cnv$ makes the category $\ghig$ into an $\oh X$-linear tensor
category~\cite{de.tc} (ACU tensor category in the terminology of
\cite{sa.tc}), and $U$ is its unit object.

\begin{definition}\label{hex.d}
Let $E_1$ and $E_2$ be objects of $\ghig$.  Then
\begin{eqnarray*}
\cH ex(E_1,E_2) &: =& \cHom_{\oh {\cG\times \cG}}(p_1^*E_1,p_2^!E_2)\\
\cH(E_1,E_2) &:= &\mu'_*\cH ex(E_1,E_2)
\end{eqnarray*}
\end{definition}

We call $\cH ex(E_1,E_2)$  the \emph{external Hom} of $E_1$ and
$E_2$. Its underlying $\oh X$-module is given by
\begin{equation*}
\begin{split}
  \cHom_{\oh X}(p_*E_1, p_*E_2 )\cong
p_*\cHom_{\oh \cG}(p^*p_*E_1, E_2) \cong
p_*\cHom_{\oh\cG}(p_{2*}p_1^*E_1,E_2) \cong  \cr p_*\cHom_{\oh
{\cG\times \cG}}(p_1^*E_1,p_2^!E_2) \cong p_*\cH ex(E_1,E_2) \cong
p_*\cH(E_1,E_2)
\end{split}
\end{equation*}
and the $\oh \cG \ot \oh \cG$ structure on $\cH ex(E_1,E_2)$ is
given by:
$$(a\ot b)h   \colon  E_1 \to E_2 \quad e_1 \mapsto bh(ae_1).$$

\begin{lemma}\label{hadj.l}
Let $E_1, E_2$, and $E_3$ be objects of $\ghig$.  Then the standard
adjunction isomorphism in the category of $\oh X$-modules
$$ \cHom_{\oh X}(E_1\ot_{\oh X}E_2,E_3) \cong \cHom_{\oh X}(E_1,\cHom_{\oh X}(E_2,E_3))$$
induces  isomorphisms
$$\Hom_{\oh \cG}(E_1\cnv E_2, E_3) \cong \Hom_{\oh \cG}(E_1,\cH(E_2,E_3))\quad \mbox{(of groups)}$$
$$\cH(E_1\cnv E_2,E_3) \cong \cH(E_1,\cH(E_2,E_3)) \quad \mbox{(in $\ghig$)}.$$
\end{lemma}
\begin{proof}
By definition,
\begin{eqnarray*}
\Hom_{\oh \cG}(E_1,\cH(E_2,E_3))
    & = &   \Hom_{\oh \cG}(E_1,\mu'_*\cHom_{\oh {\cG\times \cG}}(p_1^*E_2,p_2^!E_3))\\
    & = &   \Hom_{\oh {\cG\times \cG}}({\mu'}^*E_1,\cHom_{\oh {\cG\times \cG}}(p_1^*E_2,p_2^!E_3))\\
    & = &   \Hom_{\oh {\cG\times \cG}}({\mu'}^*E_1\ot_{\oh {\cG\times \cG}}p_1^*E_2,p_2^!E_3))
\end{eqnarray*}
Let $\alpha \colon \cG\times \cG\to \cG\times \cG$ denote the map
$(\mu',p_1)$, \ie, the map sending $(g_1,g_2)$ to
$(g_2g_1^{-1},g_1)$.  Note that $\alpha$ is an isomorphism, whose
inverse $\beta= (p_2,\mu)$ sends $(a,b)$ to $(b,ab)$.  Thus $\beta_*
= \alpha^*$, and furthermore ${\mu'}^*(E_1)\ot p_1^*E_2 =
\alpha^*(E_1\boxtimes E_2)$.  Hence
\begin{eqnarray*}
\Hom_{\oh \cG}(E_1,\cH(E_2,E_3))
    & = &   \Hom_{\oh {\cG\times \cG}}({\alpha}^*(E_1\boxtimes E_2),p_2^!E_3))\\
    & = &   \Hom_{\oh {\cG\times \cG}}({\beta}_*(E_1\boxtimes E_2),p_2^!E_3))\\
    & = &   \Hom_{\oh {\cG\times \cG}}(E_1\boxtimes E_2,\beta^!p_2^!E_3))\\
    & = &   \Hom_{\oh {\cG\times \cG}}(E_1\boxtimes E_2,\mu^!E_3))\\
    & = &   \Hom_{\oh {\cG}}(\mu_*(E_1\boxtimes E_2),E_3))\\
    & = &    \Hom_{\oh \cG}(E_1\cnv E_2,E_3)
\end{eqnarray*}
This proves the first statement.  The second statement just asserts
that the standard adjunction morphism is compatible with the
$\cG$-Higgs fields. It follows formally from the first. Indeed, it
will suffice to check that for all $E$,  the adjunction isomorphism
induces isomorphisms:
$$\Hom_{\oh G}(E, \cH(E_1\cnv E_2, E_3)) \cong
   Hom_{\oh G}(E,\cH(E_1,\cH(E_2,E_3))).$$
This follows from the first statement and the associativity of
$\cnv$.
\end{proof}

This shows that $\cH$ is the internal Hom functor of the tensor
category $(\ghig,\cnv)$ in these sense of \cite{de.tc}. As usual,
the \emph{dual} of an object $E$ of $\ghig$ is defined by
\begin{equation}\label{dualh.e}
E^\vee := \cH(E,U).
\end{equation}
The map
$$ev \colon E^\vee \cnv E \to U$$
is by definition the element of
$$Hom_{\oh G}(E^\vee \cnv E, U) = Hom_{\oh G}(E^\vee, \cH(E, U)) = Hom_{\oh G}(E^\vee, E^\vee)$$
corresponding to $\id_{E^\vee}$; it corresponds to the usual
evaluation map
$$\cHom_{\oh X}(E,\oh X) \ot E\to \oh X.$$
  For any $E_2$, one gets by functoriality maps
\begin{eqnarray*}
E_2 \cong \cH(U,E_2) & \to&\cH(E_1^\vee\cnv E_1, E_2)\\
\Hom_{\oh G}(E_2, \cH(U,E_2)) & \to  & \Hom_{\oh G}(E_2,\cH(E_1^\vee\cnv E_1, E_2))\\
\Hom_{\oh G}(E_2\cnv U,E_2)       & \to & \Hom_{\oh G}(E_2\cnv E_1^\vee \cnv E_1,E_2)\\
\cHom_{\oh G}(E_2, E_2)       & \to & \Hom_{\oh G}(E_2\cnv E^\vee_1,\cH(E_1,E_2))\\
\cHom_{\oh G}(E_2, E_2)       & \to &  \Hom_{\oh G}( E^\vee_1,\cnv E_2, \cH(E_1,E_2))\\
\end{eqnarray*}
The element of $\Hom(E_2\cnv E_1,\cH(E_1,E_2))$ corresponding to
$\id_{E_2}$ is the map
\begin{equation}\label{tto.e}
E_1^\vee \cnv E_2\to \cH(E_1,E_2)
\end{equation}
corresponding to the usual map  $E_1^\vee \ot E_2 \to \cHom_{\oh
X}(E_1,E_2)$ in the category of $\oh X$-modules. In particular it is
a  homomorphism in $\ghig$ and commutes with any endomorphism of
$E_1$ or $E_2$ in the category $\ghig$.  For example, any local
section of $\oh \cG$ defines such an endomorphism  on each $E_i$.
%Recall that such a section acts in the evident ways on the second
%factor, but the action on $E_1^\vee$ is through $\iota$.
Note that if $E_1$ is locally free and $E_1$ or $E_2$ is of finite
presentation as an $\oh X$-module, (\ref{tto.e}) an isomorphism. For
example, when $\cG$ is the cotangent space of $X$ and $\theta$ is a
Higgs field on $X$, then the Higgs field $\theta^\vee$  on $E^\vee$
is given by the usual rule, so that
$$\angles{\theta_\xi (\phi)} e +\angles \phi {\theta_\xi(e)}= 0$$
for sections $\xi$ of $T$, $\phi$ of $E^\vee$ and $e$ of $E$.

\begin{remark}\label{trhig.r}
  If $E_1$ and $E_2$ are objects of $\ghig$, the
$\oh X$-module underlying $\cH(E_1,E_2)$ is $\cHom_{\oh X}(E_1,E_2)$
and the $\oh X$-module underlying $E_1\cnv E_2$ is $E_1\ot_{\oh X}
E_2$. These $\oh X$-modules also inherit $\oh \cG$-structures by
``transport of structure'' from the $\oh \cG$-module structures of
$E_1$ and $E_2$.  When necessary we denote by $\theta_{E_i}$ the
structure coming from $E_i$ in this way and by $\theta_{tot}$ the
structure defined in (\ref{cnv.d}) and (\ref{hex.d}).  Thus
$\theta_{E_i}$ is the structure on $\cH(E_1,E_2)$ (resp. $E_1\cnv
E_2$) obtained from the structure on $\cH ex(E_1,E_2)$ (resp. $E_1
\boxtimes E_2)$ by letting $\oh \cG$ act via the morphism $p_{i*}$.
Note in particular that the $\oh \cG$-module structure on $E^\vee$
is not the structure $\theta_E$  corresponding to the action by
transport of structure on $\Hom_{\oh X}(E, \oh X)$, rather it is
given by $\iota_*\theta_E$.  Indeed, the $\oh \cG \otimes \oh
\cG$-module $\cH ex( E, \oh X)$ is annihilated by the ideal of the
graph  $\Gamma_j$ of the zero morphism $j \colon \cG \to \cG$, and
$\mu' \circ \Gamma_j = \iota$.
\end{remark}

\begin{remark}\label{cnvp.r}
A morphism $h \colon \cG' \to \cG$ of affine $X$-schemes induces a
pair of adjoint functors
$$h^* \colon \ghig \to \mbox{$\cG'$-$\hig$} \quad \mbox{and} \quad
h_*\colon    \mbox{$\cG'$-$\hig$} \to \ghig.$$ If $h$ is a
homomorphism of group schemes, these are compatible with $\cnv$ and
$\cH$. For example, let $f \colon X \to Y$ be a morphism of schemes,
let $\cH$ be a commutative affine group scheme over $Y$, and let
$f^{-1}\cH$ be its pullback to $X$. If $(E,\theta)$ is an object of
$\cH$-$\hig$, then $f^*E$ has a natural $f^{-1}\cH$-field
$f^*\theta$.  If $\cG$ is an affine group scheme over $X$ equipped
with a map $h \colon f^{-1}\cH \to \cG$, then one gets by
composition with $h$ a $\cG$-field on $f^*E$. For example, this
construction applied to the cotangent bundles, with $h$ the
differential of $f$, defines a functor $f^* \colon \hig(Y/S) \to
\hig(X/S)$.
% For another example, observe that $\mu$ and $\mu'$
% are morphisms $ \cG \times \cG \to \cG$
% and that if $E_1$ and $E_2$ are objects of $\ghig$, then the
% Higgs field $\theta$ on
% $E_1\cnv E_2$ (resp. $\cH(E_1,E_2)$) is just
% $\mu_*(E_1\boxtimes E_2)$ (resp. $\mu'_*\cH ex(E_1,E_2)$).
%   In fact, each of the two projections
% $p_i \colon \cG\times \cG \to \cG$ also allows us
% to define a Higgs field  on these modules,
% which we denote by $\theta_i$.
Finally, note that since $\iota \colon \cG \to \cG$ is a group
homomorphism, we find a canonical isomorphism
$$(E_1\cnv E_2)^\iota \cong E_1^\iota \cnv E_2^\iota.$$

%If $X/S$
%is a morphism of schemes in characteristic $p >0$, let
% $\bT^{\prime *}_{X/S}$ be the pullback of $\bT^*_{X/S}$
%by the absolute Frobenius endomorphism  of $X$.  We shall
% refer to $\bT^{\prime *}_{X/S}$-fields as \emph{$F$-Higgs fields.}
%For example, the $p$-curvature of a module with connection $(E,\nabla)$
%is an $F$-Higgs field on $E$~\cite{o.hcpc}.
% Note that a Higgs field on $E$ induces
%an $F$-Higgs field on $F_X^*(E)$, but the resulting Higgs field is just zero.
\end{remark}

\begin{remark}\label{homg.r}
Let $E_1$ and $E_2$ be object of $\ghig$. Then there is a
 natural isomorphism of sheaves of $\oh X$-modules
$$i^!\cH(E_1,E_2) \cong p_*\cHom_{\oh \cG}(E_1,E_2)  .$$
This follows from the adjointness properties of $\cH$:
\begin{eqnarray*}
i^!\cH(E_1,E_2) & :=& \cHom_{\oh \cG}\left(i_* \oh X ,\cH(E_1,E_2)\right)\\
   & \cong &\cHom_{\oh \cG}(U \cnv E_1,E_2) \cong \cHom_{\oh \cG}(E_1,E_2).
\end{eqnarray*}
We find a natural map  of $\oh X$-modules
\begin{equation}\label{gammb.e}
 \cHom_{\oh \cG}(E_1,E_2) \cong i^!\cH(E_1,E_2) \subseteq p_*\cH(E_1,E_2)
\end{equation}

This map is compatible with the actions of $\oh \cG$ induced by
transport of structure via its actions on $E_1$ and $E_2$.
\end{remark}

\subsection{Higgs transforms}
We can use an object of $\ghig$ as a kernel for what we shall call
a \emph{Higgs transform}, of which we consider the following
variants.

\begin{definition}\label{higtr.d}
Let $\cF$ be an object of $\ghig$. Define functors from $\ghig$ to
itself by:
\begin{eqnarray*}
\cH_\cF(E) &:= &\cHom_{\oh \cG}(\cF,E)  \\
\cT_\cF(E) &:= & \cF\otimes_{\oh \cG} E \\
\cT^!_\cF(E) & : = &i^!(\cF\cnv E) \\
%\cH^*_\cF(E) & : = &i^*\cH(\cF, E)
\end{eqnarray*}
\end{definition}
We view these objects as $\cG$-Higgs modules, with the Higgs field
induced by transport of structure from the field on $\cF$.  Note
that this is the same as the field induced from $E$ in the first and
second and cases, and differs by $\iota$ in the  third. Indeed,
$$i^!(\cF\cnv E) = \Hom_{\oh \cG}(i_*\oh X,\mu_*(\cF\boxtimes E))
\cong \Hom_{\oh {\cG\times \cG}}(\mu^*i_*\oh X, \cF\boxtimes E)$$
This is an $\oh{\cG\times \cG}$-module, and the action of $\oh\cG$
by transport of structure via $\cF$ corresponds to the action
induced by the first projection, while the action via $E$ is induced
by the second projection. We claim these differ by $\iota$. In fact
it is enough to check this for the $\oh {\cG\times \cG}$-module
$i_*\mu^* (\oh X)$. But this is clear, since the latter is
annihilated by the ideal of the graph of $\iota$.

The map (\ref{tto.e}) induces a  natural map of $\oh \cG$-modules
\begin{eqnarray}\label{ttoh.e}
\iota^*\cT^!_{\cF^\vee}(E) &\to & i^!(\cH(\cF, E)) \cong
\cH_{\cF}(E)
%\iota^*\cT_{\cF^\vee}(E) &\to& i^*(\cH(\cF ,E)) \cong \cH^*_{\cF}(E)
\end{eqnarray}
which is an isomorphisms if $\cF$ is finitely generated and
projective as an $\oh X$-module.  The presence of the $\iota$ is due
to the fact that $\oh \cG$ acts on $\cT^!_{\cF^\vee}(E)$ by
transport of structure via  $\cF^\vee$ and on $\cH_\cF(E)$ via
$\cF$, and these structures differ by $\iota$, as we saw in
Remark~\ref{trhig.r}.

Let $\cB$ be any object of $\ghig$. Observe that there are natural
transformations:
$$  \eta_\cB \colon  \cT_\cB \circ \cH_\cB \to \id
 \quad\mbox{and}\quad
\zeta_\cB \colon \id \to \cH_\cB\circ \cT_\cB,$$ where for any $E$,
$$\eta_{\cB,E}\colon  \cB\otimes_{\oh \cG} \cHom_{\oh \cG}(\cB,E) \to E $$
sends $b \ot h$ to $h(b)$  and
$$\zeta_{\cB,E} \colon E  \to \cHom_{\oh \cG}(\cB,\cB\otimes_{\oh \cG} E)$$
sends $e$ to the homomorphism  $b \mapsto b\ot e$. Then the
following result is immediate.
\begin{proposition}\label{ht.p}
Let $\cB$ be an object of $\ghig$ which is invertible as an $\oh
\cG$-module.  Then  $\eta_\cB$ and $\zeta_{\cB}$ are  inverse
isomorphisms.  In particular, $\cH_\cB$ and $\cT_\cB$ are
quasi-inverse equivalences of categories.
\end{proposition}\qed

In our main application, $\cB$ will be equipped with the structure
of a cocommutative coalgebra with counit in the category $\ghig$.
That is, it will be provided with morphisms
\begin{equation}\label{charm.e}
 \nu \colon \cB \to \cB\cnv \cB \quad\mbox{and} \quad
\beta \colon \cB \to U \,
\end{equation}
satisfying the usual compatibilities.  These data give the
corresponding Higgs transform additional structure.  We shall be
especially interested in the following case.

\begin{definition}\label{char.d}
A \emph{character sheaf\,}\footnote{This terminology is borrowed
from G. Lusztig} in $\ghig$ is a coalgebra $(\cB,\nu,\beta)$ such
that the map $\tilde \nu \colon \mu^* \cB \to \cB\boxtimes \cB$
corresponding to $\nu$ is an isomorphism of $\oh{ \cG\times
\cG}$-modules and the map $\tilde \beta \colon i^* \cB \to \oh X$
induced by $\beta$ is an isomorphism of $\oh X$-modules.
\end{definition}

\begin{lemma}\label{char.l}
Let $(\cB,\nu,\beta)$ be a character sheaf in $\cG$. The map $\nu$
induces an isomorphism of $\oh \cG$-modules $\cB \otimes_{\oh \cG}
\cB^\iota \cong \oh\cG$.  In particular, $\cB$ is invertible as an
$\oh \cG$-module, with inverse $\cB^\iota$.
\end{lemma}
\begin{proof}
Definition (\ref{char.d}) implies that $i^* \cB\cong \oh X$, hence
$j^* \cB \cong \oh \cG$. Let $\Gamma_\iota \colon \cG \to \cG\times
\cG$ denote the graph of $\iota$. Pulling the isomorphism  $\tilde
\nu$ back via $\Gamma_\iota$, we find an isomorphism:
$$\cB {\otimes_{\oh \cG}} \cB^\iota\cong  \Gamma_\iota^*(\cB\boxtimes \cB) \cong
\Gamma_\iota^*\mu^*(\cB) \cong j^*\cB \cong \oh \cG.$$
\end{proof}
\begin{remark}\label{char.r}
  Giving a character sheaf is equivalent to giving a commutative
extension of the group $\cG$ by the multiplicative group $\gm$.
Assume $G$ is finite and flat over $X$ and denote by $\cG^\vee$ its
Cartier dual.  Then a  $G^\vee$-torsor gives rise to a character
sheaf.  Indeed, $G^\vee$ can be identified with the group of
automorphisms of the trival extension  ${\bf E}$ of $\cG$.  Then if
${\bf L}$ is a $\cG^\vee$-torsor,
 we can form the ``twist''
${\bf L}\otimes_{\cG^\vee} {\bf E}$ of ${\bf E}$ by ${\bf L}$. This
defines an extension of $\cG$ by $\gm$ and hence a character sheaf
$\cB$.  The functor ${\bf L} \mapsto \cB$ is an equivalence from the
category of $\cG^\vee$-torsors to the category of character sheaves
which are,  locally on $X$,  isomorphic to $\oh \cG$ as character
sheaves.
\end{remark}

\begin{proposition}\label{charconv.p}
Let $(\cB,\nu,\beta)$ be a character sheaf on $\cG$. If $E_1$ and
$E_2$ are objects of $\ghig$, the standard tensor product map $\ot $
on homomorphisms and the comultiplication $\nu$ define a commutative
diagram:
\begin{diagram}
\cH_\cB(E_1) \cnv \cH_\cB(E_2)  &\rTo^\otimes &\cH_{\cB\cnv
\cB}(E_1\cnv E_2) \cr & \rdTo_\cong & \dTo_{\nu^*} \cr
  & & \cH_\cB(E_1\cnv E_2)
\end{diagram}
Thus the functor $\cH_\cB$ is an auto-equivalence of the tensor
category $\ghig$.
\end{proposition}
\begin{proof}
The diagram above can be expanded as follows.
\begin{diagram}
\mu_*\left (\cH_{\cB}(E_1)\boxtimes \cH_\cB(E_2)\right )  \cr
\uTo^{:=} &\rdTo^{:=} \cr \cH_\cB(E_1) \cnv \cH_\cB(E_2)
&&\mu_*\left (p_1^*\cHom_{\oh \cG}(\cB,E_1)\otimes_{\oh {\cG \times
\cG}} p_2^*\cHom_{\oh \cG}(\cB,E_2)\right ) \cr
 \dTo&&\dTo_\boxtimes \cr
 \cH_\cB(E_1\cnv E_2)&\lTo^{\tilde \nu^*} &\mu_*\cHom_{\oh {\cG\times \cG}}(\cB\boxtimes \cB,E_1\boxtimes E_2) \cr
\end{diagram}
The right vertical map is an isomorphism because $\cB$ is
invertible, and the map $\tilde \nu^*$ is an isomorphism because
$\tilde\nu$ is an isomorphism. The cocommutativity of the coalgebra
$\cB$ implies that $\cH_\cB$ preserves the commutativity constraint
of the tensor category $\ghig$.
\end{proof}

A change of sign allows us construct an \emph{involutive}
auto-equivalence of $\ghig$.
%This is a convenient improvement;
%more important is the fact that this auto-equivalence has the
%general form of a Riemann-Hilbert correspondence.

\begin{definition}\label{higtrp.e}
Let $\cB$ be an object of $\ghig$.  Then $\cH_\cB'$ is the functor
$\ghig \to \ghig$ sending an object $E$ to $\iota_*\cHom_{\oh
\cG}(\cB,E)$.
\end{definition}
For example, $\cH'_{\oh \cG}$ can be identified with the involutive
functor $\iota_*$.  More generally, if $(\cB,\nu,\beta)$ is a
character sheaf, then by Lemma~\ref{char.l}:
\begin{equation*}
\begin{split}
\cH_\cB'(\cH_\cB'(E)) = \iota_*\cHom_{\oh \cG}(\cB,
\iota_*\cHom_{\oh \cG}(\cB,E))
  \cong \iota_*\iota_*\cHom_{\oh \cG}(\iota^* \cB,\cHom(\cB, E)) \cr
   \cong \cHom_{\oh \cG }(\iota^*\cB\ot_{\oh \cG}\cB,E) \cong
          \cHom_{\oh \cG }(\oh \cG ,E) \cong E.
\end{split}
\end{equation*}

%Let us remark that
%$$\cH'_\cB(E) := \cHom_{\oh \cG} (\cB,\iota_*E) \cong \iota_* \cHom_{\oh \cG}(\cB^\iota,E).$$
The natural inclusion of $\oh X$-modules $\Hom_{\oh \cG}(\cB,E)
\subseteq \Hom_{\oh X}(\cB,E)$ defines morphisms of $\oh
\cG$-modules:
\begin{equation}\label{gamme.e}
\gamma_E \colon \cH'_\cB(E)\to \cH(\cB,j_*E) \quad\mbox{and}\quad
j_*\cH'_\cB(E) \to \cH(\cB,E).
\end{equation}
The morphism $\beta \colon \cB   \to  U$ and the isomorphism
$\cH(U,E) \cong E$ induce  maps of $\oh \cG$-modules:
\begin{equation}\label{beta.e}
\beta_E \colon E  \to \cH(\cB,E) \quad\mbox{and} \quad
    j_*E \to \cH(\cB,j_*E).
\end{equation}

Finally, let us consider the following diagram:
\begin{equation}\label{nue.e}
\begin{diagram}
\cH(\cB,\cH'_\cB(E)) &\rTo^{\gamma_{E*}} & \cH(\cB,\cH(\cB,j_*E))
&\quad & \cH(\cB,j_*\cH'_\cB(E)) &\rTo^{\gamma_{E*}} &
\cH(\cB,\cH(\cB,E)) \cr \dDashto^{\nu_E} && \dTo_\cong &\quad &
\dDashto^{\nu_E} && \dTo_\cong \cr \cH(\cB,j_*E)&\lTo^{\nu_{E*}}
&\cH(\cB\cnv\cB,j_*E) &\quad& \cH(\cB,E)&\lTo^{\nu_{E*}}
&\cH(\cB\cnv\cB,E)
\end{diagram}
\end{equation}
Here the top horizontal arrow is induced by the morphism $\gamma_E$
 (\ref{gamme.e}),
the right arrow is the adjunction map  of Lemma~\ref{hadj.l}, and
the bottom horizontal arrow is induced by $\nu$; the diagram defines
the arrow $\nu_E$.

\begin{theorem}\label{charac.t}
Let $(B,\nu,\beta)$ be a character sheaf for $\cG$ (\ref{char.d}).
Then for any object $E$ of $\ghig$, the arrow $\nu_E$ is an
isomorphism and induces an isomorphism $\kappa_E \colon
(\cH'_{\cB}(\cH'_{\cB}(E)) \to E$ fitting into the following
commutative diagram.
\begin{diagram}
\cH'_{\cB}(\cH'_{\cB}(E) & \rTo^{\gamma_{\cH'_\cB(E)}} &
\cH(\cB,j_*\cH'_\cB(E)) \cr \dTo^{\kappa_E} && \dTo_{\nu_E} \cr E &
\rTo^{\beta_E} & \cH(\cB,E)
\end{diagram}
\end{theorem}
\begin{proof}
Since $(\cB,\nu,\beta)$ is a character sheaf, the map $\tilde \nu$
is an isomorphism, and so induces an isomorphism of $\oh
\cG$-modules:
$$p_{2*}(\tilde \nu) \colon {p_2}_*\mu^* \cB \to {p_2}_*(\cB\boxtimes \cB).$$
Since the diagram
\begin{diagram}
G\times G &\rTo^\mu & G \cr \dTo^{p_2}&& \dTo_p \cr G & \rTo^p & X
\end{diagram}
is Cartesian, the natural map $p^*p_*\cB \to p_2^* \mu^*\cB$ is an
isomorphism.  Composing this map with ${p_2}_*(\tilde \nu)$, we find
that $\nu$ induces an isomorphism
$$ p^*p_*\cB \to {p_2}_*(\cB\boxtimes \cB).$$
%The same holds with $\cB$ in place of $\cB$.
 Taking $\cHom_{\oh \cG}(\ ,E)$
we find that the top arrow in the diagram below is an isomorphism.
\begin{diagram}
p_*\cHom_{\oh \cG}({p_2}_*(\cB\boxtimes \cB),E) &\rTo^{\cong} &
p_*\cHom_{\oh \cG}(p^*p_*\cB,E) \cr \uTo^{adj}&&\dTo_{adj}\cr
\cHom_{\oh X}(\cB,p_*\cHom_{\oh\cG}(\cB,E)) & \rTo& \cHom_{\oh
X}(\cB,p_*E)
\end{diagram}
The vertical arrows are  the adjunction isomorphisms and the bottom
horizontal arrow is $p_*(\nu_E)$.  It follows that $\nu_E$ is an
isomorphism of $\oh X$-modules, and it is compatible with the two
pairs of $\oh\cG$-module structures shown in the diagrams
\ref{nue.e}.
  Applying the functor $i^!$ to the isomorphism
$\nu_E$ in the left diagram, we find an isomorphism of $\oh
X$-modules:
$$\kappa_E \colon \cH'_\cB(\cH'_\cB(E)) \to E.$$
But $\mu_E$ is also compatible with the $\oh \cG$-module structures
in the right diagram. This implies that $\kappa_E$ is also a
homomorphism of $\oh \ccG$-modules.
\end{proof}

\begin{remark}
{If $\cB \to \cB'$ is a surjection of invertible $\oh \cG$-modules
and $E$ is an object of $\ghig$ which is annihilated by the
annihilator of $\cB'$, then the natural map $\cH_{B'}(E) \to
\cH_\cB(E)$ is an isomorphism.}
\end{remark}

It is sometimes convenient to use the dual point of view to that
taken in Theorem~\ref{charac.t}.  With the notation there, let
$\{I_n: n \in \bn\}$ denote an inverse system of  ideals of $\oh
\cG$ defining closed subschemes $\cG_n$ of $\cG$ whose support is
the zero section $\cG_0$.  Suppose further that each $\oh {\cG_n}$
is a locally free $\oh X$-module of finite rank and that for all $m,
n$,  the comultiplication map fits into commutative diagrams:
\begin{diagram}
\oh \cG &\rTo & \oh \cG \ot  \oh \cG \cr \dTo && \dTo \cr \oh
{\cG_{m+n}} &\rTo & \oh {\cG_m}\otimes \oh { \cG_n} .
\end{diagram}
Let $\cB_n := \cB/I_n\cB$ and let $\cA_n :=\cH(\cB_n,\oh X)$ Thus
the  map (\ref{ttoh.e}) induces an isomorphism:
$$\cT^!_{\cA_n}  \to \cH'_{\cB_n} := \iota_*\cH_{\cB_n}.$$
Let $\cA$ be the direct limit of the directed system $\cA_\cx$, so
that we find an injection
\begin{equation}\label{atob.e}
\cT^!_{\cA} \cong \dirlim \cT^!_{\cA_n} \to \cH'_{\cB}.
\end{equation}

The comultiplication maps on $\cB_\cx$ induce multiplication maps
$\cA_n \otimes \cA_m \to \cA_{n+m}$, and  $\cA$ inherits the
structure of an algebra in the category $\ghig$. The identity
element $1_\cA$ is the dual of  the identity section $i^\sharp$ of
$\cB$.

Let us say that an object $E$ of $\ghig$ is
\emph{$I_\cx$-continuous} if  each local section $e$ of $E$ is
annihilated by $I_n$ for some $n$, and let us denote the full
subcategory of $\ghig$ consisting of such objects by $\gchig$. If
$E$ is $I_\cx$-continuous, then any $\oh \cG$-linear homomorphism $h
\colon \cB \to E$ factors through $\cB_n$ for some $n$, so that
(\ref{atob.e})  becomes an isomorphism. Then Theorem~\ref{charac.t}
can be reformulated as follows.

%
% If $E_1$ is a direct limit of locally
%free sheaves of $\oh X$-modules of finite presentation $E_\lambda$, then
%the natural map
%$$E_1^\vee\hat\cnv E_2 \to \cH(E_1,E_2)$$
%is an isomorphism, where
%$$E_1^\vee\hat\cnv E_2  := \invlim \left (E^\vee_\lambda \cnv E_2\right).$$

\begin{theorem}\label{higtra.t}
Let $\cB_\cx$ and $\cA_\cx$ be as described above, and let $\gchig$
denote the full subcategory of $\ghig$ consisting of the
$I_\cx$-continuous objects.
\begin{enumerate}
\item{For any object $E$ of $\gchig$,
$$\alpha_E \colon E \to \cA\cnv E \quad  e \mapsto 1\ot e$$
defines a locally split injection whose image  is the annihilator of
$I$ with respect to the  Higgs field induced from the Higgs field on
$\cA$.}
\item{The functor
$$\cT^!_\cA \colon \gchig \to \gchig$$
 is an involutive equivalence.
Furthermore, the isomorphism $\cT^!_\cA\circ \cT^!_\cA \cong \id$
fits into the following commutative diagram
\begin{diagram}
\cT^!_\cA(\cT^!_\cA(E))&\rTo^{inc} & \cA \cnv \cT^!_\cA(E) &
\rTo^{\id\ot inc} & \cA \cnv \cA \cnv E \cr
 \dTo^{\kappa_E} && \dTo^{\nu_E} & \ldTo_{m\ot \id} \cr
E &\rTo^{\alpha_E} &\cA \cnv E.
\end{diagram}
Here $\kappa_E$ is compatible with the $G$-Higgs fields, and $\nu_E$
is a homomorphism with respect to the following fields:
$$\cA \cnv j_*\cT^!_\cA \to  \cA \cnv E \quad\mbox{and}\quad
\cA \cnv \cT^!_\cA \to  \cA\cnv j_*E. $$ }
\item{The functor $\cT^!_\cA$ is compatible with convolution products in the following sense.
For any two object $E_1$ and $E_2$ of $\gchig$, one has a canonical
commutative diagram
\begin{diagram}
\cT^!_\cA(E_1) \cnv \cT^!_\cA(E_2) &\rTo^\ot &\cT^!_{\cA \cnv
\cA}(E_1\cnv E_2) \cr & \rdTo_\cong & \dTo_{\mu^*} \cr &&
\cT^!_\cA(E_1\cnv E_2)
\end{diagram}}
\end{enumerate}
\end{theorem}
\begin{proof}
The first statement is clear.  The diagram in the second statement
is equivalent to the diagram in Theorem~\ref{charac.t}, and hence it
follows that the maps indicated are isomorphisms. Let us note,
however, that this can also be proved directly by d\'evissage. When
$E$ is annihilated by $\cI$ this is clear.  Since the sources and
targets of both arrows are exact functors, a d\'evissage argument
implies that the maps are isomorphisms if $E$ is nilpotent.  But any
$\cI_\cx$ continuous $E$ is, locally on $X$, a direct limit of
nilpotent objects, so the general result follows. Similarly, (3) is
a translation of Proposition~\ref{charconv.p}, and can also be
proved by d\'evissage.
\end{proof}

%\begin{remark}
%{\rm  If $\cB = \oh \cG$, then
%for any object $E$ of $\gchig$, $e \mapsto (\cdot e)$
%defines a split injection
%$$\iota_*(E) \to \cB^\vee \cnv \iota_*(E) \cong \iota_*\cB^\vee \cnv E \cong \cA \cnv E $$
%whose image is the annihilator of $I$ with the action $\theta_{tot}$.
%This follows from (\ref{htr.p}).}
%\end{remark}

\subsection{Examples and formulas}\label{exfor.ss}

Let us return to the case in which the group scheme is the group
underlying a vector bundle $\bv T$, where $T$ is a locally free
sheaf of $\oh X$-modules of finite rank, so $\oh {\bv T}$ is the
symmetric algebra $S^\cx T$. The group law in this case is given by
the unique algebra homomorphism
$$\mu^* \colon S^\cx T \to S^\cx T \ot S^\cx T$$
such that $t \mapsto 1 \ot t + t \ot 1$. Let  $ \bv_\gamma T :=
\spec \Gamma^\cx (T)$ \cite[A1]{bo.ncc} be  the divided power
envelope of the ideal of the zero section and $ \hat \bv_\gamma T :=
\spec \hat\Gamma_\cx (T)$ its completion  with respect to the
divided power filtration.  These are also group schemes, and the
group law
$$\mu^* \colon  \Gamma_\cx(T) \to \Gamma_\cx(T) \ot \Gamma_\cx(T)$$
is the unique divided power homomorphism sending $t$ to $1\ot t + t
\ot 1$.

  Let $\Omega$ be the dual of $T$ and  recall from \cite[A10]{bo.ncc}
that there is a natural isomorphism
$$\rho_n \colon S^n(\Omega) \to \cHom(\Gamma_n(T),\oh X)$$
for all $n$, and hence an isomorphism
$$S^\cx(\Omega) \to \cHom'(\Gamma_\cx(T),\oh X),$$
where the $\prime$  signifies the Matlis dual. The following
proposition is essentially contained in Theor\`eme V.1 of
\cite{roby.lplf}; we give a slightly simpler proof here for the
reader's convenience.

\begin{proposition}\label{sgd.p}
With respect to the pairing $\rho_n$ $S^n \times \Gamma_n \to \oh X$
defined above,
\begin{enumerate}
\item{The algebra multiplication
$S^i(\Omega) \ot S^j(\Omega) \to S^{i+j}(\Omega)$ is dual to the
comultiplication of $\Gamma_\cx (T)$ followed by projection:
$$\Gamma_{i+j}(T) \to \bigoplus_{a+b=i+j}\Gamma_a(T)\ot \Gamma_b(T) \to
 \Gamma_i(T) \ot \Gamma_j(T).$$}
\item{The algebra multiplication
$\Gamma_i(\Omega) \ot \Gamma_j(\Omega) \to \Gamma_{i+j}(\Omega)$ is
dual to the comultiplication followed by projection:
$$ S^{i+j}(T) \to \bigoplus_{a+b=i+j}S^a(T)\ot S^b(T) \to S^i(T) \ot S^j(T).$$}
\end{enumerate}
\end{proposition}
\begin{proof}
Let $\ell_i \colon T \to \Gamma_i(T)$ be the universal polynomial
law of degree $i$~\cite[A4]{bo.ncc}.  Then if $k = i+j$, $\ell_i\ot
\ell_j$ defines a polynomial law $T \to \Gamma_i(T) \ot
\Gamma_j(T)$, and hence a linear map $\mu'_{i,j}\colon \Gamma_{k}(T)
\to \Gamma_i(T) \ot \Gamma_j(T)$. Adding these up we find a map
$$\mu' \colon \Gamma_\cx(T) \to \Gamma_\cx(T) \ot \Gamma_\cx(T)$$
This  is the unique $\oh X$-linear map whose restriction to
$\Gamma_{k}(T)$ sends each $\ell_{k}(t)$ to $\sum_{i+j= k}
\ell_i(t)\ot \ell_j(t)$. On the other hand,  the comultiplication
$\mu^*$ of $\Gamma_\cx (T)$ is a divided power homomorphism sending
$t$ to $1\ot t + t\ot 1$, so
$$\mu^*(\ell_k(t)) = \mu^*(t^{[k]}) =
(t\ot 1 + 1 \ot t)^{[k]} = \sum_{i+j=k} t^{[i]}\ot t^{[j]}.$$ Thus,
$\mu' = \mu^*$. Now recall that $\rho_1$ is just the standard
duality map and that one deduces from $\rho_1$ a unique morphism of
algebras
$$\rho \colon S^\cx(\Omega) \to P(T,\oh X) \cong \Hom(\Gamma_\cx(T),\oh X),$$
where $P(T,\oh X)$ means the ring of polynomial laws $T \to \oh X$.
Thus if $x \in S^i(\Omega)$ and $y \in S^j(\Omega)$, $\rho_k(x y)$
is $\rho_i(x)\rho_j(y)$, where this product is taken in the algebra
$P(T,\oh X)$.  In other words, the following diagram commutes:
\begin{diagram}
\Gamma_k(T) & \rTo^{\rho_k(xy)} & \oh X \cr \dTo^{\mu'_{i,j}} &&
\uTo \cr \Gamma_i(T) \ot \Gamma_j(T) &\rTo^{\rho_i(x)\ot \rho_j(y)}&
\oh X \ot \oh X
\end{diagram}
Since $\mu'= \mu^*$ this proves (1).

Reversing the roles of $\Omega$ and $T$, we conclude from (1) that
the multiplication map $m_T$ of the algebra $S^\cx(T)$ is the dual
of the comultiplication map $\mu^*_\Omega$ of the group law of
$\Gamma_\cx(\Omega)$.  Since the latter commutes with algebra
structure $m_\Omega$ of the algebra $\Gamma_\cx(\Omega)$, it follows
that the dual $m'$ of $m_\Omega^*$ also commutes with $m_T$.  Thus
$m'$ is an algebra homomorphism $S^\cx(T) \to S^\cx(T) \ot
S^\cx(T)$. The same is true of the group law $\mu^*_T$.  We claim
that these two homomorphisms are equal, and it suffices to check
that this is true for elements of degree one.  In other words, we
have to check that the map $\mu_T^* \colon  T \to T \ot \oh X \oplus
\oh X \ot T$ is dual to the map $m_\Omega \colon \oh X \ot \Omega
\oplus \Omega \ot \oh X \to \Omega$.  But this is trivially true.
\end{proof}

It will perhaps be helpful to make everything explicit.  We shall do
this in the case $\cB = \oh {\cG}$ and $\cA = \iota_*\cB^\vee$.
Note that $\xi_i\mapsto -\xi_i$ induces an isomorphism $\cA \cong
\iota_*\cA$.

\begin{corollary}\label{expform.c}
Let $T$ and $\Omega$ be as above.  Suppose that also that $T$ is
free, and that  $(\xi_1,\ldots \xi_m)$ is a basis for $T$ and
$(\omega_1,\ldots \omega_m)$ the dual basis for $\Omega$. Let
$I_\cx$ denote the $I$-adic filtration of $S^\cx T$ (resp, the
PD-filtration of $\Gamma^\cx T$.
\begin{enumerate}
\item{If $\cG = \bv T$:
\begin{eqnarray*}
\oh \cG & = & S^\cx T; \quad \xi^I \xi^J = \xi^{I+J} \\
\iota_*\oh \cG^\vee & = & \hat  \Gamma\cx \Omega = \prod \Gamma_n
\Omega ;\quad \xi^I \omega^{[J]} =
      \omega^{[J-I]}\\
\theta_{\xi^n}(e_1\cnv e_2) & = & \sum_{i+j=n}{n \choose
i}\theta_\xi^i(e_1)\cnv  \theta_\xi^j(e_2) \quad\quad
\mbox{for}\quad
e_1\cnv e_2 \in E_1\cnv E_2\\
%\theta_\xi(h) &= &\theta_\xi \circ h - h\circ \theta_\xi \quad\quad\mbox{for}\quad h \in \cH(E_1,E_2)\\
\theta_{\xi^n}(h) &= &\sum_{i+j=n}{n\choose i}(-1)^j\theta_\xi^i \circ h\circ\theta_\xi^j \quad\quad\mbox{for}\quad h \in \cH(E_1,E_2)\\
\end{eqnarray*}}
\item{If $\cG = \hat\bv_\gamma T$:
\begin{eqnarray*}
\oh \cG & = & \hat \Gamma_\cx T; \quad \xi^{[I]} \xi^{[J]} ={(I+J)!\over I! J!} \xi^{[I+J]} \\
\iota_*\oh \cG^\vee & = & S^\cx \Omega;  \quad \xi^{[I] } \omega^{J}  =   {J!\over (J-I)!(I!)}    \omega^{J-I}\\
\theta_{\xi^{[n]}}(e_1\cnv e_2) & = &
\sum_{i+j=n}\theta_{\xi^{[i]}}(e_1)\cnv \theta_{\xi^{[j]}}(e_2)
\quad\quad \mbox{for}\quad
e_1\cnv e_2 \in E_1\cnv E_2\\
\theta_{\xi^{[n]}}(h) &= &\sum_{i+j}(-1)^j\theta_{\xi^{[i]}} \circ h
\circ \theta_{\xi^{[j]}} \quad\quad\mbox{for}\quad h \in
\cH(E_1,E_2)
\end{eqnarray*}}
\end{enumerate}
\end{corollary}

It is also of some interest to give an explicit formula for the
inverse to the map $h$ of Theorems \ref{charac.t} and
\ref{higtra.t}.  Let us do this when $G = \hat \bv_\gamma T$.
  Let $\kappa_i$ (the ``Casimir operator'') be the element of
$ S^i\Omega\otimes \Gamma^iT$ corresponding to the identity
homomorphism. In terms of the  bases above,
 $$\kappa_i := \sum_{|I|=i} \omega^I\ot \xi^{[I]} :=
 \sum_{|I|= i} \omega_1^{I_1}\ldots \omega_d^{I_d}\ot \xi_1^{[I_1]}\cdots \xi_d^{[I_d]}.$$
If $\theta$ is a $\cG$-field on $E$, $\kappa_i$ defines an
endomorphism
$$\kappa_{i,E} := \sum_{|I|= i} \omega^I \theta^{tot}_{\xi^{[I]}} \quad : \quad
S^i \Omega \ot E \to S^i \ot E$$ If $\theta$ is locally nilpotent,
so is $\theta_{tot}$, and in this case
$$\kappa_E  := \sum_{i=0}^\infty (-1)^i\kappa_{i,E}$$
is a well-defined endomorphism of $S^\cx\Omega \ot E$. The following
result is essentially classical; we shall omit the amusing and
elementary proof.

\begin{proposition}\label{casimir.p}
Let $(E,\theta)$ be a locally nilpotent $\Gamma_\cx T$-module.
\begin{enumerate}
\item The Casimir operator
$\kappa_E \colon S^\cx\Omega \ot E \to S^\cx \Omega \ot E$ defined
above is a projection operator  with image
$$\cT(E) := \left ( S^\cx\Omega \ot E\right)^{\theta_{tot}}.$$
and factors through the map $\sigma \colon S^\cx\Omega \ot E \to E$
defined by the the augmentation $S^\cx \Omega \to \oh X$.
\item The map
$$ h':= \sum_{J} \omega^J \ot( \kappa_E\circ \theta^{tot}_{t^{[J]}}) \quad :
\quad S^\cx \Omega \ot E \to S^\cx\Omega \ot \cT(E)$$ is the inverse
of  the map $h$ of Theorem~\ref{higtra.t}.
\item The map
$$\ov \kappa_E \colon E \to \cT(E)$$
induced by $\kappa$ is an isomorphism, inverse to the map
$$\eta_\tF := \sigma_\tF^* \circ j_E \colon \cT_{\lift X S}(E) \to E.$$
These isomorphisms take the field $\theta$ on $E$ to
$\iota_*\theta_{\cT}$. That is,
$$\theta_{t^{[I]}} \circ \overline \kappa_E = (-1)^I \ov \kappa_E \circ \theta_{t^{[I]}}$$
for all $I$.
\end{enumerate}
\end{proposition}

%With the same notation, let $E$ be an  object of $\bv T$-$HIG$.
%We can identify $\cH(\oh \cG, E)$ with $\Gamma_\cx(\Omega)\hat\cnv E$,
%and then write  explicit formulas:
%\begin{eqnarray*}
%\pi_E \colon \cA \ot E \to E &\quad & \sum \omega^{[I]} \ot e_I \mapsto e_0 \\
%i_E \colon E \to \cA \ot E  & \quad & e \mapsto 1\ot e \\
%\gamma_E \colon  \to \cA \ot E & \quad & e \mapsto \sum \omega^{[I]} \ot \xi^I (e_I)\\
%\xi_E \colon \cA  \ot E \to  \cA\ot E &\quad & \sum  \omega^{[I]} \ot e_I
%   \mapsto \sum_K \omega^{[K]} \left(\sum_{J\le K}{K\choose J} (-1)^|J|\xi^J(e_{K-J}) \right)
%\end{eqnarray*}

%
%
%\subsection{Higgs complexes and Higgs cohomology}
%Let us recall the meaning of the Higgs complex   Let us begin
%by recalling the definition of Higgs cohomology.  Let $X$ be a scheme,
%let $T$ be a  locally free sheaf of finite type on $X$, and let
%$E$ be a sheaf of $\oh X$-modules equipped with a $T$-Higgs field $\theta$.
%
%of sheaves of $\oh X$-modules and $\oh X$-linear maps:
%$$E \to E \ot \Omega \to E \ot \Omega^2 \to \cdots,$$
%where $\Omega^i := \Lambda^i \Omega$.  If $E$ is quasi-coherent and
%$\tilde E$ is the sheaf of $\oh{{\bf T}^*}$-modules corresponding to $E$,
%then this complex can be interpreted as $Ri^! \tilde E$, where $i \colon X
%\to {\bf T}^*$ is the zero section \cite[]{o.hcpc}.
%
%

\subsection{Azumaya algebras over group schemes}\label{aags.ss}

Recall that if $M$ is a locally free sheaf  of finite rank on a
scheme $Z$, then the (matrix) algebra $\cE:= \cEnd_{\oh Z}(M)$ is a
quasi-coherent sheaf of  $\oh Z$-algebras whose center is $\oh Z$.
Furthermore, the functor $E \to M \otimes E$ from the category of
sheaves of $\oh Z$-modules to the category of sheaves of left
$\cE$-modules is an equivalence, with quasi-inverse $\cHom_\cE(M, \
)$. A sheaf of algebras over $\oh Z$ which locally for the fppf
topology is isomorphic to $\cEnd_{\oh Z}(M)$ for some locally free
$M$ is called an  \emph{Azumaya algebra}. Note that an Azumaya
algebra, viewed  as $\oh Z$-module,  is locally free of rank $d^2$,
where $d$ is a locally constant function on $Z$. The function $d$ is
called  the \emph{index} of the Azumaya algebra. An Azumaya algebra
$\cE$ is said to be \emph{split} if it is isomorphic to $\cEnd_{\oh
Z}(M)$ for some such $M$; in this case $M$ is said to be a
\emph{splitting module}.
 If $\cA$ is an Azumaya algebra of rank $d^2$ and $M$ is a sheaf of left $\cA$-modules which
is locally free of rank $d$ over $\oh Z$, then the natural map $\cA
\to \cEnd_{\oh Z}(M)$ is necessarily an isomorphism, so $\cA$ is
split and $M$ is a splitting module~\cite{}.
%If $M$ is locally free,
%its rank is a locally constant function $d$, and $\cEnd_{\oh Z}(M)$
%is locally free of rank $d^2$.  Thus an Azumaya algebra $\cA$ is a locally
%free sheaf whose rank is a square; the square root of its rank is sometimes
%called the \emph{index} of $\cA$.

 Let us note for future reference
that if $M$  is an $\cA$-module which is locally free and of finite
rank as an $\oh Z$-module, then $M^\vee:= \cHom_{\oh Z}(M,\oh Z)$ is
a naturally a right $\cA$-module, and  for any $\cA$ module $N$
there is a natural isomorphism of $\oh Z$-modules:
\begin{equation}\label{mvee.e}
\cHom_\cA(M,N) \subseteq \cHom_{\oh Z}(M,N) \rTo^\cong M^\vee
\ot_{\oh Z} N\rTo M^\vee\ot_\cA  N
\end{equation}
We have found it convenient to use the language of 2-categories to
describe the relations among  Azumaya algebras.

\begin{definition}\label{azmorph.d}
Let $\cA$ and $\cB$ be Azumaya algebras over a scheme $Z$.
\begin{enumerate}
\item{A  \emph{1-morphism} $M: {\cal A}\to {\cal B}$ is a module $M$ over ${\cal A}^{op}\otimes  {\cal B}$ which is locally free
over ${\cal O}_Z$ of rank $(\ind{\cal A}) (\ind {\cal B})$, \ie,  a
splitting of the Azumaya algebra ${\cal A}^{op}\otimes  {\cal B}$.}
\item{If $M \colon \cA \to \cB$ and $N \colon \cB \to \cC$ are 1-morphisms, then
$N\circ M := N\ot_{\cB}M$}
\item{If $M$ and $N$ are 1-morphisms $\cA \to \cB$, a \emph{2-morphism}: $M \to N$ is an  isomorphism  of
  ${\cal A}^{op}\otimes  {\cal B}$-modules $M \to N$,  with the obvious notion of composition.}
\end{enumerate}
\end{definition}

For any object $\cA$, $\id_\cA$ is just the bimodule $\cA$.

If $M$ is a 1-morphism $\cA \to \cB$, then $M^\vee := \cHom_{\oh
Z}(M,\oh Z)$ is a ${\cB^{op}\ot \cA}$-module, \ie, a 1-morphism $\cB
\to \cA$, and the natural maps  $M^\vee \circ M \cong  End_\cB(M)
\cong \cA = \id_\cA$. Thus all 1-morphisms are equivalences (and all
2-morphisms are isomorphisms).

If $\cA$ is an Azumaya algebra over $Z$, let $\Mod(\cA)$ denote the
category of $\cA$-modules, which we may view as a stack over $Z$. A
1-morphism $M \colon \cA \to \cB$ gives rise to an equivalence of
categories:
  $$C^{-1}_M\colon   \Mod(\cA) \rTo^\sim \Mod(\cB) : E \mapsto  M\ot_\cA E,$$
together with  functorial isomorphisms
\begin{equation}\label{pr}
  C_M^{-1}(E\otimes _{{\cal O}_Z} L) \rTo^\sim C^{-1}_M( E)\otimes_{\oh Z} L
\end{equation}
for every $\oh Z$-module $L$. Conversely, every equivalence of
stacks, together with such a family of  isomorphisms (satisfying
suitable compatibilities), comes from a 1-morphism.

Recall that the category of $D_{X/S}$-modules has a tensor
structure: the tensor product of $D_{X/S}$-module $M$ and $N$ is
$M\otimes _{{\cal O}_X} N$, where the action of the vector fields on
$M\otimes _{{\cal O}_X} N$ is given by the Leibniz rule:
$\nabla_D(m\otimes n) := \nabla(m)\ot n + m\ot \nabla(n)$. We will
see that in characteristic $p$, the tensor structure on $\mic(X/S)$
comes from a more rigid structure on the Azumaya algebra $\cD_{X/S}$
which we will explain below. First we will give an abstract
definition on an Azumaya algebra over a group scheme, and then in
Example~\ref{azdt.e} we will construct a canonical tensor structure
on $\cD_{X/S}$.

\begin{definition}\label{aztens.d}
Let $G$ be a flat affine group scheme over a scheme $X$ and let
$\cA$ be an Azumaya algebra over $\oh G$. In the  notation of
(\ref{morphs.e}), a  \emph{tensor structure} on $\cA$ consists of
the following data:
\begin{enumerate}
\item  A 1-morphism
of  Azumaya algebras on $G\times G$:
  $$\delta:  \mu^*{\cal A} \to  p^*_1{\cal A}\otimes p^*_2{\cal A}. $$
\item{An associativity 2-morphism $\alpha$ as follows.
Note that $ \mu \circ (\mu,p_3) = \mu \circ (p_1,\mu)$ is the
multiplication morphism $\mu_3 \colon \cG \times \cG \times \cG \to
\cG$. Then $\alpha$ is a 2-morphism  between the 1-morphisms $f$ and
$g$ shown below:
\begin{diagram}
\mu_3^*\cA & \rTo^f & p_1^*\cA \ot p_2^*\cA \ot p_3^* \ca   \cr
&\rdTo_{(\mu,p_3)^*(\delta)} & \uTo_{\delta\ot \id} \cr
  & & \mu^*\cA \ot p_3^*\cA \cr
& \Downarrow\alpha\cr
 \mu_3^*\cA & \rTo^g & p_1^*\cA \ot p_2^*\cA \ot p_3^* \ca   & \cr
& \rdTo_{(p_1,\mu)^*(\delta)} & \uTo_{\id \ot \delta} \cr
 & & p_1^*\cA \ot \mu^*\cA
\end{diagram}
% \begin{diagram}
% \mu_3^*\cA & \rTo & p_1^*\cA \ot p_2^*\cA \ot p_3^* \ca
%   &&&& \mu_3^*\cA & \rTo & p_1^*\cA \ot p_2^*\cA \ot p_3^* \ca   & \cr
% &\rdTo_{(\mu,p_3)^*(\delta)} & \uTo_{\delta\ot \id}
% & &
%  \pile{\alpha\cr\Longrightarrow}&
%  && \rdTo_{(p_1,\mu)^*(\delta)} & \uTo_{\id \ot \delta} \cr
%   & & \mu^*\cA \ot p_3^*\cA
%  &&   && & & p_1^*\cA \ot \mu^*\cA .
% \end{diagram}
}\end{enumerate} The above data should satisfy the pentagon
condition~\cite[1.0.1]{de.tc}.
\end{definition}

If $\cG$ is commutative, then  a \emph{symmetric} tensor structure
on $\cA$ is a tensor structure as above together with an additional
datum $\gamma$. Note that ${\sigma^*(p_1^*\cA \ot p_2^*\cA)}
=(p_2^*\cA \ot p_1^*\cA)$, so that the standard commutativity
isomorphism for $\ot$ can be viewed as an isomorphism $\gamma_\ot
\colon \sigma^*(p_1^*\cA \ot p_2^*\cA) \cong p_1^* \cA \ot
p_2^*\cA$.  Since
 $\cG$ is commutative, $\mu \circ \sigma = \mu$, and hence the commutative
square below exists.  Then $\gamma$ is a 2-morphism
\begin{diagram}
\sigma^*\mu^*\cA &\rTo^{\sigma^*(\delta)} & \sigma^*(p_1^*\cA \ot
p_2^*\cA) & & \cr \dTo^= && \dTo_{\gamma_\ot}&&   \pile{\gamma \cr
\Longrightarrow}& & \Bigl(\mu^*\cA &\rTo^\delta & p_1^*\cA \ot
p_2^*\cA  \Bigr)\cr \mu^*\cA &\rDotsto & (p_1^*\cA \ot p_2^*\cA) & &
\end{diagram}
such that $\gamma^2 = \id$.  The associativity morphism $\alpha$ and
the commutativity morphism $\gamma$ should  also satisfy the hexagon
axiom \cite[1.0.2]{de.tc}.

Azumaya algebras with (resp. symmetric) tensor structure also form a
2-category: a 1-morphisms $\cA \to \cB$ is by definition a pair
$(M,\lambda)$, where $M$ is a 1-morphism of the underlying Azumaya
algebras as before and $\lambda$ is a 2-morphism
\begin{small}
\begin{diagram}
\mu^*\cA & \rTo^{\delta_A} & p_1^*\cA \ot p_2^*\cA &&& && \mu^*A\cr
&\rdTo  & \dTo_{p_1^*(M)\ot p_2^*(M)}  &&&\pile{\lambda \cr
\Longrightarrow}& &\dTo^{\mu^*(M)} & \rdTo\cr && p_1^*\cB\ot
p_2^*\cB &&&&  &\mu^*\cB &\rTo_{\delta_B}& p_1^*\cB \ot p_2^*\cB
\end{diagram}
\end{small}
compatible with the associativity (resp. and the commutativity)
2-morphisms.

Let $(\cA,\delta,\alpha)$ be an Azumaya algebra with a tensor
structure.  Then the category $\Mod(\cA)$ of modules over $\cA$ is
endowed with tensor structure: given ${\cal A}$-modules $M_1$ and
$M_2$ we define the tensor product
\begin{equation}\label{cvpr}
M_1\circledast M_2= \mu_* \left(\delta^{-1}(p_1^* M_1 \otimes  p_2^*
M_2)\right ).
\end{equation}
%We are going to see that $(\Mod(\cA),\cnv)$ has a unit object,
%and hence is a tensor category in the sense of \cite{de.tc}.
Recall from \cite{de.tc} that a unit object of a category ${\cal M}$
with a tensor structure is a pair $(U, h)$, where $U$ is an object
of ${\cal M}$ and $h$ is an isomorphism
$$h: U \simeq U \circledast U ,$$
such that the functor
 $\cnv\, U: {\cal M} \to {\cal M} $ sending an object $M$ to
$M\circledast U$ is an equivalence of categories. It is shown in
({\it loc.cit.}) that the unit object is unique up to a unique
isomorphism and that for any object $M$ there is a functorial
isomorphism $M\circledast U \simeq M$.

\begin{lemma}\label{ctsp.l}
Let $\cA$ be an Azumaya algebra over $\cG$  equipped with  a tensor
structure $(\delta,\alpha)$~(\ref{aztens.d}). Then the  restriction
$i^*{\cal A}$ of $\cA$ to the zero section has a canonical splitting
$N_0$.  Moreover, there is an isomorphism $h:  i_* N_0 \simeq i_*
N_0 \circledast i_* N_0$, and the pair $U :=(i_* N_0, h)$ is the
unit object of $\Mod(\cA) $.
%and the functor
% $\circledast i_* N_0:\Mod(\cA) \to \Mod(\cA) $ sending a $\cA$-module $M$ to
%$M\circledast i_* N_0$ is an equivalence of categories.
 %Thus  the pair $(i_*N_0, u)$ is a unit object of $\Mod(\cA)$.
\end{lemma}
 \begin{proof}
Since $\mu\circ i  = p_i\circ i = i$, the restriction $\delta_0:
i^*{\cal A}\to i^*{\cal A}\otimes _{{\cal O}_X} i^*{\cal A}$ of
$\delta$ to the zero section of $\cG$ is an   $i^*(\cA^\op\ot \cA
\ot \cA)$-module $P$. Then the $i^*\cA$-module $N_0: =
\cHom_{i^*(\cA^\op \ot \cA)} (i^*\cA, P)$ gives a splitting of
$i^*\cA$.

Next we construct the isomorphism $h$ .  By definition,
\begin{eqnarray*}
U\cnv U &=& \mu_*(\delta^{-1}(p_1^* U \ot p_2^* U)) \\
   & \cong& i_*\delta_0^{-1}(N_0\otimes N_0) \cong \cHom_{i^*(\cA\ot\cA)}(P,N_0\ot N_0)
\end{eqnarray*}
Evaluation of homomorphisms defines an isomorphism of
${i^*(A^{op}\ot A \ot A)}$-modules $i^*\cA \ot N_0 \to P$.  Thus
\begin{eqnarray*}
U\cnv U &\cong &i_*\cHom_{i^*(\cA\ot \cA)}(i^*\cA \ot N_0,N_0 \ot N_0)\\
 &\cong& i_*\cHom_{i_*\cA}(i^*\cA, N_0) \cong U
\end{eqnarray*}

Finally, we have to prove that the functor $\cnv\, U:\Mod(\cA) \to
\Mod(\cA) $  is an equivalence of categories. Let $\Gamma_j \colon
\cG \to \cG\times \cG$ be the graph of the zero section.  Since
$\mu\circ \Gamma_j = \id$, $\Gamma_j^*(\delta)$ is a 1-morphism:
$\delta_1: \cA \to \cA \otimes j^* \cA$.  If $M \in \Mod(\cA)$,
$$M\cnv U := \delta _1 ^{-1}(M\otimes p^*N_0).$$
Since $p^*N_0$ is a splitting of $j^* \cA$, the functor $\otimes
p^*N_0: \Mod(\cA) \to \Mod(\cA \otimes j^* \cA)$ is an equivalence.
Since
 $\delta_1 ^{-1}:\Mod(\cA \otimes j^* \cA)\to \Mod(\cA)$ is also an equivalence, the lemma is proved.
\end{proof}

Observe that a tensor structure on $\cA$ induces a canonical
1-morphism $\iota^*\cA \cong \cA^{op}$, obtained by pulling back
$\delta$ by the graph of $\iota$, since by Lemma~\ref{ctsp.l},
$i^*\cA$ is canonically split.  It follows that the category
$\Mod(\cA)$ has inner Homs. Let
$$\delta' \colon {\mu'}^*\cA \to p_1^*\cA^{op} \ot p_2^*\cA$$
be the pullback of $\delta$ by $(\iota,\id)$ composed with the
1-morphism $\iota^*\cA \ot \cA \to \cA^{op}\ot\cA$. If $E_1$ and
$E_2$ are objects of $\Mod(\cA)$, $\cH ex(E_1,E_2)$ is naturally a
$p^*_1\cA^{op}\ot p_2^*\cA$-module, and
$$\cH(E_1, E_2) = \mu'_*\left({\delta'}^{-1}\cH ex(E_1,E_2)\right).$$

%To simplify our notations we will write  $N_0$ instead of  $i_* N_0$.
If $\cG$ is commutative and the tensor structure is endowed with a
commutativity morphism,  then $\Mod(\cA)$ becomes a  tensor category
in the sense of \cite{de.tc}. From now on,
 we assume this to be the case.

\begin{definition}\label{tenspl.d}
A \emph{tensor splitting} of an Azumaya algebra $\cA$ equipped with
a tensor structure over $G$ is a 1-morphism (in the category of
Azumaya algebras with symmetric tensor structure): ${\cal O}_G \to
\cA$.
\end{definition}

Note that, in general,  a 1-morphism $M\colon {\cal A} \to {\cal B}$
gives rise to equivalence of  tensor categories:
$$C^{-1}_M \colon \Mod(\cA) \sim \Mod(\cB), \quad  C^{-1}_M(M_1\cnv M_2)\cong   C^{-1}_M(M_1) \cnv C^{-1}_M(M_2).$$
In particular, a tensor splitting $M$ gives an  equivalence between
the tensor categories: $C^{-1}_M \colon \Mod({\cal O}_G) \sim
\Mod(\cA)$. Observe that ${\cal O}_G$ is a commutative coalgebra
with counit in the category $\Mod({\cal O}_G)$. Thus, by ``transport
of structure'' the splitting module $M \cong C_M^{-1}(\oh G)$
becomes a commutative coalgebra with counit in $\Mod(\cA)$. In
other words, we have canonical morphisms:
$$M\to M \circledast M, \; e: M\to U .$$  In the case of the split Azumaya algebra $\cA = \oh G$
with the obvious tensor structure, the notion of a tensor splitting
boils down to the notion of a character sheaf introduced in
Definition~\ref{char.d}.

%A tensor
% In particular, a tensor splitting of ${\cal A}$
%is an equivalence
%  $$C_V: {\cal M}({\cal A})  \simeq {\cal M}({\cal O}_Y)$$
%together with functorial isomorphisms (\ref{tpr}) and
%$$ C_V(M_1\circledast M_2)\simeq   C_V (M_1) \circledast C_V(M_2)$$
%In  turn, this is equivalent to given ${\cal A}$-module $M$,
%locally free of rank $\sqrt{rk \, {\cal A}}$, together a structure
%the case, when $G$ is an affine commutative group scheme

\begin{example}\label{azdt.e}{\rm
Let us explain how the above formalism works in the case of the
Azumaya algebra of differential operators. First recall the
following lemma.

\begin{lemma}\label{pcurvet.l}
 Let  $M$ and $N$  be  objects of $\mic(X/S)$, and let $\psi_M$ and $\psi_N$ be
their $p$-curvatures. Then the p-curvature of the tensor product and
internal Hom:
 $$\psi_{M\ot N} \colon  M\ot_{\oh X} N \to  M\ot_{\oh X} N \ot F_{X/S}^*\Omega^1_{X'/S}$$
$$\psi \colon  \cHom_{\oh X}(M,N) \to  \cHom_{\oh X}(M,N) \ot F_{X/S}^*\Omega^1_{X'/S}$$
are  given respectively by the formulas
$$\psi_{M\ot N}(m\otimes n) =  \psi_M(m)  \otimes n+   m \otimes \psi_N(n) $$
$$\psi_\cH(h)   =  \psi_N\circ h -    h\circ \psi_M .$$
\end{lemma} \qed

 By definition, giving a tensor structure on  ${\cal D}_{X/S}$  amounts to giving a 1-morphism of the Azumaya algebras on $ {\bf T}^*_{X'/S}\times _{X'} {\bf T}^*_{X'/S} $
$$ \delta : \mu^* {\cal D}_{X/S} \to p^*_1 {\cal D}_{X/S} \otimes p^*_2 {\cal D}_{X/S} $$
together with the associativity and commutativity 2-morphisms
$\alpha$ and $\gamma$ . We will construct  $ \delta $ as follows.
View $\cD_{X/S}$ as a left module over itself, and endow $\Delta: =
\cD_{X/S} \otimes _{{\cal O}_X} \cD_{X/S}$ with the left
$\cD_{X/S}$-module structure defined by the Leibnitz rule above.
Note that the right action of $\cD_{X/S}$ on itself makes $\Delta$ a
right module over the algebra $\cD_{X/S} \otimes
_{F^{\cdot}_{X/S}{\cal O}_{X'}} \cD_{X/S}$ or, equivalently,  a left
module over $ \cD_{X/S}  \otimes _{F^{\cdot}_{X/S}{\cal O}_{X'}}
(\cD_{X/S} \otimes _{F^{\cdot}_{X/S}{\cal O}_{X'}} \cD_{X/S})^{op}$.
Lemma~\ref{pcurvet.l}  shows that the action factors through $
\cD_{X/S}  \otimes_{F^{\cdot}_{X/S} S^{\cdot}T_{X'/S} }   (\cD_{X/S}
\otimes _{F^{\cdot}_{X/S}{\cal O}_{X'}} \cD_{X/S})^{op}$, where the
$F^{\cdot}_{X/S} S^{\cdot}T_{X'/S} $-module structure on $\cD_{X/S}
\otimes _{F^{\cdot}_{X/S}{\cal O}_{X'}} \cD_{X/S}$ is given by
comultiplication:
$$\mu^*:F^{\cdot}_{X/S} S^{\cdot}T_{X'/S} \to F^{\cdot}_{X/S} ( S^{\cdot}T_{X'/S} \otimes_ {{\cal O}_{X'}} S^{\cdot}T_{X'/S} ) .$$
Thus, $\Delta $ gives rise to a module over the Azumaya algebra
 $\mu^* {\cal D}_{X/S}  \otimes  (p^*_1 {\cal D}_{X/S} \otimes p^*_2 {\cal D}_{X/S})^{op}$.
It is easy to see using local coordinates that as a module over the
center ${\cal  O}_{   {\bf T}^*_{X'/S}\times _{X'} {\bf T}^*_{X'/S}
}$ it is locally free of rank $p^{3 dim_S X}$ which is equal to the
index of the Azumaya algebra. Therefore we get a 1-morphism $\delta
$.

Next, let us construct the commutativity 2-morphism $\gamma$.
Consider the automorphism $\Sigma $ of the algebra
 $ \cD_{X/S}  \otimes _{F^{\cdot}_{X/S} S^{\cdot}T_{X'/S} }   (\cD_{X/S} \otimes _{F^{\cdot}_{X/S}{\cal O}_{X'}} \cD_{X/S})^{op}$
which sends $D_1 \otimes D_2 \otimes D_3$ to $D_1 \otimes D_3
\otimes D_2$.  Then, giving $\gamma$ amounts to giving  an
isomorphism
$$  \cD_{X/S} \otimes _{{\cal O}_X} \cD_{X/S} \to  (\cD_{X/S} \otimes _{{\cal O}_X} \cD_{X/S})^{\Sigma}$$
of modules over $\cD_{X/S}  \otimes _{F^{\cdot}_{X/S}
S^{\cdot}T_{X'/S} }   (\cD_{X/S} \otimes _{F^{\cdot}_{X/S}{\cal
O}_{X'}} \cD_{X/S})^{op}$ The obvious permutation does the trick.

To construct  the associativity morphism $\alpha$,  we note that
the two 1-morphisms:
$$ \mu_3^* {\cal D}_{X/S} \to p^*_1 {\cal D}_{X/S} \otimes p^*_2 {\cal D}_{X/S} \otimes  p^*_3 {\cal D}_{X/S}  $$
in Definition(\ref{aztens.d}) are given by the same module
$\cD_{X/S} \otimes _{{\cal O}_X} \cD_{X/S}\otimes _{{\cal O}_X}
\cD_{X/S}$. With this identification, we let  $\alpha$ be the
identity morphism.

It remains to show that $\alpha$ and $\gamma$  satisfy the pentagon
and hexagon axioms.  To save space, we may use the following trick.
First,  one can easily check the axioms for the restriction of
${\cal D}_{X/S}$ to the zero section $X' \to {\bf T}^*_{X'/S}$.
Furthermore, since the statement is local on $X$ and stable under a
base change $T\to S$,  we may assume
 that $S$ is reduced.
 Then any 2-morphism over ${\bf T}^*_{X'/S}\times _{X'} {\bf T}^*_{X'/S} \times \cdots \times _{X'} {\bf T}^*_{X'/S}$
is uniquely determined by its restriction to zero section. This
completes the proof.}
\end{example}

\begin{remark}
If  the base  $S$ is normal and reduced, one can prove that   the
tensor structure on ${\cal D}_{X/S}$  equipped  with an isomorphism
$\phi: N_0 \simeq F_{X/S *} {\cal  O}_X$ between the canonical
splitting module over the zero section $X' \to {\bf T}^*_{X'/S}$
described in Lemma~\ref{ctsp.l}  and the splitting $F_{X/S *} {\cal
O}_X$ of Remark~\ref{splitdes.r} is a unique (up to a unique
isomorphism).
\end{remark}

Next we shall discuss tensor Azumaya algebras over a formal and PD
completion of a group scheme  $G$ along its zero section.

Let $\hat G$ (resp. $\ccG$) denote the formal (resp. PD) completion
of $G$ along its zero section, viewed as a locally ringed space. Let
$\hat G^n$ (resp. $\ccG^n$) be the formal (resp. PD) completion of
$\cG^n$ along its zero section. Then the  definition a tensor
structure for Azumaya algebras over  $\hat G$  (resp. $\ccG$) is the
same as before, with these completed fiber products in place of
$\cG^n$.

Denote by $i_k: \cG_k \hookrightarrow \hat\cG$ (resp. $i_k: \ccG_k
\hookrightarrow \ccG$) the subscheme defined by the $k+1$ power
(resp. divided power) of the ideal of the zero section; if $M$ is an
$\oh {\hat G}$-module (resp. $\oh \ccG$-module),  $M_k :=
i_{k*}i_k^*M$. The multiplication map $\mu$ gives rise to morphisms
of schemes:  $\mu: \cG_k \times \cG_l \to \cG_{k+l} $ (resp. $\mu:
\ccG_k \times \ccG_l \to \ccG_{k+l} $). The category of
$\cA$-modules with the convolution product defined the formula
(\ref{cnv.d}) is not a tensor category, since the associativity
constraint fails in general. However, the  subcategory
$\Mod^\cx(\cA)$ of $I_\cx$-continuous modules is stable under the
convolution product and is a tensor category. (Recall that a module
is said to be $I_\cx$-continuous if each local section is locally
supported on $\cG_k$ (resp. $\ccG_k$), for some $k$.)

Let $M: {\cal O}_{\ccG}\to \cA$ be a tensor splitting.  The
splitting module $M$ inherits a structure of a topological
commutative coalgebra with counit, that is, a family of maps $M_n
\to M_k \circledast M_l$, for any $n\geq l+k$, and $e: M \to U$
satisfying the obvious compatibilities.

\begin{proposition}\label{spde.p}
Let $\cA$ be a tensor Azumaya algebra over  $\hat \cG$ (resp
$\ccG$). Then the following data are equivalent.
\begin{enumerate}
\item{A tensor splitting of $\cA$.}
\item{A splitting module $M$ for $\cA$  with the structure
of  a topological colagebra with counit.}
\item{A splitting module $M$ for $\cA$ whose topological dual $\dirlim \cH(M_n,U)$
is endowed with the structure of an algebra with unit in the tensor
category $\Mod^\cx(\cA)$.}
\end{enumerate}
\end{proposition}
\begin{proof}
 We will prove the result for $\ccG$; the formal case is similar.
The equivalence of (2) and (3) is clear, and we have already shown
how the data of (1) give the data of (2).  It remains to explain how
the data of (2) give (1).  Let $M$ be a splitting with a coalgebra
structure: $M_{k+l} \to M_k \circledast M_l$. By definition,  $M_k
\circledast M_l= \mu_* \delta ^{-1} (p_1^*M_k \otimes p_2^* M_l)$.
 Thus, using the adjointness property
of $\mu_*$ and $\mu^*$, we get a morphism: $\mu^* M_{k+l} \to \delta
^{-1} (p_1^*M_k \otimes p_2^* M_l)$. Let us consider the restriction
of the above morphism to $\ccG_k \times \ccG_l $:
\begin{equation}\label{coalgebra}
(i_k, i_l)^* \mu^* M=(i_k, i_l)^* \mu^* M_{k+l} \to  \delta ^{-1}
(p_1^*M_k \otimes p_2^* M_l)= (i_k, i_l)^* \delta ^{-1} (p_1^*M
\otimes p_2^* M)
\end{equation}

These  morphisms are compatible with change of $k$ and $l$, and we
claim that they  are all isomorphisms. Indeed, both $(i_k, i_l)^*
\mu^* M$ and \\ $(i_k, i_l)^* \delta ^{-1} {(p_1^*M \otimes p_2^*
M)}$ are splittings of the Azumaya algebra $\mu^* \cA$ over $\ccG_k
\times \ccG_l$, and, in particular, they are locally free over
 $\ccG_k \times \ccG_l$. Also, the existence of counit: $M\to M_0\simeq U$ implies that, for $k=l=0$,
 the morphism (\ref{coalgebra}) is an isomorphism. Hence (\ref{coalgebra}) is an isomorphism, for any $k$ and $l$.
Thus, we get an isomorphism:  $  \mu^* M \simeq  \delta ^{-1}
(p_1^*M \otimes p_2^* M)$.
\end{proof}

Let  ${\cal A}$ be an   Azumaya algebra  on $\ccG$ with  a symmetric
tensor  structure. Let $\tspl$ be the stack on $X$ assigning to a
scheme $U$ etale over $X$   the groupoid of    tensor splittings  of
the Azumaya algebra ${\cal A}_U$ over $\ccG\times U$ (i.e.
1-morphisms from ${\cal A}_U$ to the trivial Azumaya algebra ${\cal
O}_{\ccG\times U}$ ).  Let ${\spl_1}$ be be the stack sending $U$ to
the   groupoid whose objects are pairs $(N, e)$, where $N$ is a
splitting of the restriction   $i_1^* {\cal A}_U$ and  $e$ is  an
isomorphism $e: i^*N\simeq N_0$.

\begin{proposition}\label{tsp.p}
Let $G$ be a smooth commutative group scheme  and ${\cal A}$  be an
Azumaya algebra on $\ccG$ with  a symmetric tensor  structure, which
admits a tensor splitting etale locally on $X$. Then the  obvious
restriction functor:
$$i_1^*: {\tspl} \to {\spl_1}$$
is an equivalence of stacks.
\end{proposition}

\begin{remark}
The  stack $\spl_1$ is, in fact, a gerbe.  That is, there exists an
etale covering of $U_i$ of $X$ such that, for each $i$,
${\spl_1}_{U_i}$ is non-empty,
 and, for any etale morphism
$U\to X$, any two objects of ${\spl_1}_{U}$ are locally isomorphic.
Thus, the proposition implies that the stack $\tspl$ is also a
gerbe.
\end{remark}

\begin{proof}
Recall that, for a smooth commutative group scheme $G$, the
exponential map induces an isomorphism:
$$exp: {\cal T} \simeq \ccG$$
where ${\cal T}$ is the completed PD envelope of the zero section
$X\hookrightarrow T$ of the Lie algebra $T$. Thus, without loss of
generality we can replace $G$ by the vector group $T$.

To prove that $i_1^*$, is an equivalence we will construct the
inverse functor $exp:  \spl_1 \to \tspl$ explicitly. The reason we
call it the exponential will be clear in a moment. Given a ${\cal
A}$-module $M$ we denote by $\Gamma ^k M$ the ${\cal A}$-submodule
of the tensor power $ M^{\circledast ^{k}}$ which consists of
$S_k$-invariant sections.

Let $(N,e)$ be an object of ${\spl_1}$    over $U$. Define $$M_1=
i_{1*} N, \;  M_k= \Gamma^{k} M_1 .$$ We claim that $M_k$ is a
splitting of  $ {\cal A}_U$ over $\ccG_k$. Indeed, since the
statement is local on $U$ we may assume that ${\cal A}$ has a tensor
splitting. A tensor splitting, in turn, gives an equivalence between
the tensor category of  $I_\cdot$-continuous ${\cal A}_U$-modules
and
 $I_\cdot$-continuous ${\cal O}_{\ccG \times U}$- modules. This equivalence takes $M_1$ to a line bundle over
$\ccG_1 \times _X U$. Shrinking the base, if necessary, we may
assume that the line bundle is  trivial. Then the claim follows from
the fact that, for smooth $G$, the $S_k$-invariant multiplication
morphism
   $$\ccG_1^k= \ccG_1\times \ccG_1 \times \cdots \ccG_1 \to \ccG_k$$
induces an isomorphism $\Gamma ^k {\cal O}_{\ccG _1} \simeq {\cal
O}_{\ccG _k}$.

 The morphism $e: M_1 \to U$ gives the map:
    $$ M_k \to M_{k-1}, $$
and, moreover, $i_{k-1}^*M_k\simeq M_{k-1}$.

We define $exp(N)$ to be
         $$exp(N)= \lim _{\leftarrow} M_k $$
First of all, it is clear that $M:= exp(N)$ is a splitting of ${\cal
A}_U$. Furthermore, it is a coalgebra with counit. The coalgebra
structure comes from the canonical morphism:
$$M_{k+l} = \Gamma^{k+l} M_1 \hookrightarrow \Gamma^{k} M_1  \circledast \Gamma^{l} M_1= M_k \circledast  M_l .$$
By Proposition~\ref{spde.p} it gives a tensor splitting.

%We claim that it is, in fact, a tensor splitting. By definition,
%a tensor splitting is a splitting $N$ together with an isomorphism
%of $\mu^* {\cal A}$-modules
%$$\delta^{-1}(p_1^*N \otimes p_2^*N)\simeq \mu^* N $$
%satisfying certain compatibility conditions.
%It is enough to construct a compatible system of isomorphisms on
%$\ccG_k \times \ccG_l$:
%$\delta^{-1}(p_1^*M_k \otimes p_2^*M_l)\simeq \mu^* M_{k+l}$.
%Consider the canonical morphism:
%$$M_{k+l}\to \mu_*\delta^{-1}(p_1^*M_k \otimes p_2^*M_l)=
%M_k\circledast M_l .$$
%By adjointness, it gives rise to the map
% \begin{equation}\label{comult}
%\mu^* M_{k+l} \to
%\delta^{-1}(p_1^*M_k \otimes p_2^*M_l)
%\end{equation}

%Notice that the restriction of the  morphism (\ref{comult})
%to the zero section $X\hookrightarrow \ccG \times _X \ccG$ is an isomorphism:
%it can be identified with the identity map $N_0 \simeq N_0\otimes
%N_0= N_0$. Therefore (\ref{comult}) itself is an isomorphism.

It is clear that $i_1^*(exp(N))=N$, so it remains to construct a
canonical isomorphism
\begin{equation}\label{expis}
M\simeq exp(i_1^*(M)),
\end{equation}
 for any tensor splitting
$M$ over $U$. But this is immediate:
 the splitting $M$ gives an equivalence of tensor categories:
$C_M:\Mod({\cal A}_U) \simeq \Mod({\cal O}_{\ccG \times _X U)})$,
which takes $M$ to ${\cal O}_{\ccG \times _X U} $.
 We define (\ref{expis}) to be the morphism
corresponding under the above equivalence to  the canonical
isomorphism
$${\cal O}_{\ccG}\simeq \lim _{\leftarrow} \Gamma ^k{\cal O}_{ \ccG_{1}}.$$

%The latter comes from the comultiplication on the structure sheaf of a group %scheme:
% ${\cal O}_{\hat S_{pd}}\to {\cal O}_{\hat S_{pd}}\circledast {\cal O}_{\hat %S_{pd}}  $.
\end{proof}

Denote by $\ccG ^*$ the sheaf on $X$ whose section over  a scheme
$U$ etale over $X$ is the group of homomorphisms: $Hom_U(\ccG \times
U; \mathbb{G}_m)$.

Let ${\cal A}$  be an Azumaya algebra on $\ccG$ with  a symmetric
tensor structure, which admits a tensor splitting etale locally on
$X$. We assign to ${\cal A}$ a class $[{\cal A}]$ in $ H^2_{et}(X,
\ccG ^* )$ as follows. Choose an etale  covering of $U_i$ of $X$
together with tensor splittings $N_i$ on each $U_i$
 and tensor isomorphisms
$\phi_{ij}: N_i\simeq N_j $. Then $ \phi_{ki}\phi_{jk}\phi_{ij}$ is
a Cech cocycle with coefficients in $\ccG ^*$. Similarly, given an
Azumaya algebra over $\ccG_1$ together with a splitting over the
zero section we can construct an element of
 $H^2_{et}(X, {\cal O}^*_{1,\ccG_1 }) =  H^2_{et}(X, T^*_{G/X})$, where
${\cal O}^*_{1,\ccG_1 }$ denotes the group of invertible functions
on $\ccG_1$ equal to $1$ on $X\hookrightarrow \ccG_1$.

Recall  that for any smooth commutative group scheme $G$ over $X$ we
have the isomorphism:
  $$exp: {\cal O}^*_{1,\ccG_1 }= T^*_{G/X} \to \ccG ^*  $$
The inverse map is the restriction $i_1^*$.
  Passing to cohomologies we get a morphism:
  $$exp: H^2_{et}(X, T^*_{G/X})\simeq  H^2_{et}(X, \ccG ^*   ) \to H^2_{et}(X,{\cal O}^*_{\ccG} )$$

\begin{proposition}\label{azexp.p}
Let ${\cal A}$  be an Azumaya algebra on $\ccG$ with
 a symmetric tensor  structure, which admits a tensor splitting etale locally on $X$. Then
$$ exp [i_1^*{\cal A}]= [{\cal A}].$$
\end{proposition}
Since the restriction $i^*_1: H^2_{et}(X, \ccG ^*   ) \to
H^2_{et}(X, T^*_{G/X})$ is an isomorphism, the proposition follows
from the obvious fact: $[i_1^*{\cal A}]= i_1^*[{\cal A}]$.

%  \bibliography{all,ogus}
%  \bibliographystyle{acm}

\end{document}